\documentclass[11pt,leqno]{article}
\normalfont
\renewcommand{\baselinestretch}{1.25}

\usepackage{graphicx}
\usepackage[dvipsnames]{xcolor}

\newcommand{\IntroSection}{1}
\newcommand{\WeylSection}{2}
\newcommand{\PosetSection}{3}
\newcommand{\SplittingSection}{4}
\newcommand{\FibrousSection}{5}
\newcommand{\CrystallineSection}{6}
\newcommand{\CrystalSection}{7} 
\newcommand{\CriteriaSection}{8}
\newcommand{\GaussianSection}{9}

\newcommand{\IntroFig}{Figure 1.1}
\newcommand{\IntroTable}{Table 1.1}

\newcommand{\DiagramSymmetryFigure}{Figure 2.1}
\newcommand{\WeightRemarkOne}{Proposition 2.1}
\newcommand{\AltProp}{Proposition 2.2}
\newcommand{\WeylsDenomTheorem}{Theorem 2.3}
\newcommand{\WeylBialternantDef}{Definition 2.4}
\newcommand{\WeylBialternantProblem}{Problem 2.5}
\newcommand{\InitialTheorem}{Theorem 2.6}
\newcommand{\KostantTheorem}{Theorem 2.7}
\newcommand{\WeylCharLemma}{Lemma 2.8}
\newcommand{\FinitenessTheorem}{Theorem 2.9}
\newcommand{\FinitenessCorollary}{Corollary 2.10}
\newcommand{\FreudTheorem}{Theorem 2.11}
\newcommand{\KostkaCorollary}{Corollary 2.12}
\newcommand{\BasisProp}{Proposition 2.13}
\newcommand{\BasisTheorem}{Theorem 2.14}

\newcommand{\StarProp}{Proposition 2.16}
\newcommand{\TFAEConsequences}{Theorem 2.17}
\newcommand{\UnimodalProblem}{Problem 2.18}
\newcommand{\CaseAExample}{Example 2.19}
\newcommand{\TFAEProp}{Proposition 2.20}
\newcommand{\TFAETheorems}{Theorems 2.7 and 2.11}
\newcommand{\WeylsTheorem}{Theorem 2.21}
\newcommand{\UnimodalTheorem}{Corollary 2.22}

\newcommand{\FirstOpsLemma}{Lemma 3.1}
\newcommand{\SCDLemma}{Lemma 3.2}
\newcommand{\WGFLemma}{Lemma 3.3}
\newcommand{\WeightLemma}{Lemma 3.4}
\newcommand{\WInvariantLemma}{Lemma 3.5}
\newcommand{\WeightRemarkTwo}{Proposition 3.6}
\newcommand{\PiLemma}{Lemma 3.7}
\newcommand{\PiTheorem}{Theorem 3.8}
\newcommand{\PiSaturatedCorollary}{Corollary 3.9}
\newcommand{\PiJCorollary}{Corollary 3.10}
\newcommand{\PiRankCorollary}{Corollary 3.11}
\newcommand{\ClosingWeightRemarks}{Remarks 3.12}

\newcommand{\PiResults}{Propositions 2.1 and 3.6}

\newcommand{\InitialSplittingTheorem}{Theorem 4.1}
\newcommand{\NewSplittingRemarks}{Remarks 4.2}

\newcommand{\InitialSplittingExample}{Example 4.3}
\newcommand{\NotSturdyFigure}{Figure 4.1}
\newcommand{\InitialSplittingDefinition}{Definition 4.4}
\newcommand{\NewDefinitionRemarks}{Remarks 4.5}
\newcommand{\NewDefinitionRemark}{Remark 4.5}
\newcommand{\NotSchurPositiveFigure}{Figure 4.2}
\newcommand{\StanleyProblem}{Problem 4.6}
\newcommand{\MainCorollary}{Proposition 4.7}
\newcommand{\OperationsLemma}{Proposition 4.8}
\newcommand{\MinimalProp}{Proposition 4.9} 
\newcommand{\EdgeMinimalExamples}{Examples 4.10}
\newcommand{\EdgeMinimalQuestions}{Questions 4.11}
\newcommand{\TangledFigure}{Figure 4.3}
\newcommand{\UniqueMaximalProp}{Proposition 4.12}
\newcommand{\MinQuasiMinLemma}{Lemma 4.13}
\newcommand{\MinQuasiMinProposition}{Proposition 4.14}
\newcommand{\QuasiMinLemma}{Proposition 4.15}
\newcommand{\MinQuasiMinTheorem}{Theorem 4.16}
\newcommand{\MinQuasiMinRemark}{Remark 4.17}
\newcommand{\SupportingGraphProp}{Proposition 4.18} 
\newcommand{\SupportingGraphRemarks}{Remarks 4.19}
\newcommand{\UniqueSplittingPropOne}{Theorem 4.20}
\newcommand{\UniqueSplittingPropTwo}{Corollary 4.21}

\newcommand{\FibrousProductLemma}{Lemma 5.1}
\newcommand{\EFLemma}{Lemma 5.2}
\newcommand{\CrystalProductFigure}{Figure 5.1}
\newcommand{\CrystalProductFigures}{Figures 5.1.1 and 5.1.2}
\newcommand{\NoCommuteFigure}{Figure 5.2}
\newcommand{\CrystalExample}{Example 5.3}
\newcommand{\AssociativeCrystalLemma}{Lemma 5.4}
\newcommand{\NoncommuteExample}{Example 5.5}
\newcommand{\CommuteProblem}{Question 5.6}
\newcommand{\FibrousStructureLemma}{Lemma 5.7}
\newcommand{\WGFCorollary}{Corollary 5.8}
\newcommand{\FibrousChainStructureLemma}{Lemma 5.9}
\newcommand{\PrimaryExamples}{Examples 5.10}
\newcommand{\PrimaryNotPrimaryFigure}{Figure 5.3}
\newcommand{\TidyProductLemma}{Lemma 5.11}
\newcommand{\TidyColoringLemma}{Proposition 5.12}
\newcommand{\FibrousColoringLemma}{Theorem 5.13}

\newcommand{\FundamentalFigure}{Figure 5.4}

\newcommand{\CrystalProductStructureLemma}{Lemma 6.1}
\newcommand{\CrystalProductStructureLemmas}{Lemmas 6.1/2}
\newcommand{\FibrousInvariantLemma}{Lemma 6.2}
\newcommand{\MainColoringTheorem}{Theorem 6.3}
\newcommand{\CrystallineDefinition}{Definition 6.4}
\newcommand{\MainColoringCorollary}{Theorem 6.5}
\newcommand{\PrimaryProblem}{Problem 6.6}
\newcommand{\ProductProblem}{Problem 6.7}
\newcommand{\StembridgeUntangledTheorem}{Theorem 6.8}
\newcommand{\StembridgeUntangledCorollary}{Corollary 6.9}
\newcommand{\ExistenceCorollary}{Theorem 6.10}
\newcommand{\StembridgeExistenceTheorem}{Theorem 6.11}
\newcommand{\CrystallineExistenceTheorem}{Theorem 6.12}
\newcommand{\CrystallineReducibleProposition}{Proposition 6.13}
\newcommand{\EasyProductTheorem}{Theorem 6.14}
\newcommand{\CrystalProdWeightDiagramLemma}{Lemma 6.15}
\newcommand{\DataProp}{Proposition 6.16}
\newcommand{\CrystallineProblems}{Problems 6.17}
\newcommand{\CrystallineProblem}{Problem 6.17}
\newcommand{\CrystallineRemarks}{Remarks 6.18}
\newcommand{\CrystallineRemark}{Remark 6.18}

\newcommand{\StemLemma}{Lemma 7.1}
\newcommand{\AdmissibleTheorem}{Theorem 7.2}
\newcommand{\LittelmannRemarks}{Remarks 7.3}
\newcommand{\LittelmannRemark}{Remark 7.3}
\newcommand{\BigCrystallineTheorems}{Theorems 7.2 and 7.10}
\newcommand{\ConnectedObservation}{Proposition 7.4}
\newcommand{\CharacterObservation}{Proposition 7.5}
\newcommand{\TwoObservations}{Propositions 7.4/7.5}
\newcommand{\CrystalOpsLemma}{Proposition 7.6}
\newcommand{\CrystalProp}{Proposition 7.7}
\newcommand{\OpsProps}{Propositions 7.6/7}
\newcommand{\CharacterCorollary}{Corollary 7.8}
\newcommand{\CrystalMinQuasiMin}{Proposition 7.9}
\newcommand{\MainCrystalGraphTheorem}{Theorem 7.10}
\newcommand{\MainCrystalTheorems}{Theorems 7.2 and 7.10}
\newcommand{\WeylCorollary}{Corollary 7.11}
\newcommand{\CharacterAndWeylCorollaries}{Corollaries 7.8 and 7.11}
\newcommand{\SymmetricCorollary}{Corollary 7.12}
\newcommand{\CharacterInvariantCorollary}{Corollary 7.13}
\newcommand{\CrystalOpsCorollary}{Corollary 7.14}
\newcommand{\UniquenessCorollary}{Corollary 7.15}
\newcommand{\ConfirmSplittingTheorem}{Theorem 7.16}

\newcommand{\NewSplittingTheorem}{Theorem 8.1}
\newcommand{\ChainProductVersion}{Corollary 8.2}
\newcommand{\ApplicationTheorems}{Theorems 8.6 and 8.10}
\newcommand{\ChainProductVersions}{Corollaries 8.2.A-B}
\newcommand{\SpecialPatternFigure}{Figure 8.1}
\newcommand{\GelfandProposition}{Proposition 8.3}
\newcommand{\GelfandLemma}{Lemma 8.4}
\newcommand{\GelfandColoring}{Proposition 8.5}
\newcommand{\GelfandPrelimProps}{Propositions 8.3 and 8.5}
\newcommand{\GelfandTheorem}{Theorem 8.6}
\newcommand{\OddOrthPatternFigure}{Figure 8.2}
\newcommand{\SymplecticPatternFigure}{Figure 8.3}
\newcommand{\EvenOrthPatternFigure}{Figure 8.4}
\newcommand{\StructureProp}{Proposition 8.7}
\newcommand{\ColorLemma}{Lemma 8.8}
\newcommand{\ColorProp}{Proposition 8.9}
\newcommand{\PrelimProps}{Propositions 8.7 and 8.9}
\newcommand{\ApplicationTheorem}{Theorem 8.10}

\newcommand{\GaussianFromWGF}{Lemma 9.1}
\newcommand{\GaussianStructureProp}{Lemma 9.2}
\newcommand{\CrystalProductIdeal}{Proposition 9.3}
\newcommand{\WeightPreservingBijection}{Theorem 9.4}
\newcommand{\WeightPreservingBijectionCorollary}{Corollary 9.5}
\newcommand{\GaussianColoringTheorem}{Theorem 9.6}
\newcommand{\GaussianColoringCorollary}{Corollary 9.7}

\newcommand{\ProblemsList}{Problems/Questions 2.5, 2.18, 4.6, 4.11, 5.5, 6.6, 6.7, and 6.17}

\usepackage{amsfonts}
\usepackage{mathrsfs}
\usepackage{rotating}
\usepackage{mathdots}

\newfont{\myscbolditalics}{ecoc0500 at 11pt}

\newcommand{\myqSsc}{\mbox{\myscbolditalics s}}

\newfont{\mybolditalics}{ecbi0500 at 11pt}

\newcommand{\myqx}{\mbox{\mybolditalics x}}
\newcommand{\myqy}{\mbox{\mybolditalics y}}
\newcommand{\myqX}{\mbox{\mybolditalics X}}
\newcommand{\myqY}{\mbox{\mybolditalics Y}}
\newcommand{\myqh}{\mbox{\mybolditalics h}}

\newfont{\eulercursive}{eurm10 at 11pt}
\newcommand{\myd}{\mbox{\eulercursive d}}
\newcommand{\mya}{\mbox{\eulercursive a}}
\newcommand{\mys}{\mbox{\eulercursive s}}
\newcommand{\mym}{\mbox{\eulercursive m}}
\newcommand{\mye}{\mbox{\eulercursive e}}
\newcommand{\myl}{\mbox{\eulercursive `}}
\newfont{\smalleulercursive}{eurm10 at 9pt}
\newcommand{\mysmalll}{\mbox{\smalleulercursive `}}
\newfont{\smallereulercursive}{eurm10 at 7pt}
\newcommand{\mysmallerl}{\mbox{\smallereulercursive `}}

\newcommand{\mysmallindexM}{\mbox{\smalleulercursive M}}
\newcommand{\mysmallerindexM}{\mbox{\smallereulercursive M}}
\newfont{\smallesteulercursive}{eurm10 at 5pt}
\newcommand{\mysmallestindexM}{\mbox{\smallesteulercursive M}}

\newfont{\myslantcyrillic}{wncyi10 at 11pt}
\newcommand{\myspecialm}{\mbox{\myslantcyrillic m}}

\newcommand{\EuclideanE}{\mathfrak{E}}
\newcommand{\ASG}{\mathcal{G}_{\mbox{\tiny adm-syst}}}
\newcommand{\CG}{\mathcal{G}_{\mbox{\tiny crystal}}}

\newcommand{\CB}[2]{(\mathcal{#1},\mathcal{#2})}

\newcommand{\qbin}[3]{\left[{#1 \atop #2}\right]_{#3}}
\newcommand{\innprod}[2]{\langle #1,#2 \rangle}
\newcommand{\CharPi}{\widehat{\Pi}}
\newcommand{\sgn}{\mathrm{sgn}}

\newcommand{\QED}{\raisebox{0.5mm}{\fbox{\rule{0mm}{1.5mm}\ }}}

\newcounter{myfn}[page]
\renewcommand{\thefootnote}{\fnsymbol{footnote}}
\newcommand{\myfootnote}[1]{\setcounter{footnote}{\value{myfn}}%
    \footnote{#1}\stepcounter{myfn}}
    
\newcounter{rone}
\newcounter{rtwo}
\newcounter{rthree}
\newcounter{rfour}
\newcounter{rfive}
\newcounter{rsix}
\newcounter{rseven}
\newcommand{\myA}{\mbox{\sffamily A}}
\newcommand{\mysmallA}{\mbox{\footnotesize \sffamily A}}
\newcommand{\mytinyA}{\mbox{\tiny \sffamily A}}
\newcommand{\myB}{\mbox{\sffamily B}}
\newcommand{\mysmallB}{\mbox{\footnotesize \sffamily B}}

\newcommand{\myC}{\mbox{\sffamily C}}
\newcommand{\mysmallC}{\mbox{\footnotesize \sffamily C}}
\newcommand{\mytinyC}{\mbox{\tiny \sffamily C}}
\newcommand{\myD}{\mbox{\sffamily D}}
\newcommand{\mysmallD}{\mbox{\footnotesize \sffamily D}}

\newcommand{\myE}{\mbox{\sffamily E}}

\newcommand{\myF}{\mbox{\sffamily F}}

\newcommand{\mytinyF}{\mbox{\tiny \sffamily F}}
\newcommand{\myG}{\mbox{\sffamily G}}

\newcommand{\myK}{\mbox{\sffamily K}}

\newcommand{\mytinyK}{\mbox{\tiny \sffamily K}}
\newcommand{\myX}{\mbox{\sffamily X}}

\newcommand{\myT}{\mbox{\footnotesize \sffamily T}}

\newcommand{\aelt}{\mathbf{a}} \newcommand{\belt}{\mathbf{b}}
\newcommand{\celt}{\mathbf{c}} \newcommand{\delt}{\mathbf{d}}
\newcommand{\eelt}{\mathbf{e}} \newcommand{\felt}{\mathbf{f}}
\newcommand{\gelt}{\mathbf{g}} \newcommand{\helt}{\mathbf{h}}
\newcommand{\ielt}{\mathbf{i}} 
 
\newcommand{\melt}{\mathbf{m}} \newcommand{\nelt}{\mathbf{n}}
 
 \newcommand{\relt}{\mathbf{r}}
\newcommand{\selt}{\mathbf{s}} \newcommand{\telt}{\mathbf{t}}
\newcommand{\uelt}{\mathbf{u}} 
 \newcommand{\xelt}{\mathbf{x}}
\newcommand{\yelt}{\mathbf{y}} \newcommand{\zelt}{\mathbf{z}}

\newcommand{\ttelt}{(\telt_{1},\telt_{2})}
\newcommand{\xxelt}{(\xelt_{1},\xelt_{2})}
\newcommand{\yyelt}{(\yelt_{1},\yelt_{2})}

\newcommand{\wE}{\widetilde{E}}
\newcommand{\wF}{\widetilde{F}}

\newcommand{\Rcrystal}{R_{1} \otimes R_{2}}

\newcommand{\ecolor}{\mathbf{edgecolor}}

\newcommand{\comp}{\mathbf{comp}}

\newcommand{\weight}{\mbox{\sffamily weight}}
\newcommand{\smallweight}{\mbox{\footnotesize \sffamily weight}}
\newcommand{\tinyweight}{\mbox{\tiny \sffamily weight}}
\newcommand{\partition}{\mbox{\sffamily partition}}
\newcommand{\smallpartition}{\mbox{\footnotesize \sffamily partition}}

\newcommand{\WGF}{\mbox{\sffamily WGF}}
\newcommand{\RGF}{\mbox{\sffamily RGF}}

\newcommand{\depth}{\mbox{\sffamily depth}}

\newcommand{\Stemdel}{\overline{\delta}}

\newcommand{\myarrow}[1]{\stackrel{#1}{\rightarrow}}
\newcommand{\mybackarrow}[1]{\stackrel{\ #1}{\leftarrow}}
\newcommand{\mylongarrow}[1]{\stackrel{#1}{\longrightarrow}}

\newcommand{\disjointunion}{\setlength{\unitlength}{0.14cm}
\ 
\begin{picture}(2,2) 
\put(0,0){$\cup$}
\put(0.9,1.5){\circle*{0.5}}
\end{picture}\ }

\newcommand{\bigdisjointunion}{\setlength{\unitlength}{0.14cm}
\ 
\begin{picture}(2,2) 
\put(0,0){$\bigcup$}
\put(1.15,1.25){\circle*{0.8}}
\end{picture}\ }

\newcommand{\AoneAlphaWeights}{
\setlength{\unitlength}{1.1cm}
\begin{picture}(3,3)
\put(0.25,2.5){$\CG(\omega_{1})$}
\put(1.5,1){\circle*{0.125}}
\put(1.5,2){\circle*{0.125}}
\put(1.6,0.95){\footnotesize (-1,0)}
\put(1.6,1.95){\footnotesize (1,0)}
\put(1.5,1){\line(0,1){1}}
\put(1.45,1.45){\footnotesize {\em 1}}
\end{picture}
}

\newcommand{\AoneBetaWeights}{
\setlength{\unitlength}{1.1cm}
\begin{picture}(3,3)
\put(0.25,2.5){$\CG(\omega_{2})$}
\put(1.5,1){\circle*{0.125}}
\put(1.5,2){\circle*{0.125}}
\put(1.6,0.95){\footnotesize (0,-1)}
\put(1.6,1.95){\footnotesize (0,1)}
\put(1.5,1){\line(0,1){1}}
\put(1.45,1.45){\footnotesize {\em 2}}
\end{picture}
}

\newcommand{\AtwoAlphaWeights}{
\setlength{\unitlength}{1.1cm}
\begin{picture}(3,4)
\put(0.25,3.5){$\CG(\omega_{1})$}
\put(1.5,1){\circle*{0.125}}
\put(1.5,2){\circle*{0.125}}
\put(1.5,3){\circle*{0.125}}
\put(1.6,0.95){\footnotesize (0,-1)}
\put(1.6,1.95){\footnotesize (-1,1)}
\put(1.6,2.95){\footnotesize (1,0)}
\put(1.5,1){\line(0,1){2}}
\put(1.45,1.45){\footnotesize {\em 2}}
\put(1.45,2.45){\footnotesize {\em 1}}
\end{picture}
}

\newcommand{\AtwoBetaWeights}{
\setlength{\unitlength}{1.1cm}
\begin{picture}(3,4)
\put(0.25,3.5){$\CG(\omega_{2})$}
\put(1.5,1){\circle*{0.125}}
\put(1.5,2){\circle*{0.125}}
\put(1.5,3){\circle*{0.125}}
\put(1.6,0.95){\footnotesize (-1,0)}
\put(1.6,1.95){\footnotesize (1,-1)}
\put(1.6,2.95){\footnotesize (0,1)}
\put(1.5,1){\line(0,1){2}}
\put(1.45,1.45){\footnotesize {\em 1}}
\put(1.45,2.45){\footnotesize {\em 2}}
\end{picture}
}

\newcommand{\BtwoAlphaWeights}{
\setlength{\unitlength}{1.1cm}
\begin{picture}(3,6)
\put(0.25,5.5){$\CG(\omega_{2})$}
\put(1.5,1){\circle*{0.125}}
\put(1.5,2){\circle*{0.125}}
\put(1.5,3){\circle*{0.125}}
\put(1.5,4){\circle*{0.125}}
\put(1.5,5){\circle*{0.125}}
\put(1.6,0.95){\footnotesize (0,-1)}
\put(1.6,1.95){\footnotesize (-2,1)}
\put(1.6,2.95){\footnotesize (0,0)}
\put(1.6,3.95){\footnotesize (2,-1)}
\put(1.6,4.95){\footnotesize (0,1)}
\put(1.5,1){\line(0,1){4}}
\put(1.45,1.45){\footnotesize {\em 2}}
\put(1.45,2.45){\footnotesize {\em 1}}
\put(1.45,3.45){\footnotesize {\em 1}}
\put(1.45,4.45){\footnotesize {\em 2}}
\end{picture}
}

\newcommand{\BtwoBetaWeights}{
\setlength{\unitlength}{1.1cm}
\begin{picture}(3,5)
\put(0.25,4.5){$\CG(\omega_{1})$}
\put(1.5,1){\circle*{0.125}}
\put(1.5,2){\circle*{0.125}}
\put(1.5,3){\circle*{0.125}}
\put(1.5,4){\circle*{0.125}}
\put(1.6,0.95){\footnotesize (-1,0)}
\put(1.6,1.95){\footnotesize (1,-1)}
\put(1.6,2.95){\footnotesize (-1,1)}
\put(1.6,3.95){\footnotesize (1,0)}
\put(1.5,1){\line(0,1){3}}
\put(1.45,1.45){\footnotesize {\em 1}}
\put(1.45,2.45){\footnotesize {\em 2}}
\put(1.45,3.45){\footnotesize {\em 1}}
\end{picture}
}

\newcommand{\GtwoAlphaWeights}{
\setlength{\unitlength}{1.1cm}
\begin{picture}(3,8)
\put(0.25,7.5){$\CG(\omega_{1})$}
\put(1.5,1){\circle*{0.125}}
\put(1.5,2){\circle*{0.125}}
\put(1.5,3){\circle*{0.125}}
\put(1.5,4){\circle*{0.125}}
\put(1.5,5){\circle*{0.125}}
\put(1.5,6){\circle*{0.125}}
\put(1.5,7){\circle*{0.125}}
\put(1.6,0.95){\footnotesize (-1,0)}
\put(1.6,1.95){\footnotesize (1,-1)}
\put(1.6,2.95){\footnotesize (-2,1)}
\put(1.6,3.95){\footnotesize (0,0)}
\put(1.6,4.95){\footnotesize (2,-1)}
\put(1.6,5.95){\footnotesize (-1,1)}
\put(1.6,6.95){\footnotesize (1,0)}
\put(1.5,1){\line(0,1){6}}
\put(1.45,1.45){\footnotesize {\em 1}}
\put(1.45,2.45){\footnotesize {\em 2}}
\put(1.45,3.45){\footnotesize {\em 1}}
\put(1.45,4.45){\footnotesize {\em 1}}
\put(1.45,5.45){\footnotesize {\em 2}}
\put(1.45,6.45){\footnotesize {\em 1}}
\end{picture}
}

\newcommand{\GtwoBetaWeights}{
\setlength{\unitlength}{1.1cm}
\begin{picture}(5,12)
\put(1.75,11.5){$\CG(\omega_{2})$}
\put(4,6){\circle*{0.125}}
\put(2,5){\circle*{0.125}}
\put(2,7){\circle*{0.125}}
\put(3,1){\circle*{0.125}}
\put(3,2){\circle*{0.125}}
\put(3,3){\circle*{0.125}}
\put(3,4){\circle*{0.125}}
\put(2,6){\circle*{0.125}}
\put(3,8){\circle*{0.125}}
\put(3,9){\circle*{0.125}}
\put(3,10){\circle*{0.125}}
\put(3,11){\circle*{0.125}}
\put(4,5){\circle*{0.125}}
\put(4,7){\circle*{0.125}}
\put(3.15,10.95){\footnotesize (0,1)}
\put(3.15,9.95){\footnotesize (3,-1)}
\put(3.15,8.95){\footnotesize (1,0)}
\put(3.15,7.95){\footnotesize (-1,1)}
\put(4.15,7){\footnotesize (2,-1)}
\put(1.05,6.95){\footnotesize (-3,2)}
\put(1.15,5.95){\footnotesize (0,0)}
\put(4.2,5.95){\footnotesize (0,0)}
\put(1.05,4.95){\footnotesize (3,-2)}
\put(4.15,4.9){\footnotesize (-2,1)}
\put(3.15,3.95){\footnotesize (1,-1)}
\put(3.1,2.95){\footnotesize (-1,0)}
\put(3.1,1.95){\footnotesize (-3,1)}
\put(3.1,0.95){\footnotesize (0,-1)}
\put(3,1){\line(0,1){3}}
\put(3,8){\line(0,1){3}}
\put(3,4){\line(-1,1){1}}
\put(3,4){\line(1,1){1}}
\put(2,7){\line(1,1){1}}
\put(4,7){\line(-1,1){1}}
\put(2,5){\line(0,1){2}}
\put(4,5){\line(0,1){2}}
\put(2.95,10.45){\footnotesize {\em 2}}
\put(2.95,9.45){\footnotesize {\em 1}}
\put(2.95,8.45){\footnotesize {\em 1}}
\put(2.45,7.45){\footnotesize {\em 1}}
\put(3.45,7.45){\footnotesize {\em 2}}
\put(2.45,4.45){\footnotesize {\em 1}}
\put(3.45,4.45){\footnotesize {\em 2}}
\put(2.95,3.45){\footnotesize {\em 1}}
\put(2.95,2.45){\footnotesize {\em 1}}
\put(2.95,1.45){\footnotesize {\em 2}}
\put(1.95,5.45){\footnotesize {\em 2}}
\put(1.95,6.45){\footnotesize {\em 2}}
\put(3.95,5.45){\footnotesize {\em 1}}
\put(3.95,6.45){\footnotesize {\em 1}}
\end{picture}
}

\newcommand{\NEEdgeLabelForLatticeI}[1]{
\setlength{\unitlength}{1.5cm}
\begin{picture}(0,0)
\put(-0.25,0){
\begin{picture}(0,0)
\put(0.4,0.4){\footnotesize #1} 
\end{picture}
}
\end{picture}
}

\newcommand{\NWEdgeLabelForLatticeI}[1]{
\setlength{\unitlength}{1.5cm}
\begin{picture}(0,0)
\put(-0.25,0){
\begin{picture}(0,0)
\put(-0.525,0.4){\footnotesize #1} 
\end{picture}
}
\end{picture}
}

\newcommand{\VerticalEdgeLabelForLatticeI}[1]{
\setlength{\unitlength}{1.5cm}
\begin{picture}(0,0)
\put(-0.25,0){
\begin{picture}(0,0)
\put(-0.05,0.4){\footnotesize #1} 
\end{picture}
}
\end{picture}
}

\newcommand{\VertexTableau}[7]{
\setlength{\unitlength}{1.5cm}
\begin{picture}(0,0)
\put(-0.25,0){
\begin{picture}(0,0)
\put(0,0){\circle*{0.1}} 
\put(#6,#7){\setlength{\unitlength}{0.3cm}\begin{picture}(0,0)\put(0,0){\line(0,1){2}} \put(1,0){\line(0,1){2}} \put(2,0){\line(0,1){2}} \put(3,1){\line(0,1){1}} \put(0,0){\line(1,0){2}} \put(0,1){\line(1,0){3}} \put(0,2){\line(1,0){3}} \put(0.25,1.25){\tiny #1} \put(1.25,1.25){\tiny #2} \put(2.25,1.25){\tiny #3} \put(0.25,0.25){\tiny #4} \put(1.25,0.25){\tiny #5}\end{picture}}
\end{picture}
}
\end{picture}
}

\newcommand{\zup}[3]{z_{1}^{#1}z_{2}^{#2}z_{3}^{#3}}

\begin{document}
\pagenumbering{arabic}
\thispagestyle{empty}%
\vspace*{-0.7in}
\hfill \parbox{2in}{\hfill {\scriptsize This version: September 6, 2021} 

\vspace*{-0.1in}
\hfill{\scriptsize 1st version: November 26, 2018}}

\begin{center}
{\large \bf Poset models for Weyl group analogs of\\
symmetric functions and Schur functions} 

\vspace*{0.05in}
\renewcommand{\thefootnote}{1}
Robert G.\ Donnelly\footnote{Email: 
{\tt rob.donnelly@murraystate.edu},\ Fax: 
1-270-809-2314}

\vspace*{-0.05in}
Department of Mathematics and Statistics, Murray State
University, Murray, KY 42071
\end{center} 

\begin{abstract}
The ``Weyl symmetric functions'' studied here naturally generalize classical symmetric (polynomial) functions, and ``Weyl bialternants,'' sometimes also called Weyl characters, analogize the Schur functions. 
For this generalization, the underlying symmetry group is a finite Weyl group. 
A ``splitting poset'' for a Weyl bialternant is an edge-colored ranked poset possessing a certain structural property and a natural weighting of its elements so that the weighted sum of poset elements is the given Weyl bialternant.  
Connected such posets are of combinatorial interest in part because they are rank symmetric and rank unimodal and have nice quotient-of-product expressions for their rank 
generating functions.   
Supporting graphs of weight bases for irreducible semisimple Lie algebra representations provide one large family of examples. 
However, many splitting posets can be obtained outside of this Lie theoretic context. 

This monograph provides a tutorial on Weyl bialternants$\,$/$\,$Weyl symmetric functions and splitting posets that is largely self-contained and independent of Lie algebra representation theory.  
New results are also obtained. 
In particular, a cancelling argument of Stembridge is reworked to provide sufficient combinatorial conditions for a given poset to be splitting. 
This new splitting theorem is used to help construct what are here named crystalline splitting posets. 
The Weyl bialternants with unique splitting posets are classified, and a new combinatorial characterization of splitting posets associated with the minuscule and quasi-minuscule dominant weights is given. 
These findings are supported by other new results of a poset-structural nature. 

In addition, some new conceptual approaches are presented. 
The notion of a ``refined'' splitting poset is introduced to address not only splitting-type problems but also Littlewood--Richardson-type (i.e.\ product decomposition) and branching-type problems. 
``Crystalline splitting posets'' are introduced as a class of posets that behave much like crystal graphs from the crystal base theory of Kashiwara {\em et al} with regard to branching and decomposing products. 
A new technique called ``vertex coloring'' is used in conjunction with the new splitting theorem to produce crystalline splitting posets and, in particular, all crystal graphs.  
Via this vertex-coloring method, Stembridge's admissible systems are shown to be crystalline splitting posets. 
Crystal basis / crystal graph theory is used to demonstrate uniqueness of those crystalline splitting posets that are obtained via a certain product construction. 

\hfill {\em (continued)}

\
\vspace*{0.05in}

\begin{center}

\end{center} 
\end{abstract}

\newpage
\def\abstractname{Abstract, continued}
\begin{abstract}
This vertex-coloring technique is applied here to obtain other explicit and combinatorially interesting families of splitting posets. 
In particular, an alternate set of combinatorial criteria is developed to facilitate the study of certain ``non-crystalline'' splitting posets. 
Using these alternate criteria, some distributive lattices of Gelfand-type patterns are shown to model certain symplectic and orthogonal Weyl bialternants. 
A new proof of the Bender--Knuth (ex-)conjecture is obtained as a consequence, as well as a new symplectic analog of this enumerative identity. 
We use our methodology to produce two new proofs that minuscule posets are Gaussian.  
Further applications of this vertex-coloring technique will appear in future papers. 
\end{abstract}

\newpage
\def\abstractname{Table of contents}
\begin{abstract}
\begin{center}
\parbox{5in}{
\hspace*{-0.25in}\IntroSection. Introduction\dotfill 5\\%
\hspace*{0in}Purposes and highlights\dotfill 7\\%
\hspace*{0in}Three problems and a central formula\dotfill 10\\%
\hspace*{0in}On the current state of knowledge, some open questions,\\ \hspace*{2in} and where we hope to go from here\dotfill 11\\%
\hspace*{0in}Acknowledgements\dotfill 12\\%
\hspace*{0in}A table of case-by-case results\dotfill 13\\%
\hspace*{-0.25in}\WeylSection. Weyl symmetric functions and Weyl bialternants\dotfill 16\\%
\hspace*{0in}Roots and weights\dotfill 16\\%
\hspace*{0in}Some key quantities\dotfill 17\\%
\hspace*{0in}The Weyl group\dotfill 17\\%
\hspace*{0in}Weight diagrams\dotfill 18\\%
\hspace*{0in}Reducible root systems; root subsystems\dotfill 19\\%
\hspace*{0in}The ring of Weyl symmetric functions and some other important rings\dotfill 20\\%
\hspace*{0in}The subgroup of Weyl alternants\dotfill 20\\%
\hspace*{0in}Weyl bialternants and $\Phi$-Kostka numbers\dotfill 22\\%
\hspace*{0in}Other (potential) bases for the ring of $W$-symmetric functions\dotfill 23\\%
\hspace*{0in}Weyl bialternants, monomial $W$-symmetric functions,\\%
\hspace*{0.5in}and elementary $W$-symmetric functions for reducible root systems\dotfill 24\\%
\hspace*{0in}Finiteness of Weyl bialternants, and some consequences\dotfill 24\\%
\hspace*{0in}A recurrence for the $\Phi$-Kostka numbers\dotfill 26\\%
\hspace*{0in}Basis results\dotfill 28\\%
\hspace*{0in}Two involutions\dotfill 29\\%
\hspace*{0in}Specializations of Weyl bialternants\dotfill 30\\%
\hspace*{0in}The $\myA_n$ case\dotfill 31\\%
\hspace*{0in}Interactions with Lie algebra representation theory\dotfill 34\\%
\hspace*{-0.25in}\PosetSection. Edge-colored posets and 
${\Phi}$-structured posets\dotfill 37\\%
\hspace*{0in}Some language and notation\dotfill 37\\%
\hspace*{0in}Fibrous posets and a connection with symmetric chain decompositions\dotfill 39\\%
\hspace*{0in}Posets interacting with roots, weights, and the Weyl group\dotfill 40\\%
\hspace*{0in}Weight diagrams, again\dotfill 43\\%
\hspace*{0in}A generalized notion of weight diagram\dotfill 43\\%
\hspace*{-0.25in}\SplittingSection. Splitting posets\dotfill 50\\%
\hspace*{0in}A new splitting theorem\dotfill 50\\%
\hspace*{0in}Splitting posets and refined splitting posets\dotfill 51\\%
\hspace*{0in}A version of Stanley's problem for Weyl symmetric functions\dotfill 54\\%
\hspace*{0in}Some combinatorial properties of splitting posets\dotfill 54\\%
\hspace*{0in}Edge-minimal splitting posets\dotfill 55\\%
\hspace*{0in}The unique maximal splitting poset\dotfill 57\\%
\hspace*{0in}Minuscule and quasi-minuscule splitting posets\dotfill 57\\%
\hspace*{0in}Supporting graphs as splitting posets\dotfill 65\\%
\hspace*{0in}Weyl bialternants with unique splitting posets\dotfill 67}%
\end{center}
\end{abstract}

\newpage
\def\abstractname{Contents, continued}
\begin{abstract}
\begin{center}
\parbox{5in}{
\hspace*{-0.25in}\FibrousSection. Crystal products of fibrous posets and vertex coloring\dotfill 69\\%
\hspace*{0in}Crystal products of fibrous posets\dotfill 69\\
\hspace*{0in}Toward a more refined view of fibrous posets\dotfill 74\\ 
\hspace*{0in}Primary posets\dotfill 74\\
\hspace*{0in}Vertex coloring\dotfill 75\\%
\hspace*{-0.25in}\CrystallineSection. Crystalline splitting posets\dotfill 79\\%
\hspace*{0in}${\Phi}$-structured fibrous posets\dotfill 79\\%
\hspace*{0in}Crystalline splitting posets\dotfill 80\\%
\hspace*{0in}Constructing crystalline splitting posets\dotfill 81\\%
\hspace*{0in}When do product decomposition coefficients coincide\\ \hspace*{2in} with $\Phi$-Kostka numbers?\dotfill 83\\%
\hspace*{0in}Further consequences, questions, and comments\dotfill 86\\%
\hspace*{-0.25in}\CrystalSection. Stembridge's admissible systems and the crystal bases/graphs\\ \hspace*{3in}of Kashiwara {\em et al}\dotfill 88\\%
\hspace*{0in}Admissible systems as crystalline splitting posets\dotfill 88\\%
\hspace*{0in}Crystal graphs as crystalline splitting posets\dotfill 91\\%
\hspace*{-0.25in}\CriteriaSection. Splitting via vertex-coloring for non-fibrous posets\dotfill 100\\%
\hspace*{0in}New splitting results for certain non-fibrous posets\dotfill 100\\%
\hspace*{0in}Gelfand--Tsetlin lattices as splitting distributive lattices\dotfill 102\\%
\hspace*{0in}Splitting distributive lattices for certain symplectic\\ 
\hspace*{2.25in}and orthogonal Weyl bialternants\dotfill 106\\
\hspace*{-0.25in}\GaussianSection. Two vertex-coloring proofs that minuscule posets are Gaussian\dotfill 114\\
\hspace*{0in}Gaussian posets\dotfill 114\\ 
\hspace*{0in}Two new approaches\dotfill 114\\ 
\hspace*{0in}Some further aspects of our set-up\dotfill 115\\ 
\hspace*{0in}First new proof\dotfill 116\\ 
\hspace*{0in}Second new proof\dotfill 118\\ 
\hspace*{-0.25in}References\dotfill 120
}
\end{center}
\end{abstract}

\newpage 
\noindent 
{\Large \bf \S \IntroSection.\ Introduction.} 

Combinatorial representation theory is currently a flourishing area of research.  
Broadly speaking, the goal of this area is to advance understanding of algebraic structures and their representations using combinatorial methods, and vice-versa.   
For an excellent survey, see \cite{BR}. 
Our focus here is on one particular corner of this area: a poset theoretic study of `Weyl-group-generalized' symmetric functions and, at times, of their related semisimple Lie algebra representations.  
The purpose of this monograph is threefold. 
We recapitulate some foundational aspects of the theory of Weyl symmetric functions; we present what we believe are some new order-theoretic approaches to this Weyl group symmetric function theory; and we hope to inspire interest in these ideas by showcasing some of the beautiful objects this study has produced and by pointing out many open problems.

Now, `classical' symmetric function theory is the study of multivariate polynomials that are invariant under a certain natural action of the finite symmetric groups. 
This rich subject lies in the intersection of algebra and combinatorics, as it connects representations of the symmetric groups and general linear Lie algebras to various enumerative phenomena. 
The Weyl group generalization of classical symmetric functions opens many more vistas for observing such algebraic and combinatorial interactions, some of which we intend to explore here. 

The systematic use of posets in studying Weyl symmetric functions (and their related semisimple Lie algebra representations) originated with certain work of Richard P.\ Stanley and Robert A.\ Proctor in the late 1970's and early 1980's.  
In papers such as \cite{StanUnim}, it was clear that Stanley was aware of nice interactions between certain families of posets and Weyl characters.  
Proctor introduced the idea of semisimple Lie algebras acting on posets in papers such as \cite{PrPeck}, \cite{PrMonthly}, \cite{PrEur}, and \cite{PrGZ}.  
Since that time, there has been interest in finding combinatorial models for Weyl symmetric functions and in constructing representations using combinatorial 
methods.  
These have been topics of interest for this author (\cite{DonThesis}, \cite{DonPeck}, \cite{DonSymp}, \cite{DonSupp}, \cite{DonAdjoint}) as well as many other researchers (\cite{Alverson}, \cite{ADLPOne}, \cite{ADLMPPW}, \cite{DD}, \cite{DDgeneral}, \cite{DDW}, \cite{DLP1}, \cite{DLP2}, \cite{HL}, \cite{KN}, \cite{LS}, \cite{LP}, \cite{LitTableaux}, \cite{McClard}, \cite{Stem}, \cite{WildbergerAdv}, \cite{WildbergerEur}, \cite{WildbergerLie}, etc). 

Before we take a brief tour of the main ideas of this monograph, we indicate the key objects of our discourse.  
For formal definitions, see \S \WeylSection.  
The {\sl Weyl groups} of this monograph are viewed as symmetry groups of {\sl finite root systems}. 
Such a root system resides in a Euclidean space $\EuclideanE$, and the associated Weyl group acts on a discrete subset of $\EuclideanE$ called the {\sl weight lattice}; so-called {\sl dominant weights} are elements of the weight lattice that can be identified as nonnegative-integer tuples. 
Schur functions are certainly the most important and most interesting classical symmetric functions. 
These comprise a basis of the ring/$\mathbb{Z}$-module of symmetric functions, are directly connected to the representation theory of the general linear Lie algebras, and have various concrete characterizations as weighted sums of tableaux and as quotients of determinants. 
The {\sl Weyl bialternants} of this monograph, sometimes called Weyl characters, are the Weyl group symmetric function analogs of Schur functions. 
A Weyl bialternant can be viewed as a multivariate Laurent polynomial that is a quotient of determinant-like ``skew-invariant'' alternating sums. 
Each Weyl bialternant is uniquely associated with an irreducible representation of the associated semisimple complex Lie algebra.  
As Schur functions are naturally indexed by integer partitions, Weyl bialternants are similarly indexed by dominant weights. 

A fundamental problem in combinatorial representation theory is to find combinatorial objects whose weighted sum yields the Weyl bialternant corresponding to a given finite root system and dominant weight. 
Equivalently, one can seek objects that split the weight multiplicities for the associated irreducible representation of a semisimple Lie algebra.     
In \cite{DonSupp} \S 3 and again in \cite{ADLMPPW} \S 2, the concept of a ``splitting poset'' was introduced.  
This is an edge-colored poset with a certain structural property and possessing a natural weight function on its elements such that the weighted sum of the poset elements is a Weyl bialternant.  
As such, a splitting poset is a solution to the proposed combinatorial splitting problem.  

\begin{figure}[t]
\begin{center}
\IntroFig:  The edge-colored distributive lattice below is the Gelfand--Tsetlin lattice $L({\setlength{\unitlength}{0.125cm}\begin{picture}(3,0)\put(0,0){\line(0,1){2}} \put(1,0){\line(0,1){2}} \put(2,0){\line(0,1){2}} \put(3,1){\line(0,1){1}} \put(0,0){\line(1,0){2}} \put(0,1){\line(1,0){3}} \put(0,2){\line(1,0){3}}\end{picture}},3)$.  

\

\parbox{5.75in}{\footnotesize Following \CaseAExample, we have $n = 2$, and the partition $\mathfrak{p}$ corresponding to the shape {\setlength{\unitlength}{0.125cm}\begin{picture}(3,0)\put(0,0){\line(0,1){2}} \put(1,0){\line(0,1){2}} \put(2,0){\line(0,1){2}} \put(3,1){\line(0,1){1}} \put(0,0){\line(1,0){2}} \put(0,1){\line(1,0){3}} \put(0,2){\line(1,0){3}}\end{picture}} is $(3,2,0)$. 
Then $\displaystyle \mys_{\mathfrak{p}} = \sum_{\telt \in L}z^{wt_{GT}(\telt)}$, where $wt_{GT}(\telt)$ is the triple $(\mbox{\# of 1's},\mbox{\# of 2's},\mbox{\# of 3's})$ and $z^{wt_{GT}(\telt)}$ means $\zup{\mbox{\footnotesize \# of 1's}}{\mbox{\footnotesize \# of 2's}}{\mbox{\footnotesize \# of 3's}}$. 
Replace $z_1$ by $x_1$, $z_2$ by $x_{1}^{-1}x_{2}$, and $z_3$ by $x_{2}^{-1}$, and let $wt(\telt)$ be the pair $(\mbox{\# of 1's}-\mbox{\# of 2's},\mbox{\# of 2's}-\mbox{\# of 3's})$. 
Note that $wt(\telt)$ is the vector of integers whose $i$th coordinate is twice the rank of $\telt$ within its ``$i$-component'' less the length of this $i$-component, i.e.\ $2\rho_{i}(\telt) - l_{i}(\telt)$.  That is, a coordinate-free expression for $wt(\telt)$  is $(2\rho_{1}(\telt) - l_{1}(\telt), 2\rho_{2}(\telt) - l_{2}(\telt))$, cf.\ \S \PosetSection. 
Then the Weyl bialternant $\displaystyle \chi_{_{\omega_{1}+2\omega_{2}}}^{\mytinyA_2}$ can be expressed combinatorially as $\displaystyle \sum_{\telt \in L}x^{wt(\telt)}$.}

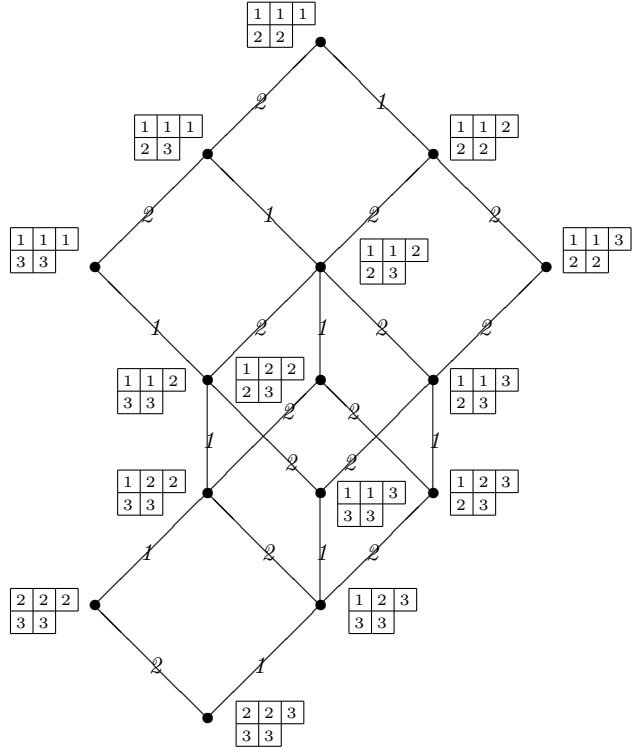
\setlength{\unitlength}{1.5cm}
\begin{picture}(4,6.5)
\put(1,0){\line(-1,1){1}}
\put(1,0){\line(1,1){2}}
\put(0,1){\line(1,1){2}}
\put(2,1){\line(-1,1){1}}
\put(2,1){\line(0,1){1}}
\put(1,2){\line(0,1){1}}
\put(2,2){\line(-1,1){2}}
\put(2,2){\line(1,1){2}}
\put(3,2){\line(-1,1){1}}
\put(3,2){\line(0,1){1}}
\put(1,3){\line(1,1){2}}
\put(2,3){\line(0,1){1}}
\put(3,3){\line(-1,1){2}}
\put(0,4){\line(1,1){2}}
\put(4,4){\line(-1,1){2}}
\put(2,6){\VertexTableau{1}{1}{1}{2}{2}{-0.65}{-0.05}}
\put(1,5){\VertexTableau{1}{1}{1}{2}{3}{-0.65}{-0.05}}
\put(3,5){\VertexTableau{1}{1}{2}{2}{2}{0.15}{-0.05}}
\put(0,4){\VertexTableau{1}{1}{1}{3}{3}{-0.75}{-0.05}}
\put(2,4){\VertexTableau{1}{1}{2}{2}{3}{0.35}{-0.15}}
\put(4,4){\VertexTableau{1}{1}{3}{2}{2}{0.15}{-0.05}}
\put(1,3){\VertexTableau{1}{1}{2}{3}{3}{-0.8}{-0.3}}
\put(2,3){\VertexTableau{1}{2}{2}{2}{3}{-0.75}{-0.2}}
\put(3,3){\VertexTableau{1}{1}{3}{2}{3}{0.15}{-0.3}}
\put(1,2){\VertexTableau{1}{2}{2}{3}{3}{-0.8}{-0.2}}
\put(2,2){\VertexTableau{1}{1}{3}{3}{3}{0.15}{-0.3}}
\put(3,2){\VertexTableau{1}{2}{3}{2}{3}{0.15}{-0.2}}
\put(0,1){\VertexTableau{2}{2}{2}{3}{3}{-0.75}{-0.25}}
\put(2,1){\VertexTableau{1}{2}{3}{3}{3}{0.25}{-0.25}}
\put(1,0){\VertexTableau{2}{2}{3}{3}{3}{0.25}{-0.25}}
\put(1,5){\NEEdgeLabelForLatticeI{{\em 2}}}
\put(3,5){\NWEdgeLabelForLatticeI{{\em 1}}}
\put(0,4){\NEEdgeLabelForLatticeI{{\em 2}}}
\put(2,4){\NWEdgeLabelForLatticeI{{\em 1}}}
\put(2,4){\NEEdgeLabelForLatticeI{{\em 2}}}
\put(4,4){\NWEdgeLabelForLatticeI{{\em 2}}}
\put(1,3){\NEEdgeLabelForLatticeI{{\em 2}}}
\put(1,3){\NWEdgeLabelForLatticeI{{\em 1}}}
\put(2,3){\VerticalEdgeLabelForLatticeI{{\em 1}}}
\put(3,3){\NWEdgeLabelForLatticeI{{\em 2}}}
\put(3,3){\NEEdgeLabelForLatticeI{{\em 2}}}
\put(1,2){\VerticalEdgeLabelForLatticeI{{\em 1}}}
\put(1.25,2.25){\NEEdgeLabelForLatticeI{{\em 2}}}
\put(2.2,1.8){\NWEdgeLabelForLatticeI{{\em 2}}}
\put(1.8,1.8){\NEEdgeLabelForLatticeI{{\em 2}}}
\put(3,2){\VerticalEdgeLabelForLatticeI{{\em 1}}}
\put(2.75,2.25){\NWEdgeLabelForLatticeI{{\em 2}}}
\put(0,1){\NEEdgeLabelForLatticeI{{\em 1}}}
\put(2,1){\VerticalEdgeLabelForLatticeI{{\em 1}}}
\put(2,1){\NWEdgeLabelForLatticeI{{\em 2}}}
\put(2,1){\NEEdgeLabelForLatticeI{{\em 2}}}
\put(1,0){\NWEdgeLabelForLatticeI{{\em 2}}}
\put(1,0){\NEEdgeLabelForLatticeI{{\em 1}}}
\end{picture}
\end{center}
\end{figure}

Let us describe a splitting poset -- a splitting distributive lattice actually -- for each Schur function. 
In particular, when $\mathfrak{p} = (p_{1} \geq p_{2} \geq \cdots \geq p_{n+1} \geq 0)$ is an integer partition and $\displaystyle \rule[-0.225in]{0mm}{0.55in} \mys_{\mathfrak{p}} := \frac{\det(z_{i}^{p_{j}+n+1-j})_{i,j=1}^{n+1}}{\det(z_{i}^{n+1-j})_{i,j=1}^{n+1}}$ is the corresponding Schur function in the $n+1$ variables $z_1,\ldots,z_{n+1}$, then one well-known combinatorial description of $\mys_{\mathfrak{p}}$ is as a weighted sum of all semistandard (i.e.\ column-strict) tableaux of shape $\mathfrak{p}$ and with entries from the set $\{1,2,\ldots,n+1\}$, as in \IntroFig: $\displaystyle \mys_{\mathfrak{p}} = \sum_{\parbox{1.75cm}{\tiny \hspace*{0.3cm}$\mathfrak{p}$-shaped\\ semistandard $T$}} z_1^{\mbox{\tiny \# of 1's in $T$}}z_2^{\mbox{\tiny \# of 2's in $T$}}{\cdots}z_{n+1}^{\mbox{\tiny \# of $n+$1's in $T$}}$, for indeterminates $z_{1},z_{2},\ldots,z_{n+1}$. 
A simple change of variables is needed in order to view $\mys_{\mathfrak{p}}$ as a ``type $\myA_{n}$'' Weyl bialternant, for details see \CaseAExample\ below. 
One splitting poset for this change-of-variables version of $\mys_{\mathfrak{p}}$ is well-known: It is the Gelfand--Tsetlin lattice of all semistandard tableaux of shape $\mathfrak{p}$ and with entries from the set $\{1,2,\ldots,n+1\}$ partially ordered by ``reverse componentwise comparison,'' i.e.\ $S \leq T$ for such tableaux $S$ and $T$ if and only if $S_{pq} \geq T_{pq}$ at all positions $(p,q)$ in the shape. 
This distributive lattice is sometimes denoted $L(\mathfrak{p},n+1)$. 
To form the weighted sum over $L(\mathfrak{p},n+1)$, there is a simple rule for weighting tableaux. 
This weighting of elements can be achieved without reference to the coordinates of the tableaux, using only an assignment of colors to the edges of the Hasse diagram for the lattice. 
(This ``coordinate-free'' approach to weighting the elements of an edge-colored poset is developed in general in \S \PosetSection.) 
The appropriate edge-coloring rule for $L(\mathfrak{p},n+1)$ is: $S \myarrow{i} T$ if for some pair $(r,s)$ we have $S_{rs}-1 = T_{rs} = i$ with $S_{pq} = T_{pq}$ for all pairs $(p,q) \not= (r,s)$. 
See \IntroFig\ for an example and \cite{StanUnim}, \cite{PrGZ}, \cite{DonSupp}, or \cite{HL} for more details. 

The problem of finding splitting posets generalizes a problem posed by Stanley in \cite{StanUnim} (Problem 3) which was at least partly inspired by the above Schur function/Gelfand--Tsetlin lattice example.  
To paraphrase, Stanley's problem is to find ranked posets -- and particularly distributive lattices -- whose rank generating functions are specializations of Weyl bialternants or, more generally, Weyl symmetric functions.  
Splitting posets provide one answer to this problem.  
One reason such posets are of combinatorial interest is that the rank generating functions for connected splitting posets of Weyl bialternants have beautiful product-over-product expressions and, applying observations originally due to Dynkin \cite{Dynkin}, are symmetric and unimodal.  
(A degree $l$ polynomial $a_{0}+a_{1}q+\cdots+a_{l}q^{l}$ in the variable $q$ is symmetric if $a_{i} = a_{l-i}$ for all $0 \leq i \leq l$ and is unimodal if for some $p$ ($0 \leq p \leq l$) it is the case that $a_{0} \leq \cdots \leq a_{p} \geq \cdots \geq a_{l}$.)  
Some well-known enumerative phenomena -- for example, unimodality of Gaussian polynomials and certain $q$-Catalan numbers, the Bender--Knuth formula for counting certain plane partitions, Gaussian-ness of minuscule posets, a combinatorial description of the Littlewood--Richardson coefficients, and many generalizations of these results -- can be obtained as applications of splitting posets. 

{\bf Purposes and highlights.} 
This monograph is a sequel of sorts to \cite{ADLMPPW} and is intended to give a readable, even browsable, introduction to the basic theory of Weyl symmetric functions$\,$/$\,$Weyl bialternants and splitting posets, 
to develop further connections with the literature, and to present new techniques and results. 
While semisimple Lie algebra representation theory  provides some of our motivation and context, knowledge of this theory is not needed to understand most of the results of this largely self-contained monograph. 
Indeed, one of our purposes is to develop as much of the story as possible without Lie algebra representation theory, in the spirit of Occam's razor but also with the idea of creating less overhead for the audience.  
The perspective is mainly combinatorial, and some combinatorial problems are noted in the discussion.  

The table of contents gives a detailed overview of the organization of the monograph. 
Here, we point out some of the highlights and give some qualitative remarks. 

In \S \WeylSection, we develop the basic theory of Weyl symmetric functions and Weyl bialternants from the foundation of finite root systems.  
There does not appear to be a similarly self-contained exposition in the literature, so in addition to providing the setting for our subsequent work, this discussion might benefit combinatorialists or others interested in this very natural generalization of symmetric function theory. 

In \S \PosetSection, we provide the basic poset theory background. 
In addition to introducing some language and notation, another goal of \S \PosetSection\ is to provide some guidance on how the structure of the posets studied here is influenced by interactions with root systems, weights, and Weyl groups. 
Many of these results are mainly auxiliary.  
However, a new theorem that characterizes a generalization of a weight diagram (see \PiTheorem) might be of independent interest, besides its uses here. 

In  \S \SplittingSection, splitting posets are formally defined and some of their basic properties are established.  
A new result provides combinatorial criteria for demonstrating that a poset is splitting (see \InitialSplittingTheorem). 
The key combinatorial data are a vertex-coloring function and a compatible bijection of certain poset elements. 
This theorem is applied three times in this monograph: First, in characterizing certain splitting posets related to minuscule and quasi-minuscule dominant weights (see \QuasiMinLemma/\MinQuasiMinTheorem); second in constructing ``crystalline'' splitting posets, which are, in a certain sense, generalizations of Stembridge's admissible systems and also of the crystal graphs of Kashiwara {\em et al} (see \MainColoringTheorem\ and \MainColoringCorollary; see also \AdmissibleTheorem\ and \MainCrystalGraphTheorem); and third, to help develop alternate combinatorial criteria for splitting (see \ChainProductVersion). 
Many examples of splitting posets are provided in \S \SplittingSection. 
These include splitting posets with extremal properties: Up to a notion of isomorphism, there are only finitely many splitting posets for a given Weyl bialternant, and among these are the unique maximal splitting poset as well as certain edge-minimal splitting posets. 
We give a new combinatorial characterization of certain splitting posets associated with the minuscule and quasi-minuscule dominant weights (\MinQuasiMinTheorem). 
We observe again (as in \S 3 of \cite{DonSupp}) that supporting graphs of weight bases for semisimple Lie algebra representations are splitting posets. 
We also determine which Weyl bialternants have unique splitting posets. 

In \S \FibrousSection, the idea of a ``crystal product'' of certain kinds of posets (``fibrous'' posets) is studied as a purely combinatorial notion. 
A vertex-coloring technique is developed that is compatible with crystal products and will allow us later to build crystalline splitting posets, cf.\ \MainColoringTheorem/\MainColoringCorollary. 

In \S \CrystallineSection, a crystalline splitting poset is constructed for each Weyl bialternant. 
This is a consequence of a vertex-coloring result (\FibrousColoringLemma), some splitting results for certain fibrous posets (\MainColoringTheorem/\MainColoringCorollary), and a result of Stembridge from \cite{Stem} that allows one to use minuscule and quasi-minuscule splitting posets as the starting point for the construction. 
Ultimately, the result is iterative in the sense that the final result is realized as the outcome of a product construction that is only carried out in principle. 
Crystalline splitting posets resolve, at least in some sense, three problems for Weyl bialternants, which we name the ``splitting'' problem, the ``product decomposition'' problem, and the ``branching'' problem. 
These three problems are discussed in more detail below.  

In \BigCrystallineTheorems\ of \S \CrystalSection, we demonstrate that Stembridge's admissible systems (see \cite{Stem}) and the crystal graphs of Kashiwara {\em et al} (see for example \cite{Kash1}, \cite{Kash2}) are instances of the crystalline splitting posets of \S \CrystallineSection. 
It follows that other constructions of crystal graphs (via Littelmann's path model \cite{Lit2}, \cite{LitPaths}, \cite{Lit4} or Lenart--Postnikov's alcove path model \cite{LP}) also yield crystalline splitting posets.  
By identifying crystal graphs as crystalline splitting posets, we also derive several useful combinatorial conclusions. 
Notably, we conclude in \UniquenessCorollary\ that the crystalline splitting poset obtained by the product construction of \S \CrystallineSection\ is independent in a certain sense of the minuscule and quasi-minuscule splitting posets used as factors in the crystal product. 

The main goal of \S \CriteriaSection\ is to develop in \NewSplittingTheorem/\ChainProductVersion.A-B a set of sufficient combinatorial conditions for certain ``non-fibrous'' posets to be splitting. 
The proof of \NewSplittingTheorem\ applies \InitialSplittingTheorem.  
For a poset that can serve as a supporting graph, it is useful to establish independently that the poset is splitting, and in certain cases \NewSplittingTheorem/\ChainProductVersion.A-B will allow us to do so (for example, \cite{DD}, \cite{DDW}). 
In general, for posets that cannot serve as or are not yet known to serve as supporting graphs, \NewSplittingTheorem/\ChainProductVersion.A-B afford us methods for obtaining splitting results (for example, some of the families of splitting distributive lattices obtained in \cite{ADLPTwo}, \cite{DDgeneral}). 

We feature some applications of our methodology in Sections \CriteriaSection\ and \GaussianSection.  
In particular, \ChainProductVersion.B is used in \ApplicationTheorems\ to show that some families of distributive lattices built from some natural Gelfand-type patterns are splitting distributive lattices for certain special linear, symplectic, and orthogonal Weyl bialternants. 
This yields another proof and some new analogs of the Bender--Knuth identity for enumerating a certain family of plane partitions. 
\ChainProductVersion.B is used in \cite{ADLPTwo} to more easily re-derive the main splitting result of \cite{ADLMPPW}. 
In \S \GaussianSection, we apply methods developed in previous sections to produce two new proofs of the Proctor--Stanley Theorem that minuscule posets are Gaussian.

Overall, we believe that the methodology developed in this monograph for building crystalline splitting posets is more straightforward combinatorially than other approaches mentioned in this introduction, at least from the viewpoint of obtaining, from scratch, a collection of crystal-like graphs that resolve our three main problems of interest for Weyl bialternants (splitting, product decomposition, branching -- these are further discussed immediately below). 
Another advantage is that we are able here and elsewhere to generalize our approach to non-fibrous posets.  
So in particular, in the spirit of Stanley's problem, we can use our approach to obtain some families of splitting distributive lattices and splitting modular lattices. 

{\bf Three problems and a central formula.} 
Three principal problems in the study of Weyl bialternants are the {\em splitting}, {\em product decomposition}, and {\em branching} problems, each of which yields to combinatorial approaches and has combinatorial consequences. 
For the discussion that follows, we need some notation that is more fully developed later. 
Beginning in \S \WeylSection, a Weyl bialternant is denoted by $\chi_{_{\lambda}}$, where $\lambda$ is a dominant weight in the weight lattice associated with a finite root system $\Phi$. 
We sometimes also write $\chi_{_{\lambda}}^{\Phi}$ to emphasize the role of the root system. 
Then $\chi_{_{\nu}}^{\Phi_J}$ denotes a ``$\Phi_{J}$-Weyl bialternant'' associated to a root subsystem $\Phi_{J} \subseteq \Phi$, where $J$ is a subset of the set $I$ that indexes the simple roots chosen for the originating root system $\Phi$ and where $\nu$ is ``$\Phi_J$-dominant.'' 

The {\em splitting problem} is the problem of explicitly realizing a Weyl bialternant or other Weyl symmetric function as a weight generating function on some combinatorially interesting set of objects.  
In this monograph, as elsewhere, we propose using certain ranked posets whose Hasse diagram edges are colored by $I$. 
One of the main goals of this monograph is to further develop this splitting poset idea, see particularly \S \SplittingSection. 

The {\em product decomposition problem} is the problem of explicitly writing a product of Weyl bialternants as a $\mathbb{Z}$-linear combination of Weyl bialternants, i.e.\ determining the  $c_{\mu}$'s in the identity $\displaystyle \chi_{_{\nu}}\chi_{_{\lambda}} = \sum_{\mbox{\footnotesize dominant }\mu}c_{\mu}\chi_{_{\mu}}$. 
Nonnegativity of the $c_{\mu}$ product decomposition coefficients can be deduced from the Weyl character formula together with the complete reducibility of modules for any semisimple Lie algebra $\mathfrak{g}$.  
In particular, if $V(\nu)$, $V(\lambda)$ etc denote irreducible $\mathfrak{g}$-modules, then $V(\nu) \otimes V(\lambda)$ can be written as $V(\mu_1) \oplus \cdots \oplus V(\mu_k)$ for some dominant weights $\mu_i$, in which case $\chi_{_{\nu}}\chi_{_{\lambda}} = \sum_{i=1}^{k}\chi_{_{\mu_i}}$. 
In this monograph we deduce nonnegativity of the $c_{\mu}$'s using the combinatorics of certain ``crystalline'' splitting posets and obtain at least an implicit description of them as counts of combinatorial objects, see \CrystallineRemark.B. 
These results are not new, but we believe our crystalline splitting poset approach is. 

The {\em branching problem} is the problem of explicitly writing a Weyl bialternant (or other Weyl symmetric function) as a $\mathbb{Z}$-linear combination of $\Phi_{J}$-Weyl bialternants when $J \subseteq I$.  That is, for the $\Phi$-Weyl bialternant $\chi_{_{\lambda}}$ and its restriction $\chi_{_{\lambda}}|_{_J}$ to a $\Phi_{J}$-Weyl symmetric function, we want to explicitly determine the  $b_{\mu}$'s in the identity $\displaystyle \chi_{_{\lambda}}|_{_J} = \sum_{\mbox{\footnotesize dominant }\mu}b_{\mu}\chi^{\Phi_J}_{_{\mu}}$.  
As in the case of the product decomposition coefficients of the previous paragraph, nonnegativity of the $b_{\mu}$ branching coefficients can be deduced from representation theory. 
We rederive these results here and obtain a combinatorial description of the restriction coefficients using our new approach of crystalline splitting posets, see \CrystallineRemark.B. 

Central to our overall perspective is the general algebraic combinatorial problem of finding solutions $R$ to the following formula for a given $J \subseteq I$ and $\Phi_{J}$-dominant weight $\nu$: 
\begin{equation}
\chi^{\Phi_J}_{_{\nu}} \cdot \WGF(R)|_{J} = \sum_{\selt \in \mathcal{S}_{J,\nu}(R)}\chi^{\Phi_J}_{_{\nu + wt^{J}(\selt)}}.
\end{equation} 
Here, we envision $R$ to be some ranked poset whose Hasse diagram edges are colored by $I$; $\WGF(R)|_{J}$ to be a certain weight generating function for $R$ that only regards those Hasse diagram edges with colors from the subset $J$ of $I$; and $\mathcal{S}_{J,\nu}(R)$ to be some specified vertex subset of $R$ (see \InitialSplittingTheorem/\InitialSplittingDefinition). 
In \NewDefinitionRemark.D we say precisely how this formula pertains to our three problems. 
Generally speaking, the $J = I$ version of the formula addresses product-decomposition-type problems, the $\nu = 0$ version addresses branching-type problems, and the version with $J=I$ and $\nu=0$ addresses purely splitting-type problems. 
An ``all-in-one'' version of the splitting, product decomposition, and/or branching problems is, then, to find posets $R$ and subsets $\mathcal{S}_{J,\nu}(R)$ satisfying the formula (1) above. 
We call $R$ a ``splitting poset'' when $J=I$ and $\nu=0$, and a ``$(J,\nu)$-splitting poset'' or ``refined splitting poset'' otherwise. 

Of course, answers to our three principal problems are known, although they have varying degrees of combinatorial explicitness. 
Given the relationship between type $\myA_{n}$-Weyl bialternants and Schur functions (see \CaseAExample), the classical ``Littlewood--Richardson rule'' (see for example \cite{StanText2}) resolves the product decomposition problem for the $\myA_{n}$ case.  
Versions of this result for the other classical irreducible root systems $\myB_{n}$, $\myC_{n}$, $\myD_{n}$ were obtained in \cite{LitTableaux} (using standard monomial theory) and \cite{Nakashima} (using crystal bases).   
Crystal graphs as well as the edge-colored directed graphs underlying Littelmann's famous path model (see \cite{Lit2}, \cite{Lit4} \S 6) provide completely general answers to the splitting, product decomposition, and branching problems.   
One consequence is that the product decomposition coefficients and the branching coefficients are nonnegative. 
Stembridge's admissible systems also answer these three basic problems and proceed from a more transparently combinatorial starting point. 
Admissible systems encompass both crystal graphs and the graphs obtained from Littelmann's path model; the latter fact is demonstrated explicitly in \cite{Stem} \S 8. 

In \S \FibrousSection\ and \S \CrystallineSection, we use a vertex-coloring technique on crystalline splitting posets to arrive at a very similar general solution to the three principal problems. 
That said, our ambition that solutions to formula (1) be combinatorially explicit and interesting leaves much room for the further study of splitting posets and their refinements in addressing our three principal problems. 

{\bf On the current state of knowledge, some open questions, and where we hope to go from here.}  
In \IntroTable, we present, to the best of our knowledge, what is currently known about splitting posets. 
This table includes all instances of splitting distributive lattices that we know of. 

Some problems and questions that are open (as far as the author knows) are explicitly identified in this monograph: See \ProblemsList. 
Some other open questions are noted elsewhere in the exposition: See the comments for items 10, 13, 14, 15 of \IntroTable; the last paragraph of \CaseAExample; and the last paragraph of \S \CrystalSection. 

Overall, this monograph and its sequels will hopefully provide some new and interesting combinatorial structures and also a useful framework for understanding them within a variety of contexts.  
It will also hopefully inspire some interest in the central problem posed here, \StanleyProblem: 
\begin{center}
\parbox{5.25in}{For families of Weyl bialternants or Weyl symmetric functions, produce splitting posets (or splitting modular/distributive lattices) with explicit and interesting combinatorial descriptions.}
\end{center} 

{\bf Acknowledgements.}  
The author thanks Wyatt Alverson, Will Atkins, Katheryn Beck, Jacob Dennerlein, Beth Donovan, Molly Dunkum, Kimmo Eriksson, Matt Gilliland, Cristian Lenart, Scott Lewis, Jordan Love, Sasha Malone, Bob Pervine, Bob Proctor, Tim Schroeder, John Stembridge, Michael Strayer, and Norman Wildberger for many helpful discussions during the development of this monograph. 
He thanks Ed Thome, Steve Cobb, and Claire Fuller of Murray State University for helping facilitate the sabbatical leaves in 2011 and 2021 that provided opportunities to more fully develop these ideas. 
There are no bounds on his devotion to his lovely and brilliant wife and collaborator, Molly Dunkum. 
He owes a particular debt of gratitude to Bob Proctor, whose work, input, perspective, and taste influenced many of the ideas presented here. 

\begin{center}
\IntroTable: Splitting posets, case by case. Entries 1 through 7 (out of 15).%

{\scriptsize
\vspace*{0.2cm} %
\begin{tabular}{|c|c|c|c|}
\hline %
\parbox{1in}{\begin{center}Irreducible root systems and dominant weights\end{center}} %
& \parbox{1in}{\begin{center}Family of splitting posets for the related Weyl bialternants\end{center}}  %
& \parbox{1in}{\begin{center}References\end{center}} %
& \parbox{2in}{\begin{center}Comments\end{center}}\\
\hline %
\hline %
\parbox{1in}{\begin{center}1. Any irreducible root system and any dominant weight\end{center}} 
& \parbox{1in}{\begin{center}Crystal graphs\end{center}}
& \parbox{1in}{\begin{center}\cite{Kash1}, \cite{Kash2}, \cite{Stem}, \S \CrystalSection\ below\end{center}} 
& \parbox{2.75in}{``Explicit'' (in particular, nonrecursive) realizations of crystal graphs can be found in \cite{KN} for types $\myA_{n}$---$\myD_{n}$; in \cite{LitTableaux} or \cite{KM} for type $\myG_{2}$; and in \cite{JS} for $\myE_{6}$.  For a general non-recursive approach, see \cite{JS}.}\\
\hline %

\parbox{1in}{\begin{center}2. Minuscule dominant weights\end{center}} 
& \parbox{1in}{\begin{center}Minuscule splitting posets\end{center}}
& \parbox{1in}{\begin{center}\cite{PrEur}, \S \SplittingSection\ below\end{center}} 
& \parbox{2.75in}{These are all distributive lattices.}\\
\hline %

\parbox{1in}{\begin{center}3. Quasi-minuscule dominant weights\end{center}} 
& \parbox{1in}{\begin{center}Quasi-minuscule splitting posets\end{center}}
& \parbox{1in}{\begin{center}\cite{Stem}, \S \SplittingSection\ below\end{center}} 
& \parbox{2.75in}{In each case, the quasi-minuscule splitting poset is the unique connected edge-minimal splitting poset.}\\
\hline %
\parbox{1in}{\begin{center}4. Highest root\end{center}} 
& \parbox{1in}{\begin{center}The $n$ extremal supporting graphs of \cite{DonAdjoint}\end{center}}
& \parbox{1in}{\begin{center}\cite{DonAdjoint}\end{center}} 
& \parbox{2.75in}{\ 

These are  indexed by the $n$ simple roots.  Each is a modular lattice.  It is distributive if and only if the irreducible root system has a branchless Dynkin diagram and the given extremal supporting graph corresponds to an end node simple root.

\

}\\
\hline %

\parbox{1in}{\begin{center}5. Highest short root\end{center}} 
& \parbox{1in}{\begin{center}The $m$ extremal supporting graphs of \cite{DonAdjoint}\end{center}}
& \parbox{1in}{\begin{center}\cite{DonAdjoint}\end{center}} 
& \parbox{2.75in}{\ 

These are indexed by the $m$ short simple roots. Each is a modular lattice and coincides with the quasi-minuscule splitting poset only in the $\myC_{2}$ and $\myG_{2}$ cases.  Such a modular lattice is distributive if and only if the given extremal supporting graph corresponds to an end node on the branchless subgraph of the Dynkin diagram induced by the short simple roots.

\

}\\
\hline %
\parbox{1in}{\begin{center}6. Rank two irreducible root systems and all dominant weights\end{center}} 
& \parbox{1in}{\begin{center}Semistandard lattices\end{center}}
& \parbox{1in}{\begin{center}\cite{ADLMPPW}, \end{center}} 
& \parbox{2.75in}{\ 

These are all distributive lattices.  See also \cite{ADLPTwo}.  It is known exactly which of these are supporting graphs, see \cite{ADLPOne}. 

\

}\\
\hline %

\parbox{1in}{\begin{center}7. Any irreducible root system and any combination of adjacency-free fundamental weights\end{center}} 
& \parbox{1in}{\begin{center}Semistandard lattices\end{center}}
& \parbox{1in}{\begin{center}\cite{DW}; for a combinatorial description of the adjacency-free fundamental weights, see \cite{DonNumbers}\end{center}} 
& \parbox{2.75in}{We obtain, in a uniform way, splitting distributive lattices for the following (root system, dominant weight) pairs: $(\myA_{n},\lambda)$, where $\lambda$ is any dominant weight; $(\myB_n,a\omega_{1}+b\omega_{n})$; $(\myC_n,a\omega_{1}+b\omega_{n})$; $(\myD_n,a\omega_{1}+b\omega_{n-1}+c\omega_{n})$; $(\myE_6,a\omega_{1}+b\omega_{6})$; $(\myE_7,a\omega_{7})$; $(\myG_2,a\omega_{1}+b\omega_{2})$.}\\
\hline %
\end{tabular}
}
\end{center}

\begin{center}
\IntroTable\ (continued): Splitting posets, case by case. Entries 8 through 13 (out of 15).%

{\scriptsize
\vspace*{0.2cm} %
\begin{tabular}{|c|c|c|c|}
\hline %
\parbox{1in}{\begin{center}Irreducible root systems and dominant weights\end{center}} %
& \parbox{1in}{\begin{center}Family of splitting posets for the related Weyl bialternants\end{center}}  %
& \parbox{1in}{\begin{center}References\end{center}} %
& \parbox{2in}{\begin{center}Comments\end{center}}\\
\hline %
\hline %
\parbox{1in}{\begin{center}8. $\myC_n$ and any fundamental weight\end{center}} 
& \parbox{1in}{\begin{center}The ``KN'' and ``De Concini'' symplectic lattices\end{center}}
& \parbox{1in}{\begin{center}Defined in \cite{DonThesis} and again in \cite{DonPeck}\end{center}} 
& \parbox{2.75in}{\ 

These are two families of distributive lattice supporting graphs for the fundamental representations of the symplectic Lie algebras.  These coincide for the $\omega_1$ fundamental weight.  For the $\omega_n$ fundamental weight, the common distributive lattice is sometimes called a ``Catalan lattice.''

\

}\\
\hline %

\parbox{1in}{\begin{center}9. $\myB_n$ with $\lambda \in \{\omega_k\}_{1 \leq k < n}$ or $\lambda = 2\omega_n$\end{center}} 
& \parbox{1in}{\begin{center}The ``KN'' and ``De Concini'' odd orthogonal lattices\end{center}}
& \parbox{1in}{\begin{center}Defined in \cite{DonThesis} and studied further in \cite{Beck}\end{center}} 
& \parbox{2.75in}{These are two families of distributive lattice supporting graphs for these (mostly) fundamental representations of the odd orthogonal Lie algebras.  These coincide for the $\omega_1$ fundamental weight.}\\
\hline %

\parbox{1in}{\begin{center}10. $\myD_n$ with $\lambda \in \{\omega_k\}_{1 < k < n-1}$\end{center}} 
& \parbox{1in}{\begin{center}{\bf ?}\end{center}}
& \parbox{1in}{\begin{center} --- \end{center}} 
& \rule[-10mm]{0mm}{21mm}\parbox{2.75in}{The problem of finding splitting modular lattices or modular lattice supporting graphs for these cases appears to be open and is, in our view, the most interesting open case.}\\
\hline %
\parbox{1in}{\begin{center}11. $\myB_n$ and $\myG_2$ with any ``one-rowed'' dominant weight \end{center}} 
& \parbox{1in}{\begin{center}The ``Molev'' and ``RS'' one-rowed odd orthogonal lattices, and the ``Molev'' and ``Littelmann'' one-rowed $\myG_2$ lattices\end{center}}
& \parbox{1in}{\begin{center}The RS lattices were obtained in \cite{RS}; the other lattices were obtained in \cite{DLP1}\end{center}} 
& \parbox{2.75in}{What we refer to as the Molev versions of these splitting distributive lattices (in fact, distributive lattice supporting graphs) coincide with the semistandard lattices of entry 7 of this table. The RS and Littelmann lattices are obtained in a different way.}\\
\hline %
\parbox{1in}{\begin{center}12. $\myC_n$ with any one-rowed dominant weight \end{center}} 
& \parbox{1in}{\begin{center}The ``Molev'' and ``RS'' one-rowed symplectic lattices\end{center}}
& \parbox{1in}{\begin{center}The  Molev lattices were observed in \cite{ADLPOne}; the RS lattices were studied in \cite{Atkins}\end{center}} 
& \parbox{2.75in}{What we refer to as the Molev versions of these splitting distributive lattices (in fact, distributive lattice supporting graphs) coincide with the semistandard lattices of entry 7 of this table. The RS lattices are obtained in a different way.}\\
\hline %
\parbox{1in}{\begin{center}13. $\myD_n$ with any one-rowed dominant weight \end{center}} 
& \parbox{1in}{\begin{center}{\bf ?}\end{center}}
& \parbox{1in}{\begin{center} --- \end{center}} 
& \parbox{2.75in}{\ 

Splitting distributive lattices for these cases are covered by entry 7 of this table and are easy analogs of the Molev lattices from entries 11 and 12.  It is an open question whether there are nice $\myD_n$ analogs of the RS lattices of entries 11 and 12.

\ 

}\\
\hline %
\end{tabular}
}
\end{center}

\newpage
\begin{center}
\IntroTable\ (continued): Splitting posets, case by case. Entries 14 and 15 (out of 15).%

{\scriptsize
\vspace*{0.2cm} %
\begin{tabular}{|c|c|c|c|}
\hline %
\parbox{1in}{\begin{center}14. $\myF_4$ and certain combinations of fundamental weights\end{center}} 
& \parbox{1in}{\begin{center}Not yet named\end{center}}
& \parbox{1in}{\begin{center}\cite{Gilliland}\end{center}} 
& \parbox{2.75in}{\ 

In \cite{Gilliland}, Gilliland produced splitting distributive lattices for $\chi_{_{2\omega_1}}^{\mytinyF_4}$ and $\chi_{_{2\omega_4}}^{\mytinyF_4}$.  The corresponding representations of the type $\myF_4$ semisimple Lie algebra have dimensions 1053 and 324 respectively.  He also showed that $\chi_{_{\omega_2}}^{\mytinyF_4}$ and $\chi_{_{\omega_3}}^{\mytinyF_4}$ do not have splitting distributive lattices.  It seems possible that the construction of \cite{Gilliland} can be extended to obtain splitting distributive lattices for $\chi_{_{a\omega_{1}+b\omega_{4}}}^{\mytinyF_4}$.

\

}\\
\hline %

\parbox{1in}{\begin{center}15. Any dominant weight in the $\myB_n$, $\myC_n$, $\myD_n$, and $\myE_6$/$\myE_7$/$\myE_8$ cases\end{center}} 
& \parbox{1in}{\begin{center}{\bf ?}\end{center}}
& \parbox{1in}{\begin{center} --- \end{center}} 
& \parbox{2.75in}{\ 

Apparently it is an open problem to find splitting distributive/modular lattices for a generic Weyl bialternant in these cases.  Supporting graphs for Molev's weight bases for irreducible representations of the classical Lie algebras are candidates (see \cite{MoOne}, \cite{MoTwo}, \cite{MoThree}).  However, it is known that some of these are not modular lattices.

\

}\\
\hline %
\end{tabular}
}
\end{center}

\newpage
\noindent 
{\Large \bf \S \WeylSection.\ Weyl symmetric functions and Weyl bialternants.} 

Our goal in this section is to explicate the basic properties of Weyl symmetric functions and Weyl bialternants primarily from a root-system/weight-lattice/Weyl-group perspective.  
There does not seem to be comparable exposition elsewhere in a single source, so the discussion here will hopefully serve as a helpful reference. 
Readers experienced with the subject of Weyl bialternants (aka Weyl characters) as a generalization of Schur functions might initially browse the results of this section as an orientation to language and notation.  
The main prerequisite is familiarity with finite root systems, which serve as our starting point. 
As our focus is on Weyl bialternants modelled by finite posets, we only consider finite root systems and not the more general root systems and characters of symmetrizable Kac--Moody theory.   
We freely use background results on finite root systems, weights, and Weyl groups from Ch.\ III of \cite{Hum}  and at times borrow from \cite{FH} \S\S 23--25, \cite{Ram} \S 5, \cite{Stem} \S 2, and \cite{StemAIM}.  
To keep the monograph as combinatorial and self-contained as possible, we limit the use of Lie algebra representation theory to the comments at the end of this section, the discussion of supporting graphs in \S \SplittingSection, the classification in \S \SplittingSection\ of Weyl bialternants with unique splitting posets, and the discussion in \S \CrystalSection\ of crystal bases/graphs.   
Since we are not at the outset assuming any connection between Weyl bialternants and the representation theory of semisimple Lie algebras (e.g.\ we will work for the most part without Weyl's character formula),  
then we have to proceed somewhat carefully. 

{\bf Roots and weights.} Let $\EuclideanE$ be an $n$-dimensional 
Euclidean space with inner product $\langle \cdot,\cdot \rangle$. At 
times we use $||v||^{2}$ to mean $\langle v,v \rangle$ when  
$v \in \mathfrak{E}$.  
Take $I$ to be an index set of cardinality $n$ (usually, $I = 
\{1,2,\ldots,n\}$). 
Let $\Phi$ be a finite rank $n$ root system in $\EuclideanE$ together 
with a choice of simple roots 
$\{\alpha_{i}\}_{i \in I}$. 
We sometimes refer to rank $n$ finite irreducible root systems 
by their classification type using 
the notation $\myX_{n}$ (where $\myX \in \{\myA, 
\myB, \myC, \myD, \myE, \myF, \myG\}$).   
Our numbering of the simple roots and associated Dynkin diagrams 
usually follows Ch.\ 11 of \cite{Hum}, with the exception that we 
begin the $\myB_{n}$ series at $n=3$ and the $\myC_{n}$ series at 
$n=2$. 
Throughout this monograph, we use ``$\Phi$'' as a prefix, superscript, or subscript in terminology/notation to emphasize dependence on the root system.  

With respect to the given choice $\{\alpha_{i}\}_{i \in I}$ of simple roots, a root $\alpha$ is positive (respectively, negative) if $\alpha = \sum_{i \in I}a_{i}\alpha_{i}$ with each $a_{i} \geq 0$ (respectively, each $a_{i} \leq 0$).  
We have $\Phi = \Phi^{+} \cup \Phi^{-}$, where $\Phi^{+}$ and $\Phi^{-}$ are respectively the sets of positive and negative roots. 
For a root $\alpha\in\Phi$, $\alpha^{\vee} = \frac{2\alpha}{\langle \alpha,\alpha \rangle}$ is its coroot.  
The set $\Phi^{\vee} = \{\alpha^{\vee}\, |\, \alpha \in \Phi\}$ is a root system -- the ``coroot'' or ``dual'' root  system'' -- with $\{\alpha_{i}^{\vee}\}_{i \in I}$ as the preferred set of simple roots. 
The Cartan matrix $M = M_{\Phi}$ has $M_{ij} = \langle \alpha_{i},\alpha_{j}^{\vee} \rangle \in \mathbb{Z}$ for all $i,j \in I$.  

The  fundamental weights $\{\omega_{i}\}_{i \in I}$ are dual 
to the 
simple coroots with respect to the inner product on $\EuclideanE$. 
The lattice 
of weights $\Lambda$ is the $\mathbb{Z}$-span of the fundamental 
weights.  
The dominant weights $\Lambda^{+}$ 
are those nonnegative integral linear 
combinations of fundamental weights.  The set $\Lambda^{++}$ denotes 
the set of {\em strongly} dominant weights, i.e.\ those 
$\lambda \in \Lambda$ such that 
$\lambda = \sum_{i \in I}a_{i}\omega_{i}$ 
with each $a_{i}$ a positive integer.  
Partially order $\Lambda$ by the rule 
$\mu \leq \nu$ if and only if $\nu - \mu = \sum_{i \in 
I}k_{i}\alpha_{i}$ with each $k_{i}$ a nonnegative integer.  

{\bf Some key quantities.} 
Certain calculations in $\Lambda$ appear often in our 
discussion.  In particular, in 
this paragraph we define Kostant's partition 
function and special elements $\varrho$ and $\varrho^{\vee}$.  
Kostant's {\em partition function} $\mathcal{P}: \Lambda 
\longrightarrow 
\mathbb{Z}_{\geq 0}$ is defined by $\mathcal{P}(\mu) := 
|\{\mbox{sequences } (k_{\alpha})_{\alpha \in \Phi^{+}}\, |\, 
\mbox{each } k_{\alpha} 
\in \mathbb{Z}_{\geq 0} \mbox{ and } \mu = \sum_{\alpha \in 
\Phi^{+}}k_{\alpha}\alpha\}|$, i.e.\ the number of 
distinct ways $\mu$ can be written as a nonnegative integer linear 
combination of positive roots. 
Let $\varrho := \sum_{i \in I}\omega_{i}$. 
It is well known that  $\varrho  = 
\frac{1}{2}\sum_{\alpha \in \Phi^{+}}\alpha$ (\cite{Hum} \S 13.3).  
Define $\varrho^{\vee} := \sum_{i \in 
I}\frac{2}{\langle \alpha_{i},\alpha_{i} \rangle}\omega_{i}$, which 
is in $\mathfrak{E}$ but not necessarily $\Lambda$.  
A simple computation shows that $\varrho^{\vee}$ is half the sum of 
the positive coroots in the dual root system $\Phi^{\vee}$. 
Observe 
that $\langle \alpha_{i},\varrho^{\vee} \rangle = 1$ for each $i \in 
I$.  
Note that for any $i \in I$, $\alpha_{i} = \sum_{j \in 
I}M_{ij}\omega_{j}$. 
From basic facts about real inner product spaces it follows that 
since the basis of simple roots is ``nonacute'' ($\langle 
\alpha_{i},\alpha_{j} \rangle \leq 0$ for all distinct $i,j \in I$), 
then the dual basis of fundamental weights is ``nonobtuse'' 
($\langle 
\omega_{i},\omega_{j} \rangle \geq 0$ for all distinct $i,j \in I$). 
If we write $\omega_{j} = \sum_{k \in I}Q_{jk}\alpha_{k}$, then 
$M^{-1} = (Q_{jk})_{j,k \in I}$, so by Cramer's rule each $Q_{jk} \in 
\mathbb{Q}$.  The fact that $\langle 
\omega_{j},\omega_{k} \rangle \geq 0$ for all $j,k \in I$ 
shows that each $Q_{jk} =  \frac{2}{\langle \alpha_{k},\alpha_{k} \rangle} \langle \omega_{j},\omega_{k} \rangle$ is nonnegative.  Observe that 
$Q_{jk} > 0$ when $\langle 
\omega_{j},\omega_{k} \rangle > 0$, so in particular 
$Q_{jj} > 0$.  So, $\langle 
\omega_{i},\varrho^{\vee} \rangle = \langle \sum_{j \in 
I}Q_{ij}\alpha_{j},\varrho^{\vee} \rangle = \sum_{j \in I}Q_{ij}$ is a 
positive rational number.  
In particular, $\langle \nu,\varrho^{\vee} \rangle > 0$ for all $\nu 
\in \Lambda^{+}$. 
Write the positive rational number $\langle \omega_{i},\varrho^{\vee} 
\rangle$ in lowest terms as 
$p_{i}/q_{i}$ for positive integers $p_{i}, q_{i}$. 
An easy calculation shows that for any $\mu \in \Lambda$, 
we have $\langle 
\mu,\varrho^{\vee} \rangle > 0$ if and only if $\langle 
\mu,\varrho^{\vee} \rangle \geq \myd_{\Phi}$, where  
$\myd_{\Phi} := \frac{1}{q_{1} \cdots q_{n}}$, which we call our {\em mesh size estimate}.  Our main use of the 
quantity $\myd_{\Phi}$ is in finiteness arguments for \FinitenessTheorem\ and \BasisProp/\BasisTheorem. 

{\bf The Weyl group.} The (finite) Weyl 
group $W = W_{\Phi} \subset GL(\EuclideanE)$ 
is generated by the simple reflections $\{s_{i}\}_{i \in 
I}$, where $s_{i}: \EuclideanE \rightarrow \EuclideanE$ 
is given by $s_{i}(v) = v - 
\langle v,\alpha_{i}^{\vee} \rangle\alpha_{i}$. Note that $W$ fixes 
the lattice $\Lambda$: We have 
$s_{i}(\omega_{j}) = \omega_{j} - \delta_{ij}\alpha_{i} \in \Lambda$ 
for all $i,j \in I$ 
if and only if $\alpha_{i} \in \Lambda$ for all $i \in I$, 
and the latter is guaranteed 
by integrality of the $M_{ij}$'s.  
A simple computation now shows that $\langle 
s_{i}(\omega_{j}),s_{i}(\omega_{k}) \rangle = 
\langle \omega_{j},\omega_{k} 
\rangle$ for all $i, j, k \in I$, and hence that $\langle \cdot,\cdot 
\rangle$ is $W$ invariant. 
It follows from \S 10.3 of \cite{Hum} that $W$ is invariant of the choice of simple roots, and that $W(\{\alpha_{i}\}_{i \in I}) = \Phi$. 

Note that for all $i, j \in I$ we have $M_{ij}M_{ji} \in \{0, 
1, 2, 3, 4\}$.  Let $\myspecialm_{ij} = 
\myspecialm_{ji}$ be the unique positive integer 
for which $M_{ij}M_{ji} = 
4\cos^{2}(\pi/\myspecialm_{ij})$, so corresponding to 
the list of numbers in the previous sentence we get 
$\myspecialm_{ij} \in \{2, 3, 
4, 6, 1\}$.  It is 
straightforward to check that the said generators for $W$ satisfy the 
following relations: for all $i,j \in I$, 
$(s_{i}s_{j})^{\myspecialm_{ij}} = 
Id$, where $Id: \mathfrak{E} \longrightarrow \mathfrak{E}$ is the 
identity transformation.  For irreducible $\Phi$, observe that the 
order of the Weyl group as determined in \S 12 of \cite{Hum} 
is the same as the order of the finite 
irreducible Coxeter group of the same type 
as determined in \S 2.11 of \cite{HumCoxeter}.  So $(W,\{s_{i}\}_{i 
\in I})$ is a Coxeter system. 

Let $w_{0}$ denote the unique longest element of the 
Weyl group (cf.\ \cite{Hum} Exercise 10.9).  It can be seen that 
$w_{0}(\Phi^{+}) = \Phi^{-}$ and that there is a unique 
permutation $\sigma_{0}: I \longrightarrow I$ such that 
for each $i \in I$, 
$w_{0}.\alpha_{i} = -\alpha_{\sigma_{0}(i)}$. 
Since $w_{0}^{2} = 
Id$ ($w_{0}^{-1}$ is also longest and therefore must be $w_{0}$), 
then $\sigma_{0}$ is an involution. 
We also get $w_{0}(\varrho) = w_{0}\big(\frac{1}{2}\sum_{\alpha \in 
\Phi^{+}}\alpha\big) = \frac{1}{2}\sum_{\beta \in \Phi^{-}}\beta = 
-\varrho$.  
Moreover, since $\langle w_{0}(\omega_{i}),\alpha_{j}^{\vee} 
\rangle = -\langle \omega_{i},\frac{2}{\langle \alpha_{j},\alpha_{j} 
\rangle}\alpha_{\sigma_{0}(j)}\rangle = 
-\langle \omega_{i},\frac{2}{\langle 
w_{0}(\alpha_{j}),w_{0}(\alpha_{j}) 
\rangle}\alpha_{\sigma_{0}(j)}\rangle = -\langle 
\omega_{i},\alpha^{\vee}_{\sigma_{0}(j)} \rangle = 
-\delta_{i,\sigma_{0}(j)}$, for all $i,j \in I$, then it must be the 
case that $w_{0}(\omega_{i}) = -\omega_{\sigma_{0}(i)}$.  
For irreducible root systems, it is well-known that $\sigma_{0}$ 
is trivial except in the cases $\myA_{n}$ $(n \geq 2)$, 
$\myD_{2k+1}$ $(k \geq 2)$, and 
$\myE_{6}$; see \DiagramSymmetryFigure.  

\begin{figure}[htb]
\begin{center}
\DiagramSymmetryFigure:  Action of the permutation $\sigma_{0}$ 
when $\sigma_{0}$ is not the identity.\\ 
{\footnotesize (In each case the permutation $\sigma_{0}$ is 
depicted as a symmetry 
of the associated Dynkin diagram.)}

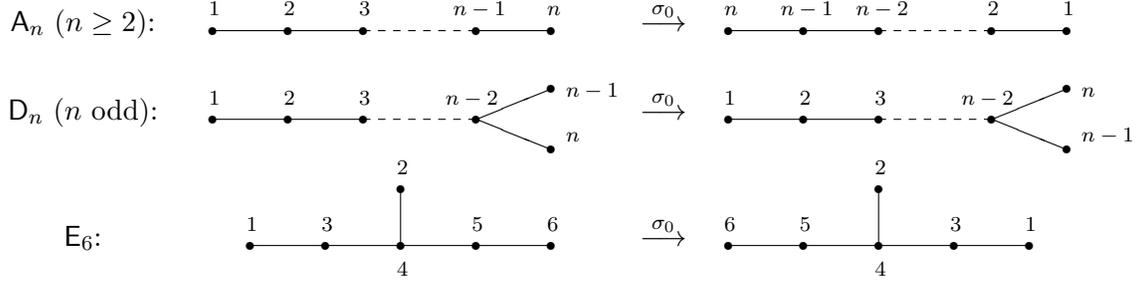
\begin{tabular}{cccc}
$\myA_{n}$ $(n \geq 2)$: \quad 
& \setlength{\unitlength}{1cm}
\begin{picture}(4.5,0.85)
\put(0,0){\circle*{0.1}}
\put(1,0){\circle*{0.1}}
\put(2,0){\circle*{0.1}}
\put(3.5,0){\circle*{0.1}}
\put(4.5,0){\circle*{0.1}}
\put(3.5,0){\line(1,0){1}}
\put(0,0){\line(1,0){2}}
\multiput(2,0)(0.2,0){8}{\line(1,0){0.1}}
\put(-0.05,0.2){\scriptsize ${1}$}
\put(0.95,0.2){\scriptsize ${2}$}
\put(1.95,0.2){\scriptsize ${3}$}
\put(3.2,0.2){\scriptsize ${n-1}$}
\put(4.45,0.2){\scriptsize ${n}$}
\end{picture}
& \hspace*{0.25in}
$\mylongarrow{\sigma_{0}}$ 
& \setlength{\unitlength}{1cm}
\begin{picture}(5,0.85)
\put(0,0){\circle*{0.1}}
\put(1,0){\circle*{0.1}}
\put(2,0){\circle*{0.1}}
\put(3.5,0){\circle*{0.1}}
\put(4.5,0){\circle*{0.1}}
\put(3.5,0){\line(1,0){1}}
\put(0,0){\line(1,0){2}}
\multiput(2,0)(0.2,0){8}{\line(1,0){0.1}}
\put(-0.05,0.2){\scriptsize ${n}$}
\put(0.7,0.2){\scriptsize ${n-1}$}
\put(1.7,0.2){\scriptsize ${n-2}$}
\put(3.45,0.2){\scriptsize ${2}$}
\put(4.45,0.2){\scriptsize ${1}$}
\end{picture}\\
$\myD_{n}$ ($n$ odd): \quad 
& \setlength{\unitlength}{1cm}
\begin{picture}(4.5,1)
\put(0,0){\circle*{0.1}}
\put(1,0){\circle*{0.1}}
\put(2,0){\circle*{0.1}}
\put(3.5,0){\circle*{0.1}}
\put(4.5,0.4){\circle*{0.1}}
\put(4.5,-0.4){\circle*{0.1}}
\put(0,0){\line(1,0){2}}
\put(3.5,0){\line(5,2){1}}
\put(3.5,0){\line(5,-2){1}}
\multiput(2,0)(0.2,0){8}{\line(1,0){0.1}}
\put(-0.05,0.2){\scriptsize ${1}$}
\put(0.95,0.2){\scriptsize ${2}$}
\put(1.95,0.2){\scriptsize ${3}$}
\put(3.1,0.2){\scriptsize ${n-2}$}
\put(4.7,0.3){\scriptsize ${n-1}$}
\put(4.7,-0.3){\scriptsize ${n}$}
\end{picture}
& \hspace*{0.25in}
$\mylongarrow{\sigma_{0}}$ 
& \setlength{\unitlength}{1cm}
\begin{picture}(5,1)
\put(0,0){\circle*{0.1}}
\put(1,0){\circle*{0.1}}
\put(2,0){\circle*{0.1}}
\put(3.5,0){\circle*{0.1}}
\put(4.5,0.4){\circle*{0.1}}
\put(4.5,-0.4){\circle*{0.1}}
\put(0,0){\line(1,0){2}}
\put(3.5,0){\line(5,2){1}}
\put(3.5,0){\line(5,-2){1}}
\multiput(2,0)(0.2,0){8}{\line(1,0){0.1}}
\put(-0.05,0.2){\scriptsize ${1}$}
\put(0.95,0.2){\scriptsize ${2}$}
\put(1.95,0.2){\scriptsize ${3}$}
\put(3.1,0.2){\scriptsize ${n-2}$}
\put(4.7,0.3){\scriptsize ${n}$}
\put(4.7,-0.3){\scriptsize ${n-1}$}
\end{picture}\\
$\myE_{6}$: \quad 
& \setlength{\unitlength}{1cm}
\hspace*{0.5cm}\begin{picture}(4,1.5)
\put(0,0){\circle*{0.1}}
\put(1,0){\circle*{0.1}}
\put(2,0){\circle*{0.1}}
\put(2,0.75){\circle*{0.1}}
\put(3,0){\circle*{0.1}}
\put(4,0){\circle*{0.1}}
\put(2,0){\line(0,1){0.75}}
\put(0,0){\line(1,0){4}}
\put(-0.05,0.2){\scriptsize ${1}$}
\put(1.95,0.95){\scriptsize ${2}$}
\put(0.95,0.2){\scriptsize ${3}$}
\put(1.95,-0.4){\scriptsize ${4}$}
\put(2.95,0.2){\scriptsize ${5}$}
\put(3.95,0.2){\scriptsize ${6}$}
\end{picture}
& \hspace*{0.25in}
$\mylongarrow{\sigma_{0}}$ 
& \setlength{\unitlength}{1cm}
\begin{picture}(5,1.5)
\put(0,0){\circle*{0.1}}
\put(1,0){\circle*{0.1}}
\put(2,0){\circle*{0.1}}
\put(2,0.75){\circle*{0.1}}
\put(3,0){\circle*{0.1}}
\put(4,0){\circle*{0.1}}
\put(2,0){\line(0,1){0.75}}
\put(0,0){\line(1,0){4}}
\put(-0.05,0.2){\scriptsize ${6}$}
\put(1.95,0.95){\scriptsize ${2}$}
\put(0.95,0.2){\scriptsize ${5}$}
\put(1.95,-0.4){\scriptsize ${4}$}
\put(2.95,0.2){\scriptsize ${3}$}
\put(3.95,0.2){\scriptsize ${1}$}
\end{picture}
\end{tabular}
\end{center}
\end{figure}

{\bf Weight diagrams.} 
For any dominant weight $\lambda$, the {\em weight diagram} 
$\Pi(\lambda)$ is the partially ordered set of weights 
$\{\mu \in \Lambda\, |\, \mu = w(\nu) 
\mbox{ for some } w \in W \mbox{ and some } \nu \in \Lambda^{+} 
\mbox{ with } \nu \leq \lambda\}$, where the partial 
order is induced from $\Lambda$.  
Weight diagrams figure prominently in our discussion.  
We record  
some combinatorial properties here and in \WeightRemarkTwo. 
A subset $\mathscr{W}$ of $\Lambda$ is {\em saturated} if for all 
$\nu \in \mathscr{W}$, $\alpha \in 
\Phi$, and $i$ between 0 and $\langle \nu,\alpha^{\vee} \rangle$ 
inclusive, we have $\nu - i\alpha \in \mathscr{W}$. 

\noindent 
{\bf \WeightRemarkOne}\ \ {\sl Let $\lambda \in \Lambda^{+}$.}  
{\sl (1)\ $\Pi(\lambda)$ is finite.}  
{\sl (2)\ $\Pi(\lambda)$ has $\lambda$ as its 
unique maximal element.}  
{\sl (3)\ The unique minimal element of $\Pi(\lambda)$ 
is $w_{0}(\lambda)$, and 
 for any $\mu \in \Pi(\lambda)$ the quantity $\langle 
\mu+\lambda,\varrho^{\vee} \rangle = \langle 
\mu-w_{0}(\lambda),\varrho^{\vee} \rangle$ is a nonnegative integer.}   
{\sl (4)\ $\Pi(\lambda)$ is the unique saturated set of weights  
whose maximal element under the induced order from $\Lambda$ is 
$\lambda$.} 
{\sl (5)\  For weights 
$\mu$ and $\nu$ in the weight diagram $\Pi(\lambda)$, 
$\mu$ is covered by $\nu$ (i.e.\ whenever $\mu \leq \xi \leq \nu$, then $\mu = \xi$ or $\nu = \xi$) if and 
only if $\mu + \alpha_{i} = \nu$ for some $i \in I$.} 

{\em Proof.} By \cite{Hum} Lemma 13.2.B, the set of 
dominant weights $\nu \leq \lambda$ is 
finite, settling claim {\sl (1)} of the lemma statement.  
By \cite{Hum} Lemma 13.2.A, we have $w.\nu \leq 
\nu$ in 
$\Lambda$ for all $w \in W$ and $\nu \in \Lambda^{+}$.  
It follows that $\Pi(\lambda)$ has $\lambda$ as its 
unique maximal element, cf.\ claim {\sl (2)}. 
With $\lambda = \sum_{i \in I}a_{i}\omega_{i}$, 
then $-w_{0}(\lambda) = \sum_{i \in I}a_{i}\omega_{\sigma_{0}(i)}$, 
which is dominant since $\lambda$ is. Also, since $-w_{0}$ is easily 
seen to permute the positive coroots, then 
$\langle -w_{0}(\lambda),\varrho^{\vee} \rangle = \langle 
\lambda,-w_{0}(\varrho^{\vee}) \rangle = \langle \lambda,\varrho^{\vee} 
\rangle$. It is also easy to see that 
$-w_{0} \in GL(\mathfrak{E})$ 
induces a poset isomorphism 
$\Pi(\lambda) \stackrel{\sim}{\rightarrow} \Pi(-w_{0}(\lambda))$: $\mu 
\leq \nu$ in $\Pi(\lambda)$ if and only if $\nu - \mu = 
\sum_{i \in I}k_{i}\alpha_{i}$ with each $k_{i} \in \mathbb{Z}_{\geq 
0}$ if and only if $w_{0}(\nu) - w_{0}(\mu) = -\sum_{i \in 
I}k_{i}\alpha_{\sigma_{0}(i)}$ if and only if $-w_{0}(\mu) \leq 
-w_{0}(\nu)$.  Since $-w_{0}(\lambda)$ is uniquely maximal in 
$\Pi(-w_{0}(\lambda))$, then $w_{0}(\lambda)$ must be uniquely 
minimal in 
$\Pi(\lambda)$. If $\mu \in \Pi(\lambda)$, then 
$w_{0}(\lambda) \leq \mu$, so $\mu-w_{0}(\lambda) = \sum_{i \in 
I}k_{i}\alpha_{i}$ where each $k_{i} \in \mathbb{Z}_{\geq 0}$.  Then, 
$\langle 
\mu-w_{0}(\lambda),\varrho^{\vee} \rangle = \sum_{i}k_{i}\langle 
\alpha_{i},\varrho^{\vee} \rangle = \sum_{i \in I}k_{i}$, clearly  
a nonnegative integer. 
This establishes {\sl (3)}.  Claim {\sl (4)} is taken from 
\S 13.4 of \cite{Hum}.  

The ``if'' part of claim {\sl (5)} is straightforward.  For the ``only if'' part, we need the 
following technical fact ({\tt *}): 
If $\mu < \nu$ in $\Lambda$ with $\nu - 
\mu = \sum_{j \in I}k_{j}\alpha_{j}$ ($k_{j} \in \mathbb{Z}_{\geq 
0}$), then for some $i \in I$ with $k_{i} > 0$ we have $\langle 
\mu,\alpha_{i}^{\vee} \rangle < 0$ or $\langle \nu,\alpha_{i}^{\vee} 
\rangle > 0$.  For otherwise, a simple computation shows that for each 
$i \in I$, the quantity $t_{i} := \langle \nu-\mu,\alpha_{i}^{\vee} 
\rangle = \sum_{j \in I}\langle \alpha_{j},\alpha_{i}^{\vee} \rangle 
k_{j}$ is nonpositive.  Relative to a total ordering of the index set 
$I$, let $[k]$ denote the column vector of $k_{i}$'s and $[t]$ the 
column vector of $t_{i}$'s.  Then $M^{\myT}[k] = [t]$, so $[k] = 
(M^{-1})^{\myT}[t]$.  Since, as argued above, the 
entries of $M^{-1}$ are nonnegative rational numbers, 
it follows that each $k_{j}$ is nonpositive. But this 
contradicts the fact that each $k_{j}$ is nonnegative and at least one 
of them is nonzero.  
So now suppose $\nu$ covers $\mu$ in 
$\Pi(\lambda)$.  Then by fact ({\tt *}), there is an $i \in I$ such 
that $k_{i} > 0$ when we write $\nu - 
\mu = \sum_{j \in I}k_{j}\alpha_{j}$ ($k_{j} \in \mathbb{Z}_{\geq 
0}$) and such that (a) $\langle 
\mu,\alpha_{i}^{\vee} \rangle < 0$ or (b) 
$\langle \nu,\alpha_{i}^{\vee} \rangle > 0$.  Since $\Pi(\lambda)$ is 
saturated, 
it follows from (a) that $\mu+\alpha_{i} \in \Pi(\lambda)$.  Since 
$k_{i} > 0$, then $\mu+\alpha_{i} \leq \nu$, and since $\mu$ is  
covered by $\nu$, then we must have $\mu+\alpha_{i} = \nu$.  
Similarly see in case (b) that $\mu = \nu-\alpha_{i}$.\hfill\QED

{\bf Reducible root systems; root subsystems.}  
Next, we detail how some of the preceding notions can be understood in the particular instances of root subsystems and reducible root systems. 

When $\Phi$ is reducible, we may write $\Phi$ as the disjoint union $\Phi_{1} \disjointunion \cdots \disjointunion \Phi_{k}$ for some irreducible root systems $\Phi_{i}$.  With $\mathfrak{E}_{i} := \mathrm{span}_{\mathbb{R}}(\Phi_{i})$, then $\mathfrak{E} = \mathfrak{E}_{1} \oplus \cdots \oplus \mathfrak{E}_{k}$, an orthogonal direct sum.  In particular, the simple roots for $\Phi$ are partitioned by the $\Phi_{i}$'s and may be ordered so that the Cartan matrix is block diagonal with blocks corresponding to the $\Phi_{i}$'s.  
We may write $\Lambda = \Lambda_{1} \oplus \cdots \oplus \Lambda_{k}$, where $\Lambda_{i} \subset \mathfrak{E}_{i}$ is the lattice of weights for $\Phi_{i}$.  So, for any $\mu \in \Lambda$ there exist unique $\mu_{i} \in \Lambda_{i}$ ($1 \leq i \leq k$) for which $\mu = \mu_{1} + \cdots \mu_{k}$, and $\mu$ is dominant (respectively, strongly dominant) in $\Lambda$ if and only if each $\mu_{i}$ is dominant (resp.\ strongly dominant) in $\Lambda_{i}$.  
We may also write $W \cong W_{1} \times \cdots \times W_{k}$, where for $1 \leq i \leq k$ we have $W_{i} := \{w|_{\mathfrak{E}_{i}}\}_{w \in  W}$ as the Weyl group associated with $\Phi_{i}$.  In particular, $W$ is a reducible Coxeter group.  Now take any $\lambda \in \Lambda^{+}$ with $\lambda = \lambda_{1} + \cdots + \lambda_{k}$ and each $\lambda_{i} \in \Lambda_{i}^{+}$.  Then the set $\Pi(\lambda)$ is identified with the Cartesian set product $\Pi(\lambda_{1}) \times \cdots \times \Pi(\lambda_{k})$. 

When $J \subseteq I$, we let $\Phi_{J}$ denote the root subsystem of $\Phi$ with simple roots $\{\alpha_{j}\}_{j \in J}$ and living in the Euclidean space $\mathfrak{E}_{J} := \mathrm{span}_{\mathbb{R}}\{\alpha_{j}\}_{j \in J} \subseteq \mathfrak{E}$. 
Some care must be taken in identifying $\Lambda_{\Phi_J}$, since $\mathrm{span}_{\mathbb{R}}\{\omega_{j}\}_{j \in J}$ need not be contained in $\mathfrak{E}_{J}$.  
In this situation, we obtain the fundamental weights $\{\omega_{j}^{J}\}_{j \in J}$ as follows: Let $M^{J}$ be the minor submatrix of the Cartan matrix $M$ consisting of those rows and columns of $M$ indexed by $J$.  Set $Q^{J} := (M^{J})^{-1}$.  
Now for each $j \in J$, set $\omega_{j}^{J} := \sum_{k \in J}Q^{J}_{jk}\alpha_{k}$.  
An easy calculation shows that $\langle \omega_{j} - \omega_{j}^{J},\alpha_{k}^{\vee} \rangle = 0$ for all $k \in J$, so $\omega_{j}^{J}$ is just the projection of $\omega_{j}$ onto the subspace $\mathfrak{E}_{J}$.  
In particular, $\langle \omega_{j}^{J},\alpha_{k}^{\vee} \rangle = \langle \omega_{j},\alpha_{k}^{\vee} \rangle = \delta_{jk}$, so $\{\omega_{j}^{J}\}_{j \in J}$ comprises the fundamental weights relative to this choice of simple roots for $\Phi_{J}$, and $\Lambda_{\Phi_J} \subset \mathfrak{E}_{J}$ is their $\mathbb{Z}$-span.  
In general, for any $v \in \mathfrak{E}$, we let $v^{J}$ denote its projection in $\mathfrak{E}_{J}$. 
Note that if we write $v = \sum_{i \in I} v_{i}\omega_{i}$, then $v^{J} = \sum_{j \in J}v_{j}\omega_{j}^{J}$. 
Now the subgroup $W_{J}$ of $W$ generated by $\{s_{j}\}_{j \in J}$ is isomorphic to the group $W_{\Phi_J}$ generated by the simple reflections $\{s_{j}|_{\mathfrak{E}_{J}}\}_{j \in J}$.  
It can be seen as follows that $(\Phi_J)^{\vee} = (\Phi^{\vee})_J$.  
If $\alpha \in \Phi_{J}$, then $\alpha = w(\alpha_{j})$ for some $j \in J$, $w \in W_{J}$. 
Then $\langle \alpha,\alpha \rangle = \langle \alpha_j,\alpha_j \rangle$, so $\alpha^{\vee} = w(\alpha_{j}^{\vee}) \in (\Phi^{\vee})_J$. 
Conversely, if $\beta = w(\alpha_{j}^{\vee}) \in (\Phi^{\vee})_J$ for some $j \in J$ and $w \in W_{J}$, then it is easy to calculate that $\beta = \alpha^{\vee}$ for $\alpha = w(\alpha_j)$, and therefore $\beta \in (\Phi_J)^{\vee}$. 
Using a Coxeter group viewpoint, it is possible to show that for any $\lambda = \sum_{i \in I}a_{i}\omega_{i} \in \Lambda^{+}$ such that $J = \{i \in I\, |\, a_{i} = 0\}$, then the stablizer $W_{\lambda}$ of $\lambda$ under the Weyl group action is just $W_{J} \cong W_{\Phi_J}$ (see for example \S 3 of \cite{DonNumbers}).  
In \PiJCorollary, we say how weight diagrams with respect to root subsystems can be understood in terms of weight diagrams with respect to the ``parent'' root system.

{\bf The ring of Weyl symmetric functions and some other important rings.} The group ring 
$\mathbb{Z}[\Lambda]$ has as its $\mathbb{Z}$-basis the formal 
symbols 
$\{e^{\mu}\, |\, \mu \in \Lambda\}$  
and multiplication rule $e^{\mu}e^{\nu} = e^{\mu + \nu}$.  The Weyl 
group ${W}$ acts on $\mathbb{Z}[\Lambda]$ by the rule 
$w.e^{\mu} := e^{w(\mu)}$.  
The {\em ring of Weyl symmetric functions} 
$\mathbb{Z}[\Lambda]^{W}$ is the ring of 
$W$-invariant elements of $\mathbb{Z}[\Lambda]$;  
elements of $\mathbb{Z}[\Lambda]^{W}$ are {\em Weyl symmetric functions} or $W$-{\em symmetric functions}.  
To emphasize dependence 
on the root system, we sometimes call these $\Phi$-Weyl symmetric functions or $W_{\Phi}$-symmetric functions. 
Now let $J \subseteq I$. For any $\chi = \sum_{\mu \in \Lambda}c_{\mu}e^{\mu} \in \mathbb{Z}[\Lambda]$, let ${\chi}|_{_J} := \sum_{\mu \in \Lambda}c_{\mu}e^{\mu^J}$, so ${\chi}|_{_J} \in \mathbb{Z}[\Lambda_{\Phi_J}]$.  It is easy to see that if $\chi$ is $W$-invariant, then ${\chi}|_{_J}$ is $W_J$-invariant. 

One way to view $W$-symmetric functions concretely as Laurent polynomials is as follows.  Fix a numbering of the simple roots, so $I = \{1,2,\ldots,n\}$.  For $i \in I$, let $z_{i} := e^{\omega_i}$, so for any weight $\mu = \sum m_{i}\omega_{i}$, then $e^{\mu}$ is the monomial $z_{1}^{m_{1}}\cdots{z_{n}}^{m_{n}}$.  Since $s_{i}(\omega_{j}) = \omega_{j} - \delta_{ij}\alpha_{i}$ and $\alpha_{i} = \sum_{k \in I}M_{ik}\omega_{k}$, then the $W$-action of a generator $s_{i}$ on any monomial $z_{1}^{m_1} \cdots z_{n}^{m_n}$ is given by 
\[s_{i}.(z_{1}^{m_1} \cdots z_{n}^{m_n}) = z_{1}^{m_{1}-M_{i,1}m_{i}} \cdots z_{i-1}^{m_{i-1}-M_{i,i-1}m_{i}} z_{i}^{-m_{i}} z_{i+1}^{m_{i+1}-M_{i,i+1}m_{i}} \cdots z_{n}^{m_{n}-M_{i,n}m_{i}},\] 
which is closely related to moves of the so-called ``numbers game'' (see e.g.\ \cite{DonNumbers} and references therein). 
With a root system / Cartan matrix / weight lattice / Weyl group in hand, one could begin a concrete theory of $W$-symmetric functions using the above numbers-game-type rule for the $W$-action on monomials as a starting point. 
Although this latter approach can be somewhat cumbersome, it analogizes the usual starting point for symmetric functions.  It should be noted that in the $\myA_{n}$ case, a change of variables is required in order to view these Laurent polynomials as classical symmetric (polynomial) functions.  See \CaseAExample. 

Although our main interest is in the ring of Weyl symmetric functions, 
we will at the outset enlarge $\mathbb{Z}[\Lambda]$ as 
follows.  
Define a linear ``height'' map $\mathrm{ht}: \mathfrak{E} 
\longrightarrow \mathbb{R}$ by 
$\mathrm{ht}(v) := \langle v,\varrho^{\vee} \rangle$ for all $v \in 
\mathfrak{E}$.  Let $\mathcal{R}$ be the ring of formal sums 
$\sum_{\mu \in \Lambda}c_{\mu}e^{\mu}$ such that $c_{\mu} \in 
\mathbb{Z}$ and for all $h \in 
\mathbb{R}$ there are only finitely many $\mu$ with $\mathrm{ht}(\mu) 
> h$ and $c_{\mu} 
\not= 0$.  
For $\varphi = \sum_{\mu \in \Lambda}c_{\mu}e^{\mu} \in \mathcal{R}$ 
and $\nu \in \Lambda$, 
we sometimes use $[e^{\nu}](\varphi)$ to denote the coefficient 
$c_{\nu}$. 
If $\varphi \not= 0$, let $\mathrm{ht}(\varphi) := 
\max\{\mathrm{ht}(\mu)\, |\, [e^{\mu}](\varphi) \not= 0\}$, 
$\mathcal{H}(\varphi) := \{\mu \in \Lambda\, |\, 
\mathrm{ht}(\mu) = 
\mathrm{ht}(\varphi) \mbox{ and } [e^{\mu}](\varphi) \not= 0\}$,  
and $\widehat{\Pi}(\varphi) := \{\mu\, |\, 
[e^{\mu}](\varphi) \not= 0\}$. 
The following observations 
will be useful. Take $\chi \in 
\mathbb{Z}[\Lambda]^{W}$ with  
$\chi \not= 0$.  Let $\mu \in \Lambda$ with $[e^{\mu}](\chi) \not= 0$.  
Then $\mu = w(\nu)$ for some $w \in W$ and $\nu 
\in \Lambda^{+}$ (Lemma 13.2.A of \cite{Hum}), so $[e^{\nu}](\chi) = 
[e^{\mu}](\chi)$. 
This shows that $\mathrm{ht}(\chi) \geq 0$ and that $\lambda \in 
\mathcal{H}(\chi)$ only if $\lambda \in \Lambda^{+}$.  

{\bf The subgroup of Weyl alternants.} 
Besides the ring of Weyl symmetric functions, the group ring contains 
another important subgroup, the subgroup 
$\mathbb{Z}[\Lambda]^{\mbox{\tiny alt}} := \{\varphi \in 
\mathbb{Z}[\Lambda]\, |\, w.\varphi = 
\det(w)\varphi \mbox{ for all } 
w \in W\}$ consisting of the 
{\em Weyl alternants} of $\mathbb{Z}[\Lambda]$. 
Note that the sum  of two Weyl alternants 
is alternating but their product is $W$-invariant. 
Say $\varphi = \sum_{\mu \in \Lambda}f_{\mu}e^{\mu}$ is a Weyl alternant.  
It is easy to see that $f_{w(\mu)} = \det(w)f_{\mu}$.  
The same 
reasoning from the end of the previous paragraph shows that if 
$\varphi \not= 0$, then $\mathrm{ht}(\varphi) \geq 0$ and that   
$\lambda \in \mathcal{H}(\varphi)$ only if $\lambda \in \Lambda^{+}$.  
We will say more about such $\varphi$ in the next result. 
Define a mapping 
$\mathcal{A}: \mathbb{Z}[\Lambda] \longrightarrow 
\mathbb{Z}[\Lambda]^{\mbox{\tiny alt}}$ 
by the rule $\mathcal{A}(\varphi) := \sum_{w 
\in W}\det(w) w.\varphi$. 
The following fundamental facts are useful for 
the subsequent discussion. 

\noindent 
{\bf \AltProp}\ \ {\sl 
(1)  Let $\lambda$ be dominant. Then $\mathcal{A}(e^{\lambda}) \not= 
0$ if and only if $\lambda$ is strongly dominant, in which case  
we have $\widehat{\Pi}(\mathcal{A}(e^{\lambda})) = W\lambda$ (the 
orbit of $\lambda$ under the $W$-action on $\Lambda$),  
$[e^{\lambda}](\mathcal{A}(e^{\lambda})) = 1$, 
$\mathrm{ht}(\mathcal{A}(e^{\lambda})) = 
\lambda$, and 
$\mathcal{H}(\mathcal{A}(e^{\lambda})) = \{\lambda\}$. 
(2)  The set   
$\{\mathcal{A}(e^{\lambda})\}_{\lambda \in \Lambda^{++}}$ is a 
$\mathbb{Z}$-basis for $\mathbb{Z}[\Lambda]^{\mbox{\tiny alt}}$. 
Moreover, if $\varphi = \sum_{\mu \in 
\Lambda}f_{\mu}e^{\mu}$ is a Weyl alternant, then $\varphi = 
\sum_{\lambda \in 
\Lambda^{++}}f_{\lambda}\mathcal{A}(e^{\lambda})$.}

{\em Proof.} For {\sl (1)}, begin with $\lambda \in \Lambda^{++}$.  
From the definition of 
$\mathcal{A}(e^{\lambda})$ it follows that 
$\widehat{\Pi}(\mathcal{A}(e^{\lambda})) 
\subseteq W\lambda$.  Now $e^{w_{1}(\lambda)} = 
e^{w_{2}(\lambda)}$ if and only if 
$w_{2}^{-1}w_{1}(\lambda) = \lambda$.  By Lemma 13.2.A of 
\cite{Hum}, we must have $w_{2}^{-1}w_{1} = Id$, i.e.\ 
$w_{1} = w_{2}$.  So, 
$[e^{w(\lambda)}](\mathcal{A}(e^{\lambda})) = \det(w)$ for 
all $w \in W$, hence $W\lambda \subseteq 
\widehat{\Pi}(\mathcal{A}(e^{\lambda}))$.  From this reasoning the 
remaining claims of the ``if'' direction of {\sl (1)} follow. 
Now assume $\lambda \in \Lambda^{+} \setminus \Lambda^{++}$. 
We will view $W$ again as a 
Coxeter group.  Then as observed earlier, the stablizer 
$W_{\lambda} 
\subseteq W$ is 
just the Coxeter subgroup $W_{J}$, where $J = \{i \in I\, |\, 
\langle \lambda,\alpha_{i}^{\vee} \rangle = 0\}$.  Since 
$\lambda \in \Lambda^{+} \setminus \Lambda^{++}$, then $J \not= 
\emptyset$. With $W^{J}$ denoting 
the corresponding set of ``minimal coset representatives'' 
(see e.g.\ \S 5.12 of 
\cite{HumCoxeter}), then each $w \in W$ is uniquely 
expressible as $w = w^{J}w_{J}$, where $w^{J} \in 
W^{J}$ and $w_{J} \in W_{J}$.  So, $\mathcal{A}(e^{\lambda}) = 
\sum_{w \in W}\det(w)e^{w(\lambda)} 
= \sum_{w_{J} \in 
W_{J}}\det(w_{J})\Big(\sum_{w^{J} \in 
W^{J}}\det(w^{J})e^{w^{J}(\lambda)}\Big) = 
\myqSsc\sum_{w_{J} \in 
W_{J}}\det(w_{J})$, where  $\myqSsc = \sum_{w^{J} \in 
W^{J}}\det(w^{J})e^{w^{J}(\lambda)}$.  Now let $\sgn_{J}: 
W_{J} \longrightarrow \{\pm 1\}$ be given by $\sgn_{J}(w_{J}) = 
\det(w_{J})$ for all $w_{J} \in W_{J}$.  Then $\sgn_{J}$ is 
a surjective group homomorphism (since $J \not= \emptyset$), hence 
$|\sgn^{-1}_{J}(-1)| = |\sgn^{-1}_{J}(1)| = \frac{1}{2}|W_{J}|$.  
Then $\myqSsc\sum_{w_{J} \in 
W_{J}}\det(w_{J}) = \myqSsc\Big(\sum_{w_{J} \in 
\sgn^{-1}_{J}(1)}(1) + \sum_{w_{J} \in 
\sgn^{-1}_{J}(-1)}(-1)\Big) = \myqSsc(\frac{1}{2}|W_{J}| - 
\frac{1}{2}|W_{J}|) = 0$.  This completes the proof of {\sl (1)}. 

For {\sl (2)}, first take $\varphi = \sum_{\lambda \in 
\Lambda^{++}}f_{\lambda}\mathcal{A}(e^{\lambda})$, 
and suppose some $f_{\nu} 
\not= 0$.  
Now $[e^{\nu}](\varphi) = \sum \det(w)f_{\lambda}$, 
where the sum is over all 
$\lambda \in \Lambda^{++}$ and 
$w \in W$ such that $w(\lambda) = \nu$.  
But by Lemma 13.2.A of \cite{Hum}, this can only happen when 
$\lambda = \nu$ and $w = Id$.  That is, 
$[e^{\nu}](\varphi) = f_{\nu} \not= 0$. It follows that the set 
$\{\mathcal{A}(e^{\lambda})\}_{\lambda \in \Lambda^{++}}$ is 
independent. 
Now say $\varphi = \sum_{\mu \in \Lambda} f_{\mu}e^{\mu}$ 
is any element of 
$\mathbb{Z}[\Lambda]^{\mbox{\tiny alt}}$.  We may write 
$\varphi = \sum_{\lambda \in \Lambda^{++}}\sum_{w \in 
W}f_{w(\lambda)}e^{w(\lambda)} + \sum f_{\mu}e^{\mu}$, 
where the latter sum is over all $\mu$ such that $\mu = w(\nu)$ 
for some $w \in W$ and $\nu \in \Lambda^{+} \setminus 
\Lambda^{++}$.  For such a $\mu = w(\nu)$, 
write $w = w^{J}w_{J}$ as in the previous paragraph, 
where $W_{\nu} = W_{J}$ is the stablizer of $\nu$.  
Let $w' := w^{J}s_{j}$ for some 
$j \in J$ (possible since $J \not= \emptyset$).  
Since $\varphi$ is alternating, we have both 
$f_{\mu} = f_{w^{J}(\nu)} = \det(w^{J})f_{\nu}$ and 
$f_{\mu} = f_{w^{J}s_{j}(\nu)} = -\det(w^{J})f_{\nu}$.  It 
follows that $f_{\nu} = 0$, hence $f_{\mu} = 0$ also.  
Thus $\varphi = \sum_{\lambda \in \Lambda^{++}}\sum_{w \in 
W}f_{w(\lambda)}e^{w(\lambda)} = \sum_{\lambda \in 
\Lambda^{++}} f_{\lambda}\mathcal{A}(e^{\lambda})$, 
which shows that the set 
$\{\mathcal{A}(e^{\lambda})\}_{\lambda \in \Lambda^{++}}$ spans 
$\mathbb{Z}[\Lambda]^{\mbox{\tiny alt}}$ and also shows that any 
$\varphi \in \mathbb{Z}[\Lambda]^{\mbox{\tiny alt}}$ 
can be written as a $\mathbb{Z}$-linear combination of the 
$\mathcal{A}(e^{\lambda})$'s as claimed.\hfill\QED 

The next result is called ``Weyl's denominator formula'' for reasons that will be made apparent shortly.  It is a generalization of the classical Vandermonde determinant.  

\noindent 
{\bf \WeylsDenomTheorem\ (Weyl's denominator formula)}\ \ {\sl  We have the following identities of Weyl alternants:} 
\begin{equation}
\mathcal{A}(e^{\varrho}) = 
e^{\varrho}\Big(\prod_{\alpha \in \Phi^{+}}(1-e^{-\alpha})\Big) = 
\prod_{\alpha \in \Phi^{+}}(e^{\alpha/2}-e^{-\alpha/2}) = 
e^{-\varrho}\Big(\prod_{\alpha 
\in \Phi^{+}}(e^{\alpha}-1)\Big).
\end{equation} 

{\em Proof.} Equality of the latter 
three expressions in equation (2) 
follows from the definitions.  Now let $\mathcal{D} := 
e^{\varrho}\Big(\prod_{\alpha \in \Phi^{+}}(1-e^{-\alpha})\Big)$.  In 
view of the fact that $s_{i}$ permutes the positive roots other 
than $\alpha_{i}$ and sends $\alpha_{i}$ to $-\alpha_{i}$ (cf.\ Lemma 
10.2.B of \cite{Hum}), 
then $s_{i}.\mathcal{D} = -\mathcal{D}$.  It 
follows that $\mathcal{D}$ is alternating.  Moreover, observe that 
$\mathcal{H}(\mathcal{D}) = \{\varrho\}$ since 
$[e^{\mu}](\mathcal{D}) \not= 0$ only if $\mu = \varrho - 
\sum_{\alpha \in S}\alpha$, where $S$ is some subset of $\Phi^{+}$. 
Also, this reasoning shows that $[e^{\varrho}](\mathcal{D}) = 1$, 
since $\varrho = \varrho - 
\sum_{\alpha \in S}\alpha$ only when $S = \emptyset$. 
Then by \AltProp.2, it follows that $\mathcal{D} = 
\mathcal{A}(e^{\varrho})$.\hfill\QED

{\bf Weyl bialternants and $\mathbf{\Phi}$-Kostka numbers.} We will shortly 
produce three distinguished bases for the ring of $W$-symmetric functions.  
The special basis of ``Weyl 
bialternants'' for $\mathbb{Z}[\Lambda]^{W}$ is the most important of 
these and is a source of much 
interesting combinatorics. 

\noindent 
{\bf \WeylBialternantDef}\ \ By Weyl's denominator formula, \[\frac{1}{\mathcal{A}(e^{\varrho})} = \frac{1}{e^{\varrho}\Big(\prod_{\alpha \in \Phi^{+}}(1-e^{-\alpha})\Big)} = e^{-\varrho}\prod_{\alpha 
\in \Phi^{+}}\sum_{k \geq 0}e^{-k\alpha}\] 
is a well-defined member of $\mathcal{R}$.  For each $\lambda \in \Lambda^{+}$, define 
\begin{equation}
\chi_{_{\lambda}} := \frac{\mathcal{A}(e^{\lambda+\varrho})}{\mathcal{A}(e^{\varrho})} =
\Bigg(\sum_{w \in 
W}\det(w)e^{w(\lambda+\varrho)}\Bigg)
\Bigg(e^{-\varrho}\prod_{\alpha 
\in \Phi^{+}}\sum_{k \geq 0}e^{-k\alpha}\Bigg), 
\end{equation} 
which is also a well-defined member of $\mathcal{R}$.  We call each $\chi_{_{\lambda}}$ {\em Weyl bialternant}.\hfill\QED  

As a quotient of Weyl alternants, it would seem obvious that $\chi_{_{\lambda}}$ is a Weyl symmetric function.  
However, we do not yet know that $\chi_{_{\lambda}}$ is in the group ring where the action of the Weyl group is defined,\myfootnote{It seems this subtle issue is not explicitly considered in the plausibility arguments of \S 24.1 of \cite{FH}.} and it is not possible to naturally extend the $W$-action on $\mathbb{Z}[\Lambda]$ to $\mathcal{R}$. 
This finiteness issue is resolved in \FinitenessTheorem\ below. 

Weyl bialternants are sometimes called Weyl characters in part because they are characters of irreducible representations of Lie groups/algebras.  
Weyl bialternants are plausibly viewed as analogs of Schur functions, as the quotient of alternating sums in (3) analogizes the classical definition of a Schur function as a quotient of determinants.  
(For more on this, see \CaseAExample\ below.) 
These somewhat ``external'' considerations give a two-fold motivation for \WeylBialternantDef. 
However, it would be nice to have a more attributive combinatorial characterization of Weyl bialternants. 
We pose this as our first problem. 

\noindent 
{\bf \WeylBialternantProblem}\ \ Find a nice set of combinatorial properties that characterize Weyl bialternants. 
One solution to this problem, stated in the language of \S\S \SplittingSection--\CrystallineSection\ of this monograph, describes $\chi_{_{\lambda}}$ as the weight generating function for a strongly untangled crystalline splitting poset whose unique maximal element has weight $\lambda$. 
But, this description requires a good bit of overhead and is not obviously intrinsic. 
So, we would like a characterizing set of properties that is commensurate with the discussion up to this point in our presentation.\hfill\QED 

Our initial results about Weyl bialternants do not require knowledge of their finiteness as members of $\mathcal{R}$.  The following fundamental observation follows immediately from the definition of Weyl bialternant. 

\noindent 
{\bf \InitialTheorem}\ \ {\sl Let $\lambda$ be a dominant weight.  Then}  
\begin{equation}
\mathcal{A}(e^{\varrho})\,\chi_{_{\lambda}} = 
\mathcal{A}(e^{\lambda+\varrho}). 
\end{equation} 
{\sl is a defining identity for $\chi_{_{\lambda}}$ in the sense that $\chi = \chi_{_{\lambda}}$ is the unique solution in $\mathcal{R}$ to the equation} $\mathcal{A}(e^{\varrho})\,\chi = \mathcal{A}(e^{\lambda+\varrho})$.\hfill\QED 

For any $\mu \in \Lambda$, let  $d_{\lambda,\mu} := [e^{\mu}](\chi_{_{\lambda}})$; then one can view $\chi_{_{\lambda}}$ as a generating function for the numbers $d_{\lambda,\mu}$.  A strong hint that these 
numbers are combinatorially interesting is that in the $\myA_{n}$ 
case, each $d_{\lambda,\mu}$ is just the Kostka number 
$K_{\mathfrak{p},\mathfrak{q}}$ for certain partitions $\mathfrak{p}, 
\mathfrak{q}$ naturally associated with $\lambda, \mu$ (see 
\CaseAExample).  So we will also refer to these Weyl bialternant expansion coefficients 
as $\Phi$-{\em Kostka numbers}.  We will shortly write down an explicit formula for the $\Phi$-Kostka numbers (\KostantTheorem) 
and a recurrence (\FreudTheorem). 
It is not obvious that the integers $d_{\lambda,\mu}$ are non-negative; this issue is settled by \KostkaCorollary\ below. 
As a formula for weight multiplicities for irreducible representations of semisimple Lie algebras, equation (5) below is due to Kostant \cite{Kostant}, and within that Lie theoretic context it is known as ``Kostant's multiplicity formula'' (KMF).  Here it is an easy consequence of the definitions. 

\noindent 
{\bf \KostantTheorem}\ \ {\sl With $\lambda \in \Lambda^{+}$ and  $\mu \in \Lambda$, then the $\Phi$-Kostka number $d_{\lambda,\mu}$ can be written as:} 
\begin{equation}
d_{\lambda,\mu} = \sum_{w \in 
W}\det(w)\mathcal{P}(w(\lambda)-\mu + 
w(\varrho)-\varrho).\end{equation} 

\vspace*{-0.35in}\hfill\QED

{\bf Other (potential) bases for the ring of $\mathbf{W}$-symmetric functions.} 
Besides the Weyl bialternants $\{\chi_{_{\lambda}}\}$,  
two other distinguished sets naturally 
indexed by dominant weights 
are the {\em monomial $W$-symmetric functions} $\{\zeta_{\lambda}\}$ 
and {\em elementary $W$-symmetric functions} $\{\psi_{\lambda}\}$.  
Fix $\lambda = 
\sum_{i \in I}a_{i}\omega_{i} \in \Lambda^{+}$.  The monomial $W$-symmetric function  
$\zeta_{\lambda} := \frac{1}{|W_{\lambda}|}\sum_{w \in 
W}e^{w(\lambda)}$. 
Observe that for any $\mu \in \Lambda$, $[e^{\mu}](\zeta_{\lambda}) = 
\frac{1}{|W_{\lambda}|}|\{w \in W\, |\, w(\lambda) = 
\mu\}|$ is nonzero if and only if $\mu \in W\lambda$ (the $W$-orbit 
of $\lambda$) if and only if $[e^{\mu}](\zeta_{\lambda}) = 1$. 
In particular, $\widehat{\Pi}(\zeta_{\lambda}) = W\lambda$. 
Now, for each $i \in I$, we call $\chi_{_{\omega_{i}}}$ a {\em fundamental 
bialternant}.  An elementary $W$-symmetric function is a monomial of fundamental bialternants.  In particular, 
we let $\psi_{\lambda} := 
\chi_{_{\omega_{1}}}^{a_{1}}{\cdots}\chi_{_{\omega_{n}}}^{a_{n}}$ 
with $\psi_{0} := e^{0}$ by convention. 
Observe that $[e^{\mu}](\psi_{\lambda}) = 
\sum\Big(\prod_{i \in I}\prod_{j=1}^{a_{i}} 
d_{\omega_{i},\mu_{i}^{(j)}}\Big)$, where this sum is over all 
ways of 
writing $\mu = \sum_{i \in I}\sum_{j=1}^{a_{i}} \mu_{i}^{(j)}$ such 
that each $d_{\omega_{i},\mu_{i}^{(j)}} \not= 0$.  
Another observation: 

\noindent 
{\bf \WeylCharLemma}\ \ {\sl Let $\lambda \in \Lambda^{+}$ and 
$\varphi \in \{\chi, \zeta, \psi\}$. If $\mu \in 
\CharPi(\varphi_{\lambda})$ then $\mu \leq \lambda$.  
Moreover, $[e^{\lambda}](\varphi_{\lambda}) = 1$, 
$\mathrm{ht}(\varphi_{\lambda}) = \mathrm{ht}(\lambda)$, and 
$\mathcal{H}(\varphi_{\lambda}) = \{\lambda\}$.} 

{\em Proof.} First take $\varphi = \chi$.  
Since for all $w \not= Id \in 
W$ we have $w(\lambda) < \lambda$ and $w(\varrho) < 
\varrho$ in $\Lambda$ (cf.\ Lemma 13.2.A of \cite{Hum}), then by 
equation (4) above 
$d_{\lambda,\lambda} = 1$. Moreover, if $d_{\lambda,\mu} \not= 0$, 
then for some nonnegative integers $k_{\alpha}$ and some $w \in 
W$ we have $\sum_{\alpha \in 
\Phi^{+}}k_{\alpha}\alpha = 
w(\lambda)-\mu+w(\varrho)-\varrho$, and hence 
$\mu \leq w(\lambda) \leq \lambda$.  So, 
$\mathrm{ht}(\chi_{_{\lambda}}) = \mathrm{ht}(\lambda)$ and 
$\mathcal{H}(\chi_{_{\lambda}}) = \{\lambda\}$.  When $\varphi = 
\zeta$, 
the conclusions of the lemma statement all follow from the 
observations about $\zeta_{\lambda}$ in the paragraph preceding the 
lemma statement, in view of the fact that $w(\lambda) \leq 
\lambda$ for all $w \in W$. Lastly take 
$\varphi = \psi$.  Suppose $[e^{\mu}](\psi_{\lambda}) \not= 0$, 
with 
$\lambda = \sum_{i \in I}a_{i}\omega_{i}$.  
In the notation of the paragraph preceding the 
lemma, we have $\prod_{i \in I}\prod_{j=1}^{a_{i}} 
d_{\omega_{i},\mu_{i}^{(j)}} \not= 0$ only if each 
$\mu_{i}^{(j)} \leq \omega_{i}$. But then 
$\mu = \sum_{i \in I}\sum_{j=1}^{a_{i}} \mu_{i}^{(j)} \leq \sum_{i \in 
I}a_{i}\omega_{i} = \lambda$.  This is an equality if and only if 
each $\mu_{i}^{(j)} = \omega_{i}$, whence 
$[e^{\lambda}](\psi_{\lambda}) = \prod_{i \in I}\prod_{j=1}^{a_{i}} 
d_{\omega_{i},\omega_{i}} = 1$. So, 
$\mathrm{ht}(\psi_{\lambda}) = \mathrm{ht}(\lambda)$ and 
$\mathcal{H}(\psi_{\lambda}) = \{\lambda\}$.\hfill\QED 

{\bf Weyl bialternants, monomial $W$-symmetric functions, and elementary $W$-symmetric functions for reducible roots systems.} 
For the record, here is what $\chi_{_{\lambda}}$, $\zeta_{\lambda}$, 
and $\psi_{\lambda}$ look like when $\Phi$ is  reducible. 
Say $\Phi = \Phi_{1} \disjointunion \cdots \disjointunion \Phi_{k}$ is reducible with each $\Phi_{i}$ irreducible, with lattice of weights $\Lambda$ 
written $\Lambda_{1} \oplus \cdots \oplus \Lambda_{k}$ and $W \cong 
W_{1} \times \cdots \times W_{k}$ as before. 
Let $\mathcal{R}_{i}$ enlarge $\mathbb{Z}[\Lambda_{i}]$ as above, and 
see that there are natural isomorphisms $\mathcal{R} \cong 
\mathcal{R}_{i} \times \cdots \times \mathcal{R}_{k}$, 
$\mathbb{Z}[\Lambda] \cong \mathbb{Z}[\Lambda_{1}] \times \cdots 
\times \mathbb{Z}[\Lambda_{k}]$, and $\mathbb{Z}[\Lambda]^{W} \cong 
\mathbb{Z}[\Lambda_{1}]^{W_{1}} \times \cdots \times 
\mathbb{Z}[\Lambda_{k}]^{W_{k}}$. 
Write $\lambda \in \Lambda^{+}$ as $\lambda = \lambda_{1} + \cdots 
\lambda_{k}$ with each $\lambda_{i} \in \Lambda_{i}^{+}$.  
It follows from the definitions that 
the $\Phi$-Weyl bialternant $\chi^{\Phi}_{_{\lambda}}$ 
is $\chi^{\Phi_{1}}_{_{\lambda_{1}}}\cdots
\chi^{\Phi_{k}}_{_{\lambda_{k}}}$, 
where $\chi^{\Phi_{i}}_{_{\lambda_{i}}}$ is the $\Phi_{i}$-Weyl 
bialternant corresponding to the dominant weight $\lambda_{i} \in 
\Lambda_{i}^{+}$. In a similar way, $\zeta_{\lambda}^{\Phi} = 
\prod_{i=1}^{k}(\zeta_{\lambda_{i}}^{\Phi_{i}})$ and 
$\psi_{\lambda}^{\Phi} = 
\prod_{i=1}^{k}(\psi_{\lambda_{i}}^{\Phi_{i}})$. 

{\bf Finiteness of Weyl bialternants, and some consequences.} 
Our definition of Weyl bialternant leaves open the question of why this member of $\mathcal{R}$ is also in the group ring $\mathbb{Z}[\Lambda]$.  
This result is an obvious consequence of the following theorem. 
Of course, this result depends on the finiteness of the root system $\Phi$.  
This is accounted for in the proof with our use of the mesh size estimate $\myd_{\Phi}$, which has the property that for all $\nu \in \Lambda$, $\langle \nu,\varrho^{\vee} \rangle > 0$ if and only if $\langle \nu,\varrho^{\vee} \rangle \geq \myd_{\Phi}$. 

\noindent 
{\bf \FinitenessTheorem\ (Finiteness Theorem)}\ \ {\sl Let $\nu \in \Lambda^{+}$.  
Then $\displaystyle \frac{1}{\mathcal{A}(e^{\nu+\varrho})}$ is a well-defined member of $\mathcal{R}$.  
Let $\phi \in \mathrm{span}_{\mathbb{Z}}\{\mathcal{A}(e^{\lambda+\nu+\varrho})\}_{\lambda \in \Lambda^{+}}$, so $\phi$ is a Weyl alternant. 
Let $\displaystyle \chi := \frac{\phi}{\mathcal{A}(e^{\nu+\varrho})}$, so $\chi$ is the unique member of $\mathcal{R}$ such that $\mathcal{A}(e^{\nu+\varrho})\, \chi = \phi$.  
Then $\chi$ is in $\mathbb{Z}[\Lambda]^{W}$.}  

{\em Proof.} In \S 2 of \cite{Stem}, Stembridge observes that the multivariate ring of formal power series $\mathbb{Z}[[e^{-\alpha_i}\, |\, i \in I]]$ is a subring of $\mathcal{R}$.  
So now consider $e^{-(\nu+\varrho)} \mathcal{A}(e^{\nu+\varrho}) = \sum_{w \in W} \det(w)e^{w(\nu+\varrho)-(\nu+\varrho)}$.  
Since $w(\nu+\varrho) \leq \nu+\varrho$ for all $w \in W$, then the latter sum of the previous sentence is in $\mathbb{Z}[[e^{-\alpha_i}\, |\, i \in I]]$.  
Also, since $w(\nu+\varrho) = \nu+\varrho$ if and only if $w$ is the identity in $W$, it follows that $[e^0](e^{-(\nu+\varrho)} \mathcal{A}(e^{\nu+\varrho})) = 1$.  
In particular, the constant term of $e^{-(\nu+\varrho)} \mathcal{A}(e^{\nu+\varrho})$ is a unit in $\mathbb{Z}$.  
By standard facts about multivariate rings of formal power series, it follows that $e^{-(\nu+\varrho)} \mathcal{A}(e^{\nu+\varrho})$ has an inverse $f$ in $\mathbb{Z}[[e^{-\alpha_i}\, |\, i \in I]]$, so $e^{\nu+\varrho}f \in \mathcal{R}$ is an inverse for $\mathcal{A}(e^{\nu+\varrho})$ in $\mathcal{R}$. 

Here is how we can conclude that $\chi$ is $W$-invariant if we can show that $\chi$ is in the group ring $\mathbb{Z}[\Lambda]$:  
Apply $w \in W$ to both sides of the equation $\mathcal{A}(e^{\nu+\varrho})\, \chi = \phi$.  
Then get $\det(w)\mathcal{A}(e^{\nu+\varrho})\, w.\chi = \det(w)\phi$, whence $\mathcal{A}(e^{\nu+\varrho})\, w.\chi = \phi$.  
Therefore $w.\chi = \chi$ for all $w \in W$.  
Note that this argument does not work if we only know $\chi \in \mathcal{R}$ since the $W$-action on the group ring does not naturally extend to $\mathcal{R}$.  

To complete the proof it suffices to show that $\phi$ is divisible by $\mathcal{A}(e^{\nu+\varrho})$.  
If $\phi = 0$, the result is trivial.  
So suppose $\phi \not= 0$.  
We will prescribe an algorithm for rewriting $\phi$ that will terminate in an expression for $\phi$ that is divisible by $\mathcal{A}(e^{\nu+\varrho})$.  

As our first step, set $\phi^{(0)} := \phi$, and write $\displaystyle \phi^{(0)} := \sum_{\lambda \in \Lambda^{+}}f_{\lambda}^{(0)} \mathcal{A}(e^{\lambda+\nu+\varrho})$ for integers $\{f_{\lambda}^{(0)}\}_{\lambda \in \Lambda^{+}}$.  
Now let \[\phi^{(1)} := \phi^{(0)} - \sum_{\mbox{\tiny $\begin{array}{c}\lambda \in \Lambda^{+} \mbox{ with}\\ \lambda+\nu+\varrho \in \mathcal{H}(\phi^{(0)})\end{array}$}}f_{\lambda}^{(0)} \zeta_{\lambda} \mathcal{A}(e^{\nu+\varrho}).\]  
If $\phi^{(1)} = 0$, then we stop, since now $\displaystyle \phi = \phi^{(0)} = \sum_{\mbox{\tiny $\begin{array}{c}\lambda \in \Lambda^{+} \mbox{ with}\\ \lambda+\nu+\varrho \in \mathcal{H}(\phi^{(0)})\end{array}$}}f_{\lambda}^{(0)} \zeta_{\lambda} \mathcal{A}(e^{\nu+\varrho})$, which is clearly divisible by $\mathcal{A}(e^{\nu+\varrho})$.  
Otherwise, if $\phi^{(1)} \not= 0$, then note that $\phi = \phi^{(0)}$ is divisible by $\mathcal{A}(e^{\nu+\varrho})$ if $\phi^{(1)}$ is.  
Now for $i = 0,1$, choose $\lambda_{i}$ dominant and such that $\lambda_{i}+\nu+\varrho \in \mathcal{H}(\phi^{(i)})$, so that $\mathrm{ht}(\phi^{(i)}) = \langle \lambda_{i}+\nu+\varrho,\varrho^{\vee}\rangle > 0$.  
Now clearly $\mathrm{ht}(\phi^{(0)}) > \mathrm{ht}(\phi^{(1)})$, so therefore $\mathrm{ht}(\phi^{(0)}) - \mathrm{ht}(\phi^{(1)}) = \langle \lambda_{0}-\lambda_{1},\varrho^{\vee}\rangle > 0$.  
But this means that $\langle \lambda_{0}+\lambda_{1},\varrho^{\vee}\rangle \geq \myd_{\Phi}$.  
That is, $\mathrm{ht}(\phi^{(0)}) \geq \mathrm{ht}(\phi^{(1)}) + \myd_{\Phi}$. 

Repeat this procedure, where at the $k$th step ($k \geq 1$) we have \[\phi^{(k)} = \phi^{(k-1)} - \sum_{\mbox{\tiny $\begin{array}{c}\lambda \in \Lambda^{+} \mbox{ with}\\ \lambda+\nu+\varrho \in \mathcal{H}(\phi^{(k-1)})\end{array}$}}f_{\lambda}^{(k-1)} \zeta_{\lambda} \mathcal{A}(e^{\nu+\varrho}).\]  
The procedure terminates at this step if $\phi^{(k)} = 0$, because then $\phi^{(k-1)}$ is divisible by $\mathcal{A}(e^{\nu+\varrho})$ and so, in turn, are each of $\phi^{(k-2)}, \ldots , \phi^{(1)},$ and $\phi^{(0)} = \phi$.  
On the other hand, if $\phi^{(k)} \not= 0$, then as before we obtain that $\mathrm{ht}(\phi^{(k-1)}) \geq \mathrm{ht}(\phi^{(k)}) + \myd_{\Phi}$ and that $\phi^{(k-1)}$ is divisible by $\mathcal{A}(e^{\nu+\varrho})$ if $\phi^{(k)}$ is.  
Then we repeat the procedure for a $(k+1)$st step.  At the $k$th step we know that $\mathrm{ht}(\phi) = \mathrm{ht}(\phi^{(0)}) \geq \mathrm{ht}(\phi^{(1)}) + \myd_{\Phi} \geq \mathrm{ht}(\phi^{(2)}) + 2\myd_{\Phi} \geq \cdots \geq \mathrm{ht}(\phi^{(k)}) + k\myd_{\Phi}$, so therefore the procedure must terminate in no more than $\lfloor \mathrm{ht}(\phi)/\myd_{\Phi} \rfloor$ steps.\hfill\QED  

\noindent 
{\bf \FinitenessCorollary}\ \ {\sl Let $\lambda \in \Lambda^{+}$.  (1) Then we have $\chi_{_{\lambda}} \in \mathbb{Z}[\Lambda]^{W}$, so in particular for all $\nu \in \Lambda$ and $w \in W$ we have $d_{\lambda,w(\nu)} = d_{\lambda,\nu}$.  (2) Also, $\displaystyle \chi_{_{\lambda}} = \sum_{\mbox{\tiny $\begin{array}{c}\nu \in \Lambda^{+}\\ \mbox{with } \nu \in \Pi(\lambda)\end{array}$}} d_{\lambda,\nu} \zeta_{\nu}$, and $\CharPi(\chi_{_{\lambda}}) \subseteq \Pi(\lambda)$.  In particular, $d_{\lambda,\mu} \not= 0$ only if $\mu \in \Pi(\lambda)$.}

{\em Proof.} In view of \FinitenessTheorem, {\sl (1)} is obvious.  For {\sl (2)}, note that for each 
$\mu \in \Lambda$, 
there is a unique 
$\nu \in \Lambda^{+}$ for which $\mu = w(\nu)$ 
(Lemma 13.2.A of 
\cite{Hum}). From this fact and the fact that $d_{\lambda,\nu} = 
d_{\lambda,w(\nu)}$ for all $w \in W$ and $\nu \in 
\Lambda^{+}$ we obtain: $\displaystyle 
\chi_{_{\lambda}} = \sum_{\mu \in 
\Lambda}d_{\lambda,\mu}e^{\mu} = \sum_{\nu \in 
\Lambda^{+}}d_{\lambda,\nu}\Bigg(\frac{1}{|W_{\nu}|}\sum_{w 
\in W}e^{w(\nu)}\Bigg) = \sum_{\nu \in 
\Lambda^{+}}d_{\lambda,\nu}\zeta_{\nu}$. By \WeylCharLemma, 
$d_{\lambda,\nu} \not= 0$ only if $\nu \leq \mu$. Then 
$\displaystyle \chi_{_{\lambda}} = \sum_{\nu 
\in \Lambda^{+},\nu \leq 
\lambda}d_{\lambda,\nu}\zeta_{\nu}$. Therefore, 
$\displaystyle \CharPi(\chi_{_{\lambda}}) \subseteq \bigcup_{\nu \in \Lambda^{+},\nu \leq \lambda} \CharPi(\zeta_{\nu}) = \bigcup_{\nu \in \Lambda^{+},\nu \leq \lambda} W\nu = \Pi(\lambda)$, from which it follows that $d_{\lambda,\mu} \not= 0$ only if $\mu \in \Pi(\lambda)$.\hfill\QED

{\bf A recurrence for the $\mathbf{\Phi}$-Kostka numbers.} 
The following recurrence for the 
$\Phi$-Kostka numbers as weight multiplicities for semisimple Lie 
algebra representations is due to Freudenthal, and in that context it is known as ``Freudenthal's multiplicity formula'' (FMF). 
We derive it here without the context of Lie algebra 
representation theory.   

\noindent 
{\bf \FreudTheorem}\ \ {\sl Let $\lambda \in \Lambda^{+}$.  Then for 
all $\mu \in \Lambda$, we have:} 
\begin{equation}
(||\lambda+\varrho||^{2} - ||\mu+\varrho||^{2}) d_{\lambda,\mu} = 
2\sum_{\alpha \in \Phi^{+}}\sum_{k \geq 1}\langle 
\alpha,\mu+k\alpha \rangle d_{\lambda,\mu+k\alpha}. 
\end{equation}
{\sl The numbers $\{d_{\lambda,\mu}\}_{\mu \in \Lambda}$ are uniquely determined by the recurrence relation (6) above together with the conditions that $d_{\lambda,\mu} = 0$ when $\mu \not\in \Pi(\lambda)$ and that $d_{\lambda,\lambda} = 1$.}

The proof we give below closely follows the proof of a $q$-analog 
of equation (6) derived by Lansky in \cite{Lansky}, although 
at $q=1$ the 
argument simplifies a bit.  As an alternative, one can simply take 
$q=1$ in Lansky's $q$-analog, but some work is required to 
simplify that result to get equation (6). 
The argument given below 
nicely illustrates a 
typical generating function 
technique, namely ``differentiating'' both sides of the identity (4) 
that defines $\chi_{_{\lambda}}$.  

Before giving the proof, 
we note precisely in what sense the above 
formula recursively determines the $\Phi$-Kostka numbers.  
By \WeylCharLemma, $d_{\lambda,\mu} = 0$ if $\mu \not\leq \lambda$, 
so on the right-hand side of (6) there are at most a finite number of 
weights of the form $\mu+k\alpha$ such that $d_{\lambda,\mu+k\alpha} 
\not= 0$. Also, for any $\mu \in \Pi(\lambda)$ 
with $\mu \not=\lambda$, 
it is the case that 
$||\lambda+\varrho||^{2} - ||\mu+\varrho||^{2} > 0$ by Lemma 
13.4.C of \cite{Hum}.  
Moreover, by \WeightRemarkOne.2, we have $\mu < \lambda$, so $\lambda - \mu = \sum_{i \in I} k_{i}\alpha_{i}$ for some nonnegative integers $k_{i}$ which are not all zero.  For such $\mu$, declare the ``depth'' of $\mu$ in $\Pi(\lambda)$ to be the quantity 
\[\mathrm{depth}(\mu) := \sum_{i \in I} k_{i} = \langle \lambda-\mu,\varrho^{\vee} \rangle = \mathrm{ht}(\lambda) - \mathrm{ht}(\mu),\] 
which is a positive integer since $\mu < \lambda$.  Of course, for each $\alpha \in \Phi^{+}$ and each integer $k \geq 1$, if $\mu + k\alpha \in \Pi(\lambda)$, then $\mathrm{depth}(\mu + k\alpha) < \mathrm{depth}(\mu)$.  Then from the starting point of $d_{\lambda,\lambda} = 1$ (cf.\ \WeylCharLemma), we can therefore induct on $\mathrm{depth}(\mu)$ and use (6) to determine $d_{\lambda,\mu}$ for all $\mu \in \Pi(\lambda)$. So, (6) uniquely determines those 
$d_{\lambda,\mu}$  for which $\mu \in \Pi(\lambda)$ and $\mu 
\not=\lambda$.   
However, the recursive procedure breaks down when $\mu < \lambda$ 
with 
$||\lambda+\varrho||^{2} - ||\mu+\varrho||^{2} = 0$.  Indeed, it is 
easy to check that this problematic situation occurs 
(for example) when $\mu = 
w_{0}(\lambda+\varrho) - \varrho$.  
So, by itself equation (6) falls short of uniquely determining all 
the $\Phi$-Kostka numbers.  But in view of \FinitenessCorollary, we have $d_{\lambda,\mu} = 0$ if $\mu \not\in \Pi(\lambda)$.  In \KostkaCorollary\ below we show that $d_{\lambda,\mu} \not= 0$ if and only if $\mu \in \Pi(\lambda)$, and in this case we must have $d_{\lambda,\mu} > 0$.  

{\em Proof of \FreudTheorem.} 
Consider $\mathcal{Z}_{\mathbb{R}} := 
\mathbb{Z}[\Lambda] \otimes_{\mathbb{Z}} \mathbb{R}$.  
We define an operator $\Delta: \mathcal{Z}_{\mathbb{R}} 
\longrightarrow \mathcal{Z}_{\mathbb{R}}$ by the rule 
$\Delta(\sum_{\mu \in \Lambda}c_{\mu}e^{\mu}) = \sum_{\mu \in 
\Lambda}||\mu||^{2}c_{\mu}e^{\mu}$.  Check that 
$\Delta(\mathcal{A}(e^{\lambda+\varrho})) = 
||\lambda+\varrho||^{2}\mathcal{A}(e^{\lambda+\varrho})$ for all 
$\lambda \in \Lambda^{+}$.  Also note that 
$\Delta(\chi_{_{\lambda}}) = \sum_{\mu \in 
\Lambda}||\mu||^{2}d_{\lambda,\mu}e^{\mu}$.  In applying $\Delta$ to 
equation (4), we need to know how $\Delta$ acts on products.  

Let $\nabla: \mathcal{Z}_{\mathbb{R}} \longrightarrow \mathfrak{E} 
\otimes \mathcal{Z}_{\mathbb{R}}$ be given by $\nabla(\sum_{\mu 
\in \Lambda}c_{\mu}e^{\mu}) = \sum_{\mu \in 
\Lambda}c_{\mu}\mu \otimes e^{\mu}$.  
Define the pairing $(\cdot,\cdot): \mathfrak{E} 
\otimes \mathcal{Z}_{\mathbb{R}} \times \mathfrak{E} 
\otimes \mathcal{Z}_{\mathbb{R}} \longrightarrow 
\mathcal{Z}_{\mathbb{R}}$ by the rule $(\lambda_{1} \otimes 
e^{\mu_{1}}, \lambda_{2} \otimes e^{\mu_{2}}) := \langle 
\lambda_{1},\lambda_{2} \rangle e^{\mu_{1}+\mu_{2}}$. 
Now confirm the following product rules applied to $f, g \in 
\mathcal{Z}_{\mathbb{R}}$: 
\begin{eqnarray*}
\nabla(fg) & = & (\nabla f)\, (g) + (f)\, (\nabla g)\\
\Delta(fg) & = & (\Delta f)\, (g) + 2(\nabla f, \nabla g) + 
(f)\, (\Delta g)
\end{eqnarray*}

Let $f := \mathcal{A}(e^{\varrho})$ and $g := \chi_{_{\lambda}}$.  
Then $\Delta(fg) = 
\Delta(\mathcal{A}(e^{\lambda+\varrho})) = 
||\lambda+\varrho||^{2}\mathcal{A}(e^{\lambda+\varrho}) = 
||\lambda+\varrho||^{2}(fg)$, with the latter equality following from equation (3). 
By the product rule, 
$\Delta(fg)  =  
(\Delta f)\, (g) + 2(\nabla f, \nabla g) + (f)\, (\Delta 
g) = ||\varrho||^{2}(fg) + 2(\nabla f, \nabla g) + f\sum_{\mu \in 
\Lambda}||\mu||^{2}d_{\lambda,\mu}e^{\mu}$. Thus:  
$\sum_{\mu \in \Lambda}(||\lambda+\varrho||^{2} - ||\mu||^{2} - 
||\varrho||^{2})d_{\lambda,\mu}e^{\mu} = 
2f^{-1}(\nabla f,\nabla g)$.  We need to simplify the right-hand side 
of the latter identity. 

For any $S \subseteq \Phi^{+}$, define $h_{S} := \prod_{\beta \in 
\Phi^{+}\setminus{S}}(1-e^{-\beta})$.  By Weyl's denominator 
formula, $f = e^{\varrho}h_{\emptyset}$.  Also, 
$h_{\emptyset}^{-1}h_{S} = \prod_{\beta \in S}(1-e^{-\beta})^{-1}$. 
Using the product rule applied to $k$-fold products, we get $\nabla 
f = \nabla(e^{\varrho}\prod_{\alpha \in \Phi^{+}}(1-e^{-\alpha})) = 
\varrho \otimes e^{\varrho}h_{\emptyset} + \sum_{\alpha \in 
\Phi^{+}}\alpha \otimes e^{\varrho-\alpha}h_{\{\alpha\}}$. So:  
\begin{eqnarray*} 
2f^{-1}(\nabla f,\nabla g) & = & 
2e^{-\varrho}h_{\emptyset}^{-1}\Big( 
\varrho \otimes e^{\varrho}h_{\emptyset} + \sum_{\alpha \in 
\Phi^{+}}\alpha \otimes e^{\varrho-\alpha}h_{\{\alpha\}},\sum_{\mu \in 
\Lambda}d_{\lambda,\mu}\mu \otimes e^{\mu} \Big)\\
 & = & (2\varrho \otimes e^{0},\sum_{\mu \in 
\Lambda}d_{\lambda,\mu}\mu \otimes e^{\mu}) + 2\Big( \sum_{\alpha \in 
\Phi^{+}}\alpha \otimes e^{-\alpha}(1-e^{-\alpha})^{-1},\sum_{\mu \in 
\Lambda}d_{\lambda,\mu}\mu \otimes e^{\mu} \Big)\\
 & = & \sum_{\mu \in \Lambda}2\langle \varrho,\mu \rangle 
 d_{\lambda,\mu}e^{\mu} + 2\Big( \sum_{\alpha \in 
\Phi^{+}}\sum_{k \geq 1}\alpha \otimes e^{-k\alpha},\sum_{\mu \in 
\Lambda}d_{\lambda,\mu}\mu \otimes e^{\mu} \Big)\\ 
 & = & \sum_{\mu \in \Lambda}2\langle \varrho,\mu \rangle 
 d_{\lambda,\mu}e^{\mu} + 2\sum_{\mu \in \Lambda}\sum_{\alpha \in 
\Phi^{+}}\sum_{k \geq 1}\langle \alpha,\mu \rangle d_{\lambda,\mu} 
e^{\mu-k\alpha}\\ 
 & = & \sum_{\mu \in \Lambda}2\langle \varrho,\mu \rangle 
 d_{\lambda,\mu}e^{\mu} + 2\sum_{\mu \in \Lambda}\sum_{\alpha \in 
\Phi^{+}}\sum_{k \geq 1}\langle 
\alpha,\mu+k\alpha \rangle d_{\lambda,\mu+k\alpha} 
e^{\mu},
 \end{eqnarray*}  
where the latter is obtained by reindexing.  Combine this last 
expression with the identity from the previous paragraph to get 
\[\sum_{\mu \in \Lambda}(||\lambda+\varrho||^{2} - ||\mu||^{2} - 
2\langle \varrho,\mu \rangle - 
||\varrho||^{2})d_{\lambda,\mu}e^{\mu} 
= 2\sum_{\mu \in \Lambda}\sum_{\alpha \in 
\Phi^{+}}\sum_{k \geq 1}\langle 
\alpha,\mu+k\alpha \rangle d_{\lambda,\mu+k\alpha} 
e^{\mu}.\] Equate 
coefficients of the $e^{\mu}$'s to get equation (6).  That the recurrence (6) together with the conditions in the theorem statement uniquely the numbers $\{d_{\lambda,\mu}\}$ follows from the discussion in the paragraph preceding the proof.\hfill\QED 

\noindent
{\bf \KostkaCorollary}\ \ {\sl Let $\lambda \in \Lambda^{+}$.  Then for all $\mu \in \Lambda$ we have $d_{\lambda,\mu} \geq 0$, and moreover $d_{\lambda,\mu} > 0$ if and only if $\mu \in \Pi(\lambda)$.  In particular, $\CharPi(\chi_{_{\lambda}}) = 
\Pi(\lambda)$.} 

{\em Proof.} In view of \FinitenessCorollary, it suffices to show that for each $\nu \in \Pi(\lambda) \cap 
\Lambda^{+}$, 
$d_{\lambda,\nu} > 0$. We do so by inducting on the quantity 
$\mathrm{depth}(\nu) = \mathrm{ht}(\lambda) - \mathrm{ht}(\nu)$, which 
is a nonnegative integer. For $\nu \in \Pi(\lambda) \cap 
\Lambda^{+}$ with $\mathrm{depth}(\nu) = 0$, then necessarily $\nu = \lambda$.  
In this case $d_{\lambda,\nu} = 1$.  Now say $d_{\lambda,\eta} > 0$ 
for all $\eta \in \Pi(\lambda) \cap 
\Lambda^{+}$ with $\mathrm{depth}(\eta) < \mathrm{depth}(\nu)$. Apply \FreudTheorem\ to $\nu$:  
\[(||\lambda+\varrho||^{2} - ||\nu+\varrho||^{2}) d_{\lambda,\nu} = 
2\sum_{\alpha \in \Phi^{+}}\sum_{k \geq 1}\langle \alpha,\nu+k\alpha 
\rangle d_{\lambda,\nu+k\alpha}.\]  On the left-hand side, the 
quantity $||\lambda+\varrho||^{2} - ||\nu+\varrho||^{2}$ is positive 
by Lemma 13.4.C of \cite{Hum}.  On the right hand side, it is easy to 
see that each $\langle \alpha,\nu+k\alpha 
\rangle > 0$ since $\nu \in \Lambda^{+}$.  Now suppose $\nu+k\alpha 
\in \Pi(\lambda)$ for some $k > 0$ and $\alpha \in \Phi^{+}$.  Then 
$\nu+k\alpha = w(\eta)$ for some $w \in W$ and $\eta \in 
\Pi(\lambda) \cap \Lambda^{+}$.  So $\nu < \nu+k\alpha = w(\eta) 
\leq \eta$.  This implies that $d_{\lambda,\eta} > 0$ (by the 
induction hypothesis) and that $d_{\lambda,\nu+k\alpha} = 
d_{\lambda,\eta}$ (since $\chi_{_{\lambda}}$ is $W$-invariant).  So we 
only need to show that $\nu+k\alpha \in \Pi(\lambda)$ for some $k > 0$ 
and $\alpha \in \Phi^{+}$.  
Since $\nu < \lambda$, then $\nu$ is not maximal in 
$\Pi(\lambda)$ by \WeightRemarkOne.2.  But then by 
\WeightRemarkOne.5, it must be the case that for some $i \in I$ we 
have $\nu + \alpha_{i} = \xi \in \Pi(\lambda)$.   Thus 
$d_{\lambda,\nu} > 0$, completing the induction argument.\hfill\QED

{\bf Basis results.} 
The following proposition gives sufficient conditions for a set 
$\{\phi_{\lambda}\}_{\lambda \in \Lambda^{+}}$ to be a basis for the 
ring of $W$-symmetric functions.  To show that 
$\{\phi_{_{\lambda}}\}$ is a spanning set, we prescribe an algorithm 
for writing any nonzero $W$-symmetric function $\chi$ as a $\mathbb{Z}$-linear 
combination of the $\phi_{\lambda}$'s that  
essentially inducts on the height $\mathrm{ht}$ of $W$-symmetric functions.  
Such induction is possible since if $\mathrm{ht}(\chi) > 
\mathrm{ht}(\chi')$ then $\mathrm{ht}(\chi) - 
\mathrm{ht}(\chi') \geq \myd_{\Phi}$ (our mesh size estimate).  

\noindent 
{\bf \BasisProp}\ \ {\sl Suppose 
$\{\phi_{\lambda}\}_{\lambda \in \Lambda^{+}} \subset 
\mathbb{Z}[\Lambda]^{W}$ satisfies: (1) $\phi_{0} = e^{0}$ and 
for all $\lambda \in \Lambda^{+}$, 
$[e^{\lambda}](\phi_{\lambda}) = 1$; and (2) For 
all $\lambda, \nu \in \Lambda^{+}$, if $\mathrm{ht}(\nu) \geq 
\mathrm{ht}(\lambda)$ and $\nu \in \CharPi(\phi_{\lambda})$, then $\nu = 
\lambda$. 
Then  
$\{\phi_{\lambda}\}_{\lambda \in \Lambda^{+}}$ is a 
$\mathbb{Z}$-basis for 
$\mathbb{Z}[\Lambda]^{W}$.  Moreover, 
if for some $W$-symmetric function $\chi$ we have 
$\chi = \sum_{\lambda \in 
\Lambda^{+}}c_{\lambda}\phi_{\lambda}$ 
and if for some $\nu \in \Lambda^{+}$ 
we have $c_{\nu} \not= 0$, then  
$\mathrm{ht}(\nu) \leq \mathrm{ht}(\chi)$.} 

{\em Proof.}  Say $\varphi = \sum_{\lambda \in 
\Lambda^{+}}c_{\lambda}\phi_{\lambda} = \sum_{\mu \in 
\Lambda}f_{\mu}e^{\mu}$ with at least one but at 
most finitely many 
$c_{\lambda}$'s nonzero. Pick $\nu \in \Lambda^{+}$ with $c_{\nu} 
\not= 0$ and such that 
$\mathrm{ht}(\nu) = \max\{\mathrm{ht}(\lambda)\, |\, 
c_{\lambda} \not= 
0\}$.  Then by {\sl (2)}, $\nu \in \CharPi(\phi_{\lambda})$ for some 
$\lambda \in 
\Lambda^{+}$ only if $\lambda = \nu$.  It follows from {\sl (1)} 
that $f_{\nu} = 
c_{\nu}[e^{\nu}](\phi_{\nu}) = c_{\nu}$, so $\varphi \not= 0$.   
This shows 
independence of the set $\{\phi_{\lambda}\}_{\lambda \in 
\Lambda^{+}}$. 

To show that the $\mathbb{Z}$-span of 
$\{\phi_{\lambda}\}_{\lambda \in \Lambda^{+}}$ is all of 
$\mathbb{Z}[\Lambda]^{W}$, we first establish the 
following claim ({\tt *}): If 
$\chi = \sum_{\mu \in \Lambda}b_{\mu}e^{\mu}$ 
and $\chi' := \chi - \sum_{\lambda 
\in \mathcal{H}(\chi)}b_{\lambda}\phi_{\lambda} = \sum_{\mu \in 
\Lambda}  
b'_{\mu}e^{\mu}$ are nonzero $W$-symmetric functions, then 
$\mathrm{ht}(\chi') < \mathrm{ht}(\chi)$ and moreover  
$b_{\nu}' \not= 0$ for $\nu \in \Lambda^{+}$ only if $\nu \in 
\CharPi(\chi)$. 
Indeed, suppose $b'_{\nu} \not= 0$ 
for some $\nu \in \Lambda^{+}$. Then 
(i) $b_{\nu} \not= 0$ or 
(ii)  $[e^{\nu}](\phi_{\lambda}) \not= 0$ 
for some $\lambda \in \mathcal{H}(\chi)$.  
In case (i), then $\nu \in \CharPi(\chi)$ so by definition 
$\mathrm{ht}(\nu) \leq 
\mathrm{ht}(\chi)$.  Supposing $\nu \in \mathcal{H}(\chi)$, then $\nu 
\not\in \CharPi(\phi_{\nu'})$ for all $\nu' \not= \nu$ in 
$\mathcal{H}(\chi)$ by condition {\sl (2)}.  But then 
$[e^{\nu}](\phi_{\nu}) = 1$ by condition {\sl (1)}, which would give us 
$b'_{\nu} = b_{\nu} - b_{\nu}[e^{\nu}](\phi_{\nu}) = 0$.  
Therefore $\nu \not\in 
\mathcal{H}(\chi)$, so $\mathrm{ht}(\nu) < 
\mathrm{ht}(\chi)$. 
In case (ii), $\nu \in \CharPi(\phi_{\lambda})$ and hence 
$\mathrm{ht}(\nu) 
\leq \mathrm{ht}(\lambda)$.  If $\mathrm{ht}(\nu) = 
\mathrm{ht}(\lambda)$, then $\nu = \lambda$ by condition {\sl (2)}, 
but then as in case (i) before 
we would get  
$b'_{\nu} = b_{\nu} - b_{\nu}[e^{\nu}](\phi_{\nu}) = 0$.  Therefore 
$\mathrm{ht}(\nu) 
< \mathrm{ht}(\lambda) = \mathrm{ht}(\chi)$. Based on the results of 
these two cases, we conclude that $\mathrm{ht}(\nu) < 
\mathrm{ht}(\chi)$ for each $\nu \in \Lambda^{+}$ such that $b_{\nu}' 
\not= 0$, and therefore that 
$\mathrm{ht}(\chi') < \mathrm{ht}(\chi)$ as desired.  

Now suppose $\chi = \sum_{\mu \in \Lambda} 
b_{\mu}e^{\mu}$ is a nonzero $W$-symmetric function. 
We provide an algorithm for expressing $\chi$ 
as $\sum_{\lambda \in \Lambda^{+}} 
c_{\lambda}\phi_{\lambda}$ for some integers $c_{\lambda}$ with 
only finitely many nonzero.  Let $\chi^{(0)} := \chi$ with 
$b^{(0)}_{\mu} := b_{\mu}$ and set 
$\chi^{(1)} := 
\chi^{(0)}-\sum_{\lambda \in 
\mathcal{H}(\chi)} b^{(0)}_{\lambda}\phi_{\lambda} = \sum_{\mu \in 
\Lambda}  
b^{(1)}_{\mu}e^{\mu}$ for some integers $b^{(1)}_{\mu}$. If 
$\chi^{(1)} = 0$, then we are done. Else, by ({\tt *}) we have 
$\mathrm{ht}(\chi^{(1)}) < \mathrm{ht}(\chi^{(0)})$ and $\nu \in 
\CharPi(\chi^{(0)})$ if $b^{(1)}_{\nu} \not= 0$ for some 
$\nu \in \Lambda^{+}$.  Now  
$\mathrm{ht}(\chi^{(0)}) - \mathrm{ht}(\chi^{(1)})$ is at least 
$\myd_{\Phi} > 0$.  So as we continue this 
procedure, we will reach some $\chi^{(k)} = 
\chi^{(k-1)}-\sum_{\lambda \in 
\mathcal{H}(\chi^{(k-1)})}b^{(k-1)}_{\lambda}\phi_{\lambda}$ with 
$\chi^{(k)} = b_{0}^{(k)}e^{0}$ for some integer $b^{(k)}_{0}$.  Now 
$e^{0} = \phi_{0}$ by condition {\sl (1)}, so we now have $\chi$ as a 
$\mathbb{Z}$-linear combination of the $\phi_{\lambda}$'s.  
Also, for each $i$ with $1 \leq i \leq k$ and   
each nonzero $b_{\nu}^{(i)}$ with $\nu \in \Lambda^{+}$ we have 
$\mathrm{ht}(\nu) < \mathrm{ht}(\chi)$.  Since $b_{\nu}^{(0)} \not= 0$ 
for some $\nu \in \Lambda^{+}$ only if $\nu \in \mathcal{H}(\chi)$, 
then $\mathrm{ht}(\nu) = \mathrm{ht}(\chi)$.  
This shows that any nonzero $W$-symmetric function $\chi$ can be written as 
$\sum_{\lambda \in \Lambda^{+}}c_{\lambda}\phi_{\lambda}$ with 
$c_{\lambda} \not= 0$ only if $\mathrm{ht}(\lambda) \leq 
\mathrm{ht}(\chi)$.\hfill\QED  

An {\em amenable basis} for $\mathbb{Z}[\Lambda]^{W}$  is 
any $\mathbb{Z}$-basis 
$\{\phi_{\lambda}\}_{\lambda \in \Lambda^{+}}$ 
that satisfies the two  
conditions of \BasisProp. 
It is easy to see that the monomial 
$W$-symmetric functions are an amenable basis for the ring of $W$-symmetric functions.  
In \BasisTheorem\ we argue that 
the Weyl bialternants and the elementary $W$-symmetric functions are amenable bases as well.  

\noindent
{\bf \BasisTheorem\ (Basis Theorem)}\ \ 
{\sl The Weyl bialternants} $\{\chi_{_{\lambda}}\}_{\lambda \in \Lambda^{+}}$, {\sl the monomial $W$-symmetric functions} 
$\{\zeta_{\lambda}\}_{\lambda \in \Lambda^{+}}$, {\sl and the elementary $W$-symmetric functions} $\{\psi_{\lambda}\}_{\lambda \in \Lambda^{+}}$ {\sl are  
amenable bases for} $\mathbb{Z}[\Lambda]^{W}$. {\sl The latter result means that the fundamental $W$-symmetric functions} $\{\chi_{_{\omega_{i}}}\}_{_{i \in I}}$ {\sl are an algebraic basis 
for $\mathbb{Z}[\Lambda]^{W}$, so 
$\mathbb{Z}[\Lambda]^{W}$ is isomorphic to the polynomial 
ring in $n$ variables over $\mathbb{Z}$.}

{\em Proof.} It suffices to    
check conditions {\sl (1)} and {\sl (2)} from \BasisProp\ 
for each of the sets 
$\{\zeta_{\lambda}\}_{\lambda \in \Lambda^{+}}$, 
$\{\chi_{_{\lambda}}\}_{\lambda \in \Lambda^{+}}$, and 
$\{\psi_{\lambda}\}_{\lambda \in \Lambda^{+}}$. 
First consider $\{\zeta_{\lambda}\}$. For condition {\sl (1)}, 
$\zeta_{0} = e^{0}$ follows from the definitions, and 
$[e^{\lambda}](\zeta_{\lambda}) = 1$ follows from \WeylCharLemma.  
For condition {\sl (2)}, say $\nu \in \CharPi(\zeta_{\lambda})$ with 
$\mathrm{ht}(\nu) \geq 
\mathrm{ht}(\lambda)$.  Now 
$\nu \in \CharPi(\zeta_{\lambda})$ means $\nu = 
w(\lambda)$ for some $w \in W$, hence $\nu \in 
\Pi(\lambda)$.  Then $\nu \leq \lambda$, hence $\mathrm{ht}(\nu) \leq 
\mathrm{ht}(\lambda)$.  So $\mathrm{ht}(\nu) = \mathrm{ht}(\lambda)$.  
The only way this can happen in $\Pi(\lambda)$ is if $\nu = \lambda$. 
Second consider $\{\chi_{_{\lambda}}\}$. 
For condition {\sl (1)}, $d_{\lambda,\lambda}=1$ by 
\WeylCharLemma.  
Since $\Pi(0) = \{0\}$, then by \KostkaCorollary, it 
follows that $\chi_{_{0}} = 
\sum_{\mu \in \Pi(0)}d_{0,\mu}e^{\mu} = 
d_{0,0}e^{0} = e^{0}$. 
For condition {\sl (2)}, suppose 
$\mathrm{ht}(\nu) \geq \mathrm{ht}(\lambda)$ and $\nu 
\in \CharPi(\chi_{_{\lambda}})$ for 
some $\nu, \lambda \in \Lambda^{+}$.  
Note that $\CharPi(\chi_{_{\lambda}}) = \Pi(\lambda)$, so 
\WeightRemarkOne\ applies.  Then 
$\nu \in \Pi(\lambda)$ implies that $\nu \leq \lambda$, and 
therefore $\mathrm{ht}(\nu) \leq \mathrm{ht}(\lambda)$.  So 
$\mathrm{ht}(\nu) = \mathrm{ht}(\lambda)$.  The only way this can 
happen in $\Pi(\lambda)$ is if $\nu = \lambda$. 
Next consider $\{\psi_{\lambda}\}$. Now 
$[e^{\lambda}](\psi_{\lambda}) = 1$ by \WeylCharLemma, and 
$\psi_{0} = e^{0}$ by convention. 
This establishes condition {\sl (1)}.  For {\sl (2)},  
suppose $\mathrm{ht}(\nu) \geq \mathrm{ht}(\lambda)$ and $\nu 
\in \CharPi(\psi_{\lambda})$ for some $\nu, \lambda \in \Lambda^{+}$. 
If $\lambda = \sum_{i \in 
I}a_{i}\omega_{i}$, then we have $\nu = \sum_{i \in 
I}\sum_{j=1}^{a_{i}}\nu_{i}^{(j)}$, where each $\nu_{i}^{(j)} \in 
\Pi(\omega_{i})$. 
Then $\mathrm{ht}(\nu) = \sum_{i \in 
I}\sum_{j=1}^{a_{i}}\mathrm{ht}(\nu_{i}^{(j)}) \leq \sum_{i \in 
I}\sum_{j=1}^{a_{i}}\mathrm{ht}(\omega_{i}) = \mathrm{ht}(\lambda)$. 
So $\mathrm{ht}(\nu) = \mathrm{ht}(\lambda)$. 
It follows that $\mathrm{ht}(\nu_{i}^{(j)}) = \mathrm{ht}(\omega_{i})$ 
for each $i \in I$ and $1 \leq j \leq a_{i}$. Therefore 
we must have $\nu_{i}^{(j)} = \omega_{i}$, 
so $\nu = \lambda$.\hfill\QED   

{\bf Two involutions.}  The following involutions of the group ring are sometimes useful.  Let $^{*}$ be the involution of $\mathbb{Z}[\Lambda]$ induced by 
$(e^{\mu})^{*} = e^{-\mu}$. Similarly, define $^{\bowtie}$  
to be the involution of $\mathbb{Z}[\Lambda]$ induced by 
$(e^{\mu})^{\bowtie} = e^{w_{0}(\mu)}$. 
These are the {\em dualizing} and {\em bow tie} involutions, respectively. 
Note that while $\mathbb{Z}[\Lambda]$ is closed under $^{*}$ and 
$^{\bowtie}$, $\mathcal{R}$ is not. The next result says how these involutions interact with the ring of Weyl symmetric functions. 

\noindent 
{\bf \StarProp}\ \ {\sl Let $\lambda \in \Lambda^{+}$.  Then $\chi_{_{\lambda}}^{*} = 
\chi_{_{-w_{0}(\lambda)}}$ and $\chi_{_{\lambda}}^{\bowtie} = 
\chi_{_{\lambda}}$.  In particular, $^{*}$ is an involution of $\mathbb{Z}[\Lambda]^{W}$ and $^{\bowtie}$ restricts to the identity on $\mathbb{Z}[\Lambda]^{W}$.} 

{\em Proof.} For any $\nu \in \Lambda^{+}$, we have 
$\mathcal{A}(e^{\nu+\varrho})^{*} = 
\sum\det(w)e^{-w(\nu+\varrho)} = \sum 
\det(w \cdot w_{0})e^{-w \cdot w_{0}(\nu+\varrho)} = \sum 
\det(w \cdot w_{0})e^{w(-w_{0}(\nu)+\varrho)} = 
\det(w_{0})\sum \det(w)e^{w(-w_{0}(\nu)+\varrho)} = 
\det(w_{0})\mathcal{A}(e^{-w_{0}(\nu)+\varrho})$,  
where each sum is over all 
$w \in W$. Then applying $*$ to both sides of 
$\mathcal{A}(e^{\varrho})\, \chi_{_{\lambda}} = 
\mathcal{A}(e^{\lambda+\varrho})$ we get: 
$\det(w_{0})\mathcal{A}(e^{\varrho})\, \chi_{_{\lambda}}^{*} = 
\det(w_{0})\mathcal{A}(e^{-w_{0}(\lambda)+\varrho})$, or simply   
$\mathcal{A}(e^{\varrho})\, \chi_{_{\lambda}}^{*} = 
\mathcal{A}(e^{-w_{0}(\lambda)+\varrho})$. 
Therefore $\chi_{_{\lambda}}^{*} = 
\chi_{_{-w_{0}(\lambda)}}$.   Similarly, $\mathcal{A}(e^{\nu+\varrho})^{\bowtie} = \det(w_{0})\mathcal{A}(e^{\nu+\varrho})$ for any $\nu \in \Lambda^{+}$.  Then we have $\mathcal{A}(e^{\nu+\varrho}) = \det(w_{0})\mathcal{A}(e^{\nu+\varrho})^{\bowtie} = \det(w_{0})\mathcal{A}(e^{\varrho})^{\bowtie} \chi_{_{\lambda}}^{\bowtie} = \mathcal{A}(e^{\nu+\varrho})\chi_{_{\lambda}}^{\bowtie}$, so $\chi_{_{\lambda}}^{\bowtie} = \chi_{_{\lambda}}$.  The claims of the last sentence of the proposition follow from the fact that $\{\chi_{_{\lambda}}\}_{\lambda \in \Lambda^{+}}$ is a basis for $\mathbb{Z}[\Lambda]^{W}$.\hfill\QED 

{\bf Specializations of Weyl bialternants.} 
The next result, and particularly its part {\sl (1)}, has important 
consequences for the combinatorics of splitting posets, cf.\ 
\MainCorollary. 

\noindent 
{\bf \TFAEConsequences\ (Specializations)}\ \ 
{\sl Let $\lambda \in \Lambda^{+}$. Then:}\\
\hspace*{0.25in}{\sl (1)}\ \ \parbox[t]{5.75in}{{\sl The three quantities of the following identity are indeed equal and are a symmetric polynomial whose 
degree is the nonnegative integer $2\langle \lambda,\varrho^{\vee} \rangle$:} 
\[q^{\langle \lambda,\varrho^{\vee} 
\rangle}\sum_{\mu \in \Lambda}d_{\lambda,\mu}q^{\langle 
\mu,\varrho^{\vee} 
\rangle} = 
\sum_{\mu \in \Lambda}d_{\lambda,\mu}q^{\langle 
\mu+\lambda,\varrho^{\vee} 
\rangle} = 
\frac{\prod_{\alpha \in \Phi^{+}}(1-q^{\langle 
\lambda+\varrho,\alpha^{\vee} \rangle})}{\prod_{\alpha \in 
\Phi^{+}}(1-q^{\langle 
\varrho,\alpha^{\vee} \rangle})}.\]}\\
\hspace*{0.25in}{\sl (2)}\ \ \parbox[t]{5in}{{\sl We have:} 
$\displaystyle \sum_{\mu \in \Lambda}d_{\lambda,\mu} = 
\frac{\prod_{\alpha \in \Phi^{+}}\langle 
\lambda+\varrho,\alpha^{\vee} \rangle}{\prod_{\alpha \in 
\Phi^{+}}\langle 
\varrho,\alpha^{\vee} \rangle}.$}

\vspace{0.1in}
The proof is below. 
The formula of part {\sl (2)} of \TFAEConsequences\ is known as the ``Weyl dimension formula'' since it records the dimension of a related irreducible semisimple Lie algebra representation. 
This Lie theory connection is discussed in more detail at the end of this section.  
The quantity $\sum_{\mu \in \Lambda}d_{\lambda,\mu}q^{\langle \mu,\varrho^{\vee} \rangle}$ from {\sl (1)} can be viewed as $\chi_{_{\lambda}} = \sum_{\mu \in \Lambda}d_{\lambda,\mu}e^{\mu}$ specialized by $e^{\omega_i} \mapsto q^{\langle \omega_i,\varrho^{\vee} \rangle}$. 
In \cite{Panyushev}, the $q$-polynomials of 
{\sl (1)} are termed ``Dynkin 
polynomials.'' 
In \cite{Ram}, the latter identity of {\sl (1)} is called the ``quantum dimension formula.''  

\noindent 
{\bf \UnimodalProblem}\ \ For families of Weyl bialternants, provide a combinatorial proof of the 
unimodality of their associated 
Dynkin polynomials.  This seems to be a difficult problem, even in very 
special cases.  For example, a combinatorial 
proof of the unimodality of the Dynkin polynomial associated with 
the fundamental $\myA_{n}$-bialternant $\chi_{_{\omega_{k}}}$  
-- which is just a 
$q$-binomial coefficient 
-- was only 
first obtained in the late 1980's 
by O'Hara (\cite{OHara}; see also \cite{Zeilberger}). A proof of 
unimodality using Lie algebra representation theory is noted 
below.\hfill\QED  

{\em Proof of \TFAEConsequences.} 
For {\sl (1)}, note that by by parts {\sl (2)} and {\sl (3)} of \WeightRemarkOne, we have $\langle 
\mu+\lambda,\varrho^{\vee} \rangle \in \{0,1,2,\ldots,l\}$ for all 
$\mu \in \Pi(\lambda)$, where 
$l = \langle \lambda+\lambda,\varrho^{\vee} \rangle$ and $0 = \langle 
w_{0}(\lambda)+\lambda,\varrho^{\vee} \rangle$.  
So, the left-hand side is a polynomial in $q$ of the claimed degree.  
For $0 \leq k \leq l$, let $\mathcal{D}_{k} := \{\mu \in 
\Pi(\lambda)\, |\, \langle \mu+\lambda,\varrho^{\vee} \rangle = k\}$.  
So we can write 
$\displaystyle \sum_{\mu \in \Pi(\lambda)}d_{\lambda,\mu}q^{\langle 
\mu+\lambda,\varrho^{\vee} \rangle} = \sum_{k=0}^{l}D_{k}q^{k}$, 
with $\displaystyle D_{k} := \sum_{\mu \in 
\mathcal{D}_{k}}d_{\lambda,\mu}$.  We know $d_{\lambda,\mu} = 
d_{\lambda,w_{0}(\mu)}$, so to show symmetry of the polynomial 
coefficients (i.e.\ $D_{k} = D_{l-k}$ for all $0 \leq k \leq l$) it 
suffices to show that $w_{0}(\mathcal{D}_{k}) = \mathcal{D}_{l-k}$ 
for all $0 \leq k \leq l$.  Now for all $\nu \in \Lambda$, it is easy 
to check that $\langle \nu,\varrho^{\vee} \rangle = \langle 
-w_{0}(\nu),\varrho^{\vee} \rangle$.  Then for $\mu \in 
\mathcal{D}_{k}$ we calculate: 
$\langle w_{0}(\mu)+\lambda,\varrho^{\vee} \rangle = 
\langle w_{0}(\mu)-w_{0}(\lambda),\varrho^{\vee} \rangle = 
\langle w_{0}(\mu-\lambda),\varrho^{\vee} \rangle = 
\langle \lambda-\mu,\varrho^{\vee} \rangle = 
2\langle \lambda,\varrho^{\vee} \rangle - 
\langle \mu+\lambda,\varrho^{\vee} \rangle= l-k$.  So, 
$w_{0}(\mathcal{D}_{k}) \subseteq \mathcal{D}_{l-k}$.  Since this is 
true for all $0 \leq k \leq l$, then 
$w_{0}(\mathcal{D}_{l-k}) \subseteq \mathcal{D}_{k}$ as well, so 
$w_{0}(\mathcal{D}_{k}) = \mathcal{D}_{l-k}$. 
Then symmetry of the $q$-polynomial follows. 
To complete the proof of {\sl (2)}, 
we need to work with the group ring $\mathbb{Z}[\mathfrak{E}]$, where 
$W$ acts as before and we define $\mathcal{A}(e^{u}) := \sum_{w \in 
W}\det(w)e^{w(u)}$ for any $u \in \mathfrak{E}$.  
For any $v \in \mathfrak{E}$, 
the mapping $\Psi_{v}: \mathbb{Z}[\mathfrak{E}] \longrightarrow 
\mathbb{R}[[q]]$ with $\Psi_{v}(e^{u}) = q^{\langle u,v 
\rangle}$ extends uniquely to a ring homomorphism.  It is easy to 
check that for all $u, v \in \mathfrak{E}$, 
$\Psi_{v}(\mathcal{A}(e^{u})) = \Psi_{u}(\mathcal{A}(e^{v}))$. 
Since $\varrho^{\vee}$ is half the sum of 
the positive coroots in the dual root system $\Phi^{\vee}$, then 
$\displaystyle \mathcal{A}(e^{\varrho^{\vee}}) = 
\prod_{\alpha \in \Phi^{+}}
(e^{\alpha^{\vee}/2}-e^{-\alpha^{\vee}/2})$ by \WeylsDenomTheorem. 
Then $\sum_{\mu \in \Lambda}d_{\lambda,\mu}q^{\langle \mu,\varrho^{\vee} 
\rangle} = 
\Psi_{\varrho^{\vee}}(\chi_{_{\lambda}}) = 
\Psi_{\varrho^{\vee}}(\mathcal{A}(e^{\lambda+\varrho}))/
\Psi_{\varrho^{\vee}}(\mathcal{A}(e^{\varrho})) = 
\Psi_{\lambda+\varrho}(\mathcal{A}(e^{\varrho^{\vee}}))/
\Psi_{\varrho}(\mathcal{A}(e^{\varrho^{\vee}})) = 
\prod_{\alpha \in \Phi^{+}}
(q^{\langle \lambda+\varrho,\alpha^{\vee}/2 
\rangle} - q^{-\langle \lambda+\varrho,\alpha^{\vee}/2 
\rangle})/
\prod_{\alpha \in \Phi^{+}}
(q^{\langle \varrho,\alpha^{\vee}/2 
\rangle} - q^{-\langle \varrho,\alpha^{\vee}/2 
\rangle})$.  Simplify to obtain  
the identity in the proposition statement.  For {\sl (2)}, evaluate 
{\sl (1)} at $q=1$.\hfill\QED

\noindent 
{\bf \CaseAExample: The} $\myA_{n}$ {\bf case.} In this extended example we explicitly show how in case $\myA_{n}$ the bases $\{\chi_{_{\lambda}}\}$, $\{\zeta_{\lambda}\}$, and $\{\psi_{\lambda}\}$ (indexed by dominant weights) for the $W$-symmetric function ring are related respectively to the bases of Schur functions,  monomial symmetric functions, and elementary symmetric functions for the ring of symmetric functions in $n+1$ variables.  Throughout the example, we take $I = \{1,\ldots,n\}$ as our index set for simple roots and fundamental weights.  From the generators-and-relations description of the associated Weyl group $W$, it is easily seen that $W$ coincides with the symmetric group $\mathfrak{S}_{n+1}$ on $n+1$ letters when we identify the Weyl group generator $s_{i}$ ($i \in I$) with the symmetric group transposition $(i,i+1)$. 

Many quantities of interest will be indexed by sequences of nonnegative integers.  If $\mathfrak{p} = (p_{i})_{i=1}^{n+1}$ and $\mathfrak{q} = (q_{i})_{i=1}^{n+1}$ are two such sequences, then $\mathfrak{p}+\mathfrak{q}$ is their componentwise sum in the usual way.  Call $\mathfrak{p}$ a {\em partition} if $p_{1} \geq p_{2} \cdots \geq p_{n+1} \geq 0$, and set $|\mathfrak{p}| := \sum_{i=1}^{n+1} p_{i}$.  If $\mathfrak{q}$ is another partition, we write $\mathfrak{q} \triangleleft \mathfrak{p}$ if $q_{1} + \cdots + q_{i} \leq p_{1} + \cdots + p_{i}$ for $1 \leq i \leq n+1$.   
For $1 \leq k \leq n+1$, set $\mathbf{1}_{k} := (1,\ldots,1,0,\ldots,0)$, a partition with $k$ leading $1$'s. 
Then for any nonnegative integer $a$ we have $a\mathbf{1}_{k} = (a,\ldots,a,0,\ldots,0)$. 
The notation $\mathbf{1}$ means $\mathbf{1}_{n+1}$. 

The symmetric group $W \cong \mathfrak{S}_{n+1}$ acts naturally on the polynomial ring $\mathbb{Z}[z_{1},\ldots,z_{n+1}]$ via the rule that generator $s_{i}$ transposes the variables $z_{i}$ and $z_{i+1}$. The subring $\mathbb{Z}[z_{1},\ldots,z_{n+1}]^{W}$ of $W$-invariant polynomials is just the ring of symmetric functions in $n+1$ variables.  For a sequence $\mathfrak{p} = (p_{1},\ldots,p_{n+1})$ of nonnegative integers, the monomial $z^{\mathfrak{p}}$ is defined to be $z_{1}^{p_{1}}\cdots{z}_{n+1}^{p_{n+1}}$.  
Set $\weight(\mathfrak{p}) := \sum_{i \in I}(p_{i}-p_{i+1})\omega_{i}$.  For the ring homomorphism $\Theta: \mathbb{Z}[z_{1},\ldots,z_{n+1}] \longrightarrow \mathbb{Z}[\Lambda]$ induced by $z_{i} \stackrel{\Theta}{\mapsto} e^{\omega_{i}-\omega_{i-1}}$ (where $\omega_{0} := 0 =: \omega_{n+1}$), observe that $\Theta(z^{\mathfrak{p}}) = e^{\smallweight(\mathfrak{p})}$.  
The mapping $\Theta$ is easily seen to be surjective.  

Obviously the ideal $\langle z_{1}z_{2}\cdots{z}_{n+1} - 1 \rangle$ resides in $\ker(\Theta)$.  
We will argue for equality of these sets, but first a little bookkeeping. 
Any $f \in \mathbb{Z}[z_1,\cdots,z_{n+1}]$ can be written uniquely as $f = \sum_{\mathfrak{p}}c_{\mathfrak{p}}z^{\mathfrak{p}}$ for some integers $c_{\mathfrak{p}}$, where the sum is over all sequences $\mathfrak{p} = (p_{1},\ldots,p_{n+1})$ of nonnegative integers. 
It is easy to see that for two such sequences, $\weight(\mathfrak{q}) = \weight(\mathfrak{q}')$ if and only if there is some integer $k \geq \max\{-q_1,\cdots,-q_{n+1}\}$ such that $q'_{j} = q_{j} + k$ for each $1 \leq j \leq n+1$. 
Moreover, amongst the $\mathfrak{q}'$ sequences, exactly one has the property that some $q'_{j} = 0$. 
Let $\mathscr{Z}$ be the set of all nonnegative integer sequences $\mathfrak{p} = (p_1,\cdots,p_{n+1})$ such that some $p_j = 0$. 
Now observe that $\weight: \mathscr{Z} \longrightarrow \Lambda$ is a bijection. 
Then, we can write $f = \sum_{\mathfrak{p} \in \mathscr{Z}}z^{\mathfrak{p}}(\sum_{k \geq 0} c_{\mathfrak{p}+k\mathbf{1}}z^{k\mathbf{1}})$, where at most a finite number of the $c_{\mathfrak{p}+k\mathbf{1}}$'s are nonzero. 
So, $\Theta(f) = \sum_{\mu \in \Lambda}(\sum_{k \geq 0}c_{\mathfrak{p}+k\mathbf{1}})e^{\mu}$, where in the second summation $\mathfrak{p}$ is the unique element of  $\mathscr{Z}$ such that $\weight(\mathfrak{p}) = \mu$. 
Suppose now that $\Theta(f) = 0$.  
Then for each $\mathfrak{p} \in \mathscr{Z}$, we have $\sum_{k \geq 0}c_{\mathfrak{p}+k\mathbf{1}} = 0$. 
So for each such $\mathfrak{p}$, it follows that $1$ is a root of the polynomial $\sum_{k \geq 0}c_{\mathfrak{p}+k\mathbf{1}}Z^k$, where $Z$ is some indeterminate, and therefore that $Z-1$ is a factor. 
In particular, $z_{1}z_{2}\cdots{z}_{n+1} - 1$ is a factor of each $\sum_{k \geq 0} c_{\mathfrak{p}+k\mathbf{1}}z^{k\mathbf{1}}$ in the expression $f = \sum_{\mathfrak{p} \in \mathscr{Z}}z^{\mathfrak{p}}(\sum_{k \geq 0} c_{\mathfrak{p}+k\mathbf{1}}z^{k\mathbf{1}})$.  
That is, $f \in \langle z_{1}z_{2}\cdots{z}_{n+1} - 1 \rangle$. 

It is easy to see that the mapping $\Theta$ preserves the $W$-action.  
Therefore $\Theta(\mathbb{Z}[z_{1},\ldots,z_{n+1}]^{W}) \subseteq \mathbb{Z}[\Lambda]^{W}$.  
We wish to show that the latter subset containment is an equality.  
To this end, define a mapping $\mya$ on $\mathbb{Z}[z_{1},\ldots,z_{n+1}]$ by the rule $\mya(f) := \sum_{w \in W}\mathrm{sgn}(w)w.f$, where $\mathrm{sgn}(w) = \det(w)$ is $1$ (respectively $-1$) if as a  permutation $w$ is even (resp.\ odd).  
Then for any partition $\mathfrak{p} = (p_{1},\ldots,p_{n+1})$, we have \[\mya(z^{\mathfrak{p}}) = \det(z_{i}^{p_{j}})_{i,j=1}^{n+1}.\] 
Let $\mathfrak{r} := (n,n-1,\ldots,1,0)$, so $\Theta(z^{\mathfrak{r}}) = e^{\varrho}$.  
The classical definition of Schur functions sets \[\mys_{\mathfrak{p}} := \frac{\mya(z^{\mathfrak{p}+\mathfrak{r}})}{\mya(z^{\mathfrak{r}})}.\] 
See for example \cite{StanText2} \S 7.15 for a justification that $\mys_{\mathfrak{p}} \in \mathbb{Z}[z_{1},\ldots,z_{n+1}]^{W}$. 
Since $\Theta(\mya(z^{\mathfrak{p}+\mathfrak{r}})) = \mathcal{A}(\Theta(z^{\mathfrak{p}+\mathfrak{r}})) = \mathcal{A}(e^{\smallweight(\mathfrak{p})+\varrho})$, then it follows that $\Theta(\mys_{\mathfrak{p}}) = \chi_{_{\tinyweight(\mathfrak{p})}}$. 
Since  the set $\{\chi_{_{\tinyweight(\mathfrak{p})}}\}$ (indexed by partitions $\mathfrak{p}$) is the full set of Weyl bialternants in $\mathbb{Z}[\Lambda]^{W}$, then $\Theta(\mathbb{Z}[z_{1},\ldots,z_{n+1}]^{W}) = \mathbb{Z}[\Lambda]^{W}$.  
Let $\overline{\Theta} := \Theta|_{\mathbb{Z}[z_{1},\ldots,z_{n+1}]^{W}}$.  
So $\overline{\Theta}$ is a surjective ring homomorphism from $\mathbb{Z}[z_{1},\ldots,z_{n+1}]^{W}$ onto $\mathbb{Z}[\Lambda]^{W}$ and has kernel $\langle z_{1}\cdots{z}_{n+1} - 1 \rangle \cap \mathbb{Z}[z_{1},\ldots,z_{n+1}]^{W}$.  

For the partition $\mathfrak{p} = (p_{1},\ldots,p_{n+1})$, 
the monomial symmetric function  
$\mym_{\mathfrak{p}}$ is $\sum 
z^{\mathfrak{q}}$, where the sum is over all distinct permutations 
$\mathfrak{q}$ 
of the $(n+1)$-element sequence $\mathfrak{p}$.  It is easy to see that 
$\overline{\Theta}(\mym_{\mathfrak{p}}) = 
\zeta_{\mathrm{weight}(\mathfrak{p})}$.  
See \cite{StanText2} Ch.\ 7 for arguments that 
the Schur functions $\{\mys_{\mathfrak{p}}\}$ and monomial 
symmetric functions $\{\mym_{\mathfrak{p}}\}$ 
are bases for 
$\mathbb{Z}[z_{1},\ldots,z_{n+1}]^{W}$.  
Indeed, there it is shown that 
we may write $\mys_{\mathfrak{p}} = 
\sum_{\mathfrak{q}}K_{\mathfrak{p},\mathfrak{q}} 
\mym_{\mathfrak{q}}$, where the sum is over all partitions 
$\mathfrak{q} = (q_{1},\ldots,q_{n+1})$ and the ``Kostka number'' 
$K_{\mathfrak{p},\mathfrak{q}}$ is the number of semistandard 
(i.e.\ column-strict) tableaux 
of shape $\mathfrak{p}$ with entries from 
the set $\{1,2,\ldots,n+1\}$ and having $q_{1}$ 1's, $q_{2}$ 2's, 
etc. 

For any integer $1 \leq k 
\leq n+1$, let $\displaystyle 
\mye_{k} := \sum_{1 \leq i_{1} < \cdots < i_{k} \leq 
n+1}z_{i_{1}}\cdots{z}_{i_{k}}$.  
A consequence of the above description of Schur functions as 
$\mathbb{Z}$-linear combinations of monomial symmetric functions is 
that $\mye_{k} = \mys_{\mathbf{1}_{k}}$. It follows that 
$\overline{\Theta}(\mye_{k}) = \chi_{_{\omega_{k}}}$ for $k \in I$. 
(Of course, $\overline{\Theta}(e_{n+1}) = 
\overline{\Theta}(z_{1}\cdots{z}_{n+1}) = 1$.)  
For any sequence 
$\mathfrak{a} = (a_{1},\ldots,a_{n+1})$ of nonnegative integers 
set $\mye^{\mathfrak{a}} := 
\mye_{1}^{a_{1}}\cdots\mye_{n+1}^{a_{n+1}}$, an elementary symmetric function.  Then we 
have $\overline{\Theta}(\mye^{\mathfrak{a}}) = \psi_{\lambda}$, 
where $\lambda = \sum_{i \in I}a_{i}\omega_{i}$.  
That the set 
$\{\mye^{\mathfrak{a}}\}$ 
of elementary symmetric functions  
is a basis for  
$\mathbb{Z}[z_{1},\ldots,z_{n+1}]^{W}$ is sometimes called the 
Fundamental Theorem of Symmetric Functions, cf.\ Theorem 7.4.4 of 
\cite{StanText2}. 

For any $\lambda = \sum_{i \in I}a_{i}\omega_{i} 
\in \Lambda$, let $\widehat{\lambda}$ be the $(n+1)$-element sequence  
$(\sum_{j=i}^{n}a_{j})_{i=1}^{n+1}$, so in particular 
$\widehat{\lambda}_{n+1} = 0$. 
Let $\mathcal{S} := 
\mathrm{span}_{\mathbb{Z}}\{\mys_{\widehat{\lambda}}\}_{\lambda \in 
\Lambda^{+}}$, i.e.\ the integer linear span of Schur functions 
indexed by all partitions $\mathfrak{q} = (q_{1},\ldots,q_{n},q_{n+1})$ 
where $q_{n+1} = 0$.  We claim that $\mathcal{S}$ is a subring of 
$\mathbb{Z}[z_{1},\ldots,z_{n+1}]^{W}$ and that 
$\overline{\Theta}|_{_{\mathcal{S}}}$ is a ring isomorphism from 
$\mathcal{S}$ onto $\mathbb{Z}[\Lambda]^{W}$. 
Indeed, consider the product 
$\mys_{\widehat{\mu}}\mys_{\widehat{\nu}} = 
\sum_{\mathfrak{p}}
c_{\widehat{\mu},\widehat{\nu}}^{\mathfrak{p}}\mys_{\mathfrak{p}}$, 
where the sum is over all partitions $\mathfrak{p} = 
(p_{1},\ldots,p_{n+1})$.  The integers 
$c_{\widehat{\mu},\widehat{\nu}}^{\mathfrak{p}}$ are the 
Littlewood-Richardson coefficients, and 
by Theorem A1.3.3 of \cite{StanText2}, 
$c_{\widehat{\mu},\widehat{\nu}}^{\mathfrak{p}}$ counts certain 
semistandard tableaux of skew shape 
$\mathfrak{p}/\widehat{\mu}$ and type $\widehat{\nu}$.  In particular, 
the number of $n+1$'s appearing in any of those tableaux must be 
$\widehat{\nu}_{n+1} = 0$.  This means that the shape corresponding 
to $\mathfrak{p}$ can have no more than $n$ parts, i.e.\ $p_{n+1} = 
0$.  Thus the product $\mys_{\widehat{\mu}}\mys_{\widehat{\nu}}$ is 
in $\mathcal{S}$.  Observe that 
$\mathcal{S} \cap \langle z_{1}\cdots{z}_{n+1} - 1 \rangle = \{0\}$ (we cannot factor $z_{1}{\cdots}{z}_{n+1}$ out of any $\mys_{\widehat{\lambda}}$), 
so that $\overline{\Theta}|_{_{\mathcal{S}}}$ is one-to-one.  Since it is 
clearly onto as well, then  $\overline{\Theta}|_{_{\mathcal{S}}}$ is a 
ring isomorphism from 
$\mathcal{S}$ onto $\mathbb{Z}[\Lambda]^{W}$, as claimed.

For any dominant weight $\lambda \sum_{i \in I}a_{i}\omega_{i}$, let $\mathfrak{a}(\lambda) = (a_{1},a_{2},\ldots,a_{n},0)$.  
The subring generated by $\{\mye^{\mathfrak{a}(\lambda)}\}_{\lambda \in \Lambda^{+}}$ is in $\mathcal{S}$, since each elementary symmetric function is a product of Schur functions $\mye_{k} = \mys_{\mathbf{1}_{k}}$ where $1 \leq k \leq n$.  
Since $\overline{\Theta}|_{_{\mathcal{S}}}(\mye^{\mathfrak{a}(\lambda)}) = \psi_{\lambda}$ and $\{\psi_{\lambda}\}_{\lambda \in \Lambda^{+}}$ is a $\mathbb{Z}$-basis for $\mathbb{Z}[\Lambda]^{W}$, then $\{\mye^{\mathfrak{a}(\lambda)}\}_{\lambda \in \Lambda^{+}}$ is a $\mathbb{Z}$-basis for $\mathcal{S}$.  
However, $\{\mym_{\widehat{\lambda}}\}_{\lambda \in \Lambda^{+}}$ does not in general span $\mathcal{S}$.  
For example, when $n=1$, we can write $\mys_{(2,0)} = \mym_{(2,0)}+\mym_{(1,1)}$, but $\mym_{(1,1)}$ is not expressible as a $\mathbb{Z}$-linear combination of $\{\mym_{(k,0)}\}_{k \geq 0}$. 

We would like to understand the relationship between Kostka numbers in the usual sense and the $\myA_{n}$-Kostka numbers.  
For any dominant weight $\lambda = \sum_{i \in I}a_{i}\omega_{i}$, let $\mathfrak{p} := \widehat{\lambda}$.  
Let $\mathscr{P}(\mathfrak{p}) := \{\mbox{partitions } \mathfrak{q} \mbox{ with } K_{\mathfrak{p},\mathfrak{q}} \not= 0\}$.  
It is well-known that $K_{\mathfrak{p},\mathfrak{q}} \not= 0$ for the partition $\mathfrak{q}$ if and only $|\mathfrak{q}| = |\mathfrak{p}|$ and $\mathfrak{q} \triangleleft \mathfrak{p}$.  
We also set $\mathscr{W}(\lambda) := \{\nu \in \Lambda^{+}\, |\, \nu \leq \lambda\}$. 
Now $\displaystyle \overline{\Theta}(\mys_{\mathfrak{p}}) = \overline{\Theta}\left(\sum_{\mathfrak{q} \in \mathscr{P}(\mathfrak{p})}K_{\mathfrak{p},\mathfrak{q}} \mym_{\mathfrak{q}}\right) = \sum_{\mathfrak{q} \in \mathscr{P}(\mathfrak{p})}K_{\mathfrak{p},\mathfrak{q}} \zeta_{\smallweight(\mathfrak{q})} = \sum_{\nu \in \mathscr{W}(\lambda)} d_{\lambda,\nu}\zeta_{\nu}$.  
In particular, $\weight: \mathscr{P}(\mathfrak{p}) \longrightarrow \mathscr{W}(\lambda)$ is a surjective set mapping.  
It is easy to see that when $\weight(\mathfrak{q}) = \weight(\mathfrak{q}')$ with $|\mathfrak{q}| = |\mathfrak{q}'|$, then $\mathfrak{q} = \mathfrak{q}'$, and hence that this set mapping is injective.  
To describe the inverse $(\weight|_{\mathscr{P}(\mathfrak{p})})^{-1}$, we begin with an observation that involves straightforward computations: A weight $\nu = \sum b_{i}\omega_i$ is in $\mathscr{W}(\lambda)$ if and only if (1) $b_i \geq 0$ for all $i \in I$, (2) $\sum_{i=1}^{j}i(n+1-j)(a_i - b_i) + \sum_{i=j+1}^{n}j(n+1-i)(a_i - b_i) \geq 0$ for all $j \in I$, and (3) $\sum_{i=1}^{j}i(n+1-j)(a_i - b_i) + \sum_{i=j+1}^{n}j(n+1-i)(a_i - b_i)$ is divisible by $n+1$ for some $j \in I$ (in which case it can be shown that this quantity is divisible by $n+1$ for all $j \in I$). 
In this case, we have $\lambda - \nu = \sum_{i \in I}k_i(\nu)\alpha_i$ where the nonnegative integers $k_i(\nu)$ are given by $\frac{1}{n+1}\left[\sum_{i=1}^{j}i(n+1-j)(a_i - b_i) + \sum_{i=j+1}^{n}j(n+1-i)(a_i - b_i)\right]$. 
With this in mind, we define $\partition_{\lambda}: \mathscr{W}(\lambda) \longrightarrow \mathscr{P}(\mathfrak{p})$ by the rule $\partition_{\lambda}(\nu) = \widehat{\nu} + k_{n}(\nu)\mathbf{1}$. 
Using conditions (1), (2), and (3), it is straightforward to argue that $\partition_{\lambda}(\nu) \in \mathscr{P}(\mathfrak{p})$ for each $\nu \in \mathscr{W}(\lambda)$ and that $\partition_{\lambda} = (\weight|_{\mathscr{P}(\mathfrak{p})})^{-1}$. 
(More generally, one can see that a weight $\mu = \sum b_{i}\omega_i$ is in $\Pi(\lambda)$ if and only if the above conditions (2) and (3) hold. 
Then $\lambda - \mu = \sum_{i \in I}k_i(\nu)\alpha_i$ where the nonnegative integers $k_i(\mu)$ are given by the same formula as above. 
Then $\partition_{\lambda}$ is a bijective mapping from $\Pi(\lambda)$ to the set of nonnegative integer sequences $\mathfrak{q} = (q_{1},\ldots,q_{n+1})$ for which the coefficient of $z^{\mathfrak{q}}$ in $\mys_{\mathfrak{p}}$ is nonzero, and its inverse is the $\weight$ mapping.) 
For all $\mathfrak{q} \in \mathscr{P}(\mathfrak{p})$ and $\nu \in \mathscr{W}(\lambda)$, we therefore obtain that $K_{\mathfrak{p},\mathfrak{q}} = d_{\lambda,\nu}$ if and only if $\weight(\mathfrak{q}) = \nu$ if and only $\partition_{\lambda}(\nu) = \mathfrak{q}$. 
In this way, we have an interpretation of the nonnegative integers $d_{\lambda,\nu}$ as a count of combinatorial objects. 

Next we give some remarks on inverse images under $\Theta$.  
The following facts about the preimages of $\chi_{_{\lambda}}$, $\zeta_{\lambda}$, and $\psi_{\lambda}$ under $\overline{\Theta}$ are straightforward. 
\[
\overline{\Theta}^{-1}(\chi_{_{\lambda}}) \cap \{\mys_{\mathfrak{p}}|\mbox{ \small partitions } \mathfrak{p}\}  =  \{\mys_{\widehat{\lambda}+k\mathbf{1}}\}_{k \geq 0}  = \{\mym_{\mathbf{1}}^{k} \mys_{\widehat{\lambda}}\}_{k \geq 0}\]
\[\overline{\Theta}^{-1}(\zeta_{\lambda}) \cap \{\mym_{\mathfrak{p}}|\mbox{ \small partitions } \mathfrak{p}\} = \{\mym_{\smallpartition_{\lambda'}(\lambda)}\, |\, \lambda \leq \lambda' \in \Lambda^{+}\}  =  \{\mym_{\mathbf{1}}^{k} \mym_{\widehat{\lambda}}\}_{k \geq 0}\]
\[\overline{\Theta}^{-1}(\psi_{\lambda}) \cap \{\mye^{\mathfrak{a}}|\mbox{ \small nonnegative integer sequences } \mathfrak{a}\}  =  \{\mye^{\mathfrak{a}(\lambda)} e_{n+1}^{k}\}_{k \geq 0} = \{\mym_{\mathbf{1}}^{k} \mye^{\mathfrak{a}(\lambda)}\}_{k \geq 0}
\]
Next, for $\lambda \in \Lambda^{+}$, consider the subgroup $\mathscr{L}_\lambda := \{\phi \in \mathbb{Z}[\Lambda]\, |\, \widehat{\Pi}(\phi) \subseteq \Pi(\lambda)\}$ of $\mathbb{Z}[\Lambda]$. 
Using an inductive argument as in the proof of \BasisProp, one can see that its subgroup $\mathscr{L}_\lambda^{W}$ of $W$-invariants is spanned by those monomial $W$-symmetric functions $\zeta_{\nu}$ for which $\nu \in \mathscr{W}(\lambda)$. 
We define a group homomorphism $\Psi_{\lambda}: \mathscr{L}_\lambda \longrightarrow \mathbb{Z}[z_1,\ldots,z_{n+1}]$ by requiring that $\Psi_{\lambda}(e^{\mu}) := z^{\smallpartition_{\lambda}(\mu)}$ for each $\mu \in \Pi(\lambda)$ and extending $\mathbb{Z}$-linearly.  
Since the $\partition_{\lambda}$ mapping is injective, then $\Psi_\lambda$ is as well. 
Its inverse is $\Theta|_{\Psi_\lambda(\mathscr{L}_\lambda)}$. 
One can also check that $\Psi_\lambda$ preserves the $W$-action, so that $\Psi_\lambda(\mathscr{L}_{\lambda}^{W}) \subset \mathbb{Z}[z_1,\ldots,z_{n+1}]^{W}$. 
Of course, $\Psi_\lambda(\zeta_{\nu}) = \mym_{\smallpartition_{\lambda}(\nu)}$ for all $\nu \in \mathscr{W}(\lambda)$. 
In summary, the subring $\mathcal{S}$ of $\mathbb{Z}[z_1,\ldots,z_{n+1}]^{W}$ explored above nicely correlates certain Schur functions with the $\myA_{n}$-Weyl bialternants and certain elementary symmetric functions with the elementary $\myA_{n}$-Weyl symmetric functions.  On the other hand, the subring $\mathcal{S}$ does not so nicely correlate monomial symmetric functions with monomial $\myA_{n}$-Weyl symmetric functions; however, the subgroups $\Psi_\lambda(\mathscr{L}_{\lambda}^{W})$ do. 

We close this discussion with some remarks about Weyl symmetric function analogs of other standard symmetric functions. 
We do not know of a general Weyl symmetric function analog of the complete homogeneous symmetric functions.  
In classical symmetric function theory, a so-called ``skew'' Schur function $\mys_{\mathfrak{p}/\mathfrak{q}}$ is defined as $\sum c_{\mathfrak{q},\mathfrak{n}}^{\mathfrak{p}} \mys_{\mathfrak{n}}$, where the sum is over all partitions $\mathfrak{n}$ and $c_{\mathfrak{q},\mathfrak{n}}^{\mathfrak{p}}$ is a Littlewood-Richardson coefficient. 
An obvious analog of skew Schur functions in the general Weyl symmetric function setting would be to define $\chi_{_{\lambda/\nu}} := \sum_{\mu \in \Lambda^{+}} c_{\mu,\nu}^{\lambda} \chi_{_{\mu}}$, where the $c_{\mu,\nu}^{\lambda}$'s are the product decomposition coefficients from the identity $\chi_{_{\lambda}} \chi_{_{\mu}} = \sum_{\nu \in \Lambda^{+}} c_{\mu,\nu}^{\lambda} \chi_{_{\nu}}$. 
We do not know how well these or related objects have been studied.\hfill\QED

{\bf Interactions with Lie algebra representation theory.} 
So far we have required no Lie algebra representation theory. However, some combinatorial questions about Weyl bialternants and $W$-symmetric functions are naturally 
addressed in that context.  For example, while a combinatorial proof of unimodality of all Dynkin polynomials has not yet been discovered (cf.\  \UnimodalProblem), there is an algebraic proof, see \UnimodalTheorem\ below. 

Weyl's character formula (WCF), as stated in \TFAEProp\ below, is a crucial point of contact between $W$-symmetric function theory and Lie algebra representation theory.  
\TFAETheorems\ allow us to conclude almost immediately that Weyl's character formula is equivalent to Kostant's multiplicity formula (KMF) for weight space dimensions and Freudenthal's multiplicity formula (FMF) for weight space dimensions, a fact we record below as \TFAEProp.  
The terminology and notation we use to set up and state the next several results borrow from  the discussion of supporting graphs in \S \SplittingSection\ (see also \cite{DonSupp}). Let $\mathfrak{g}$ be a complex finite-dimensional semisimple Lie algebra associated with the root system $\Phi$ and the given choice of simple roots. Let $V$ be a complex finite-dimensional $\mathfrak{g}$-module.  The character of the representation $V$ is defined to be $\mbox{char}(V) := \sum_{\mu \in \Lambda}(\dim V_{\mu})e^{\mu}$.  

\noindent 
{\bf \TFAEProp}\ \ {\sl Fix $\lambda \in \Lambda^{+}$.  Let $V(\lambda)$ be an irreducible finite-dimensional $\mathfrak{g}$-module with highest weight $\lambda$. Then the following are equivalent:}\\ 
\hspace*{0.25in}(WCF)\ \ $\mbox{char}(V(\lambda)) = \chi_{_{\lambda}}$, {\sl so for all $\mu \in \Lambda$, $\dim V(\lambda)_{\mu} = d_{\lambda,\mu}$.}\\
\hspace*{0.25in}(KMF)\ \ {\sl For all $\mu \in \Lambda$,} $\dim V(\lambda)_{\mu} = \sum_{w \in W}\det(w)\mathcal{P}(w(\lambda)-\mu + w(\varrho)-\varrho)$.\\
\hspace*{0.25in}(FMF)\ \ {\sl For all weights $\mu$ not in $\Pi(\lambda)$, we have $\dim V(\lambda)_{\mu} = 0$.  We have $\dim V(\lambda)_{\lambda} = 1$,\\ \hspace*{0.76in} 
and for all $\mu \in \Pi(\lambda)$ with $\mu < \lambda$ we have} \[\hspace*{0.25in}(||\lambda+\varrho||^{2} - ||\mu+\varrho||^{2}) \dim V(\lambda)_{\mu} = 2\sum_{\alpha \in \Phi^{+}}\sum_{k \geq 1}\langle \alpha,\mu+k\alpha, \rangle \dim V(\lambda)_{\mu+k\alpha}.\]

{\em Proof.} \KostantTheorem\ shows WCF $\iff$ KMF.  \FreudTheorem\ shows WCF $\iff$ FMF.\hfill\QED

So once any one of the Lie theoretic assertions WCF, KMF, or FMF from \TFAEProp\ has been established, we automatically have all of them.  Of course, these results are classical and well-known: 

\noindent 
{\bf \WeylsTheorem}\ \ {\sl Each of the equivalent statements of \TFAEProp\ is true.}\hfill\QED 

For a proof of WCF, see \S 24 of \cite{Hum} or \S 25 of \cite{FH}. Weyl's original analytic proof is recapitulated in \S 26 of \cite{FH}.  In \cite{Stem}, Stembridge uses combinatorial methods to construct certain graphs related to crystal bases.  Weight generating functions for these graphs are Weyl bialternants on the one hand, but on the other hand, since they are related to crystal bases, they are characters for the associated irreducible $\mathfrak{g}$-module.  WCF is an immediate  consequence of this reasoning.  We give a new version of this result in \S\CrystalSection.  The equivalence of WCF and KMF appears to have been first observed by Cartier \cite{Cartier}.  Kostant's original Lie-representation-theoretic proof did not use WCF, see \cite{Kostant}.  FMF was first derived using Lie representation theoretic (as opposed to $W$-symmetric function) techniques, for such a proof see for example \S 22 of \cite{Hum}.  

The following unimodality corollary of WCF is attributed to Dynkin \cite{Dynkin}, although the particular argument given below, which makes use of a ``principal three-dimensional subalgebra,'' is due to Proctor \cite{PrEur}. 

\noindent 
{\bf \UnimodalTheorem\ (Unimodality Theorem)}\ \ 
{\sl The symmetric Dynkin polynomial of 
\TFAEConsequences.1 is unimodal.}

{\em Proof.}  Our notation borrows from the discussion of supporting graphs in \S 
\SplittingSection. Let 
$V(\lambda)$ be an irreducible $\mathfrak{g}$-module with maximal 
vector of (dominant) weight $\lambda$.  The 
Chevalley generators 
$\{\myqx_{i},\myqy_{i},\myqh_{i}\}_{i \in I}$ for 
$\mathfrak{g}$ respect the given choice of simple roots. For each $i 
\in I$, let $c_{i} := 2\langle \omega_{i},\varrho^{\vee} \rangle$, and 
set $\myqx := \sum c_{i}\myqx_{i}$, $\myqy := \sum \myqy_{i}$, and 
$\myqh := \sum c_{i}\myqh_{i}$, with each sum over all $i \in I$. 
It follows readily that $[\myqx,\myqy] = \myqh$.  From the identity 
$\sum_{i \in I}c_{i}M_{ji} = 2$ it follows that 
$[\myqh,\myqx] = 2\myqx$  and 
$[\myqh,\myqy] = -2\myqy$.  Therefore  
the Lie subalgebra $\mathfrak{s}$ of $\mathfrak{g}$ spanned 
by $\{\myqx, \myqy, \myqh\}$ is isomorphic to 
$\mathfrak{sl}(2,\mathbb{C})$.  

We use the notation $V'$ when 
regarding $V(\lambda)$ to be an 
$\mathfrak{s}$-module under the induced action. Then by basic 
$\mathfrak{sl}(2,\mathbb{C})$ representation theory (cf \S 7 of 
\cite{Hum}), we get that $V' = \bigoplus_{j=1}^{l}V'(k_{j})$ for some 
positive integer $l$ and some nonnegative integers $k_{j}$ ($1 \leq j 
\leq l)$, 
where $V'(k_{j})$ is an irreducible 
$\mathfrak{sl}(2,\mathbb{C})$-module with highest weight $k_{j} \in 
\mathbb{Z}_{\geq 0}$.  From Theorem 7.2 of \cite{Hum}, $m \in 
\Pi(k_{j})$ if and only if $m \in 
\{-k_{j},-k_{j}+2,\ldots,k_{j}-2,k_{j}\}$, in which case the 
$m$-weight space for $V'(k_{j})$ is one-dimensional. 
Let $m \in \mathbb{Z}$, and let $V'_{m}$ be the $m$-weight space 
for the $\mathfrak{s}$-module $V'$.  From the structure of 
$\mathfrak{sl}(2,\mathbb{C})$-modules as described in \S 7 of 
\cite{Hum}, we see that $\dim(V'_{m}) = |\{j\, |\, 1 \leq j \leq l, m 
\in \Pi(k_{j})\}|$. Therefore 
$\dim(V'_{m-2}) \leq 
\dim(V'_{m})$ if $m \leq 0$ and $\dim(V'_{m}) \geq 
\dim(V'_{m+2})$ if $m \geq 0$.  

Any weight basis for the $\mathfrak{g}$-module 
$V(\lambda)$ is also a weight basis for the $\mathfrak{s}$-module 
$V'$. 
Moreover, for any $\mu \in \Lambda$, check that if $v 
\in V_{\mu}$, then $\myqh.v = 2\langle \mu,\varrho^{\vee}\rangle v$.  
Therefore $\dim(V'_{m}) = \sum \dim(V(\lambda)_{\mu}) = 
\sum d_{\lambda,\mu}$, where the two summations are  
over all $\mu$ for which $2\langle \mu,\varrho^{\vee} \rangle = m$, 
and where we invoke WCF in that latter equality.  
Say for some such $\mu$ that $d_{\lambda,\mu} \not= 0$. Then 
$2\langle \mu,\varrho^{\vee} \rangle$ is an integer.  Moreover,  $\mu 
\in \Pi(\lambda)$ means that $w_{0}(\lambda) \leq 
\mu \leq \lambda$.  In particular, see that 
$-2\langle \lambda,\varrho^{\vee} 
\rangle \leq 2\langle \mu,\varrho^{\vee} \rangle \leq 2\langle 
\lambda,\varrho^{\vee} \rangle$ and that 
$2\langle \mu,\varrho^{\vee} \rangle$ has the same 
parity as $2\langle \lambda,\varrho^{\vee} \rangle$.   
Then, \[\dim(V') =  
\dim(V'_{-2\langle \lambda,\varrho^{\vee} \rangle}) + 
\dim(V'_{-2\langle \lambda,\varrho^{\vee} \rangle + 2}) + 
\cdots + \dim(V'_{2\langle \lambda,\varrho^{\vee} \rangle - 2}) + 
\dim(V'_{2\langle \lambda,\varrho^{\vee} \rangle}),\]
where the terms on the right-hand side form a unimodal sequence.  
For any $p \in \{0,1,2,\ldots,2\langle \lambda,\varrho^{\vee} 
\rangle\}$, then $\dim(V'_{-2\langle \lambda,\varrho^{\vee} \rangle + 
2p}) = \sum d_{\lambda,\mu}$, where the sum on the right-hand side is 
over all $\mu \in \Lambda$ with $2\langle \mu,\varrho^{\vee} 
\rangle = -2\langle \lambda,\varrho^{\vee} \rangle + 2p$, i.e.\ over 
all $\mu \in \Lambda$ with $\langle \mu+\lambda,\varrho^{\vee} 
\rangle = p$.  But this is precisely the coefficient of $q^{p}$ in 
the Dynkin polynomial of \TFAEConsequences.1.  Therefore the said 
Dynkin polynomial is unimodal.\hfill\QED

\newpage
\noindent 
{\Large \bf \S \PosetSection.\ Edge-colored posets and 
${\Phi}$-structured posets.} 

This section furnishes the combinatorial language and poset prerequisite results needed for our discussion of splitting posets in subsequent sections. 
At the outset, none of the $W$-symmetric function theory from \S \WeylSection, beyond the initial definitions, is used.  
Later, we begin making connections between the combinatorial structures of this section and the material of the previous section. 
It seems safe for readers to browse the definitions and results given here and refer back to this section as needed later on. 
From here on, all graphs and posets are assumed to be finite unless stated otherwise.

{\bf Some language and notation.} 
Our combinatorial conventions here largely 
follow  \S 2 of \cite{ADLMPPW} and \cite{DonDistributive}. 
Let $I$ be any set.  An {\em edge-colored directed graph with
edges colored by the set $I$} is a directed graph $R$ with vertex set 
$\mathcal{V}(R)$ and directed-edge set $\mathcal{E}(R)$ together
with a function $\ecolor_{R}\, :\, \mathcal{E}(R) \longrightarrow I$ 
assigning to each edge of $R$ an element (``color'') from the
set $I$.  If an edge $\selt \rightarrow \telt$ in $R$ is assigned 
color $i \in I$, we write $\selt \myarrow{i} \telt$.  
For $i \in I$, we let 
$\mathcal{E}_{i}(R)$ denote the set of edges in $R$ of color $i$, 
so $\mathcal{E}_{i}(R) = \ecolor_{R}^{-1}(i)$.  All such directed graphs in this monograph will be taken to have no multiple edges between nodes and no loops. 
A {\em path} in $R$ from $\xelt$ to $\yelt$ is a sequence $\mathcal{P} = (\xelt = \xelt_{0}, \xelt_{1}, \ldots, \xelt_{k} = \yelt)$ wtih colors $(i_{j})_{j=1}^{k}$ and such that 
for $1 \leq j \leq k$ we have either $\xelt_{j-1} \myarrow{i_j} \xelt_{j}$ or $\xelt_{j} \myarrow{i_j} \xelt_{j-1}$.  This path has {\em length} $k$, and we allow paths to have length $0$.  For any $i \in I$, we let $a_{i}(\mathcal{P}) := |\{j \in \{1,2,\ldots,k\}\, |\, \xelt_{j-1} \myarrow{i_j} \xelt_{j} \mbox{ in $\mathcal{P}$ and } i_j = i\}|$ (a count of ``ascending'' edges of color $i$ in the path) and $d_{i}(\mathcal{P}) := |\{j \in \{1,2,\ldots,k\}\, |\, \xelt_{j} \myarrow{i_j} \xelt_{j-1} \mbox{ in $\mathcal{P}$ and } i_j = i\}|$ (a count of ``descending'' edges of color $i$). 

If $J$ is a subset of $I$,
remove all edges from $R$ whose colors are not in $J$; connected
components of the resulting edge-colored directed graph are called
{\em J-components} of $R$. 
For any $\telt$ in $R$ and any $J \subseteq I$, we let 
$\mathbf{comp}_{J}(\telt)$ denote the $J$-component of $R$ 
containing $\telt$.  
The {\em dual} $R^{*}$ is the edge-colored directed graph whose vertex 
set $\mathcal{V}(R^{*})$ is the set of symbols 
$\{\telt^{*}\}_{\telt{\in}R}$ 
together with colored edges 
$\mathcal{E}_{i}(R^{*}) := 
\{\telt^{*} \myarrow{i} \selt^{*}\, |\, 
\selt \myarrow{i} \telt \in \mathcal{E}_{i}(R)\}$ 
for each $i \in I$. Let $Q$ be another edge-colored
directed graph with edge colors from $I$. 
Regarding $R$ and $Q$ to have disjoint 
vertex sets, then the {\em disjoint sum} $R
\oplus Q$ is the edge-colored directed graph whose vertex set is 
the disjoint union $\mathcal{V}(R) \disjointunion \mathcal{V}(Q)$ 
and whose 
colored edges are $\mathcal{E}_{i}(R) \disjointunion 
\mathcal{E}_{i}(Q)$ 
for each $i \in I$.  
We denote by $R \cup Q$ the edge-colored directed graph with $\mathcal{V}(R \cup Q) := \mathcal{V}(R) \cup \mathcal{V}(Q)$ and $\mathcal{E}_{i}(R \cup Q) := \mathcal{E}_{i}(R) \cup \mathcal{E}_{i}(Q)$ for each $i \in I$. 
If $\mathcal{V}(Q) 
\subseteq \mathcal{V}(R)$ and $\mathcal{E}_{i}(Q) \subseteq 
\mathcal{E}_{i}(R)$ for each $i \in
I$, then $Q$ is an {\em edge-colored subgraph} of $R$. 
Let $R \times Q$ denote the edge-colored directed graph whose vertex 
set is the
Cartesian product $\{(\selt,\telt)|\selt \in R,\telt \in Q\}$ and with
colored edges $(\selt_{1},\telt_{1}) \myarrow{i}
(\selt_{2},\telt_{2})$ if and only if $\selt_{1} = \selt_{2}$ in
$R$ with $\telt_{1} \myarrow{i} \telt_{2}$ in $Q$ or $\selt_{1}
\myarrow{i} \selt_{2}$ in $R$ with $\telt_{1} = \telt_{2}$ in $Q$.
Two edge-colored directed graphs are {\em isomorphic} if there is a
bijection between their vertex sets that preserves edges and edge
colors.  If $R$ is an edge-colored directed graph with edges 
colored by the 
set $I$, and if $\sigma\, :\, I \longrightarrow I'$ is a mapping of 
sets, then we let $R^{\sigma}$ be the edge-colored directed graph with 
edge color function $\ecolor_{R^{\sigma}} := \sigma \circ 
\ecolor_{R}$. 
We call $R^{\sigma}$ 
a {\em recoloring} of $R$. Observe that $(R^{*})^{\sigma} \cong 
(R^{\sigma})^{*}$. Suppose $I$ indexes a set of simple roots 
for a root 
system $\Phi$ as in \S \WeylSection.  The permutation $\sigma_{0}$ 
discussed there is related to the longest Weyl group element $w_{0}$.  
The $\sigma_{0}$-{\em recolored dual}  
$R^{\bowtie}$ is the edge-colored directed graph 
$(R^{\sigma_{0}})^{*} \cong (R^{*})^{\sigma_{0}}$. 
More colloquially, we sometimes call $R^{\bowtie}$ the {\em bow tie} of $R$. 
(Elsewhere, we have denoted this poset by the notation $R^{\triangle}$ rather than $R^{\bowtie}$.  
We now prefer the bow tie symbol ``$\bowtie$'' because its $\mathbb{Z}_{2} \oplus \mathbb{Z}_{2}$ symmetry group is suggestive of the two independent involutive processes being used to obtain $R^{\bowtie}$ from $R$.) 
The following observations are straightforward. 

\noindent 
{\bf \FirstOpsLemma}\ \ {\sl Up to poset isomorphism, the Cartesian product $\times$ and the disjoint sum $\oplus$ are commutative and associative binary operations, and $\times$ distributes over $\oplus$. 
Both $\times$ and $\oplus$ interact with $^{*}$, $^{\sigma}$, and $^{\bowtie}$ in the natural way (i.e.\ $(R_{1} \times R_{2})^{*} \cong R_{1}^{*} \times R_{2}^{*}$, $(R_{1} \oplus R_{2})^{\sigma} \cong R_{1}^{\sigma} \oplus R_{2}^{\sigma}$, etc). 
The dual $^{*}$ and the bow tie $^{\bowtie}$ operations are involutions.}\hfill\QED

Identify a partially ordered set $R$ with its {\em Hasse diagram}, 
that is, the 
directed graph whose edges depict the {\em covering relations} for the 
poset: for elements $\selt$ and $\telt$ in $R$ the directed edge 
$\selt \rightarrow \telt$ means that $\selt < \telt$ and 
if $\selt \leq \xelt \leq \telt$ 
then $\selt = \xelt$ or $\xelt = \telt$. 
Via the Hasse diagram, terminology that applies to  
directed graphs will be applied to posets (connected, edge-colored, 
dual, etc). 
We say $R$  
is {\em ranked} if there exists a surjective function $\rho : R 
\longrightarrow \{0,1,\ldots,l\}$ such that $\rho(\selt) + 1 = 
\rho(\telt)$ whenever $\selt \rightarrow \telt$; in this case  
$\rho$ is a {\em rank function}, 
$\rho(\xelt)$ is the {\em rank} of any element $\xelt$ in $R$, and 
the number $l$ is the {\em length} of $R$ with respect to $\rho$.  
This ranked poset is {\em rank symmetric} if $|\rho^{-1}(i)| = 
|\rho^{-1}(l-i)|$ for all $i \in \{0,1,\ldots,l\}$. 
It is {\em rank unimodal} if there is an $m$ such that
$|\rho^{-1}(0)| \leq |\rho^{-1}(1)| \leq \cdots \leq
|\rho^{-1}(m)| \geq |\rho^{-1}(m+1)| \geq \cdots \geq
|\rho^{-1}(l)|$. 
It has a {\em symmetric chain 
decomposition} or {\em SCD} if there exist chains 
$R_{1}, \ldots, R_{k}$ in $R$ 
such that (1) as a set $R = R_{1} \disjointunion \cdots 
\disjointunion 
R_{k}$ (disjoint 
union), and (2) for $1 \leq i \leq k$, if 
$\xelt_{i}$ and $\yelt_{i}$ are respectively the minimal and 
maximal elements of the chain $R_{i}$, then $\rho(\yelt_{i}) + 
\rho(\xelt_{i}) = l$ and $\rho(\yelt_{i}) - \rho(\xelt_{i})$ is the 
length of $R_{i}$.  Each $R_{i}$ is a {\em symmetric chain} 
for the SCD. 
If $R$ has an SCD, then one can see that 
$R$ is rank symmetric and rank unimodal. 
Define the {\em rank generating function} for a ranked poset $R$ to be $\RGF(R,q) := \sum_{i=0}^{l}|\rho^{-1}(i)|q^{i} = \sum_{\xelt \in R} q^{\rho(\xelt)}$.  

A {\em lattice} $L$ is a poset for which any two given elements $\xelt$ and $\yelt$ of $L$ have a (unique) least upper bound, denoted $\xelt \vee \yelt$ and called their {\em join}, and a (unique) greatest lower bound, denoted $\xelt \wedge \yelt$ and called their {\em meet}. 
Such a lattice is necessarily connected and has a unique maximal element $\mathbf{max}(L)$ and a unique minimal element $\mathbf{min}(L)$. 
This lattice is {\em modular} by definition if and only if $L$ is ranked and $\rho(\xelt \wedge \yelt) + \rho(\xelt \vee \yelt) = \rho(\xelt) + \rho(\yelt)$ for any $\xelt, \yelt \in L$. 
The lattice $L$ is {\em distributive} if and only if meets distribute over joins and vice-versa;  
that is, $\xelt \wedge (\yelt \vee \zelt) = (\xelt \wedge \yelt) \vee (\xelt \wedge \zelt)$ and $\xelt \vee (\yelt \wedge \zelt) = (\xelt \vee \yelt) \wedge (\xelt \vee \zelt)$ for any given $\xelt, \yelt, \zelt \in L$. 
Any distributive lattice is modular, but not all modular lattices are distributive. 

In a distributive lattice $L$, an element $\xelt$ is a {\em join irreducible} if $\xelt$ covers exactly one other element of $L$. 
Let $\mathbf{j}(L)$ be the set of join irreducible elements of $L$ with partial order induced by $L$; we call $\mathbf{j}(L)$ the {\em poset of join irreducibles} of $L$. 
Now let $P$ be a poset.  
A subset $\mathcal{X}$ from $P$ is an {\em order ideal} if, for any $\xelt \in \mathcal{X}$ and any $\xelt'$ in $P$ with $\xelt' \leq \xelt$, we have $\xelt' \in \mathcal{X}$. 
Let $\mathbf{J}(P)$ be the set of order ideals from $P$ partially ordered by subset containment. 
Since meets and joins in $\mathbf{J}(P)$ are (respectively) just intersections and unions of sets, it is easy to see that $\mathbf{J}(P)$ is a distributive lattice. 
What is sometimes called The Fundamental Theorem of Finite Distributive Lattices asserts that for any distributive lattice $L$ and any poset $P$ we have:
\[\mathbf{J}\left(\rule[-2.5mm]{0mm}{7mm}\mathbf{j}(L)\right)\ \mbox{\sl is isomorphic to}\ L,\ \mbox{\sl and}\ \ \mathbf{j}\left(\rule[-2.5mm]{0mm}{7mm}\mathbf{J}(P)\right)\ \mbox{\sl is isomorphic to}\ P.\] 

Now suppose $L$ is a distributive lattice whose edges are colored by $I$. 
Assign color $i$ to a join irreducible $\xelt$ if $i$ is the color of the edge beneath $\xelt$. 
We denote the resulting vertex-colored poset of join irreducibles by $\mathbf{j}_{color}(L)$. 
Now, for any poset $P$ that is vertex-colored by $I$, note that there is a covering relation $\mathcal{X} \rightarrow \mathcal{Y}$ in $\mathbf{J}(P)$ if and only if $\mathcal{X} \subset \mathcal{Y}$ and there exists some $\yelt \in \mathcal{Y} \setminus \mathcal{X}$ such that $\mathcal{Y}=\mathcal{X}\cup\{\yelt\}$. 
If $j$ is the color of the vertex $\yelt$, then assign color $j$ to the edge $\mathcal{X} \rightarrow \mathcal{Y}$ in the distributive lattice $\mathbf{J}(P)$. 
We denote the resulting edge-colored distributive lattice by $\mathbf{J}_{color}(P)$. 
Now, $\mathbf{J}_{color}(P)$ has the property that in any ``diamond'' of colored edges, parallel edges have the same colors; that is, if \parbox{1.1cm}{\begin{center}
\setlength{\unitlength}{0.2cm}
\begin{picture}(4,3.5)
\put(2,0){\circle*{0.5}} \put(0,2){\circle*{0.5}}
\put(2,4){\circle*{0.5}} \put(4,2){\circle*{0.5}}
\put(0,2){\line(1,1){2}} \put(2,0){\line(-1,1){2}}
\put(4,2){\line(-1,1){2}} \put(2,0){\line(1,1){2}}
\put(0.75,0.55){\em \small k} \put(2.7,0.7){\em \small l}
\put(0.7,2.7){\em \small i} \put(2.75,2.55){\em \small j}
\end{picture} \end{center}} is an edge-colored subgraph of the order  
diagram of $\mathbf{J}_{color}(P)$, then $i = l$ and $j = k$.  
We say the distributive lattice $\mathbf{J}_{color}(P)$ is {\em diamond-colored}. 
The following is a colorized version of The Fundamental Theorem of Finite Distributive Lattices: Assume that a poset $P$ is vertex-colored by $I$, that a distributive lattice $L$ is edge-colored by $I$, and that $L$ is diamond-colored. 
Then:
\[\mathbf{J}_{color}\left(\rule[-2.5mm]{0mm}{7mm}\mathbf{j}_{color}(L)\right)\ \mbox{\sl is isomorphic to}\ L,\ \mbox{\sl and}\ \ \mathbf{j}_{color}\left(\rule[-2.5mm]{0mm}{7mm}\mathbf{J}_{color}(P)\right)\ \mbox{\sl is isomorphic to}\ P,\] 
where these isomorphisms preserve colors as well as poset structure.

For the remainder of this section, $R$ is a ranked poset whose Hasse diagram edges are colored by a set $I$. 
We let $\rho_{i}(\xelt)$ denote the rank of an element $\xelt$ in $R$ within its $i$-component $\mathbf{comp}_{i}(\xelt)$ and $l_{i}(\xelt)$ denote the length of $\mathbf{comp}_{i}(\xelt)$.  
We use $\delta_{i}(\xelt)$ to denote the quantity $l_{i}(\xelt) - \rho_{i}(\xelt)$, and we set $m_{i}(\xelt) := 2\rho_{i}(\xelt)-l_{i}(\xelt) = \rho_{i}(\xelt) - \delta_{i}(\xelt)$. 
For $J \subseteq I$, $\xelt$ is $J$-{\em prominent} 
if $\delta_{j}(\xelt) = 0$ for all $j \in J$, 
and $\xelt$ is $J$-{\em maximal} (respectively $J$-{\em minimal}) if for any edge $\xelt \myarrow{i} \yelt$ (resp.\ $\yelt \myarrow{i} \xelt$) in $R$ we have $i \not\in J$. 
When $J = I$, we simply call such $\xelt$ ``prominent,'' ``maximal'' or etc. 

For this paragraph,  
the elements of $R$ will be denoted by 
$v_{1},\ldots,v_{n}$, so $n = |R|$.  For an integer $k \geq 0$, 
let $\bigwedge^{k}(R)$ denote the set of all $k$-element subsets of 
the vertex set of $R$.  If $k > n$, then $\bigwedge^{k}(R) = 
\emptyset$.  If $k = 0$ or $k = n$ then $\bigwedge^{k}(R)$ is a set 
with one element.  For $\selt, \telt \in \bigwedge^{k}(R)$, write 
$\selt \myarrow{i} \telt$ if and only if $\selt$ and $\telt$ differ by 
exactly one element in the sense that $(\selt - \telt, \telt - 
\selt) = (\{v_{p}\},\{v_{q}\})$ and $v_{p} \myarrow{i} v_{q}$ in $R$.  
Use the notation $\bigwedge^{k}(R)$ to refer to this edge-colored 
directed graph, which we call the $k$th {\em exterior power} of $R$.  
Similarly let $\mathbb{S}^{k}(R)$ denote the set of all $k$-element 
multisubsets of the vertex set of $R$ and define colored, directed 
edges $\selt \myarrow{i} \telt$ between elements of $\mathbb{S}^{k}(R)$.  
Call $\mathbb{S}^{k}(R)$ the $k$th {\em symmetric power} of $R$. 
It can be shown that $\bigwedge^{k}(R)$ and 
$\mathbb{S}^{k}(R)$ are ranked posets whose covering relations are the 
colored, directed edges prescribed in this paragraph. 

{\bf Fibrous posets and a connection with symmetric chain decompositions.} 
If for some $i \in I$, all $i$-components of $R$ are chains, then we say $R$ is $i$-{\em fibrous}. 
If $R$ is $i$-fibrous for all $i \in I$, then we say $R$ is {\em fibrous}. 
The following lemma straightforwardly connects the fibrous property with symmetric chain decompositions. 

\noindent 
{\bf \SCDLemma}\ \ {\sl With $R$ as above (but not necessarily fibrous), fix an edge color $i \in I$. 
Suppose $R'$ is ranked and edge-colored by $I$ as well, and is $i$-fibrous. 
Let $\rho_{i}'$ denote the color $i$ rank function on $R'$, and similarly use the notation $\delta_{i}'$, $l_{i}'$. 
Set $m_{i}' := \rho_{i}' - \delta_{i}' = 2\rho_{i}' - l_{i}'$. 
Suppose that $R = R'$ as sets. 
Suppose that for all $\xelt, \yelt$ in $R$ we have $\xelt \myarrow{i} \yelt$ in $R$ if $\xelt \myarrow{i} \yelt$ in $R'$. 
Also suppose that for all $\xelt \in R$ we have $m_{i}(\xelt) = m_{i}'(\xelt)$. 
Let $\mathcal{C}$ be the $i$-component $\comp_{i}(\xelt)$ of some $\xelt \in R$. 
Let $\mathcal{C}_1,\ldots,\mathcal{C}_k$ be the collection of $i$-components in $R'$ which contain some element of $\mathcal{C}$. 
Then $\mathcal{C}_1,\ldots,\mathcal{C}_k$ are the symmetric chains for a symmetric chain decomposition of $\mathcal{C}$.}

{\em Proof.} Obviously each $\mathcal{C}_j$ is a chain in $R'$. 
Moreover $\mathcal{C} \supseteq \mathcal{C}_1 \disjointunion \cdots \disjointunion \mathcal{C}_k$ since edges in $R'$ are also edges in $R$. 
Clearly $\mathcal{C} \subseteq \mathcal{C}_1 \disjointunion \cdots \disjointunion \mathcal{C}_k$. 
Now let $\xelt_j$ be the minimum element and $\yelt_j$ the maximum element of $\mathcal{C}_j$ in $R'$. 
Since edges in $R'$ are also edges in $R$, then $l_{i}'(\xelt_j) = \rho_i(\yelt_j) - \rho_i(\xelt_j)$, which is also $\delta_{i}'(\xelt_j)$. 
Then $\rho_i(\xelt_j) - \delta_i(\xelt_j) = m_i(\xelt_j) = m_{i}'(\xelt_j) = \rho_{i}'(\xelt_j) - \delta_{i}'(\xelt_j)$. 
So $\rho_i(\xelt_j) + \rho_{i}(\yelt_j) = \delta_i(\xelt_j) - \delta_{i}'(\xelt_j) + \rho_{i}(\yelt_j) = \delta_i(\xelt_j) - \rho_i(\yelt_j) + \rho_i(\xelt_j) + \rho_i(\yelt_j) = l_i(\xelt_j)$, which is the length of the $i$-component $\mathcal{C}$ in $R$. 
Therefore we have a symmetric chain decomposition $\mathcal{C} = \mathcal{C}_1 \disjointunion \cdots \disjointunion \mathcal{C}_k$.\hfill\QED

So, one way to demonstrate that a ranked poset $P$ has a symmetric chain decomposition is to realize $P$ as an $i$-component of some ranked and edge-colored poset that contains an $i$-fibrous subposet in the sense of the preceding lemma. 
This might not be as far-fetched as it sounds: There are many ``splitting posets'' (this term is defined in \S\SplittingSection\ below) that naturally contain crystal graphs (which are indeed fibrous) as subposets in the manner prescribed above. 
That said, for all of the examples of this phenomenon that we know of, the $i$-components of the ``parent'' splitting poset can be seen much more directly to possess symmetric chain decompositions. 
For instance, in many such cases the $i$-components are products of chains. 

{\bf Posets interacting with roots, weights, and the Weyl group.} 
For the remainder of this section, $R$ is a ranked poset with edges colored by $I$, an index set of cardinality $n$ for a choice of simple roots for a finite root system $\Phi$ as in \S \WeylSection.  
For any $J \subseteq I$, we let $wt^{J}: R \longrightarrow \Lambda_{\Phi_J}$ be the function given by $wt^{J}(\xelt) = \sum_{j \in J}m_{j}(\xelt)\omega_{j}^{J}$ for all $\xelt \in R$. 
Call $wt^{J}(\xelt)$ the $J$-{\em weight} of $\xelt$. 
When $J=I$, we write $wt$ for $wt^{I}$.  
Note that $wt^{J}(\xelt)$ is just the projection $(wt(\xelt))^{J}$ of $wt(\xelt)$ onto $\Lambda_{\Phi_J}$. 

The {\em weight generating function for} $R$ is $\WGF(R) := \displaystyle \sum_{\xelt \in R} e^{wt(\xelt)}$, an element of the group ring $\mathbb{Z}[\Lambda]$. 
By collecting like terms, we can rewrite $\WGF(R)$ as $\displaystyle \sum_{\mu \in \Lambda} \mathscr{C}_{\mu}(R)e^{\mu}$, where $\mathscr{C}_{\mu}(R)$ is the nonnegative integer size of the set $\{\xelt \in R\, |\, wt(\xelt) = \mu\}$. 
We observe the following. 

\noindent 
{\bf \WGFLemma}\ \ {\sl Let each of $R$, $R_1$, and $R_2$ be ranked posets with edges colored by the index set $I$ for our chosen basis of simple roots for $\Phi$. 
Then} $\WGF(R_1 \oplus R_2) = \WGF(R_1) + \WGF(R_2)$, $\WGF(R_1 \times R_2) = \WGF(R_1)\, \WGF(R_2)$, $\WGF(R^{*}) = \WGF(R)^{*}$, {\sl and} $\WGF(R^{\bowtie}) = \WGF(R)^{\bowtie}$.

{\em Proof.} It follows easily from the definitions that $\WGF(R_1 \oplus R_2) = \WGF(R_1) + \WGF(R_2)$. 
In $R_1 \times R_2$, note that $\rho_{i}\xxelt = \rho_i(\xelt_1)+\rho_i(\xelt_2)$ and $\delta_i\xxelt = \delta_i(\xelt_1) + \delta_i(\xelt_2)$. 
So $m_i\xxelt = m_i(\xelt_1)+m_i(\xelt_2)$, and therefore $wt\xxelt = wt(\xelt_1) + wt(\xelt_2)$. 
Then for any $\mu \in \Lambda$, $\mathscr{C}_{\mu}(R_1 \times R_2) = \sum(\mathscr{C}_{\nu_1}(R_1)\mathscr{C}_{\nu_2}(R_2))$, where the latter sum is over all weights $\nu_1,\nu_2$ such that $\nu_1 + \nu_2 = \mu$. 
The identity $\WGF(R_1 \times R_2) = \WGF(R_1)\, \WGF(R_2)$ follows. 
Since we have $wt(\xelt^{*}) = -wt(\xelt)$ for all $\xelt \in R$, then $\mathscr{C}_{\mu}(R^{*}) = \mathscr{C}_{-\mu}(R)$. 
It follows that $\WGF(R^{*}) = \WGF(R)^{*}$. 
Setwise, we can write $R^{\bowtie} = \{\xelt^{\bowtie}\, |\, \xelt \in R\}$. 
Then for all $\xelt \in R$, we can see that $m_{i}(\xelt^{\bowtie}) = -m_{\sigma_0(i)}(\xelt)$, so that $wt(\xelt^{\bowtie}) = w_0(wt(\xelt))$. 
Then $\mathscr{C}_{\mu}(R^{\bowtie}) = \mathscr{C}_{w_0(\mu)}(R)$, so therefore $\WGF(R^{\bowtie}) = \WGF(R)^{\bowtie}$.\hfill\QED

Given our choice of simple roots for the root system $\Phi$, we say $R$ is ${\Phi}$-{\em structured} if it has the following property: 
$wt(\selt) + \alpha_{i} = wt(\telt)$ whenever $\selt \myarrow{i} \telt$ in $R$, or equivalently (by \WeightLemma) for all $j \in I$, $m_{j}(\selt) + M_{ij} = m_{j}(\telt)$ whenever $\selt \myarrow{i} \telt$ in $R$. 
If $R$ is ${\Phi}$-structured, then $R$ is ${\Phi_J}$-structured for each $J \subseteq I$; note that $\WGF(R)|_{J} = \displaystyle \sum_{\xelt \in R} e^{wt^J(\xelt)} \in \mathbb{Z}[\Lambda_{\Phi_J}]$.

\noindent 
{\bf \WeightLemma}\ \ {\sl Let $R$ be a ranked poset with edges colored by the index set $I$ for our chosen basis of simple roots for $\Phi$.  Let $\xelt, \yelt \in R$.  
(1) If $J \subseteq I$, then for each $j \in J$, $m_{j}(\xelt) = \langle wt^{J}(\xelt),\alpha_{j}^{\vee} \rangle$. 
(2) If $R$ is connected and ${\Phi}$-structured, and if $wt(\xelt) = wt(\yelt)$, then $\xelt$ and $\yelt$ have the same rank in $R$. 
(3) Say $Q$ is another ranked poset with edges colored by $I$, and suppose both $Q$ and $R$ are ${\Phi}$-structured.  Then $R \oplus Q$, $R \times Q$, $R^{*}$, $R^{\bowtie}$ are ${\Phi}$-structured as well. (4) For any nonnegative integer $k$, 
$\bigwedge^{k}(R)$  
and $\mathbb{S}^{k}(R)$ are ${\Phi}$-structured.}

{\em Proof.}  Part {\sl (1)} is an immediate consequence of the definitions.  
For {\sl (2)}, since $R$ is connected, we can take a path $\mathcal{P} = (\xelt = \xelt_{0}, \xelt_{1}, \ldots, \xelt_{k} = \yelt)$ such that for $1 \leq j \leq k$ we have either $\xelt_{j-1} \myarrow{i_j} \xelt_{j}$ or  $\xelt_{j} \myarrow{i_j} \xelt_{j-1}$.  
Now $\sum_{i \in I}a_i(\mathcal{P})$ counts the number of ``upward'' steps in the path from $\xelt$ to $\yelt$, and $\sum_{i \in I}d_i(\mathcal{P})$ counts the number of ``downward'' steps.  
Since $R$ is ${\Phi}$-structured, then $wt(\xelt) + \sum_{i \in I}(a_i(\mathcal{P}) - d_i(\mathcal{P}))\omega_i = wt(\yelt)$, hence $\sum_{i \in I}(a_i(\mathcal{P}) - d_i(\mathcal{P}))\omega_i = 0$.  
Then $a_i(\mathcal{P}) = d_i(\mathcal{P})$ for each $i \in I$. 
So $\sum_{i \in I}a_i(\mathcal{P}) = \sum_{i \in I}d_i(\mathcal{P})$, and hence $\xelt$ and $\yelt$ have the same rank. 
The results in part {\sl (3)} of the lemma statement follow from the proof of \WGFLemma.

For (4), we only show that  
$\bigwedge^{k}(R)$ has the $\Phi$-structure property;   
since the argument that $\mathbb{S}^{k}(R)$ 
is $\Phi$-structured is similar to the the argument for the $k$th 
exterior power, we omit the details of that proof.   For any $\selt \in \bigwedge^{k}(R)$, 
define $\mu_{i}(\selt) := \sum_{v_{j} \in \selt}m_{i}(v_{j})$, where 
$m_{i}(v_{j}) = \rho_{i}(v_{j}) - \delta_{i}(v_{j})$ is calculated in 
$R$ for each $v_{j}$.  Then define $\mu(\selt) := (\mu_{i}(\selt))_{i 
\in I_{n}}$.  First, note that if $\selt \myarrow{i} \telt$ in 
$\bigwedge^{k}(R)$, then $(\selt - \telt,\telt - \selt) = 
(\{v_{p}\},\{v_{q}\})$ with $v_{p} \myarrow{i} v_{q}$ in $R$.  Then 
$m_{j}(v_{p}) + M_{ij} = m_{j}(v_{q})$ in $R$.  It now follows that 
$\mu_{j}(\selt) + M_{ij} = \left(\sum_{v_{r} \in 
\selt}m_{j}(v_{r})\right)+M_{ij} = m_{j}(v_{p}) + M_{ij} + 
\sum_{v_{r} \not= v_{p}}m_{j}(v_{r}) = m_{j}(v_{q}) + \sum_{v_{r} 
\not= v_{p}}m_{j}(v_{r}) = \sum_{v_{r} \in \telt}m_{j}(v_{r}) = 
\mu_{j}(\telt)$.  From this it follows that $\mu(\selt) + 
M^{(i)} = \mu(\telt)$.  

Now suppose $\selt = \relt_{0} \myarrow{i_{1}} \relt_{1} 
\myarrow{i_{2}} \relt_{2} \myarrow{i_{3}} \cdots \myarrow{i_{p}} 
\relt_{p} = \selt$ in $\bigwedge^{k}(R)$. Then $\mu(\selt) = 
\mu(\selt) + \sum_{i \in I_{n}}a_{i}M^{(i)}$, where $a_{i}$ counts the 
number of times there is an edge of color $i$ in our given path from 
$\selt$ to itself.  So, $\sum_{i \in I_{n}}a_{i}M^{(i)} = 0$.  
Since $M$ is nonsingular, then the $M^{(i)}$'s are linearly 
independent. So each $a_{i} = 0$.  Hence $\bigwedge^{k}(R)$ is 
acyclic, so we may define a partial order on $\bigwedge^{k}(R)$ in 
the following way: $\selt \leq \telt$ if and only if there is an 
`ascending' path $\selt = \relt_{0} \myarrow{i_{1}} \relt_{1} 
\myarrow{i_{2}} \relt_{2} \myarrow{i_{3}} \cdots \myarrow{i_{p}} 
\relt_{p} = \telt$ from $\selt$ to $\telt$ in $\bigwedge^{k}(R)$.  
Suppose $\selt \myarrow{i} \telt$ and that $\selt \leq \xelt \leq 
\telt$.  So then we have an ascending path $\selt = \relt_{0} \myarrow{i_{1}} 
\relt_{1} 
\myarrow{i_{2}} \cdots \myarrow{i_{q-1}} \relt_{q-1} 
\myarrow{i_{q}} \relt_{q} = \xelt 
\myarrow{i_{q+1}} \relt_{q+1} \myarrow{i_{q+2}} \cdots \myarrow{i_{p}} 
\relt_{p} = \telt$.  In this case we get $M^{(i)} = \sum_{i \in 
I_{n}}a_{i}M^{(i)}$, where $a_{i}$ is as before.  Then $a_{j} = 
\delta_{ij}$, from which we see that $\xelt = \selt$ or $\xelt = 
\telt$.  Hence each $\selt \myarrow{i} \telt$ is a covering relation 
for the partial order on $\bigwedge^{k}(R)$.  

Finally, we show that $\bigwedge^{k}(R)$ is ranked. It suffices to 
show this on each connected component of $\bigwedge^{k}(R)$.  So let 
$\mathcal{C}$ be such a connected component.  An ordered pair of elements 
$(\xelt, \yelt)$ from $\bigwedge^{k}(R)$ is {\em ascending of color} $i$ 
if $\xelt \myarrow{i} \yelt$ and {\em descending of color} $i$ if 
$\yelt \myarrow{i} \xelt$.  For a path 
$\mathcal{P} = (\selt = \relt_{0}, \relt_{1}, \ldots, \relt_{p} = 
\telt)$, let $a_{i}(\mathcal{P})$ count the number of ascending pairs 
of color $i$ in the path $\mathcal{P}$ and $d_{i}(\mathcal{P})$ count 
the number of descending pairs of color $i$.  Let 
$\sigma(\mathcal{P}) := \sum_{i \in I_{n}}(a_{i}-d_{i})$. Call this 
quantity the {\em signed length} of the path $\mathcal{P}$.  Define a 
new relation $<_{\mathcal{C}}$ on $\mathcal{C}$ by 
declaring $\selt <_{\mathcal{C}} \telt$ if and only if there is a 
path $\mathcal{P}$ from $\selt$ to $\telt$ such that 
$\sigma(\mathcal{P}) > 0$.  Then define $\leq_{\mathcal{C}}$ by the 
rule that $\selt \leq_{\mathcal{C}} \telt$ if and only if $\selt 
<_{\mathcal{C}} \telt$ or $\selt = \telt$.  We claim that 
$\leq_{\mathcal{C}}$ is a partial order on $\mathcal{C}$.  Clearly 
$\leq_{\mathcal{C}}$ is reflexive.  Use concatenation of paths to see 
that $\leq_{\mathcal{C}}$ is transitive.  Finally, we check that 
$\leq_{\mathcal{C}}$ is asymmetric.  
Suppose $\selt \leq_{\mathcal{C}} \telt$ and $\telt \leq_{\mathcal{C}} 
\selt$.  If $\selt \not= \telt$, then $\selt <_{\mathcal{C}} \telt$ 
and $\telt <_{\mathcal{C}} \selt$.  So there is a path $\mathcal{P}$ 
from $\selt$ to $\telt$ for which $\sigma(\mathcal{P}) > 0$ and a 
path $\mathcal{P}'$ from $\telt$ to $\selt$ for which 
$\sigma(\mathcal{P}') > 0$.  But now we can see that $\mu(\selt) + 
\sum_{i \in I_{n}}(a_{i}(\mathcal{P})-d_{i}(\mathcal{P}))M^{(i)} = 
\mu(\telt) = \mu(\selt) - \sum_{i \in 
I_{n}}(a_{i}(\mathcal{P}')-d_{i}(\mathcal{P}'))M^{(i)}$. By linear 
independence of the $M^{(i)}$'s we conclude that 
$(a_{i}(\mathcal{P})-d_{i}(\mathcal{P})) + 
(a_{i}(\mathcal{P}')-d_{i}(\mathcal{P}')) = 0$ for all $i \in 
I_{n}$.  But then $\sigma(\mathcal{P}) = \sum_{i \in 
I_{n}}(a_{i}(\mathcal{P})-d_{i}(\mathcal{P})) = 
-\sum_{i \in I_{n}}(a_{i}(\mathcal{P}')-d_{i}(\mathcal{P}')) = 
-\sigma(\mathcal{P}')$.  So $\sigma(\mathcal{P})$ and 
$\sigma(\mathcal{P}')$ cannot both be positive.  From this 
contradiction we conclude that $\selt = \telt$.  Then 
$\leq_{\mathcal{C}}$ is a partial order on $\mathcal{C}$. 

Now choose $\xelt$ to be a minimal element of $\mathcal{C}$ with 
respect to this partial order.  For any $\selt \in \mathcal{C}$, we 
declare $\rho_{\mathcal{C}}(\selt) := \sigma(\mathcal{P})$, where 
$\mathcal{P}$ is any path from $\xelt$ to $\selt$.  We claim that 
$\rho_{\mathcal{C}}(\selt)$ does not depend on the choice of path 
from $\xelt$ to $\selt$.  To see this, suppose $\mathcal{Q}$ is 
another path from $\xelt$ to $\selt$.  Then from the facts that 
$\mu(\xelt) + \sum_{i \in 
I_{n}}(a_{i}(\mathcal{P})-d_{i}(\mathcal{P}))M^{(i)}$ and 
$\mu(\xelt) + \sum_{i \in 
I_{n}}(a_{i}(\mathcal{Q})-d_{i}(\mathcal{Q}))M^{(i)}$, 
we deduce that $a_{i}(\mathcal{P}) - d_{i}(\mathcal{P}) = 
a_{i}(\mathcal{Q}) - d_{i}(\mathcal{Q})$ for all $i \in I_{n}$.  
Hence $\sigma(\mathcal{P}) = \sigma(\mathcal{Q})$.  Since $\xelt$ is 
minimal with respect to the partial order $\leq_{\mathcal{C}}$ on 
$\mathcal{C}$, it must be the case that $\rho_{\mathcal{C}}(\selt) = 
0$ for all $\selt \in \mathcal{C}$.  Finally, suppose $\selt 
\myarrow{i} \telt$ is a covering relation in $\bigwedge^{k}(R)$ for 
elements $\selt$ and $\telt$ in $\mathcal{C}$.  Then any path 
$\mathcal{P}$ from $\xelt$ to $\selt$ may be extended via $\selt 
\myarrow{i} \telt$ to a path 
$\mathcal{Q}$ from $\xelt$ to $\telt$.  Then $\sigma(\mathcal{Q}) = 
\sigma(\mathcal{P}) + 1$, and hence $\rho_{\mathcal{C}}(\telt) = 
\rho_{\mathcal{C}}(\selt) + 1$.  Then $\rho_{\mathcal{C}}$ is a rank 
function for $\mathcal{C}$.  It follows that $\bigwedge^{k}(R)$ is 
ranked.\hfill\QED 

The weight generating function $\WGF(R)$ is in the ring of Weyl symmetric functions $\mathbb{Z}[\Lambda]^{W}$ (i.e.\ is $W$-invariant) if and only if  $w.\WGF(R) = \WGF(R)$ for all $w \in W$ if and only if $\mathscr{C}_{\mu}(R) = \mathscr{C}_{s_{i}(\mu)}(R)$ for all $i \in I$ and $\mu \in \Lambda$, since the $s_{i}$'s generate $W$.  
The $W$-invariance of $\WGF(R)$ is implied by a number of combinatorial conditions on $R$, such as the one developed next. 

\noindent 
{\bf \WInvariantLemma}\ \ {\sl Let $R$ be an ${\Phi}$-structured poset such that for each $i \in I$, each $i$-component of $R$ is rank symmetric.  
Then for all $J \subseteq I$,} $\WGF(R)|_{J}$ {\sl is $W_J$-invariant.} 

{\em Proof.} For each $i \in I$, construct an involutive correspondence between elements of $R$ as follows: in each $i$-component, let $\tau_{i}$ be a matching of the elements of correspondingly symmetric ranks, so each $\xelt$ of rank $\rho_{i}(\xelt)$ in its $i$-component is matched with some $\tau_{i}(\xelt)$ of rank $\rho_{i}(\tau_{i}(\xelt)) = l_{i}(\xelt)-\rho_{i}(\xelt)$ in such a way that $\tau_{i}(\tau_{i}(\xelt)) = \xelt$. 
This is possible since each $i$-component of $R$ is rank symmetric.  
Since $R$ is ${\Phi}$-structured, one can see that $wt(\tau_{i}(\xelt)) = wt(\xelt) + (l_{i}(\xelt)-2\rho_{i}(\xelt))\alpha_{i} = wt(\xelt) - m_{i}(\xelt)\alpha_{i}$.  But by \WeightLemma, $wt(\xelt) - m_{i}(\xelt)\alpha_{i} = wt(\xelt) - \langle wt(\xelt),\alpha_{i}^{\vee} \rangle\alpha_{i}$, which is just $s_{i}(wt(\xelt))$.  Then $s_{i}(wt(\xelt)) = wt(\tau_{i}(\xelt))$.  
Then for all $\mu \in \Lambda$, $\mathscr{C}_{\mu}(R) = |\{\xelt \in R|wt(\xelt) = \mu\}| = |\{\xelt \in R|wt(\tau_{i}(\xelt)) = \mu - m_{i}(\xelt)\alpha_{i}\}| = |\{\yelt \in R|wt(\yelt)=s_{i}(\mu)\}| = \mathscr{C}_{s_{i}(\mu)}(R)$.  
Since $\WGF(R)$ is $W$-invariant, then $\WGF(R)|_J$ is $W_J$-invariant for all $J \subseteq I$.\hfill\QED

{\bf Weight diagrams, again.} 
The following proposition is a continuation of \WeightRemarkOne\ and establishes several more poset-theoretic properties of weight diagrams. 
Of course, the proposition applies when the saturated set of weights $\mathscr{W}$ is $\Pi(\lambda)$ for some dominant $\lambda$. 

\noindent 
{\bf \WeightRemarkTwo}\ \ {\sl Let $\mathscr{W} \subset \Lambda$ be a finite saturated set of weights.  
Regard $\mathscr{W}$ to be a subposet of $\Lambda$ in the induced order.  
(1) For $\mu, \nu \in \mathscr{W}$, we have $\mu \rightarrow \nu$ in the poset 
$\mathscr{W}$ if and only if for some $i \in I$ we have $\mu + 
\alpha_{i} = \nu$.  
(2) Regarding $\mathscr{W}$ to be edge-colored by the rule from part (1), then the $i$-components of $\mathscr{W}$ are chains, for each $i \in I$.  
(3) $\mathscr{W}$ is ranked.  
(4) For any $\mu \in \mathscr{W}$, we have $wt(\mu) = \mu$, and hence $\mathscr{W}$ is an ${\Phi}$-structured poset.} 

{\em Proof.} 
For {\sl (1)}, the proof of \WeightRemarkOne.5 works here, since that proof only depended on the property of being saturated.  
{\sl (2)} amounts to an observation. 
For {\sl (3)}, note that it is enough to show that for any two paths $\mathcal{P}$ and $\mathcal{Q}$ from $\mu$ to $\nu$ in the edge-colored poset $\mathscr{W}$, then $\sum_{i \in I}(a_i(\mathcal{P})-d_i(\mathcal{P})) = \sum_{i \in I}(a_i(\mathcal{Q})-d_i(\mathcal{Q}))$.  
This follows from the following computation: $\sum_{i \in I}(a_i(\mathcal{P})-d_i(\mathcal{P}))\alpha_i  = \nu - \mu = \sum_{i \in I}(a_i(\mathcal{Q})-d_i(\mathcal{Q}))\alpha_i$.  
For {\sl (4)}, we only need to argue that $m_i(\mu) = \langle \mu,\alpha_{i}^{\vee} \rangle$ for each $\mu \in \mathscr{W}$ and $i \in I$. 
Let 
\[\mu_{0} \myarrow{i} \cdots \myarrow{i} \mu \myarrow{i} \cdots 
\myarrow{i} \mu_{1}\] 
be the $i$-component $\mathbf{comp}_{i}(\mu)$.  
Since $\mathscr{W}$ is saturated, then $\langle \mu_{0},\alpha_{i}^{\vee} \rangle \geq -l_{i}(\mu)$ and $\langle \mu_{1},\alpha_{i}^{\vee} \rangle \leq l_{i}(\mu)$.  
But $\langle \mu_{1},\alpha_{i}^{\vee} \rangle = \langle \mu_{0} + l_{i}(\mu)\alpha_{i},\alpha_{i}^{\vee} \rangle \geq -l_{i}(\mu) + 2l_{i}(\mu) = l_{i}(\mu)$.  
An easy calculation now shows that $\langle \mu,\alpha_{i}^{\vee} \rangle = m_{i}(\mu)$.\hfill\QED

From now on, we regard any finite, saturated set of weights to be edge-colored and ${\Phi}$-structured as in \WeightRemarkTwo. 

{\bf A generalized notion of weight diagram.} 
Our next theorem regards a sort-of generalization of weight diagrams for ${\Phi}$-structured posets, here denoted ``$\Pi(R)$.'' 
This theorem and its corollaries will be applied in proofs of some poset-structural results (\MainCorollary, \MinQuasiMinTheorem, \EasyProductTheorem/\DataProp), will give a weight-theoretic method for calculating ranks within connected and ${\Phi}$-structured posets, and will aid in understanding how weight diagrams behave under restriction to some $J \subseteq I$.  
So, we view this theorem as a very basic and sometimes useful result about ${\Phi}$-structured posets. 
Some set-up: For an ${\Phi}$-structured poset $R$, denote by $\Pi(R)$ the edge-colored directed graph whose elements are $\{wt(\selt)\}_{\selt \in R}$ and with colored directed edges given by $\mu \myarrow{i} \nu$ whenever there exists an edge $\xelt \myarrow{i} \yelt$ in $R$ for which $wt(\xelt) = \mu$ and $wt(\yelt) = \nu$.  
Call $\Pi(R)$ the {\em weight diagram for} $R$ or, generically, a {\em generalized weight diagram}. 
We record some basic features of $\Pi(R)$ in the following lemma.  

\noindent 
{\bf \PiLemma}\ \ {\sl With $R$ as above, then $\Pi(R)$ is is the edge-colored Hasse diagram for a ranked poset.  
If $R$ is connected, then so is $\Pi(R)$.} 

{\em Proof.} That $\Pi(R)$ is connected when $R$ is follows from the definitions. 
To see that $\Pi(R)$ is the Hasse diagram for a ranked poset, we let $\mu \leq \nu$ in $\Pi(R)$ if and only if there exists a path $\mathcal{P} = (\mu = \mu_{0}, \ldots, \mu_{k} = \nu)$ and colors $(i_j)_{j=1}^{k}$ such that for $1 \leq j \leq k$ we have $\mu_{j-1} \myarrow{i_j} \mu_{j}$. 
It is easy to see that $\leq$ is reflexive and transitive. 
If $\mu \leq \nu$ and $\nu \leq \mu$, then we would have a path from $\mu$ to $\nu$ as above and another path $\mathcal{Q} = (\nu = \nu_{0}, \ldots, \nu_{l} = \mu)$ from $\nu$ to $\mu$ such that $\nu_{j-1} \myarrow{i'_j} \nu_{j}$ for $1 \leq j \leq l$.  
Since $R$ is ${\Phi}$-structured, then $\mu = \sum_{i \in I}a_i(\mathcal{Q})\omega_i + \nu = \sum_{i \in I}a_i(\mathcal{Q})\omega_i + \sum_{i \in I}a_i(\mathcal{P})\omega_i + \mu$. 
It follows that for each $i$, the nonnegative integer sum $a_i(\mathcal{P}) + a_i(\mathcal{Q})$ is zero, hence $a_i(\mathcal{P}) = a_i(\mathcal{Q}) = 0$. 
Therefore $\mu = \nu$. 

To show that $\Pi(R)$ is ranked, it suffices to assume that $\Pi(R)$ is connected.  
Let $\eta \in \Pi(R)$ be any weight such that $\mathrm{ht}(\eta) \leq \mathrm{ht}(wt(\telt))$ for all $\telt \in R$.  
Then for any $\mu \in \Pi(R)$, let $\rho(\mu) := \mathrm{ht}(\mu) - \mathrm{ht}(\eta) = \langle \mu - \eta,\varrho^{\vee} \rangle$. 
Connectedness of $\Pi(R)$ guarantees that $\rho(\telt)$ is a nonnegative integer. 
It is easy now to check that $\rho$ is a rank function for $\Pi(R)$.\hfill\QED

We more fully investigate the relationship between $R$ and $\Pi(R)$ in \PiTheorem\ and its corollaries below. 
A nonempty set $\mathscr{D}$ of prominent elements of $R$ is {\em indomitable} in $R$ if for any prominent $\xelt \in R$, there is some $\delt \in \mathscr{D}$ such that $wt(\xelt) \leq wt(\delt)$.  
It is {\em minimally indomitable} in $R$ if no proper subset of $\mathscr{D}$ is indomitable. 
We note that any indomitable set $\mathscr{D}'$ contains a minimally indomitable subset. 
To see this, we simply pare $\mathscr{D}'$ down as follows.  
First, form $\mathscr{D}''$ from $\mathscr{D}'$ by removing all $\delt' \in \mathscr{D}'$ such that $wt(\delt') < wt(\delt'')$ for some $\delt'' \in \mathscr{D}'$.  
Then note that $\disjointunion_{\mu \in \Pi(R)} \{\delt \in \mathscr{D}''\, |\, wt(\delt) = \mu\}$ is a set partition of $\mathscr{D}''$.  
Now form a set $\mathscr{D}$ from $\mathscr{D}''$ by choosing precisely one $\delt$ from each nonempty subset in the partition.  
It is easy to see now that $\mathscr{D} \subseteq \mathscr{D}'$ is minimally indomitable in $R$.
It is also easy to see that if $\mathscr{D}_{1}$ and $\mathscr{D}_{2}$ are two minimally indomitable sets in $R$, then $|\mathscr{D}_{1}| = |\mathscr{D}_{1}|$. 

\noindent 
{\bf \PiTheorem}\ \  {\sl Let $R$ be ${\Phi}$-structured.  
Suppose $\mathscr{N}$ is any set of elements of $R$ such that $wt(\nelt) \in \Lambda^{+}$ for each $\nelt \in \mathscr{N}$.  
Then $\displaystyle \Pi(R) = \bigcup_{\nelt \in \mathscr{N}} \Pi(wt(\nelt))$ (an equality of edge-colored directed graphs) if and only if  $\displaystyle \Pi(R) = \bigcup_{\nelt \in \mathscr{N}} \Pi(wt(\nelt))$ (an equality of sets) if and only if $\mathscr{D} \subseteq \mathscr{N}$ for some minimally indomitable set $\mathscr{D}$ of prominent elements in $R$.} 

{\em Proof.}  We first argue that  it suffices to prove the result when $\Pi(R)$ is connected.  
Indeed, if $\Pi(R) = \Pi_{1} \oplus \cdots \oplus \Pi_{k}$ is a disjoint sum of the connected components of $\Pi(R)$, then we can set $R_{i} := \{\xelt \in R\, |\, wt(\xelt) \in \Pi_{i}\}$ for $1 \leq i \leq k$.  
Then $\Pi(R_{i}) = \Pi_{i}$, each $R_{i}$ corresponds a disjoint sum of connected components of $R$, and $R = R_{1} \oplus \cdots \oplus R_{k}$.  
Suppose we have $\mathscr{D} \subseteq \mathscr{N}$ as in the theorem statement.  Let $\mathscr{D}_{i} := \mathscr{D} \cap R_{i}$ and $\mathscr{N}_{i} := \mathscr{N} \cap R_{i}$. 
In particular, each $\mathscr{D}_{i}$ is minimally indomitable in $R_{i}$, and $\mathscr{D}_{i} \subseteq \mathscr{N}_{i}$. 
The hypothesis that the theorem statement holds when $\Pi(R)$ is connected would now apply to $\Pi(R_{i})$ and allow us to conclude that 
$\displaystyle \Pi(R_{i}) = \bigcup_{\nelt \in \mathscr{N}_{i}} \Pi(wt(\nelt))$ is an equality of edge-colored directed graphs.  It follows easily that $\displaystyle \Pi(R) = \bigcup_{\nelt \in \mathscr{N}} \Pi(wt(\nelt))$ is also an equality of edge-colored directed graphs, and hence an equality of sets. 
Now, if $\displaystyle \Pi(R) = \bigcup_{\nelt \in \mathscr{N}} \Pi(wt(\nelt))$ is an equality of sets, where $\mathscr{N}$ is as in the theorem statement, then we can let $\mathscr{N}_{i} := \mathscr{N} \cap R_{i}$.  Then $\displaystyle \Pi(R_{i}) = \bigcup_{\nelt \in \mathscr{N}_{i}} \Pi(wt(\nelt))$ (set equality), so there is a minimally indomitable set $\mathscr{D}_{i}$ in $R_{i}$ such that $\mathscr{D}_{i} \subseteq \mathscr{N}_{i}$.  Then $\mathscr{D} := \mathscr{D}_{1} \disjointunion \cdots \disjointunion \mathscr{D}_{k}$ is a minimally indomitable set in $R$ with $\mathscr{D} \subseteq \mathscr{N}$. 

So from here on, we assume that $\Pi(R)$ is connected.  We first assume that $\mathscr{D} \subseteq \mathscr{N}$ as in the theorem statement. 
Pick $\lambda \in \Pi(R)$ such that $\mathrm{ht}(\lambda) \geq \mathrm{ht}(wt(\telt))$ for all $\telt \in R$.  
Then set $\delta(\telt) := \mathrm{ht}(\lambda)-\mathrm{ht}(wt(\telt)) = \langle \lambda-wt(\telt),\varrho^{\vee} \rangle$.  
The connectedness hypothesis for $\Pi(R)$ implies that $\delta(\telt)$ is a nonnegative integer for all $\telt \in R$.  We will induct on $\delta(\telt)$ to show that $\displaystyle \Pi(R) \subseteq \bigcup_{\nelt \in \mathscr{N}} \Pi(wt(\nelt))$ (a containment of subsets) and for each edge $\selt \myarrow{i} \telt$ in $R$ we have $wt(\selt) \myarrow{i} wt(\telt)$ as an edge in $\displaystyle \bigcup_{\nelt \in \mathscr{N}} \Pi(wt(\nelt))$. 

Say $\selt \myarrow{i} \telt$ with  $\delta(\telt) = 0$. 
Then $wt(\telt) = wt(\delt)$ where $\delt \in \mathscr{D} \subseteq \mathscr{N}$. 
Applying \WeightLemma, $\langle wt(\telt),\alpha_{i}^{\vee} \rangle = m_{i}(\telt) = \rho_{i}(\telt) > 0$.  
Since $0 < 1 \leq \langle wt(\telt),\alpha_{i}^{\vee} \rangle = \langle \lambda,\alpha_{i}^{\vee} \rangle$, then $\mu := \lambda - \alpha_{i} \in \Pi(wt(\delt))$ since $\Pi(wt(\delt))$ is saturated, cf.\ \WeightRemarkOne.4. 
Since $\mu + \alpha_{i} = \lambda$, then $\mu \myarrow{i} \lambda$ in $\displaystyle \bigcup_{\nelt \in \mathscr{N}} \Pi(wt(\nelt))$.  
That is, $wt(\selt) \myarrow{i} wt(\telt)$ in $\displaystyle \bigcup_{\nelt \in \mathscr{N}} \Pi(wt(\nelt))$. 
Now say for all $\yelt \in R$ with $\delta(\yelt) \leq k$, it is the case that $wt(\xelt) \myarrow{i} wt(\yelt)$ is an edge in $\displaystyle \bigcup_{\nelt \in \mathscr{N}} \Pi(wt(\nelt))$ whenever $\xelt \myarrow{i} \yelt$ is an edge in $R$.  

So now suppose that $\selt \myarrow{i} \telt$ in $R$ with $\delta(\telt) = k+1$.  
If $\telt$ is prominent, then the same reasoning we used in the case $\delta(\telt) = 0$ shows that $wt(\selt) \myarrow{i} wt(\telt)$ is an edge in $\displaystyle \bigcup_{\nelt \in \mathscr{N}} \Pi(wt(\nelt))$. 
Now suppose $\telt$ is not prominent, and take any $j \in I$ such that $\delta_{j}(\telt) > 0$ in $\mathbf{comp}_{j}(\telt)$. 
Let $\uelt$ be an element of $\mathbf{comp}_{j}(\telt)$ for which $\delta_{j}(\uelt) = 0$, and let $\xelt$ be any element of $\comp_{j}(\telt)$ with $\delta_j(\xelt) > 0$.   
Then $wt(\xelt) = wt(\uelt) - \delta_j(\xelt)\alpha_{j}$, since $R$ is 
${\Phi}$-structured.  
Since $\uelt$ is maximal in $\mathbf{comp}_{j}(\telt)$, then $\uelt' \myarrow{j} 
\uelt$ for some $\uelt' \in R$.  
In particular, the induction hypothesis applies to this edge, and we conclude that  
$wt(\uelt) \in \displaystyle \bigcup_{\nelt \in \mathscr{N}} \Pi(wt(\nelt))$. 
Now  $\langle wt(\uelt),\alpha_{j}^{\vee} \rangle = \rho_{j}(\uelt)-\delta_{j}(\uelt) = \rho_{j}(\uelt)$ by \WeightLemma.  
We now know that $0 < \delta_j(\xelt) \leq \rho_j(\uelt) = \langle wt(\uelt),\alpha_{j}^{\vee} \rangle$. 
In particular, we have $wt(\xelt) = wt(\uelt) - \delta_{j}(\xelt)\alpha_{j}$ with 
$0 < \delta_{j}(\xelt) \leq \langle wt(\uelt),\alpha_{j}^{\vee} \rangle$.  
Since $\displaystyle \bigcup_{\nelt \in \mathscr{N}} \Pi(wt(\nelt))$ 
is saturated, then $wt(\xelt) \in \displaystyle \bigcup_{\nelt \in \mathscr{N}} \Pi(wt(\nelt))$.  
In particular, $wt(\telt)$ is in this set.  
Suppose that $\delta_i(\telt) > 0$. 
Then the above reasoning applies when $j=i$ to show that $wt(\selt)$ is in $\displaystyle \bigcup_{\nelt \in \mathscr{N}} \Pi(wt(\nelt))$ as well. 
Since $\selt \myarrow{i} \telt$ and $R$ is ${\Phi}$-structured, then $wt(\selt) + \alpha_{i} = wt(\telt)$, and hence $wt(\selt) \myarrow{i} wt(\telt)$ is an edge in $\displaystyle \bigcup_{\nelt \in \mathscr{N}} \Pi(wt(\nelt))$.  
On the other hand suppose that $\delta_{i}(\telt) = 0$. 
Then $\langle wt(\telt),\alpha_{i}^{\vee} \rangle = 
\rho_{i}(\telt)-\delta_{i}(\telt) = \rho_{i}(\telt)$ by \WeightLemma.  
So $0 < 1 \leq \rho_{i}(\telt) = \langle wt(\telt),\alpha_{i}^{\vee} \rangle$, and since $wt(\telt) \in \displaystyle \bigcup_{\nelt \in \mathscr{N}} \Pi(wt(\nelt))$ and $\displaystyle \bigcup_{\nelt \in \mathscr{N}} \Pi(wt(\nelt))$ is saturated, it follows that $wt(\selt) = wt(\telt) - \alpha_{i} \in \displaystyle \bigcup_{\nelt \in \mathscr{N}} \Pi(wt(\nelt))$ also. 
So, $wt(\selt) \myarrow{i} wt(\telt)$ is an edge in $\displaystyle \bigcup_{\nelt \in \mathscr{N}} \Pi(wt(\nelt))$. 
This completes the induction argument. 

So $\displaystyle \Pi(R) \subseteq \bigcup_{\nelt \in \mathscr{N}} \Pi(wt(\nelt))$ (a containment of subsets) and for each edge $\selt \myarrow{i} \telt$ in $R$ we have $wt(\selt) \myarrow{i} wt(\telt)$ in $\Pi(R)$.  
We now show that $\displaystyle \Pi(R) \supseteq \bigcup_{\nelt \in \mathscr{N}} \Pi(wt(\nelt))$ and that if $\mu \myarrow{i} \nu$ is an edge in $\Pi(R)$ if it is an edge in $\displaystyle \bigcup_{\nelt \in \mathscr{N}} \Pi(wt(\nelt))$. 

Before doing so, we show that $\displaystyle \bigcup_{\nelt \in \mathscr{N}} \Pi(wt(\nelt)) \subseteq \bigcup_{\delt \in \mathscr{D}'} \Pi(wt(\delt))$ as sets for any minimally indomitable set $\mathscr{D}'$.  To do this, it suffices to show that for any $\nelt \in \mathscr{N}$, there is some $\delt \in \mathscr{D}'$ for which $wt(\nelt) \leq wt(\delt)$.  To this end, set $\nelt_0 := \nelt$.  If $\nelt$ is not prominent, then choose $i_1$ such that $\delta_{i_1}(\nelt) > 0$.  Then choose $\nelt_{1} \in \comp_{i_1}(\nelt)$ such that $\delta_{i_1}(\nelt_1) = 0$.  Then $wt(\nelt) \leq wt(\nelt_1)$, and the rank of $\nelt$ is strictly less than the rank of $\nelt_1$ in $\comp_{I}(\nelt)$.  If $\nelt_1$ is prominent then since $\mathscr{D}'$ is minimally indomitable, we must have $wt(\nelt_1) = wt(\delt)$ for some $\delt \in \mathscr{D}'$.  If not, then repeat this process, replacing $\nelt$ with $\nelt_1$.  Thus we get a sequence $\nelt = \nelt_0, \nelt_1, \ldots$ where the ranks of these elements are strictly increasing.  Since the sequence must be finite, it must terminate with some $\nelt_k$ that is prominent.  Then we reason as before that $wt(\nelt_k) = wt(\delt)$ for some $\delt \in \mathscr{D}'$.  Since $wt(\nelt) = wt(\nelt_0) \leq wt(\nelt_1) \leq \cdots \leq wt(\nelt_k)$, then we get $wt(\nelt) \leq wt(\delt)$, as desired.  

Since $\mathscr{D} \subseteq \mathscr{N}$, it follows from the previous paragraph that $\displaystyle \bigcup_{\nelt \in \mathscr{N}} \Pi(wt(\nelt)) = \bigcup_{\delt \in \mathscr{D}} \Pi(wt(\delt))$ is an equality of sets and therefore of edge-colored directed graphs. 
This means that we only need to show that $\displaystyle \Pi(R) \supseteq \bigcup_{\delt \in \mathscr{D}} \Pi(wt(\delt))$ and that $\mu \myarrow{i} \nu$ is an edge in $\Pi(R)$ if it is an edge in $\displaystyle \bigcup_{\delt \in \mathscr{D}} \Pi(wt(\delt))$.
So we let $\delt \in \mathscr{D}$ and set $\lambda := wt(\delt)$. 
For any $\mu \in \Pi(\lambda)$, let $r_{i}(\mu)$ and $\myl_{i}(\mu)$ denote respectively the rank of $\mu$ in its $i$-component of $\Pi(\lambda)$ and the length of this $i$-component. 
Write $\lambda - \mu = \sum_{i\in{I}}k_{i}(\mu)\alpha_{i}$, where for each $i \in I$ we have $k_{i}(\mu) \geq 0$ since $\mu \leq \lambda$. 
Let $\depth(\mu) := \sum_{i\in{I}} k_{i}(\mu) \geq 0$.  
We induct on $\depth(\nu)$ to show that if $\mu \myarrow{i} \nu$ in $\Pi(\lambda)$, then $\mu = wt(\selt)$ and $\nu = wt(\telt)$ for some edge $\selt \myarrow{i} \telt$ in $R$.  

For the basis step of this induction argument, say $\mu \myarrow{i} \nu$ with 
$\depth(\nu) = 0$.  
Then $\nu = \lambda$, and $\lambda \in \Pi(R)$ since $\lambda = wt(\delt)$.  
So, $l_{i}(\delt) = \rho_{i}(\delt) = 2\rho_{i}(\delt) - l_{i}(\delt) = m_{i}(\delt) = \langle wt(\delt),\alpha_{i}^{\vee} \rangle = \langle \lambda,\alpha_{i}^{\vee} \rangle = 2r_{i}(\lambda)-\myl_{i}(\lambda) = \myl_{i}(\lambda) > 0$.  
Since $\delt$ must be maximal in its $i$-component in $R$ and since $l_{i}(\delt) > 0$, then there must be $\selt \in R$ with $\selt \myarrow{i} \delt$.  
Then $wt(\selt) = \lambda - \alpha_{i} = \mu$, completing the basis step of the induction argument.  
Suppose now that for some $k \geq 0$, it is that case that whenever $\pi \myarrow{i} \xi$ in $\Pi(\lambda)$ with $\depth(\xi) \leq k$, then there is an edge $\xelt \myarrow{i} \yelt$ in $R$ with $wt(\xelt) = \pi$ and $wt(\yelt) = \xi$.  
Now suppose that $\mathrm{depth}(\nu) = k+1$ for some $\nu \in \Pi(\lambda)$ and that $\mu \myarrow{i} \nu$ in $\Pi(\lambda)$.  
We consider two cases: (1) $\nu$ is not maximal in its $i$-component in $\Pi(\lambda)$, and (2) $\nu$ is maximal in its $i$-component. 
In case (1), since any monochromatic component of $\Pi(\lambda)$ is necessarily a chain, then we have $\mu \myarrow{i} \nu \myarrow{i} \cdots \myarrow{i} \pi \myarrow{i} \xi$, where $\xi$ is maximal in this $i$-component of $\Pi(\lambda)$, $\nu < \xi$, and possibly $\pi = \nu$. 
Applying the induction hypothesis to the edge $\pi \myarrow{i} \xi$, there is an edge $\xelt \myarrow{i} \yelt$ in $R$ with $wt(\xelt) = \pi$ and $wt(\yelt) = \xi$. 
Then $\rho_i(\yelt)-\delta_i(\yelt) = \langle wt(\yelt),\varrho^{\vee} \rangle = \langle \xi,\varrho^{\vee} \rangle = \myl_i(\xi) = r_i(\xi)$. 
Also, $r_i(\xi) - r_i(\nu) > 0$ since $\nu$ is not the least element in its $i$-component. 
It follows that $\rho_{i}(\yelt) > r_i(\xi) - r_i(\nu)$. 
Therefore there must be some $\telt$ in $\comp_{i}(\yelt)$ with $\rho_i(\telt) = \rho_i(\yelt) - (r_i(\xi) - r_i(\nu))$ and such that $\telt$ is not minimal in $\comp_{i}(\yelt)$, i.e.\ there is an edge $\selt \myarrow{i} \telt$. 
Then $wt(\telt) = wt(\yelt) - (r_i(\xi) - r_i(\nu))\alpha_{i} = \xi - (r_i(\xi) - r_i(\nu))\alpha_{i} = \nu$, and therefore $wt(\selt) = \mu$. 
Now there are $r_i(\xi) - r_i(\nu) > 0$ steps from $\xi$ down to $\nu$. 
In case (2), since $\nu$ is not maximal in $\Pi(wt(\delt))$, there is some $\nu' \in 
\Pi(wt(\delt))$ for which $\nu \myarrow{j} \nu'$.  
By the induction hypothesis, we have $\xelt \myarrow{j} \yelt$ in $R$ with $wt(\xelt) = \nu$ and $wt(\yelt) = \nu'$.  
With $\nu$ maximal in its $i$-component of $\Pi(wt(\delt))$, it follows that $0 < \myl_{i}(\nu) = r_{i}(\nu) = \langle \nu,\alpha_{i}^{\vee} \rangle = \langle 
wt(\xelt),\alpha_{i}^{\vee} \rangle = 2\rho_{i}(\xelt) - 
l_{i}(\xelt)$. 
In particular, $\rho_{i}(\xelt) > 0$, so in $\mathbf{comp}_{i}(\xelt)$, there must be an edge $\selt \myarrow{i} \telt$ with $\rho_{i}(\telt) = \rho_{i}(\xelt)$.  
Then since $R$ is ${\Phi}$-structured, we have $wt(\telt) = wt(\xelt) = \nu$.  
Then $wt(\selt) = wt(\telt) - \alpha_{i} = \mu$.  
This completes the induction step.

So we have shown that if $\mathscr{D} \subseteq \mathscr{N}$ for some minimally indomitable set $\mathscr{D}$ in $R$, then $\displaystyle \Pi(R) = \bigcup_{\nelt \in \mathscr{N}} \Pi(wt(\nelt))$ is an equality of edge-colored directed graphs.  Obviously if the latter holds, then we also have $\displaystyle \Pi(R) = \bigcup_{\nelt \in \mathscr{N}} \Pi(wt(\nelt))$ is an equality of sets.  

To finish the proof of the theorem, we assume the preceding equality of sets and show that $\mathscr{D} \subseteq \mathscr{N}$ for some minimally indomitable set $\mathscr{D}$ in $R$. 
By the observations from the paragraph preceding the theorem statement, it suffices to show that $\mathscr{D}' := \{\nelt \in \mathscr{N}\, |\, \nelt \mbox{ is prominent}\}$ is indomitable.  
So let $\xelt \in R$ be prominent.  
We will show that  there is some $\delt \in \mathscr{D}'$ such that $wt(\xelt) \leq wt(\delt)$.  
If $\xelt$ is in $\mathscr{N}$ then we are done, as we can take $\delt = \xelt$.  
Otherwise, since $\displaystyle \Pi(R) = \bigcup_{\melt \in \mathscr{N}} \Pi(wt(\melt))$, then we must have $wt(\xelt) \leq wt(\nelt_0)$ for some $\nelt_0 \in \mathscr{N}$. 
If $\nelt_0$ is prominent, then we are done. 
Otherwise, then do what was done four paragraphs back to find some prominent $\delt_1$ such that $wt(\nelt_0) \leq wt(\delt_1)$.  
If $wt(\delt_1) = wt(\melt)$ for some $\melt \in \mathscr{N}$, then we are done.  
Otherwise, since $\displaystyle \Pi(R) = \bigcup_{\melt \in \mathscr{N}} \Pi(wt(\melt))$, we can find $\nelt_1 \in \mathscr{N}$ for which $wt(\delt_1) < wt(\nelt_1)$.  
If $\nelt_1$ is prominent, then we are done.  
Otherwise, continue by finding $\delt_2$, $\nelt_2$, $\delt_3$, $\nelt_3$, etc.  
The finiteness of $R$, and therefore of $\Pi(R)$, ensures that this process will terminate, in which case we obtain $\delt \in \mathscr{D}'$ such that $wt(\xelt) \leq wt(\delt)$.\hfill\QED

Before we derive some corollaries, we give some remarks. 
In the statement of \PiTheorem, one can take $\mathscr{N}$ to be the set of all maximal elements of $R$ whose weights are dominant, or the set of all prominent elements of $R$. 
Also, it follows from this theorem that $\Pi(R) = \Pi(\lambda)$ for some dominant $\lambda$ if and only if each minimally indomitable set $\mathscr{D}$ in $R$ is a singleton $\mathscr{D} = \{\delt\}$ for some $\delt \in R$ with $wt(\delt) = \lambda$. 
The following is used in the proof of \PiJCorollary. 

\noindent 
{\bf \PiSaturatedCorollary}\ \ {\sl Let $\mathscr{W} \subset \Lambda$ be finite and saturated, and regard $\mathscr{W}$ to be an edge-colored Hasse diagram as in \WeightRemarkTwo.  Let $\mathscr{D}$ be the set of maximal elements in $\mathscr{W}$.  Then $\mathscr{D}$ is precisely the set of prominent elements in $\mathscr{W}$ and is minimally indomitable, and $\mathscr{W} = \displaystyle \bigcup_{\lambda \in \mathscr{D}} \Pi(\lambda)$, an equality of edge-colored directed graphs.}

{\em Proof.} First, observe that $\mathscr{W} = \Pi(\mathscr{W})$ is an equality of edge-colored directed graphs.  
Now suppose $\lambda \in \mathscr{W}$ is maximal, and let $i \in I$.   
Since the $i$-component of $\lambda$ is a chain, then $\lambda$ must be at the top of the chain.  
It follows that $\lambda$ is prominent. 
Since any prominent $\lambda$ is necessarily maximal, it follows that $\mathscr{D}$ is precisely the set of prominent elements in $\mathscr{W}$.  
If $\lambda$ and $\lambda'$ are both prominent, then $wt(\lambda) \leq wt(\lambda')$ means $\lambda \leq \lambda'$, in which case maximality of $\lambda$ means that $\lambda = \lambda'$.  
In particular, the set of prominent elements of $\mathscr{W}$ is minimally indomitable. 
The final claim of the corollary statement now follows from \PiTheorem.\hfill\QED

\noindent
{\bf \PiJCorollary}\ \ {\sl Let $\mathscr{W} \subset \Lambda$ be finite and saturated, and let $J \subseteq I$.  Let $\mathcal{C}$ be some $J$-component of $\mathscr{W}$.  Let $\mathscr{D}$ be the set of all maximal elements of $\mathcal{C}$, i.e.\ the set of all $J$-maximal elements of $\mathscr{W}$ that are in $\mathcal{C}$.  Then $\mathcal{C} \cong \displaystyle \bigcup_{\lambda \in \mathscr{D}} \Pi(\lambda^{J})$.}

{\em Proof.} It follows from the definitions that for all $\mu \in \mathcal{C}$, $wt^{J}(\mu)$ is just the projection $\mu^{J}$.  We now argue that $wt^{J}: \mathcal{C} \longrightarrow \Lambda_{\Phi_J}$ is injective.  Indeed, if $\mu^{J} = \nu^{J}$ for $\mu,\nu \in \mathcal{C}$, then take a path $\mathcal{P}$ from $\mu$ to $\nu$ in $\mathcal{C}$.  Then $\mu + \sum_{j \in J}(a_j(\mathcal{P}) - d_j(\mathcal{P}))\alpha_j = \nu$.  Then $\mu^{J} + \sum_{j \in J}(a_j(\mathcal{P}) - d_j(\mathcal{P}))\alpha_j = \nu^{J}$, so $\sum_{j \in J}(a_j(\mathcal{P}) - d_j(\mathcal{P}))\alpha_j = 0$.  So, $\mu = \nu$.  A similar argument shows that $\mu \myarrow{j} \nu$ in $\mathcal{C}$ if and only if $\mu^{J} + \alpha_j = \nu^{J}$.  So if we regard $wt^{J}(\mathcal{C})$ to be an edge-colored directed graph via the rule $\mu^{J} \myarrow{j} \nu^{J}$ if and only if $\mu^{J} + \alpha_j = \nu^{J}$, then $\mathcal{C} \cong wt^{J}(\mathcal{C})$. Observe now that $\lambda \in \mathcal{C}$ is maximal if and only if $\lambda^{J}$ is maximal in $wt^{J}(\mathcal{C})$. 

We claim that $wt^{J}(\mathcal{C})$ is saturated as a subset of $\Lambda_{\Phi_J}$.  To see this, we let $\mu = \sum_{i \in I}m_i\omega_i$ be any weight in $\mathcal{C}$ and let $\alpha \in \Phi_J$.    
Since $\mathscr{W}$ is saturated, it follows that $\mu-k\alpha \in \mathscr{W}$ if $k$ is any integer between $0$ and $\langle \mu,\alpha^{\vee} \rangle$ (inclusive). 
Since $\mu$ and $\mu-k\alpha$ are comparable in $\mathscr{W}$, they are in the same connected component of $\mathscr{W}$. 
Since there is a path from $\mu$ to $\mu-k\alpha$ in $\mathscr{W}$, and since $\alpha \in \Phi_J$ is a sum of $\{\pm\alpha_{j}\}_{j \in J}$, then there is a path in $\mathscr{W}$ from $\mu$ to $\mu-k\alpha$ along edges that are $J$-colored. 
In particular, $\mu-k\alpha \in \mathcal{C}$, so $\mu^{J}-k\alpha \in wt^{J}(\mathcal{C})$. 
Since $\alpha \in \Phi_J$, then $\alpha^{\vee} \in (\Phi_J)^{\vee} = (\Phi^{\vee})_J$, so $\alpha^{\vee} = \sum_{j \in J}k_j\alpha_{j}^{\vee}$.  
In particular, $\langle \mu,\alpha^{\vee} \rangle = \langle \sum_{i \in I}m_i\omega_i,\sum_{j \in J}k_j\alpha_{j}^{\vee} \rangle = \langle \sum_{j \in J}m_j\omega_{j}^{J},\sum_{j \in J}k_j\alpha_{j}^{\vee} \rangle = \langle \mu^{J},\alpha^{\vee} \rangle$. 
So now we know that if $k$ is any integer between $0$ and $\langle \mu^{J},\alpha^{\vee} \rangle$ (inclusive), then $\mu^{J}-k\alpha \in wt^{J}(\mathcal{C})$. 

We are now in the situation of \PiSaturatedCorollary, so it follows that $\mathcal{C} \cong wt^{J}(\mathcal{C}) = \displaystyle \bigcup_{\lambda \in \mathscr{D}} \Pi(\lambda^{J})$.\hfill\QED

\noindent 
{\bf \PiRankCorollary}\ \ {\sl Let $R$ be ${\Phi}$-structured and connected. 
Let $\mathscr{L} = \{\lambda \in \Lambda\, |\, \lambda \in \Pi(R) \mbox{ and } \mathrm{ht}(\lambda) \geq \mathrm{ht}(wt(\xelt)) \mbox{ for all } \xelt \in R\}$.  
Then $\selt \in R$ has rank zero if and only if $w_{0}(wt(\selt)) \in \mathscr{L}$, and $\uelt \in R$ has maximum rank if and only if $wt(\uelt) \in \mathscr{L}$. 
Now take any $\lambda \in \mathscr{L}$.  
Then for any $\telt \in R$, the rank of $\telt$ is $\langle wt(\telt) - w_{0}(\lambda), \varrho^{\vee} \rangle = \langle wt(\telt) + \lambda, \varrho^{\vee} \rangle$, and its depth is $\langle \lambda - wt(\telt), \varrho^{\vee} \rangle = \langle - w_{0}(\lambda) - wt(\telt), \varrho^{\vee} \rangle$. 
Moreover, $R$ has length $\langle \lambda - w_{0}(\lambda), \varrho^{\vee} \rangle = 2\langle \lambda, \varrho^{\vee} \rangle$.} 

{\em Proof.} Let $\mathscr{D}$ be any minimally indomitable set in $R$.  Suppose $\selt \in R$ has rank zero.  
Since $wt(\selt) \in \Pi(R) = \displaystyle  \bigcup_{\delt \in \mathscr{D}} \Pi(wt(\delt))$, then $w_0(wt(\selt)) = wt(\uelt)$ for some $\uelt \in R$.  
Suppose $\mathrm{ht}(wt(\xelt)) > \mathrm{ht}(w_0(wt(\selt)))$ for some $\xelt \in R$.  
Take any path $\mathcal{P}$ in $R$ from $\uelt$ to $\xelt$, which is possible since $R$ is connected.  
Then $wt(\xelt) - wt(\uelt) = \sum_{i \in I}(a_i(\mathcal{P}) - d_i(\mathcal{P}))\alpha_i$, and $0 < \mathrm{ht}(wt(\xelt)) - \mathrm{ht}(wt(\uelt)) = \langle wt(\xelt)-wt(\uelt),\varrho^{\vee} \rangle = \sum_{i \in I}(a_i(\mathcal{P}) - d_i(\mathcal{P}))$. 
Since $wt(\xelt) \in \Pi(R) = \displaystyle  \bigcup_{\delt \in \mathscr{D}} \Pi(wt(\delt))$, then $w_0(wt(\xelt)) = wt(\relt)$ for some $\relt \in R$.  
So, $\mathrm{ht}(wt(\selt)) - \mathrm{ht}(wt(\relt)) = \langle wt(\selt)-wt(\relt),\varrho^{\vee} \rangle = -\langle w_0(wt(\xelt)-wt(\uelt)),\varrho^{\vee} \rangle = -\langle \sum_{i \in I}(a_i(\mathcal{P}) - d_i(\mathcal{P}))w_0(\alpha_i),\varrho^{\vee} \rangle = \sum_{i \in I}(a_i(\mathcal{P}) - d_i(\mathcal{P})) > 0$. 
Now take a path $\mathcal{Q}$ from $\relt$ to $\selt$ in $R$. 
A computation similar to those we have done already will show that $wt(\selt) - wt(\relt) = \sum_{i \in I}(a_i(\mathcal{Q}) - d_i(\mathcal{Q}))\alpha_i$, and that $\mathrm{ht}(wt(\selt)) - \mathrm{ht}(wt(\relt)) = \sum_{i \in I}(a_i(\mathcal{Q}) - d_i(\mathcal{Q}))$.  
Therefore $\sum_{i \in I}(a_i(\mathcal{Q}) - d_i(\mathcal{Q}))$ is positive. 
But this means that the path $\mathcal{Q}$ from $\relt$ to $\selt$ has more ``upward'' edges than ``downward'' edges, and hence $\relt$ has rank less than zero, a contradiction. 
So we conclude that $\mathrm{ht}(wt(\xelt)) \leq \mathrm{ht}(w_0(wt(\selt)))$ for all $\xelt \in R$. 
Since $w_0(wt(\selt)) = wt(\uelt) \in \Pi(R)$, then it follows that $w_0(wt(\selt)) \in \mathscr{L}$. 
  
Now suppose $w_0(wt(\selt)) \in \mathscr{L}$.  
Let $\telt \in R$, and let $\mathcal{P}$ be any path in $R$ from $\selt$ to $\telt$.  
To show that $\selt$ has rank zero, it suffices to show that $\sum_{i \in I} (a_i(\mathcal{P}) - d_i(\mathcal{P})) \geq 0$. 
The fact that $\Pi(R) = \displaystyle  \bigcup_{\delt \in \mathscr{D}} \Pi(wt(\delt))$ means that there exist $\uelt, \xelt \in R$ such that $w_0(wt(\selt)) = wt(\uelt)$ and $w_0(wt(\telt)) = wt(\xelt)$. 
So $0 \leq \mathrm{ht}(wt(\uelt)) - \mathrm{ht}(wt(\xelt)) = \langle w_0(wt(\selt)-wt(\telt)),\varrho^{\vee} \rangle = \langle wt(\telt)-wt(\selt),\varrho^{\vee} \rangle = \sum_{i \in I} (a_i(\mathcal{P}) - d_i(\mathcal{P}))$. 

The argument that $\uelt \in R$ has maximum rank if and only if $wt(\uelt) \in \mathscr{L}$ is similar.  
So now let $\lambda \in \mathscr{L}$, and let $\telt \in R$.  
We show that the rank of $\telt$ can be calculated as claimed. 
So, take $\selt \in R$ with $w_0(wt(\selt)) = \lambda$. 
Then $\selt$ has rank zero. 
If $\mathcal{P}$ is any path from $\selt$ to $\telt$, then the rank of $\telt$ is $\sum_{i \in I} (a_i(\mathcal{P}) - d_i(\mathcal{P})) = \langle wt(\telt)-wt(\selt),\varrho^{\vee} \rangle = \langle wt(\telt) - w_0(\lambda),\varrho^{\vee} \rangle = \langle wt(\telt) + \lambda,\varrho^{\vee} \rangle$.  
The latter equality follows from the fact that $\langle \lambda,\varrho^{\vee} \rangle = \langle -w_0(\lambda),\varrho^{\vee} \rangle$, cf.\ \WeightRemarkOne.3. 
The remaining claims of the corollary follow from similar computations.\hfill\QED

\noindent 
{\bf \ClosingWeightRemarks}\ \ We close this section with the following remarks about weight diagrams and generalized weight diagrams. 

(A) If $\Pi(R)$ is connected, then $R$ need not be: Consider that $\Pi(S \oplus S) = \Pi(S)$ when $S$ is connected and ${\Phi}$-structured. 

(B) Weight diagrams can be ``tangled'' in the sense of the following examples.   
In case $\myG_{2}$, the set union of weight diagrams $\Pi(5\omega_{1}) \cup \Pi(3\omega_{2})$ with the induced order from $\Lambda$ is connected and ${\myG_{2}}$-structured, and its two maximal elements $5\omega_{1}$ and $3\omega_{2}$ have the same rank but different weights.  
For the connected and ${\myG_{2}}$-structured poset $\Pi(7\omega_{1}) \cup \Pi(4\omega_{2})$, the two maximal elements have different ranks. 

(C) Let $\lambda \in \Lambda^{+}$. 
By \PiRankCorollary, the unique rank function for $\Pi(\lambda)$ is given by $\rho(\mu) = \langle \mu - w_{0}(\lambda), \varrho^{\vee} \rangle = \langle \mu + \lambda, \varrho^{\vee} \rangle$.\hfill\QED

\newpage
\noindent 
{\Large \bf \S \SplittingSection.\ Splitting posets.}  

In this section, we bring together the poset perspective of \S \PosetSection\ and the theory of Weyl symmetric functions in \S \WeylSection\ and initiate our combinatorial study of Weyl bialternants. 
We begin with a new splitting theorem (\InitialSplittingTheorem) and a new splitting definition (\InitialSplittingDefinition) that provide some specific combinatorial and algebraic context for the kinds of problems we are interested in. 
The purpose of the remainder of this section is to begin developing some basic consequences of our splitting poset definition, to provide some examples, and to begin connecting splitting posets to other objects in the literature. 
Before we proceed, we offer a comment on the logical dependence of the results of this section and those of \S \WeylSection: {\em We do not assume Weyl's character formula (\WeylsTheorem) unless explicitly stated otherwise.}   

{\bf A new splitting theorem.} 
The first result of this section was directly inspired by Theorem 2.4 of \cite{Stem}.  
However, for finite systems, our set-up is more general. 
In a certain sense, this result gives sufficient conditions for a poset $R$ to split weight multiplicities.  
This result motivates the definition of ``$(J,\nu)$-splitting poset'' a.k.a.\ ``refined splitting poset'' in \InitialSplittingDefinition. 

\noindent 
{\bf \InitialSplittingTheorem}\ \ {\sl Let $R$ be a ranked poset with edges colored by an index set $I$ for a choice of simple roots for $\Phi$.  
Let $J \subseteq I$, and let $\nu = \sum_{j \in J}\nu_{j}\omega_{j}^{J} \in \Lambda_{\Phi_J}^{+}$. 
Suppose that} $\WGF(R)|_{J}$ {\sl is $W_{J}$-invariant.  
Suppose that there is a subset $\mathcal{S}_{J,\nu}(R)$ of $R$ such that $\nu + wt^{J}(\selt) \in \Lambda_{\Phi_J}^{+}$ for all $\selt$ in $\mathcal{S}_{J,\nu}(R)$. 
Further suppose that there is a bijection $\tau : R \setminus \mathcal{S}_{J,\nu}(R) \longrightarrow R \setminus \mathcal{S}_{J,\nu}(R)$ and a function $\kappa : R \setminus \mathcal{S}_{J,\nu}(R) \longrightarrow J$ such that for all $\xelt \in R \setminus \mathcal{S}_{J,\nu}(R)$, we have $wt^{J}(\tau(\xelt)) = wt^{J}(\xelt) - (1 + \nu_{\kappa(\xelt)} + m_{\kappa(\xelt)}(\xelt))\alpha_{\kappa(\xelt)}$. Then,}
\[\chi^{\Phi_J}_{_{\nu}} \cdot \WGF(R)|_{J} = \sum_{\selt \in \mathcal{S}_{J,\nu}(R)}\chi^{\Phi_J}_{_{\nu + wt^{J}(\selt)}}.\] 

\noindent 
{\bf \NewSplittingRemarks.}\ \ Before we give the proof, first some observations/comments.  

(A) The key identity at the end of the theorem statement pertains to each of the splitting, product decomposition, and branching problems.  
For more on this, see \InitialSplittingDefinition/\NewDefinitionRemark.D below. 

(B) We call $\kappa$ a {\em vertex-coloring function}.  
While the bijection $\tau$ and the vertex-coloring function $\kappa$ are both needed in order to make the preceding theorem work, in many situations, some of which are explored in later sections, only a vertex-coloring function is needed. 
So, vertex-coloring will be a crucial feature of our poset-theoretic approach to studying Weyl symmetric functions. 

(C) One way to guarantee that $\WGF(R)|_{J}$ is $W_{J}$-invariant is for $R$ to be ${\Phi_J}$-structured and for the $j$-components of $R$ to be rank symmetric for each $j \in J$, cf. \WInvariantLemma.  

(D) While the theorem can be useful in the generality with which it is stated for subsets $J$ of $I$, for purposes of the proof we lose no generality in assuming that $J = I$, since the hypotheses of the theorem allow us to simply ignore all edges with colors in $I \setminus J$. 

(E) The hypothesis that $R$ is an edge-colored, ranked poset can be relaxed to just the assumption that $R$ is a set together with a weight-assigning function $\weight^{J} : R \longrightarrow \Lambda_{\Phi_J}$ such that the group ring element $\sum_{\xelt \in R} e^{\smallweight^{J}(\xelt)}$ is $W_J$-invariant; the proof given below would not have to change much at all. 
However, it is more convenient for our purposes to maintain the poset-oriented language of the theorem statement.

(F) It is sometimes helpful for purposes of specialization to view a given $\chi_{_{\lambda}}$ as a function in the variables $z_{k} := e^{\omega_{k}}$, for $k \in I$. Then we will use $\chi_{_{\lambda}}|_{\{z_{k} := Q_{k}\}_{k \in I}}$ to notate the algebraic expression obtained by replacing each occurrence of each $z_{k} := e^{\omega_{k}}$ with a given quantity $Q_{k}$.\hfill\QED 

{\em Proof of \InitialSplittingTheorem.} 
As observed above, we lose no generality in assuming that $J=I$, since otherwise we can simply remove from the picture all edges with colors not in the set $J$. 
Let $\mu \in \Lambda$ and $i \in I$.  
Then $\displaystyle \mathcal{A}(e^{\mu}) = \sum_{w \in W} \det(w)e^{w(\mu)} = \sum_{w \in W} \det(w \cdot s_i)e^{(w \cdot s_i)(\mu)} = -\sum_{w \in W} \det(w)e^{w(\mu - \langle \mu,\alpha_{i}^{\vee} \rangle \alpha_i)} = -\mathcal{A}(e^{\mu - \langle \mu,\alpha_{i}^{\vee} \rangle \alpha_i})$.  
Taking $\mu = \varrho + \nu + wt(\xelt)$ and combining the previous computation with the hypothesis that $wt(\tau(\xelt)) = wt(\xelt) - (1 + \nu_{\kappa(\xelt)} + m_{\kappa(\xelt)}(\xelt))\alpha_{\kappa(\xelt)}$, we obtain that 
$\mathcal{A}(e^{\varrho + \nu + wt(\xelt)}) = -\mathcal{A}(e^{\varrho + \nu + wt(\tau(\xelt))})$.
Then using the pairing $\tau$, we compute that $\sum_{\xelt \in R \setminus \mathcal{S}_{I,\nu}(R)}\mathcal{A}(e^{\varrho+\nu+wt(\xelt)}) = \sum_{\xelt \in R \setminus \mathcal{S}_{I,\nu}(R)}\mathcal{A}(e^{\varrho+\nu+wt(\tau(\xelt))}) = -\sum_{\xelt \in R \setminus \mathcal{S}_{I,\nu}(R)}\mathcal{A}(e^{\varrho+\nu+wt(\xelt)})$, which yields the desired cancelling of certain terms we will need shortly: $\sum_{\xelt \in R \setminus \mathcal{S}_{I,\nu}(R)}\mathcal{A}(e^{\varrho+\nu+wt(\xelt)}) = 0$. 

In what follows, note how the $W$-invariance of $\WGF(R)$ is used in the initial computation: 
\begin{eqnarray*}
\mathcal{A}(e^{\varrho+\nu})\WGF(R) & = & \sum_{w \in 
W}\left(\det(w)e^{w(\varrho+\nu)}\WGF(R)\right)\\ 
& = & \sum_{w \in 
W}\left(\det(w)e^{w(\varrho+\nu)}\left(\sum_{\xelt \in 
R}e^{w(wt(\xelt))}\right)\right)\\  
& = & \sum_{\xelt \in R}\left(\sum_{w \in 
W}\det(w)e^{w(\varrho+\nu+wt(\xelt))}\right)\\ 
& = & \sum_{\xelt \in R}\mathcal{A}(e^{\varrho+\nu+wt(\xelt)})\\ 
& = & \sum_{\selt \in \mathcal{S}_{I,\nu}(R)} \mathcal{A}(e^{\varrho + \nu + wt(\selt)}) + \sum_{\xelt \in R \setminus 
\mathcal{S}_{I,\nu}(R)}\mathcal{A}(e^{\varrho+\nu+wt(\xelt)})\\ & = & 
\sum_{\selt \in \mathcal{S}_{I,\nu}(R)} \mathcal{A}(e^{\varrho + \nu + wt(\selt)}).
\end{eqnarray*}
Divide both sides by $\mathcal{A}(e^{\varrho})$ to get the identity $\displaystyle \chi_{_{\nu}} \cdot \WGF(R) = \sum_{\selt \in \mathcal{S}_{I,\nu}(R)}\chi_{_{\nu + wt(\selt)}}$.\hfill\QED

\noindent 
{\bf \InitialSplittingExample}\ \ In this example, we take $\Phi = \myC_{2}$. 
See \NotSturdyFigure\ for a picture of the poset $R$ and some further data. 
We take $J = I = \{1,2\}$, and $\nu = 0$. 
The given poset $R$ is not ${\myC_2}$-structured. 
With vertices as labelled in the figure, we let $\mathcal{S}_{I,0}(R) := \{\mathbf{max}\}$, and define the functions $\tau : R \setminus \mathcal{S}_{I,0}(R) \longrightarrow R \setminus \mathcal{S}_{I,0}(R)$ and $\kappa : R \setminus \mathcal{S}_{I,0}(R) \longrightarrow I$ by the table presented in the figure. 
The data from that table confirms the main requirements of \InitialSplittingTheorem.  
It follows that $\WGF(R) = \chi^{\mytinyC_2}_{_{2\omega_{1}}}$.\hfill\QED

\begin{figure}[ht]
\begin{center}
{\small {\NotSturdyFigure:} A non-$\myC_{2}$-structured poset $R$ that nonetheless generates the Weyl bialternant $\WGF(R) = \chi^{\mytinyC_2}_{_{2\omega_{1}}}$.}\\
{\footnotesize (See \InitialSplittingExample\ for further details about this figure.)}

\vspace*{0.1in}
\setlength{\unitlength}{0.7cm}
\begin{picture}(15,10)
\put(-2.45,9){\Large $R$}
\put(-1.7,9.05){\vector(3,-1){1.4}}
\thicklines
\multiput(4.5,0.5)(-1.5,1.5){3}{\circle*{0.22}}
\multiput(4.5,3.5)(-1.5,1.5){4}{\circle*{0.22}}
\multiput(4.5,6.5)(-1.5,1.5){3}{\circle*{0.22}}
\put(1.7,9.4){\scriptsize $\mathbf{max}$}
\put(3.2,7.9){\scriptsize $\mathbf{b}$}
\put(4.7,6.4){\scriptsize $\mathbf{d}$}
\put(3.2,4.9){\scriptsize $\mathbf{e}$}
\put(4.7,3.4){\scriptsize $\mathbf{g}$}
\put(-0.4,7.9){\scriptsize $\mathbf{a}$}
\put(1.1,6.4){\scriptsize $\mathbf{c}$}
\put(1.1,3.4){\scriptsize $\mathbf{f}$}
\put(2.6,1.9){\scriptsize $\mathbf{h}$}
\put(4.1,0.4){\scriptsize $\mathbf{i}$}
%
%
\put(4.4,0.6){\color{Red}{\line(-1,1){1.3}}}
\put(3.6,1.1){\color{Red}{\em 1}}
\put(2.9,2.1){\color{Red}{\line(-1,1){1.3}}}
\put(2.1,2.6){\color{Red}{\em 1}}
%
\put(4.4,3.6){\color{Red}{\line(-1,1){1.3}}}
\put(3.6,4.1){\color{Red}{\em 1}}
\put(2.9,5.1){\color{Cyan}{\line(-1,1){1.3}}}
\put(2.1,5.6){\color{Cyan}{\em 2}}
\put(1.4,6.6){\color{Red}{\line(-1,1){1.3}}}
\put(0.6,7.1){\color{Red}{\em 1}}
\put(4.4,6.6){\color{Cyan}{\line(-1,1){1.3}}}
\put(3.6,7.1){\color{Cyan}{\em 2}}
\put(2.9,8.1){\color{Red}{\line(-1,1){1.3}}}
\put(2.1,8.6){\color{Red}{\em 1}}
%
%
\put(3.1,2.1){\color{Cyan}{\line(1,1){1.3}}}
\put(3.6,2.6){\color{Cyan}{\em 2}}
\put(1.6,3.6){\color{Cyan}{\line(1,1){1.3}}}
\put(2.1,4.1){\color{Cyan}{\em 2}}
\put(3.1,5.1){\color{Red}{\line(1,1){1.3}}}
\put(3.6,5.6){\color{Red}{\em 1}}
%
\put(1.6,6.6){\color{Red}{\line(1,1){1.3}}}
\put(2.1,7.1){\color{Red}{\em 1}}
\put(0.1,8.1){\color{Red}{\line(1,1){1.3}}}
\put(0.6,8.6){\color{Red}{\em 1}}
\put(7,5.25){\parbox[t]{3.5in}{\footnotesize 
\begin{tabular}{|c||c|c|c|}
\hline
Vertex $\xelt$ & $\kappa(\xelt)$ & $\tau(\xelt)$ & $wt(\xelt)$\\
\hline
\hline
$\aelt$ & {\color{Red}$1$} & $\gelt$ & $0$\\
\hline
$\belt$ & {\color{Red}$1$} & $\celt$ & ${\color{Red}\alpha_{1}}+{\color{Cyan}\alpha_{2}}$\\
\hline
$\celt$ & {\color{Red}$1$} & $\belt$ & ${\color{Cyan}\alpha_{2}}$\\
\hline
$\delt$ & {\color{Cyan}$2$} & $\delt$ & ${\color{Red}\alpha_{1}}$\\
\hline
$\eelt$ & {\color{Cyan}$2$} & $\felt$ & $0$\\
\hline
$\felt$ & {\color{Cyan}$2$} & $\eelt$ & $-{\color{Cyan}\alpha_{2}}$\\
\hline
$\gelt$ & {\color{Red}$1$} & $\aelt$ & $-{\color{Red}\alpha_{1}}$\\
\hline
$\helt$ & {\color{Red}$1$} & $\ielt$ & $-{\color{Red}\alpha_{1}}-{\color{Cyan}\alpha_{2}}$\\
\hline
$\ielt$ & {\color{Red}$1$} & $\helt$ & $-2{\color{Red}\alpha_{1}}-{\color{Cyan}\alpha_{2}}$\\
\hline
\end{tabular}}}
\end{picture}

\end{center}
\end{figure}

{\bf Splitting posets and refined splitting posets.}  
In view of the preceding theorem, we make the following definition.  

\noindent 
{\bf \InitialSplittingDefinition}\ \ 
We say an ${\Phi}$-structured poset $R$ is a {\em splitting poset for a Weyl symmetric function} $\chi$ if $\WGF(R) = \chi$. 
A version of this definition that pertains to the product decomposition and branching problems is as follows: Let $J \subseteq I$ and $\nu \in \Lambda_{\Phi_J}^{+}$.  
A $(J,\nu)${\em -splitting poset} is an ${\Phi}$-structured poset $R$ together with a set $\mathcal{S}_{J,\nu}(R)$ for which $\nu + wt^{J}(\selt) \in \Lambda_{\Phi_J}^{+}$ for all $\selt \in \mathcal{S}_{J,\nu}(R)$ and 
\begin{equation}\chi^{\Phi_J}_{_{\nu}} \cdot \WGF(R)|_{J} = \sum_{\selt \in \mathcal{S}_{J,\nu}(R)}\chi^{\Phi_J}_{_{\nu + wt^{J}(\selt)}}.\end{equation} 
A $(J,\nu)$-splitting poset will sometimes more generically be referred to as a {\em refined splitting poset}.\hfill\QED 


\noindent 
{\bf \NewDefinitionRemarks.}\ \  
Some initial comments on the definition: 

(A) We will sometimes use the phrase ``splitting poset'' as a catch-all for all the kinds of posets featured in our definition.  
But most of the time, the phrase ``splitting poset'' by itself and without the ``$(J,\nu)$'' or ``refined'' prefix will refer in particular to the $J=I$ and $\nu = 0$ situation of the first part of the definition. 

(B) The explicit meaning of equation (6) depends on the explicitness of the descriptions of $R$ and of $\mathcal{S}_{J,\nu}(R)$. 

(C) It is a consequence of \FinitenessTheorem\ that $\WGF(R)|_{J}$ is a $W_{J}$-symmetric function if $R$ is $(J,\nu)$-splitting. 
However, $\WGF(R)|_{J}$ need not be a \underline{nonnegative} integer linear combination of Weyl bialternants. 
For example, if $R$ is as pictured in \NotSchurPositiveFigure, then $\WGF(R) = \chi^{\mytinyA_2}_{_{\omega_1 + \omega_2}} - \chi^{\mytinyA_2}_{_{0}}$. 
Taking $J=I$, $\nu = \omega_1$, and $\mathcal{S}_{I,\omega_1}(R) = \{\aelt,\belt\}$, then $\chi^{\mytinyA_2}_{_{\omega_1}} \cdot \WGF(R) = \sum_{\selt \in \mathcal{S}_{I,\omega_1}(R)}\chi^{\mytinyA_2}_{_{\omega_1 + wt(\selt)}} = \chi^{\mytinyA_2}_{_{2\omega_1 + \omega_2}} + \chi^{\mytinyA_2}_{_{2\omega_2}}$.

\begin{figure}[ht]
\begin{center}
{{\NotSchurPositiveFigure:} An $\myA_{2}$-structured poset $R$ for which $\WGF(R) = \chi^{\mytinyA_2}_{_{\omega_1 + \omega_2}} - \chi^{\mytinyA_2}_{_{0}}$.}\\
{\small (See \NewDefinitionRemark.C for further details about this figure.)}

\vspace*{0.1in}
\setlength{\unitlength}{0.7cm}
\begin{picture}(5.5,7)
\put(-1.35,6.3){\Large $R$}
\put(-0.6,6.45){\vector(3,-1){1.6}}
\thicklines
\multiput(3,0.5)(-1.5,1.5){2}{\circle*{0.22}}
\multiput(4.5,2)(-1.5,1.5){3}{\circle*{0.22}}
\multiput(4.5,5)(-1.5,1.5){2}{\circle*{0.22}}
%
%
\put(2.9,0.6){\color{Red}{\line(-1,1){1.3}}}
\put(2.1,1.1){\color{Red}{\em 1}}
\put(4.4,2.1){\color{Red}{\line(-1,1){1.3}}}
\put(3.6,2.6){\color{Red}{\em 1}}
\put(2.9,3.6){\color{Cyan}{\line(-1,1){1.3}}}
\put(2.1,4.1){\color{Cyan}{\em 2}}
\put(4.4,5.1){\color{Cyan}{\line(-1,1){1.3}}}
\put(3.6,5.6){\color{Cyan}{\em 2}}
%
%
\put(3.1,0.6){\color{Cyan}{\line(1,1){1.3}}}
\put(3.6,1.1){\color{Cyan}{\em 2}}
\put(1.6,2.1){\color{Cyan}{\line(1,1){1.3}}}
\put(2.1,2.6){\color{Cyan}{\em 2}}
\put(3.1,3.6){\color{Red}{\line(1,1){1.3}}}
\put(3.6,4.1){\color{Red}{\em 1}}
\put(1.6,5.1){\color{Red}{\line(1,1){1.3}}}
\put(2.1,5.6){\color{Red}{\em 1}}
\end{picture}
\end{center}
\end{figure}

(D) Obviously our interest in such objects is motivated in part by the special cases when $(J,\nu)$ is one of the pairs $(I,\nu)$, $(J,0)$, or $(I,0)$, and how these special cases pertain to the aforementioned product decomposition, branching, and splitting problems. 
For $(I,\nu)$, equation (6) becomes 
\[\chi_{_{\nu}} \cdot \WGF(R) = \sum_{\selt \in \mathcal{S}_{I,\nu}(R)}\chi_{_{\nu + wt(\selt)}},\] 
which is a sort of Littlewood--Richardson rule for the decomposition of the product $\chi_{_{\nu}} \cdot \WGF(R)$ as a nonnegative-integer linear combination of Weyl bialternants.  
For $(J,0)$, equation (6) becomes 
\[\WGF(R)|_{J} = \sum_{\selt \in \mathcal{S}_{J,0}(R)}\chi^{\Phi_J}_{_{wt^{J}(\selt)}},\] 
which is a sort of branching rule for $\WGF(R)$ when we view it as a $W_{J}$-symmetric function. 
Finally, for $(I,0)$ equation (6) becomes 
\[\WGF(R) = \sum_{\selt \in \mathcal{S}_{I,0}(R)}\chi_{_{wt(\selt)}}.\] 
In this case, $R$ is a splitting poset for the $W$-symmetric function $\WGF(R)$.  
Of course, when $\mathcal{S}_{I,0}(R) = \{\selt\}$, then $\WGF(R) = \chi_{_{wt(\selt)}}$, which is a case of especial interest because of the combinatorial consequences for $R$, see \MainCorollary\ below.

(E) We have several reasons for including the ${\Phi}$-structure property as one of the requirements for a poset to be $(J,\nu)$-splitting. 
First, this constraint is natural in that it is shared with crystal graphs, admissible systems, and supporting graphs. 
Splitting posets are generalizations of these structures: We will eventually see that crystal graphs and admissible systems are $(J,\nu)$-splitting posets for all pairs $(J,\nu)$ and that supporting graphs are $(I,0)$-splitting posets. 
Second, this feature is crucial in order to guarantee a certain combinatorial niceness of splitting posets. 
In particular, a connected $(I,0)$-splitting poset is rank symmetric and rank unimodal, and its rank generating function can be calculated via root and weight related computations (see \MainCorollary\ below). 
Obviously these outcomes can fail in the absence of the ${\Phi}$-structure property, see \InitialSplittingExample\ above. 
Also, some interesting classification results are possible when this constraint is imposed. 
For example, we have been able to demonstrate elsewhere that if a modular lattice $L$ is `$A$-structured' by the rows of some real $n \times n$ matrix $A = (A_{ij})_{i,j \in I}$ such that $\ecolor_{L}: \mathcal{E}(L) \longrightarrow I$ is surjective, then $A$ is a Cartan matrix. 
And third, the ${\Phi}$-structure property is useful.   
For example, it is used to deduce $W$-invariance of the weight generating function $\WGF(R)$ in \WInvariantLemma. 
Moreover, in \MainColoringTheorem, it subsumes the bijection $\tau$ so that the refined splitting problem reduces to finding the special vertex-coloring function $\kappa$.\hfill\QED  

{\bf A version of Stanley's problem for Weyl symmetric functions.} 

\noindent 
{\bf \StanleyProblem}\ \ Our variation of Stanley's problem (Problem 3 of \cite{StanUnim}) is: For (families of) Weyl bialternants or, more generally, Weyl symmetric functions, produce splitting posets (or splitting modular/distributive lattices) with explicit and interesting combinatorial descriptions.\hfill\QED 

Overall, our main programmatic interest is in the $(I,0)$ version of this splitting problem and sometimes the $(J,0)$ version as well. 
That said, in later sections of this monograph we will work within our poset-theoretic framework to produce, from scratch, a ``crystalline'' splitting poset for any given $\chi_{_{\lambda}}$ that is $(J,\nu)$-splitting for all $J \subseteq I$ and $\nu \in \Lambda_{\Phi_J}^{+}$. 

{\bf Some combinatorial properties of splitting posets.}  
Splitting posets for Weyl bialternants have certain salient combinatorial features, as stated in the next result. 
This result and the succeeding result (\OperationsLemma) are, respectively,  restatements of Proposition 2.4 and Lemma 2.2 of \cite{ADLMPPW} but, with the exception of the unimodality claim of \MainCorollary\ below,  are demonstrated here with no dependence on Lie algebra representation theory. 
For a more explicit version of the $\RGF$ formula in \MainCorollary, see Theorem 10.6.3 of \cite{DonDistributive}. 

\noindent 
{\bf \MainCorollary}\ \ {\sl 
Let $\lambda \in \Lambda^{+}$ and suppose $R$ is a splitting poset for $\chi_{_{\lambda}}$. 
Let $\rho: R \longrightarrow \mathbb{Z}$ be the function  $\xelt \stackrel{\rho}{\mapsto} \langle wt(\xelt)+\lambda,\varrho^{\vee} \rangle$.  
Then there is a unique element $\melt \in R$ such that $\rho(\melt)$ is largest, and $\lambda = wt(\melt)$; there is a unique element $\nelt \in R$ such that $\rho(\nelt) = 0$, and $w_{0}(\lambda) = wt(\nelt)$; and $\rho$ is a rank function with $\rho(R) = \{0,1,\ldots,2\langle\lambda,\varrho^{\vee}\rangle\}$.  
Relative to $\rho$, $R$ is rank symmetric 
and}
\[\RGF(R,q) = \frac{\mbox{$\displaystyle \prod_{\alpha \in 
\Phi^{+}}$}(1-q^{\langle \lambda+\varrho,\alpha^{\vee} 
\rangle})}{\mbox{$\displaystyle \prod_{\alpha \in 
\Phi^{+}}(1-q^{\langle \varrho,\alpha^{\vee} \rangle})$}} = q^{-\langle w_{0}.\lambda,\varrho^{\vee} \rangle} \chi_{_{\lambda}}|_{\{z_{k} := q^{\langle \omega_{k},\varrho^{\vee} \rangle}\}_{k \in I}}.\]
{\sl Assuming Weyl's character formula (\WeylsTheorem), then} $\RGF(R,q)$ {\sl is unimodal as well.} 

{\em Proof.} Since $d_{\lambda,\lambda} = 1$, then there is a unique $\melt \in R$ with $wt(\melt) = \lambda$. 
Since by \WeightRemarkOne\ $\lambda$ is the unique maximal element of $\Pi(\lambda)$, then $\langle wt(\xelt),\varrho^{\vee} \rangle < \langle wt(\melt),\varrho^{\vee} \rangle$ for all $\xelt \not= \melt$ in $R$, and hence $\rho(\xelt) < \rho(\melt)$.  
Similarly obtain the unique element $\nelt$.  
Now let $R^{0}$ be the connected component of $R$ containing $\melt$.  
By \PiRankCorollary, $\rho|_{R^{0}}$ is the unique rank function for $R^{0}$, and $\rho(R^{0}) = \{0,1,\ldots,2\langle\lambda,\varrho^{\vee}\rangle\}$. 
From \PiTheorem, it follows that $\Pi(\lambda) \subseteq \Pi(R^{0})$, a setwise inclusion. 
Since $\Pi(R^{0}) \subseteq \Pi(R) = \Pi(\lambda)$, then $\Pi(R^{0}) = \Pi(\lambda)$ is an equality of sets. 
So $\rho(R) = \{0,1,\ldots,2\langle \lambda,\varrho^{\vee} \rangle\}$ as well. 
Now suppose $\xelt \myarrow{i} \yelt$ in $R$.  Then since $wt(\xelt) + \alpha_{i} = wt(\yelt)$, it follows that $\rho(\xelt) + 1 = \rho(\yelt)$.  
So $\rho$ is a rank function for $R$.  
Now, $\RGF(R,q) = \sum_{\xelt \in R}q^{\langle wt(\xelt)+\lambda,\varrho^{\vee}\rangle} = \sum_{\mu \in \Pi(\lambda)}d_{\lambda,\mu}q^{\langle \mu+\lambda,\varrho^{\vee}\rangle}$, which becomes the right-hand side of the identity of the proposition statement by 
\TFAEConsequences.1.  
Rank symmetry and rank unimodality of $R$ now follow from \TFAEConsequences.1 and \UnimodalTheorem.\hfill\QED 

\begin{figure}[t]
\begin{center}
{\TangledFigure.1:} In the language of \S \FibrousSection, the $\myA_1 \oplus \myA_1$-structured posets below are primary. 

{\footnotesize $\big($These are connected and edge-minimal splitting posets for $2\chi_{_{\omega_{1}+\omega_{2}}}^{\mytinyA_{1} \oplus \mytinyA_{1}}$ (left picture) and $3\chi_{_{\omega_{1}+\omega_{2}}}^{\mytinyA_{1} \oplus \mytinyA_{1}}$ (right picture).$\big)$}

\setlength{\unitlength}{0.7cm}
\begin{picture}(8,8.5)
\thicklines
\put(4,0){\circle*{0.22}}
\put(2,4){\circle*{0.22}}
\put(0,8){\circle*{0.22}}
\put(6,4){\circle*{0.22}}
\put(8,8){\circle*{0.22}}
\put(2,7){\circle*{0.22}}
\put(4,6){\circle*{0.22}}
\put(6,7){\circle*{0.22}}
\put(3.925,0.15){\color{Red}\line(-1,2){1.85}}
\put(2.9,1.7){\color{Red}{\scriptsize \em 1}}
\put(4.075,0.15){\color{Cyan}\line(1,2){1.85}}
\put(4.85,1.7){\color{Cyan}{\scriptsize \em 2}}
\put(1.925,4.15){\color{Cyan}\line(-1,2){1.85}}
\put(0.9,5.7){\color{Cyan}{\scriptsize \em 2}}
\put(6.075,4.15){\color{Red}\line(1,2){1.85}}
\put(6.85,5.7){\color{Red}{\scriptsize \em 1}}
\put(1.875,7.05){\color{Red}\line(-2,1){1.75}}
\put(0.9,7.35){\color{Red}{\scriptsize \em 1}}
\put(3.875,6.05){\color{Cyan}\line(-2,1){1.75}}
\put(2.9,6.35){\color{Cyan}{\scriptsize \em 2}}
\put(4.125,6.05){\color{Red}\line(2,1){1.75}}
\put(4.9,6.35){\color{Red}{\scriptsize \em 1}}
\put(6.125,7.05){\color{Cyan}\line(2,1){1.75}}
\put(6.9,7.35){\color{Cyan}{\scriptsize \em 2}}
\end{picture}
\hspace*{1in}
\begin{picture}(8,8.5)
\thicklines
\put(4,0){\circle*{0.22}}
\put(2,2){\circle*{0.22}}
\put(0,8){\circle*{0.22}}
\put(6,2){\circle*{0.22}}
\put(8,8){\circle*{0.22}}
\put(1,6){\circle*{0.22}}
\put(2,4){\circle*{0.22}}
\put(3,6){\circle*{0.22}}
\put(4,8){\circle*{0.22}}
\put(5,6){\circle*{0.22}}
\put(6,4){\circle*{0.22}}
\put(7,6){\circle*{0.22}}
\put(3.9,0.1){\color{Red}\line(-1,1){1.8}}
\put(2.8,0.95){\color{Red}{\scriptsize \em 1}}
\put(4.1,0.1){\color{Cyan}\line(1,1){1.8}}
\put(4.85,0.95){\color{Cyan}{\scriptsize \em 2}}
\put(1.9,2.1){\color{Cyan}\line(-1,3){1.92}}
\put(0.9,4.3){\color{Cyan}{\scriptsize \em 2}}
\put(6.1,2.1){\color{Red}\line(1,3){1.92}}
\put(6.825,4.3){\color{Red}{\scriptsize \em 1}}
\put(0.95,6.125){\color{Red}\line(-1,2){0.875}}
\put(0.4,7){\color{Red}{\scriptsize \em 1}}
\put(1.95,4.125){\color{Cyan}\line(-1,2){0.875}}
\put(1.4,5){\color{Cyan}{\scriptsize \em 2}}
\put(2.05,4.125){\color{Red}\line(1,2){0.875}}
\put(2.35,5){\color{Red}{\scriptsize \em 1}}
\put(3.05,6.125){\color{Cyan}\line(1,2){0.875}}
\put(3.35,7){\color{Cyan}{\scriptsize \em 2}}
\put(4.95,6.125){\color{Red}\line(-1,2){0.875}}
\put(4.4,7){\color{Red}{\scriptsize \em 1}}
\put(5.95,4.125){\color{Cyan}\line(-1,2){0.875}}
\put(5.4,5){\color{Cyan}{\scriptsize \em 2}}
\put(6.05,4.125){\color{Red}\line(1,2){0.875}}
\put(6.35,5){\color{Red}{\scriptsize \em 1}}
\put(7.05,6.125){\color{Cyan}\line(1,2){0.875}}
\put(7.35,7){\color{Cyan}{\scriptsize \em 2}}
\end{picture}
\end{center}
\end{figure}

Various operations on splitting posets yield new splitting posets. 

\noindent 
{\bf \OperationsLemma}\ \ {\sl (1) Let $\lambda_{1},\ldots,\lambda_{k}$ 
be dominant weights.  
Let $R$ be a splitting poset for the Weyl 
symmetric function $\chi := \sum \chi_{_{\lambda_{i}}}$.  
Then $R^{*}$ is a splitting 
poset for $\chi^{*} = \sum \chi_{_{-\omega_{0}(\lambda_{i})}}$.  
The $\sigma_{0}$-recolored dual $R^{\bowtie}$ is a 
splitting poset for $\chi^{\bowtie} = \chi$.}  
{\sl (2) If $R$ and $Q$ are splitting posets for 
Weyl symmetric functions $\chi$ and  
$\chi'$ respectively, then the edge-colored poset 
$Q \times R$ is a splitting poset for 
$\chi'\, \chi$, 
and the edge-colored poset $Q \oplus R$ is a splitting 
poset for $\chi' + \chi$.  (3) Suppose $\Phi = \Phi_{1} \disjointunion 
\Phi_{2}$ for root systems $\Phi_{1}$ and $\Phi_{2}$.  
For $i=1,2$ let 
$\chi_{_{i}}$ be a $\Phi_{i}$-Weyl symmetric function with splitting 
poset $R_{i}$.  
Viewing each $\chi_{_{i}}$ as 
a $\Phi$-Weyl symmetric function, the edge-colored poset $R_{1} \times 
R_{2}$ is a splitting poset for the $\Phi$-Weyl symmetric function
$\chi_{_{1}}\, \chi_{_{2}}$.}

{\em Proof.} {\sl (2)} and {\sl (3)} are routine and follow from the 
definitions.  The claims of {\sl (1)} easily reduce to knowing that 
for each $i$, 
$\chi_{_{\lambda_{i}}}^{*} = \chi_{_{-w_{0}(\lambda_{i})}}$ and that 
$\chi_{_{\lambda_{i}}}^{\bowtie} = \chi_{_{\lambda_{i}}}$, cf.\ \StarProp.\hfill\QED 

{\bf Edge-minimal splitting posets.} 
A splitting poset $R$ for a $W$-symmetric function $\chi$ is {\em edge-minimal} if no other splitting poset for $\chi$ is isomorphic to a proper edge-colored subgraph of $R$.  

\noindent 
{\bf \MinimalProp}\ \ {\sl Let $R$ be a splitting poset for a Weyl symmetric function $\chi$.  
If $R$ is fibrous, then $R$ is edge-minimal.} 

{\em Proof.} We may 
write $\chi = \sum_{\mu \in \Lambda}d_{\mu}e^{\mu}$ 
for certain nonnegative integers $d_{\mu}$.  Let $Q$ be any 
splitting poset for $\chi$, and let 
$Q_{\mu} := \{\telt \in Q\, |\, wt(\telt) = \mu\}$.  For any integer 
$m$ and any $i \in I$, let $Q_{m,i} := \{\telt \in Q\, |\, m_{i}(\telt) = 
m\}$.  Then, $Q_{m,i} = \disjointunion_{\!\!\!\mu \in \Lambda, 
\langle \mu,\alpha_{i}^{\vee} \rangle = m\,}Q_{\mu}$, a disjoint union.  
Since $|Q_{\mu}| = d_{\mu}$, 
then $|Q_{m,i}| = \sum_{\mu \in \Lambda, 
\langle \mu,\alpha_{i}^{\vee} \rangle = m}d_{\mu}$. Now suppose $R'$ is a 
splitting poset for $\chi$ that is isomorphic to a proper edge-colored 
subgraph of $R$.  Let $m$ be the largest integer for which there is 
an $i \in I$ with an $i$-component $\mathbf{comp}_{i}(\telt)$ in $R$ of 
length $m$ such that the subgraph of $R$ isomorphic to $R'$ does not 
retain all of the edges of $\mathbf{comp}_{i}(\telt)$.  Now, 
$\mathbf{comp}_{i}(\telt)$ is a chain, and without loss of generality 
we may assume that $\telt$ is the maximal element of this chain.  
So in $R'$, $m'_{i}(\telt)$ is strictly less than $m_{i}(\telt) = m$ 
in $R$.  Given our choice for $m$ and $i$, we cannot by removing 
edges from $R$ obtain any other element $\selt$ for which 
$m'_{i}(\selt) = m$.  That is, $|R'_{m,i}| < |R_{m,i}|$.  But 
this contradicts that fact that we must have $|R'_{m,i}| = 
|R_{m,i}| = \sum_{\mu \in \Lambda, 
\langle \mu,\alpha_{i}^{\vee} \rangle = m}d_{\mu}$.\hfill\QED

\vspace{0.1in}
\noindent
{\bf \EdgeMinimalExamples}\ \ The splitting poset 
\parbox{1.7cm}{
\setlength{\unitlength}{0.35cm}
\begin{picture}(0,2)
\put(0,0){\circle*{0.35}}
\put(2,0){\circle*{0.35}}
\put(1,1){\circle*{0.35}}
\put(4,1){\circle*{0.35}}
\put(2.6,0.75){$\oplus$}
\put(0,2){\circle*{0.35}}
\put(2,2){\circle*{0.35}}
\put(0,0){\line(1,1){2}}
\put(2,0){\line(-1,1){2}}
\end{picture}} 
for the $A_{1}$-Weyl symmetric function 
$2\chi_{_{2\omega_{1}}}$ is non-fibrous and edge-minimal. 
The edge-minimal splitting posets depicted in \TangledFigure\ are ``tangled'' in the sense that they are connected but each has more than one maximal element. 
This is a feature we will seek to avoid, see \MinQuasiMinTheorem\ and \StembridgeUntangledTheorem.\hfill\QED

\noindent
{\bf \EdgeMinimalQuestions}\ \ 
An open question is whether for two edge-minimal splitting posets $Q$ and $R$ for some Weyl symmetric function, it is that case that $|\mathcal{E}_{i}(Q)| = |\mathcal{E}_{i}(R)|$ for all $i \in I$. 
A related question is: If $R$ is an edge-minimal splitting poset for some Weyl bialternant, must $R$ be fibrous?\hfill\QED 

\begin{figure}[t]
\begin{center}
{\TangledFigure.2:} In the language of \S \FibrousSection, the $\myA_2$-structured posets below are primary. 

{\footnotesize $\big($These are connected and edge-minimal splitting posets for $2\chi_{_{\omega_{1}+\omega_{2}}}^{\mytinyA_{2}}$ (top picture) and $3\chi_{_{\omega_{1}+\omega_{2}}}^{\mytinyA_{2}}$ (bottom picture).$\big)$}

\setlength{\unitlength}{0.7cm}
\begin{picture}(12,10.5)
\thicklines
\put(2,0){\circle*{0.22}}
\put(0,2){\circle*{0.22}}
\put(0,4){\circle*{0.22}}
\put(0,6){\circle*{0.22}}
\put(6,10){\circle*{0.22}}
\put(12,6){\circle*{0.22}}
\put(12,4){\circle*{0.22}}
\put(12,2){\circle*{0.22}}
\put(10,0){\circle*{0.22}}
\put(8,2){\circle*{0.22}}
\put(8,4){\circle*{0.22}}
\put(8,6){\circle*{0.22}}
\put(6,8){\circle*{0.22}}
\put(4,6){\circle*{0.22}}
\put(4,4){\circle*{0.22}}
\put(4,2){\circle*{0.22}}
\put(1.9,0.1){\color{Red}\line(-1,1){1.8}}
\put(0.85,0.9){\color{Red}{\scriptsize \em 1}}
\put(0,2.125){\color{Cyan}\line(0,1){1.75}}
\put(-0.2,2.9){\color{Cyan}{\scriptsize \em 2}}
\put(0,4.125){\color{Cyan}\line(0,1){1.75}}
\put(-0.2,4.9){\color{Cyan}{\scriptsize \em 2}}
\put(0.1,6.1){\color{Red}\line(3,2){5.775}}
\put(2.8,7.9){\color{Red}{\scriptsize \em 1}}
\put(11.9,6.1){\color{Cyan}\line(-3,2){5.775}}
\put(9,7.9){\color{Cyan}{\scriptsize \em 2}}
\put(12,2.125){\color{Red}\line(0,1){1.75}}
\put(11.8,2.9){\color{Red}{\scriptsize \em 1}}
\put(12,4.125){\color{Red}\line(0,1){1.75}}
\put(11.8,4.9){\color{Red}{\scriptsize \em 1}}
\put(10.1,0.1){\color{Cyan}\line(1,1){1.8}}
\put(10.85,0.9){\color{Cyan}{\scriptsize \em 2}}
\put(9.9,0.1){\color{Red}\line(-1,1){1.8}}
\put(8.85,0.9){\color{Red}{\scriptsize \em 1}}
\put(8,2.125){\color{Cyan}\line(0,1){1.75}}
\put(7.8,2.9){\color{Cyan}{\scriptsize \em 2}}
\put(8,4.125){\color{Cyan}\line(0,1){1.75}}
\put(7.8,4.9){\color{Cyan}{\scriptsize \em 2}}
\put(7.9,6.1){\color{Red}\line(-1,1){1.8}}
\put(6.85,6.9){\color{Red}{\scriptsize \em 1}}
\put(4.1,6.1){\color{Cyan}\line(1,1){1.8}}
\put(4.85,6.9){\color{Cyan}{\scriptsize \em 2}}
\put(4,2.125){\color{Red}\line(0,1){1.75}}
\put(3.8,2.9){\color{Red}{\scriptsize \em 1}}
\put(4,4.125){\color{Red}\line(0,1){1.75}}
\put(3.8,4.9){\color{Red}{\scriptsize \em 1}}
\put(2.1,0.1){\color{Cyan}\line(1,1){1.8}}
\put(2.85,0.9){\color{Cyan}{\scriptsize \em 2}}
\end{picture}

\vspace*{-0.15in}
\setlength{\unitlength}{0.7cm}
\begin{picture}(20,10.5)
\thicklines
\put(2,0){\circle*{0.22}}
\put(0,2){\circle*{0.22}}
\put(0,4){\circle*{0.22}}
\put(0,6){\circle*{0.22}}
\put(10,10){\circle*{0.22}}
\put(20,6){\circle*{0.22}}
\put(20,4){\circle*{0.22}}
\put(20,2){\circle*{0.22}}
\put(18,0){\circle*{0.22}}
\put(16,2){\circle*{0.22}}
\put(16,4){\circle*{0.22}}
\put(16,6){\circle*{0.22}}
\put(14,8){\circle*{0.22}}
\put(12,6){\circle*{0.22}}
\put(12,4){\circle*{0.22}}
\put(12,2){\circle*{0.22}}
\put(10,0){\circle*{0.22}}
\put(8,2){\circle*{0.22}}
\put(8,4){\circle*{0.22}}
\put(8,6){\circle*{0.22}}
\put(6,8){\circle*{0.22}}
\put(4,6){\circle*{0.22}}
\put(4,4){\circle*{0.22}}
\put(4,2){\circle*{0.22}}
\put(1.9,0.1){\color{Red}\line(-1,1){1.8}}
\put(0.85,0.9){\color{Red}{\scriptsize \em 1}}
\put(0,2.125){\color{Cyan}\line(0,1){1.75}}
\put(-0.2,2.9){\color{Cyan}{\scriptsize \em 2}}
\put(0,4.125){\color{Cyan}\line(0,1){1.75}}
\put(-0.2,4.9){\color{Cyan}{\scriptsize \em 2}}
\put(0.075,6.025){\color{Red}\qbezier(0.05,0.025)(4.9,2)(9.75,3.975)}
\put(7.8,9.1){\color{Red}{\scriptsize \em 1}}
\put(19.925,6.025){\color{Cyan}\qbezier(-0.05,0.025)(-4.9,2)(-9.75,3.975)}
\put(12,9.1){\color{Cyan}{\scriptsize \em 2}}
\put(20,2.125){\color{Red}\line(0,1){1.75}}
\put(19.8,2.9){\color{Red}{\scriptsize \em 1}}
\put(20,4.125){\color{Red}\line(0,1){1.75}}
\put(19.8,4.9){\color{Red}{\scriptsize \em 1}}
\put(18.1,0.1){\color{Cyan}\line(1,1){1.8}}
\put(18.85,0.9){\color{Cyan}{\scriptsize \em 2}}
\put(17.9,0.1){\color{Red}\line(-1,1){1.8}}
\put(16.85,0.9){\color{Red}{\scriptsize \em 1}}
\put(16,2.125){\color{Cyan}\line(0,1){1.75}}
\put(15.8,2.9){\color{Cyan}{\scriptsize \em 2}}
\put(16,4.125){\color{Cyan}\line(0,1){1.75}}
\put(15.8,4.9){\color{Cyan}{\scriptsize \em 2}}
\put(15.9,6.1){\color{Red}\line(-1,1){1.8}}
\put(14.85,6.9){\color{Red}{\scriptsize \em 1}}
\put(12.1,6.1){\color{Cyan}\line(1,1){1.8}}
\put(12.85,6.9){\color{Cyan}{\scriptsize \em 2}}
\put(12,2.125){\color{Red}\line(0,1){1.75}}
\put(11.8,2.9){\color{Red}{\scriptsize \em 1}}
\put(12,4.125){\color{Red}\line(0,1){1.75}}
\put(11.8,4.9){\color{Red}{\scriptsize \em 1}}
\put(10.1,0.1){\color{Cyan}\line(1,1){1.8}}
\put(10.85,0.9){\color{Cyan}{\scriptsize \em 2}}
\put(9.9,0.1){\color{Red}\line(-1,1){1.8}}
\put(8.85,0.9){\color{Red}{\scriptsize \em 1}}
\put(8,2.125){\color{Cyan}\line(0,1){1.75}}
\put(7.8,2.9){\color{Cyan}{\scriptsize \em 2}}
\put(8,4.125){\color{Cyan}\line(0,1){1.75}}
\put(7.8,4.9){\color{Cyan}{\scriptsize \em 2}}
\put(7.9,6.1){\color{Red}\line(-1,1){1.8}}
\put(6.85,6.9){\color{Red}{\scriptsize \em 1}}
\put(4.1,6.1){\color{Cyan}\line(1,1){1.8}}
\put(4.85,6.9){\color{Cyan}{\scriptsize \em 2}}
\put(4,2.125){\color{Red}\line(0,1){1.75}}
\put(3.8,2.9){\color{Red}{\scriptsize \em 1}}
\put(4,4.125){\color{Red}\line(0,1){1.75}}
\put(3.8,4.9){\color{Red}{\scriptsize \em 1}}
\put(2.1,0.1){\color{Cyan}\line(1,1){1.8}}
\put(2.85,0.9){\color{Cyan}{\scriptsize \em 2}}
\end{picture}
\end{center}
\end{figure}

{\bf The unique maximal splitting poset.}  
Given a Weyl bialternant $\chi_{_{\lambda}}$, we aim to build a connected splitting poset.  
So take $\lambda \in \Lambda^{+}$, and for each $\mu \in \Pi(\lambda)$, let $U_{\mu}$ be the set of symbols $\{\uelt_{1}^{(\mu)}, \uelt_{2}^{(\mu)}, \ldots, \uelt_{d_{\lambda,\mu}}^{(\mu)}\}$.  
If $\mu \myarrow{i} \nu$ in $\Pi(\lambda)$, construct edges $\uelt_{p}^{(\mu)} \myarrow{i} \uelt_{q}^{(\nu)}$ for all $1 \leq p \leq d_{\lambda,\mu}$ and $1 \leq q \leq d_{\lambda,\nu}$.  
With $\disjointunion_{\!\!\!\mu \in \Pi(\lambda)}U_{\mu}$ as the underlying vertex set, denote by $U(\lambda)$ the resulting edge-colored directed graph.  
For reasons made apparent by the following proposition, we call $U(\lambda)$ the {\em unique maximal splitting poset} for $\chi_{_{\lambda}}$. 

\noindent 
{\bf \UniqueMaximalProp}\ \ {\sl Let $\lambda \in \Lambda^{+}$. 
Then $U(\lambda)$ is the Hasse diagram for a connected splitting poset for $\chi_{_{\lambda}}$.  
Moreover, any splitting poset for $\chi_{_{\lambda}}$ is isomorphic as an edge-colored poset to some edge-colored subgraph of $U(\lambda)$.} 

{\em Proof.}  Connectedness of $U(\lambda)$ follows from connectedness of $\Pi(\lambda)$.  
Clearly $U(\lambda)$ is acyclic (no directed cycles), since whenever we have a path $\uelt_{p_{0}}^{\mu_{0}} \myarrow{i_{1}} \uelt_{p_{1}}^{\mu_{1}} \myarrow{i_{2}} \cdots \myarrow{i_{k}} \uelt_{p_{k}}^{\mu_{k}}$ in $U(\lambda)$ it is necessarily the case that $\mu_{0} + \alpha_{i_{1}} + \cdots + \alpha_{i_{k}} = \mu_{k}$.  
So, the 
transitive closure of the edge relations for $U(\lambda)$ uniquely 
determines a partial ordering of the elements of $U(\lambda)$ with 
respect to which the covering relations are precisely the defining 
edge relations $\uelt_{p}^{(\mu)} \myarrow{i} \uelt_{q}^{(\nu)}$.  That 
$U(\lambda)$ is ranked follows from the fact that $\Pi(\lambda)$ is 
ranked, cf.\ \WeightRemarkTwo.  
It also follows that for any $\uelt_{p}^{\mu} \in U(\lambda)$, 
$\rho_{i}(\uelt_{p}^{\mu})$ is the same as the 
rank of $\mu$ in its $i$-component in $\Pi(\lambda)$ and that 
$l_{i}(\uelt_{p}^{\mu})$ is the length of the $i$-component for $\mu$ 
in $\Pi(\lambda)$.  Therefore,  
$m_{i}(\uelt_{p}^{\mu}) = \langle \mu,\alpha_{i}^{\vee} \rangle$, so 
$wt(\uelt_{p}^{\mu}) = \mu$.  So, $U(\lambda)$ is 
${\Phi}$-structured.  Since 
$|\{\uelt_{1}^{(\mu)}, \uelt_{2}^{(\mu)}, \ldots, 
\uelt_{d_{\lambda,\mu}}^{(\mu)}\}| = 
d_{\lambda,\mu}$, then $\WGF(U(\lambda)) = \chi_{_{\lambda}}$ and $U(\lambda)$ is a splitting poset for $\chi_{_{\lambda}}$. 
Now say $R$ is any splitting 
poset for $\chi_{_{\lambda}}$.  For each $\mu \in \Pi(\lambda)$, let 
$R_{\mu} := \{\selt\, |\, wt(\selt) = \mu\}$; then $|R_{\mu}| = 
d_{\lambda,\mu}$. Let $\phi: R \longrightarrow U(\lambda)$ be any 
function such that for each $\mu \in \Pi(\lambda)$, $\phi|_{R_{\mu}}$ 
puts the sets 
$R_{\mu}$ and $\{\uelt_{1}^{\mu}, \ldots, 
\uelt_{d_{\lambda,\mu}}^{\mu}\}$ 
in one-to-one correspondence. Clearly if $\selt \myarrow{i} \telt$ 
in $R$, then $\phi(\selt) \myarrow{i} \phi(\telt)$ in $U(\lambda)$.  
So the image $\phi(R)$ is an edge-colored subgraph of $U(\lambda)$ 
for which $R \cong \phi(R)$ as edge-colored posets.\hfill\QED 

{\bf Minuscule and quasi-minuscule splitting posets.} The special splitting posets described next will be used later on as building blocks for constructing other splitting posets. 
Our discussion of these objects is patterned after \cite{StemDom} and \cite{Stem}. 
The minuscule splitting posets are well-studied (see e.g.\ \cite{Green}, \cite{PrEur}) and are also called minuscule lattices. 
These are all diamond-colored distributive lattices whose vertex-colored posets of join irreducibles are uniformly characterized in \cite{StrayerThesis} and \cite{StrayerArxiv} as part of a more general program that also accounts for their Kac--Moody generalizations.
Throughout this subsection, we assume that $\Phi$ is irreducible, a hypothesis we reiterate for emphasis in lemma and proposition statements. 

\noindent 
{\bf \MinQuasiMinLemma}\ \ {\sl Assume that $\Phi$ is irreducible. 
The connected components of the partially ordered set $\Lambda$ are precisely the cosets of the subgroup $\mathbb{Z}\Phi$ in $\Lambda$. 
The number of such connected components is finite. 
Let $\mathcal{D}$ be the set of dominant weights within any given connected component of $\Lambda$, and give $\mathcal{D}$ the induced partial ordering from $\Lambda$.  
Then $\mathcal{D}$ has a unique minimal element. 
Moreover, if the zero weight is in $\mathcal{D}$, then the highest short root is the unique minimal element of $\mathcal{D} \setminus \{0\}$.} 

{\em Proof.} It follows from the definitions that the cosets $\Lambda / \mathbb{Z}\Phi$ are the connected components of $\Lambda$. 
The relations matrix for the finitely-generated $\mathbb{Z}$-module $\Lambda / \mathbb{Z}\Phi$ is just the Cartan matrix $M_{\Phi}$. 
Its Smith normal form has no zeros on the diagonal since $M_{\Phi}$ is invertible. 
It follows from the structure theorem for finitely-generated $\mathbb{Z}$-modules that $\Lambda / \mathbb{Z}\Phi$ is a finite group, and therefore $\Lambda$ has only finitely many connected components. 
That $\mathcal{D}$ has a unique minimal element follows from Corollary 1.4 of \cite{StemDom}. 
Now say $0 \in \mathcal{D}$. 
That $\mathcal{D} \setminus \{0\}$ has the highest short root of $\Phi$ as its unique minimal element follows from Proposition 2.1 of \cite{StemDom}.\hfill\QED 

The pairs $(\Phi,\lambda)$ for which $\Phi$ is irreducible and $\lambda$ is nonzero and minimal have the following well-known classification: 
$(\myA_{n},\omega_{k})$ ($1 \leq k \leq n$);  
$(\myB_{n},\omega_{n})$;  
$(\myC_{n},\omega_{1})$; 
$(\myD_{n},\omega_{1})$; $(\myD_{n},\omega_{n-1})$;$(\myD_{n},\omega_{n})$; 
$(\myE_{6},\omega_{1})$; $(\myE_{6},\omega_{6})$; 
and $(\myE_{7},\omega_{7})$. 
The pairs $(\Phi,\lambda)$ for which $\Phi$ is irreducible and $\lambda$ is the highest short root are: 
$(\myA_{n},\omega_{1}+\omega_{n})$;  
$(\myB_{n},\omega_{1})$;  
$(\myC_{n},\omega_{2})$; 
$(\myD_{n},\omega_{2})$; 
$(\myE_{6},\omega_{2})$; 
$(\myE_{7},\omega_{1})$; 
$(\myE_{8},\omega_{8})$; 
$(\myF_{4},\omega_{4})$; 
and $(\myG_{2},\omega_{1})$. 

\noindent 
{\bf \MinQuasiMinProposition}\ \ {\sl Assume that $\Phi$ is irreducible. 
Let $\lambda \in \Lambda^{+}$ with $\lambda \not= 0$.  Then the following (M1)--(M4) are equivalent:\\ 
\hspace*{0.15in}(M1) $\{\langle \lambda,\alpha^{\vee} \rangle\}_{\alpha \in \Phi} \subseteq \{0,\pm 1\}$.\\
\hspace*{0.15in}(M2) $\Pi(\lambda) = W\lambda$.\\
\hspace*{0.15in}(M3) For all $\nu \in \Lambda^{+}$, we have $\nu \leq \lambda$ if and only if $\nu = \lambda$.\\
\hspace*{0.15in}(M4) $\lambda$ is the unique minimal element amongst the dominant weights in some connected component of $\Lambda$, cf.\ \MinQuasiMinLemma.

Moreover, the following (QM1)--(QM4) are equivalent:\\
\hspace*{0.15in}(QM1) $\{\langle \lambda,\alpha^{\vee} \rangle\}_{\alpha \in \Phi} \subseteq \{0,\pm 1,\pm 2\}$ and there is a unique $\beta \in \Phi$ such that $\langle \lambda,\beta^{\vee} \rangle = 2$.\\
\hspace*{0.15in}(QM2) $\Pi(\lambda) = W\lambda \bigdisjointunion \{0\}$.\\
\hspace*{0.15in}(QM3) For all $\nu \in \Lambda^{+}$, we have $\nu \leq \lambda$ if and only if $\nu = 0$ or $\nu = \lambda$.\\
\hspace*{0.15in}(QM4) $\lambda$ is the highest short root, which by \MinQuasiMinLemma\ means that within the connected component of $\Lambda$ that contains the zero weight, $\lambda$ is the unique minimal element amongst the nonzero dominant weights.}

{\em Proof.} For (M1) $\Longrightarrow$ (M2), we take advantage of our knowledge of $\Pi(\lambda)$ as a ranked poset.  
In particular, it is ${\Phi}$-structured and has $\lambda$ as its unique maximal element.  
We will induct on the depth $\mathrm{depth}(\mu)$ of weights $\mu \in \Pi(\lambda)$.  If $\mathrm{depth}(\mu) = 0$, then $\mu = \lambda$, which is obviously in $W\lambda$. 
Our induction hypothesis is that for some nonnegative integer $k$ and any weight $\mu \in \Pi(\lambda)$ with $\mathrm{depth}(\mu) \leq k$, then $\mu \in W\lambda$. 
Suppose now that $\mathrm{depth}(\mu) = k+1$ for some $\mu \in \Pi(\lambda)$. 
Since $\mu \not= \lambda$, then by \WeightRemarkOne.5, $\mu \myarrow{i} \nu$ for some $i \in I$ and $\nu \in \Pi(\lambda)$. 
In fact, we have 
\[\mu' \myarrow{i} \cdots \myarrow{i} \mu \myarrow{i} \nu \myarrow{i} \cdots \myarrow{i} \nu',\]
where $\mu'$ and $\nu'$ are, respectively, the minimal and maximal elements of the $i$-component of $\mu$. 
Then $\langle \nu',\alpha_{i}^{\vee} \rangle > 0$.  
But $\nu' \in W\lambda$ by the induction hypothesis, so we can write $\nu' = w(\lambda)$ for some $w \in W$. 
Then $\langle \nu',\alpha_{i}^{\vee} \rangle = \langle \lambda,w^{-1}(\alpha_{i}^{\vee}) \rangle \in \{0,\pm 1\}$. 
Since the quantity in question is positive, it follows that $\langle \nu',\alpha_{i}^{\vee} \rangle = 1$. 
Then $\nu' = \nu$ and $\mu' = \mu$. 
In particular, $\mu = s_{i}(\nu) = (s_{i}w)(\lambda) \in W\lambda$. 

(M2) $\Longrightarrow$ (M3) is a consequence of our definition of $\Pi(\lambda)$. 
For (M3) $\Longrightarrow$ (M4), it follows directly from (M3) that $\lambda$ is minimal amongst the dominant weights within its connected component of $\Lambda$, and applying \MinQuasiMinLemma\ we conclude that $\lambda$ is uniquely minimal. 
(M3) $\Longrightarrow$ (M4) follows from Proposition 1.12 of \cite{StemDom}. 

For (QM1) $\Longrightarrow$ (QM2), we can apply Lemma 4.6 of \cite{Stem} to conclude that $\lambda$ is the highest short root.  
It follows that $W\lambda$ is the set of short roots, cf.\ Lemma 10.4.C of \cite{Hum}. 
Then by Proposition 2.1/Remark 2.2 of \cite{StemDom}, it follows that $\Pi(\lambda) = W\lambda \bigdisjointunion \{0\}$. 
That (QM2) $\Longrightarrow$ (QM3) is a consequence of our definition of $\Pi(\lambda)$. 
Assuming (QM3), then we can apply \MinQuasiMinLemma\ to get (QM4).  
Assuming (QM4), then (QM1) follows from Lemma 2.3 of \cite{StemDom}. 
\hfill\QED

The special properties of the weights described in the preceding proposition facilitate the construction of splitting posets for the companion Weyl bialternants. 
The following language follows \cite{Stem}.  
Let $\mu \in \Lambda$. 
Say $\mu$ is {\em minuscule} if $\{\langle \mu,\alpha^{\vee} \rangle\}_{\alpha \in \Phi} = \{0,\pm 1\}$. 
It is easy to see that $\mu$ is minuscule if and only if $\mu \in W\lambda$, where $\lambda$ is a nonzero dominant weight meeting any one of the equivalent conditions (M1)--(M4) of \MinQuasiMinProposition. 
We say $\mu$ is {\em quasi-minuscule} if $\{\langle \mu,\alpha^{\vee} \rangle\}_{\alpha \in \Phi} = \{0,\pm 1,\pm 2\}$ and there is a unique $\beta \in \Phi$ such that $\langle \mu,\beta^{\vee} \rangle = 2$. 
It is easy to see that $\mu$ is quasi-minuscule if and only if $\mu \in W\lambda$, where $\lambda$ is a nonzero dominant weight meeting any one of the equivalent conditions (QM1)--(QM4) of \MinQuasiMinProposition. 

Suppose $\lambda$ is dominant and minuscule. 
Let $R(\lambda)$ be the ${\Phi}$-structured poset $\Pi(\lambda) = W\lambda$. 
Now $d_{\lambda,\lambda} = 1$ by \WeylCharLemma. 
For $\mu \in \Pi(\lambda) = W\lambda$ we have $\mu = w(\lambda)$ for some $w \in W$. 
So by the $W$-invariance of $\chi_{_{\lambda}}$, we have $d_{\lambda,\mu} = d_{\lambda,w(\lambda)} = d_{\lambda,\lambda} = 1$. 
So the unique maximal splitting poset for $\chi_{_{\lambda}}$ coincides with $\Pi(\lambda)$ and is therefore edge-minimal. 
That is, $R(\lambda)$ is the unique splitting poset for $\chi_{_{\lambda}}$.  
This fact is recorded in \QuasiMinLemma\ below. 
We call $R(\lambda)$ the {\em minuscule splitting poset} for the dominant minuscule weight $\lambda$. 

Now suppose $\lambda$ is dominant and quasi-minuscule, so $\lambda$ is the highest short root.  
All of the short roots in $\Phi$ are $W$-conjugate, cf.\ Lemma 10.4.C of \cite{Hum}. 
In particular, all of the short simple roots are in $W\lambda$. 
We let $R(\lambda)$ be the set $W\lambda \bigdisjointunion \{\overline{\alpha_i}\, |\, \alpha_{i} \mbox{ is short and simple}\}$, where $\{\overline{\alpha_i}\}$ is regarded as a set of symbols. 
Place colored directed edges between elements of $R(\lambda)$ according to the following rule: $\alpha \myarrow{i} \beta$ if and only if (1) $\alpha$ and $\beta$ are short roots and $\alpha \myarrow{i} \beta$ in $\Pi(\lambda)$, (2) $\alpha = -\alpha_{i}$ and $\beta = \overline{\alpha_i}$, or (3) $\alpha = \overline{\alpha_i}$ and $\beta = \alpha_i$. 
Observe that $R(\lambda)$ is the Hasse diagram for a ranked poset and that for all $\xelt \in R(\lambda)$ we have $wt(\xelt) = \alpha$ if $\xelt = \alpha$ is a short root in $\Phi$, and we have $wt(\xelt) = 0$ if $\xelt = \overline{\alpha_i}$ for some simple root $\alpha_i$. 
It is now easy to see that $R(\lambda)$ is ${\Phi}$-structured with $\lambda$ as its unique maximal element. 
Since its $i$-components are chains, it follows from \WInvariantLemma\ that $\WGF(R(\lambda))$ is $W$-invariant. 
Let $\mathcal{S}_{I,0}(R(\lambda)) = \{\lambda\}$. 
For all $\xelt \in R(\lambda) \setminus \{\lambda\}$, define 
\[\tau(\xelt) = \left\{\begin{array}{ll}
\overline{\alpha_i} & \hspace*{0.25in} \mbox{if } \xelt = -\alpha_i\\
-\alpha_i & \hspace*{0.25in} \mbox{if } \xelt = \overline{\alpha_i}\\
\xelt & \hspace*{0.25in} \mbox{otherwise}
\end{array}\right.
\hspace*{0.25in}\mbox{and}\hspace*{0.25in}
\kappa(\xelt) = \left\{\begin{array}{ll}
i & \hspace*{0.25in} \mbox{if } \xelt = -\alpha_i\\
i & \hspace*{0.25in} \mbox{if } \xelt = \overline{\alpha_i}\\
k & \hspace*{0.25in} \mbox{\parbox{2in}{\footnotesize otherwise, where we choose $k$ to be any color so that $\xelt \myarrow{k} \yelt$ is an edge in $R(\lambda)$}}
\end{array}\right.
\]
It is easy now to check that $R(\lambda)$ meets the requirements of \InitialSplittingTheorem, so $\WGF(R(\lambda)) = \chi_{_{\lambda}}$. 
We call $R(\lambda)$ the {\em quasi-minuscule splitting poset} for the dominant quasi-minuscule weight $\lambda$. 
In the next result we argue that $R(\lambda)$ enjoys a certain uniqueness property.  

\noindent 
{\bf \QuasiMinLemma}\ \ {\sl Keep the notation of the preceding paragraphs, and assume that $\Phi$ is irreducible. 
If $\lambda$ is dominant and minuscule, then the minuscule splitting poset $R(\lambda)$ is the unique splitting poset for $\chi_{_{\lambda}}$, and it is connected and edge-minimal. 
If $\lambda$ is dominant and quasi-minuscule, then the quasi-minuscule splitting poset $R(\lambda)$ is the only connected edge-minimal splitting poset for $\chi_{_{\lambda}}$.}

{\em Proof.} For the minuscule case, the claims follow from the second paragraph above the proposition statement together with the observation that $R(\lambda) = \Pi(\lambda)$ is connected. 
Now consider the case that $\lambda$ is dominant and quasi-minuscule. 
Since the $i$-components of $R(\lambda)$ are chains, then $R(\lambda)$ is an edge-minimal splitting poset for $\chi_{_{\lambda}}$. 
Connectedness of $R(\lambda)$ is easily deduced from the connectedness of $\Pi(\lambda)$. 
We note that since $\WGF(R(\lambda)) = \chi_{_{\lambda}}$, then $d_{\lambda,0}$ just counts the number of short simple roots. 

Now suppose $R$ is any connected and edge-minimal splitting poset for $\chi_{_{\lambda}}$. 
We use the following language in our argument that $R \cong R(\lambda)$ as edge-colored directed graphs. 
An edge $\xelt \myarrow{i} \yelt$ in $R$ is a color $i$ edge ``above'' $\xelt$ and ``below'' $\yelt$. 
We call any $\zelt \in R$ with $wt(\zelt) = 0$ a ``middle rank'' element of $R$. 
The non-middle-rank elements of $R$ can be identified with the short roots. 
Now suppose $\alpha_{i} \in R$ is short and simple. 
There are no color $i$ edges above $\alpha_{i}$ since $2\alpha_{i}$ is not a root. 
Also, any color $i$ edge that is above a middle rank element must also be below $\alpha_{i}$. 
From these two facts it follows that $\alpha_{i}$ is the unique highest rank element in its $i$-component. 
Similarly, $-\alpha_{i}$ is the unique lowest rank element in this $i$-component. 
Further, since $\langle \alpha_{i},\alpha_{i}^{\vee} \rangle = 2$, it follows then that there is some middle rank element $\zelt_{i}$ such that $-\alpha_{i} \myarrow{i} \zelt_{i} \myarrow{i} \alpha_{i}$. 
We could remove any other color $i$ edges incident with $\alpha_{i}$ or $-\alpha_{i}$ without violating the splitting property. 
Since $R$ is edge-minimal, it follows that the only color $i$ edges incident with $\alpha_{i}$ and $-\alpha_{i}$ are those in the chain $-\alpha_{i} \myarrow{i} \zelt_{i} \myarrow{i} \alpha_{i}$, which is therefore the $i$-component of $\alpha_{i}$. 
Take another short simple root $\alpha_{j}$ with $i \not= j$. 
If $\zelt_{j} = \zelt_{i}$, then there would be some middle rank element $\zelt$ that is not connected to any short simple roots or their negatives, i.e.\ $R$ would be disconnected. 
However, $R$ is connected, so $\zelt_{i}$ and $\zelt_{j}$ must be distinct, i.e.\ $|\{\zelt_{i}\, |\, \alpha_{i} \mbox{ is short and simple}\}| = d_{\lambda,0}$. 
So each middle rank element is part of exactly one nontrivial color $i$ component, which is a length two chain. 
It is evident now that the mapping $R \longrightarrow R(\lambda)$ which identifies short roots with themselves and sends $\zelt_{i} \mapsto \overline{\alpha_{i}}$ is an isomorphism of edge-colored directed graphs.\hfill\QED 

We close this subsection by offering a purely combinatorial characterization of minuscule and quasi-minuscule splitting posets whose proof only uses general principles and, in particular, does not depend on the classification of minuscule and quasi-minuscule dominant weights by root system type. 
In \S \GaussianSection, we consider other combinatorial-structural aspects of minuscule splitting posets, and in particular the fact that they are diamond-colored distributive lattices, cf.\ Theorem 12.6 of \cite{DonDistributive}. 

\noindent 
{\bf \MinQuasiMinTheorem}\ \ {\sl Let $\Phi$ be irreducible, and let $R$ be a ranked poset with edges colored by the index set $I$ for our choice of simple roots. 
\underline{Minuscule case:} $R$ is (isomorphic to) a minuscule splitting poset if and only if it is connected, fibrous, and ${\Phi}$-structured; has at least one edge; and has the following properties: 
(1) $l_i(\xelt) \leq 1$ for all $i \in I$ and $\xelt \in R$, and 
(2) whenever $\langle \alpha_{i},\alpha_{j}^{\vee} \rangle = 0 = \langle \alpha_{j},\alpha_{i}^{\vee} \rangle$ then any $\{i,j\}$-component of $R$ has a unique maximal element.  
In this case, $R$ has a unique maximal element $\melt$, $\lambda := wt(\melt)$ is dominant and minuscule, and $R \cong R(\lambda)$.  
\underline{Quasi-minuscule case:} $R$ is (isomorphic to) a quasi-minuscule splitting poset if and only if it is connected, fibrous, and ${\Phi}$-structured; has an $i$-component of length at least two for some $i \in I$; and has the following properties: (1) If $l_i(\xelt) \geq 2$ for some $i \in I$ and $\xelt \in R$, then $\comp_i(\xelt)$ is a chain $\xelt_0 \myarrow{i} \xelt_1 \myarrow{i} \xelt_2$ such that $\delta_j(\xelt_0) = 0 = \rho_j(\xelt_2)$ and $\comp_j(\xelt_1) = \{\xelt_1\}$ for all $j \not= i$; (2) is the same as condition (2) in the minuscule case; and (3) whenever $\langle \alpha_{i},\alpha_{j}^{\vee} \rangle = -1 = \langle \alpha_{j},\alpha_{i}^{\vee} \rangle$ then any $\{i,j\}$-component of $R$ has a unique maximal element. 
In this case, $R$ has a unique maximal element $\melt$, $\lambda := wt(\melt)$ is dominant and quasi-minuscule, and $R \cong R(\lambda)$.} 

{\em Proof.} 
First, we demonstrate the necessity of the stated combinatorial conditions.  
So, take $R = R(\lambda)$ and assume $\lambda$ is dominant and minuscule. 
That $R = \Pi(\lambda)$ is connected, fibrous, and ${\Phi}$-structured follows from properties of $\Pi(\lambda)$, cf.\ \PiResults. 
Since $\lambda \not= 0$, then $\langle \lambda,\alpha_{i}^{\vee} \rangle > 0$ for some $i \in I$, hence $\nu \myarrow{i} \lambda$ in $\Pi(\lambda)$, so $R = \Pi(\lambda)$ has at least one edge. 
Since for any $\mu \in \Pi(\lambda) = W\lambda$ we can write $\mu = \sigma(\lambda)$ for some $\sigma \in W$, then $\langle \mu,\alpha_{i}^{\vee} \rangle = \langle \lambda,\sigma^{-1}(\alpha_{i}^{\vee}) \rangle \in \{0,\pm{1}\}$. 
So, $l_{i}(\mu) \leq 1$, establishing property (1) in the minuscule case. 
Now pick any $\mu \in R = \Pi(\lambda)$ and consider $\comp_{\{i,j\}}(\mu)$, where $i,j \in I$ with $\langle \alpha_{i},\alpha_{j}^{\vee} \rangle = 0 = \langle \alpha_{j},\alpha_{i}^{\vee} \rangle$ and necessarily $i \not= j$. 
We aim to show that $\comp_{\{i,j\}}(\mu)$ has a unique maximal element. 
We do so by describing all possible configurations for $\comp_{\{i,j\}}(\mu)$.  
Of course, if $\comp_{\{i,j\}}(\mu)$ consists only of the vertex $\mu$, or of a single edge having $\mu$ as one of its vertices, then we are done. 
So now suppose that $\mu$ is incident with edges of colors $i$ and $j$. 
We must have one of the following cases, where $\nu$ and $\pi$ are some elements of the component: 
(a) $\nu \myarrow{i} \mu \mybackarrow{j} \pi$, 
(b) $\nu \myarrow{i} \mu \myarrow{j} \pi$, 
(c) $\nu \myarrow{j} \mu \myarrow{i} \pi$, or
(d) $\nu \mybackarrow{i} \mu \myarrow{j} \pi$. 
In case (a), we have $\pi = \nu + \alpha_{i} - \alpha_{j}$, hence $\langle \pi,\alpha_{i}^{\vee} \rangle = \langle \nu,\alpha_{i}^{\vee} \rangle + \langle \alpha_{i},\alpha_{i}^{\vee} \rangle + \langle \alpha_{j},\alpha_{i}^{\vee} \rangle = 1$. 
Hence, $\nu \mybackarrow{j} (\pi - \alpha_{i}) \myarrow{i} \pi$ in $R = \Pi(\lambda)$. 
Since each one-color component in $R = \Pi(\lambda)$ has at most one edge, it follows that the the four edges identified so far in the component $\comp_{\{i,j\}}(\mu)$ must be all the edges.  
That is, $\comp_{\{i,j\}}(\mu)$ has the form 
\parbox{1.1cm}{\begin{center}
\setlength{\unitlength}{0.2cm}
\begin{picture}(4,3)
\put(2,0){\circle*{0.5}} \put(0,2){\circle*{0.5}}
\put(2,4){\circle*{0.5}} \put(4,2){\circle*{0.5}}
\put(0,2){\line(1,1){2}} \put(2,0){\line(-1,1){2}}
\put(4,2){\line(-1,1){2}} \put(2,0){\line(1,1){2}}
\put(0.75,0.55){\em \small j} \put(2.7,0.7){\em \small i}
\put(0.7,2.7){\em \small i} \put(2.75,2.55){\em \small j}
\end{picture} \end{center}}. 
Clearly this component has a unique maximal element. 
Similar analyses in cases (b), (c), and (d) show that $\comp_{\{i,j\}}(\mu)$ has the form
\parbox{1.1cm}{\begin{center}
\setlength{\unitlength}{0.2cm}
\begin{picture}(4,3)
\put(2,0){\circle*{0.5}} \put(0,2){\circle*{0.5}}
\put(2,4){\circle*{0.5}} \put(4,2){\circle*{0.5}}
\put(0,2){\line(1,1){2}} \put(2,0){\line(-1,1){2}}
\put(4,2){\line(-1,1){2}} \put(2,0){\line(1,1){2}}
\put(0.75,0.55){\em \small j} \put(2.7,0.7){\em \small i}
\put(0.7,2.7){\em \small i} \put(2.75,2.55){\em \small j}
\end{picture} \end{center}}. 
This completes the proof that the stated combinatorial conditions are necessary when $R = R(\lambda)$ for a dominant and minuscule $\lambda$. 
 
Now assume $R = R(\lambda)$ with $\lambda$ dominant and quasi-minuscule. 
It follows straightforwardly from the definition of $R(\lambda)$ that $R$ is connected, fibrous, and ${\Phi}$-structured and that there is an edge below $\lambda$ in $R$. 
Now suppose we have $\xelt \in R$ with $\delta_{i}(\xelt) = 0$ and $l_{i}(\xelt) \geq 2$. 
Then $\xelt \not= \overline{\alpha_{j}}$ for all $j$, since $l_{i}(\overline{\alpha_{j}}) = 0$ when $j \not= i$ and $l_{i}(\overline{\alpha_{i}}) = 1$. 
So $\xelt = \mu$ for some $\mu \in W\lambda$. 
Therefore $\langle \mu,\alpha_{i}^{\vee} \rangle = l_{i}(\mu) \geq 2$. 
Now since $\mu \in W\lambda$, then $\langle \mu,\alpha^{\vee} \rangle \in \{0, \pm{1}, \pm{2}\}$ for all $\alpha \in \Phi$ and there is a unique $\beta \in \Phi$ for which $\langle \mu,\beta^{\vee} \rangle = 2$. 
It follows that $\langle \mu,\alpha_{i}^{\vee} \rangle = 2$. 
Writing $\mu = \sigma(\lambda)$, we get $2 = \langle \mu,\alpha_{i}^{\vee} \rangle = \langle \lambda,\sigma^{-1}(\alpha_{i}^{\vee}) \rangle$. 
Since there is only one coroot $\gamma^{\vee}$ such that $\langle \lambda,\gamma^{\vee} \rangle = 2$, namely $\gamma^{\vee} = \lambda^{\vee}$, then $\sigma^{-1}(\alpha_{i}^{\vee}) = \lambda^{\vee}$. 
In particular, it follows that $\alpha_{i}$ is in the $W$-orbit of $\lambda$. 
Then $\alpha_{i}$ has the same length as $\mu$, hence $\langle \mu,\alpha_{i}^{\vee} \rangle = \langle \alpha_{i},\mu^{\vee} \rangle = 2$. 
Therefore $\mu^{\vee} = \alpha_{i}^{\vee}$, so $\mu = \alpha_{i}$. 
Then the $i$-component of $\xelt = \mu$ is $-\alpha_{i} \myarrow{i} \overline{\alpha_{i}} \myarrow{i} \alpha_{i}$. 
Since $\alpha_{i} - \alpha_{j}$ is not a root for $j \not= i$, then $\rho_{j}(\alpha_{i}) = 0$. 
Similarly, $\delta_{j}(-\alpha_{i}) = 0$ when $j \not= i$. 
Then we have established property (1) in the quasi-minuscule case. 

It remains to be checked that any $\{i,j\}$-component of $R = R(\lambda)$ has a unique maximal element whenever $\langle \alpha_{i},\alpha_{j}^{\vee} \rangle = 0 = \langle \alpha_{j},\alpha_{i}^{\vee} \rangle$ or when $\langle \alpha_{i},\alpha_{j}^{\vee} \rangle = -1 = \langle \alpha_{j},\alpha_{i}^{\vee} \rangle$. 
First take $\langle \alpha_{i},\alpha_{j}^{\vee} \rangle = 0 = \langle \alpha_{j},\alpha_{i}^{\vee} \rangle$, and suppose $\mathcal{C}$ is an $\{i,j\}$-component of $R$. 
If $\mathcal{C}$ has no $i$- or $j$-components of length more than one, then every element of $\mathcal{C}$ is in $\Pi(\lambda)$, and the analysis of the first paragraph of the proof can be used to conclude that $\mathcal{C}$ has a unique maximal element. 
So now suppose that $\mathcal{C}$ has an $i$- or $j$-component of length at least two. 
Without loss of generality, we assume that this length-at-least-two component has color $i$. 
By the definition of $R(\lambda)$, this $i$-component has length exactly two and is the chain $-\alpha_{i} \myarrow{i} \overline{\alpha_i} \myarrow{i} \alpha_{i}$.  
Since $\langle \alpha_{i},\alpha_{j}^{\vee} \rangle = 0 = \langle \alpha_{j},\alpha_{i}^{\vee} \rangle$, then none of the vectors $\pm \alpha_{i} \pm \alpha_{j}$ is a root and therefore none of them is in $\Pi(\lambda)$.  
We conclude that $-\alpha_{i}$, $\overline{\alpha_i}$, and $\alpha_{i}$ have no incident color $j$ edges in $R$.  
Therefore the given component $\mathcal{C}$ consists precisely of the length two color $i$ chain $-\alpha_{i} \myarrow{i} \overline{\alpha_i} \myarrow{i} \alpha_{i}$, which obviously has a unique maximal element. 

Next take $\langle \alpha_{i},\alpha_{j}^{\vee} \rangle = -1 = \langle \alpha_{j},\alpha_{i}^{\vee} \rangle$, and let $\mathcal{C}$ be an $\{i,j\}$-component of $R$. 
Suppose $\mathcal{C}$ has no $i$- or $j$-components of length more than one, so $\mathcal{C}$ is contained (setwise) entirely within the $\Pi(\lambda)$ part of $R$. 
If $\mathcal{C}$ consists only of a single element or only of a single edge, then clearly $\mathcal{C}$ has a unique maximal element. 
Suppose now $\mathcal{C}$ has multiple edges, and suppose $\mu \in \mathcal{C}$ is incident with at least two edges. 
Then $\mu$ is incident with exactly two edges in $\mathcal{C}$, which are necessarily of different colors for two reasons: in the $\Pi(\lambda)$ part of $R$ we cannot have 
\parbox{1.1cm}{\begin{center}
\setlength{\unitlength}{0.2cm}
\begin{picture}(4,1.25)
\put(2,0){\circle*{0.5}} 
\put(0,2){\circle*{0.5}}
\put(4,2){\circle*{0.5}}
\put(2,0){\line(-1,1){2}}
\put(2,0){\line(1,1){2}}
\put(0.75,0.55){\em \small k} 
\put(2.7,0.7){\em \small k}
\end{picture} \end{center}} 
or 
\parbox{1.1cm}{\begin{center}
\setlength{\unitlength}{0.2cm}
\begin{picture}(4,1.25)
\put(0,0){\circle*{0.5}}
\put(2,2){\circle*{0.5}} 
\put(4,0){\circle*{0.5}}
\put(0,0){\line(1,1){2}} 
\put(4,0){\line(-1,1){2}} 
\put(0.7,0.7){\em \small k} 
\put(2.75,0.55){\em \small k}
\end{picture} \end{center}} 
and moreover the length of any one-color component is at most one.   
The incident edges with $\mu$ form one of the four configurations (a), (b), (c), or (d) from the first paragraph of the proof. 
In case (a), one can see that $1 = \langle \mu,\alpha_{i}^{\vee} \rangle = \langle \pi+\alpha_{j},\alpha_{i}^{\vee} \rangle = \langle \pi,\alpha_{i}^{\vee} \rangle - 1$, whence $\langle \pi,\alpha_{i}^{\vee} \rangle = 2$.  
That is, $m_{i}(\pi) = 2$, so $l_{i}(\pi) \geq 2$. 
But this violates the assumption that no $i$- or $j$-component of $\mathcal{C}$ has length more than one. 
So we rule out case (a), and similarly we rule out case (d). 
Now consider case (b). 
Here, $\langle \pi,\alpha_{i}^{\vee} \rangle = \langle \mu+\alpha_{j},\alpha_{i}^{\vee} \rangle = 1 - 1 = 0$, hence $\pi$ has no incident color $i$ edges. 
Similarly, $\langle \nu,\alpha_{j}^{\vee} \rangle = 0$, and so $\nu$ has no incident color $j$ edges.  
That is, in case (b), $\mathcal{C}$ consists only of the elements $\{\mu,\nu,\pi\}$. 
Therefore $\mathcal{C}$ has a unique maximal element. 
In case (c), a similar argument also shows that (setwise) $\mathcal{C} = \{\mu,\nu,\pi\}$, so $\mathcal{C}$ has a unique maximal element. 

It remains to show that if $\langle \alpha_{i},\alpha_{j}^{\vee} \rangle = -1 = \langle \alpha_{j},\alpha_{i}^{\vee} \rangle$ and if $\mathcal{C}$ is an $\{i,j\}$-component of $R$ with an $i$- or $j$-component of length more than one, then $\mathcal{C}$ has a unique maximal element. 
A $k$-component of length more than one will, by the definition of quasi-minuscule splitting poset, have length exactly two.  
So without loss of generality, we assume that the length of an $i$-component in $\mathcal{C}$ is two. 
This component has the form $-\alpha_{i} \myarrow{i} \overline{\alpha_i} \myarrow{i} \alpha_{i}$. 
Using $\langle \alpha_{i},\alpha_{j}^{\vee} \rangle = -1 = \langle \alpha_{j},\alpha_{i}^{\vee} \rangle$, we see that $\alpha_{i}+\alpha_{j}$ and $-\alpha_{i}-\alpha_{j}$ (and therefore $\alpha_{j}$, $\overline{\alpha_j}$, and $-\alpha_{j}$) are in $\mathcal{C}$. 
Now $2\alpha_{i}+\alpha_{j}$, $\alpha_{i}+2\alpha_{j}$, $-2\alpha_{i}-\alpha_{j}$, and $-\alpha_{i}-2\alpha_{j}$ are not roots, so setwise $\mathcal{C}$ is $\{-\alpha_{i}-\alpha_{j}, -\alpha_{i}, \overline{\alpha_i}, \alpha_{i}, -\alpha_{j}, \overline{\alpha_j}, \alpha_{j}, \alpha_{i}+\alpha_{j}\}$.  Then 
\parbox{1.1cm}{\begin{center}
\setlength{\unitlength}{0.2cm}
\begin{picture}(4,10)
\put(2,0){\circle*{0.5}} 
\put(0,2){\circle*{0.5}}
\put(0,5){\circle*{0.5}} 
\put(0,8){\circle*{0.5}} 
\put(4,2){\circle*{0.5}}
\put(4,5){\circle*{0.5}} 
\put(4,8){\circle*{0.5}} 
\put(2,10){\circle*{0.5}} 
\put(2,0){\line(1,1){2}}
\put(2,0){\line(-1,1){2}}
\put(4,2){\line(0,1){6}} 
\put(0,2){\line(0,1){6}} 
\put(0,8){\line(1,1){2}}
\put(4,8){\line(-1,1){2}}
\put(0.75,0.55){\em \small j} 
\put(2.7,0.6){\em \small i}
\put(0.2,2.75){\em \small i} 
\put(3,2.75){\em \small j}
\put(0.2,6.25){\em \small i} 
\put(3,6.25){\em \small j}
\put(0.75,8.7){\em \small j} 
\put(2.7,8.65){\em \small i}
\end{picture} \end{center}} is edge-color isomorphic to the $\{i,j\}$-component $\mathcal{C}$, so $\mathcal{C}$ has a unique maximal element.  

Now we show that the stated combinatorial conditions are sufficient. 
We begin with the minuscule case and argue that any such $R$ must have a unique maximal element. 
Indeed suppose that we have $\selt \leftarrow \relt \rightarrow \telt$ for some $\relt$, $\selt$, and $\telt$ in $R$. 
Since $R$ is fibrous, we have $\selt \mybackarrow{i} \relt \myarrow{j} \telt$ for some colors $i \not = j$ in $I$. 
Since the length of any $j$-component is at most one in $R$, then $-1 \leq m_{j}(\selt) = \langle wt(\selt),\alpha_{j}^{\vee} \rangle = \langle wt(\relt),\alpha_{j}^{\vee} \rangle + \langle \alpha_{i},\alpha_{j}^{\vee} \rangle = -1 + \langle \alpha_{i},\alpha_{j}^{\vee} \rangle$, and therefore $0 \leq \langle \alpha_{i},\alpha_{j}^{\vee} \rangle$.  
Since $i$ and $j$ are distinct, we have $\langle \alpha_{i},\alpha_{j}^{\vee} \rangle \leq 0$. 
It follows that $\langle \alpha_{i},\alpha_{j}^{\vee} \rangle = 0 = \langle \alpha_{j},\alpha_{i}^{\vee} \rangle$. 
Let $\mathcal{C}$ be the $\{i,j\}$-component containing $\relt$, $\selt$, and $\telt$. 
By hypothesis, $\mathcal{C}$ has a unique maximal element. 
It follows from this reasoning that any two elements of $R$ share a common upper bound. 
Finiteness of $R$ means, then, that $R$ has a unique maximal element. 

Let $\melt$ be the unique maximal element of $R$. 
Since $\melt$ is the unique maximal element of $R$ then $\melt$ is prominent, hence $\lambda := wt(\melt)$ is dominant. 
We will show that $\lambda$ is minuscule. 
Pick a coroot $\alpha \in \Phi$, so $\alpha^{\vee} = \sigma(\alpha_{i}^{\vee})$ for some simple coroot $\alpha_{i}$ and some Weyl group element $\sigma$. 
Then $\langle \lambda,\alpha^{\vee} \rangle = \langle \sigma^{-1}(\lambda),\alpha_{i} \rangle$. 
Since $\Pi(R) = \Pi(\lambda)$ by \PiTheorem, then $\sigma^{-1}(\lambda) = wt(\xelt)$ for some $\xelt \in R$. 
Since all $i$-components of $R$ have length at most one, then $\langle wt(\xelt),\alpha_{i}^{\vee} \rangle \in \{0,\pm 1\}$, so $\langle \lambda,\alpha^{\vee} \rangle \in \{0,\pm 1\}$. 
Since $\lambda$ meets the conditions of (M1) of \MinQuasiMinProposition, $\lambda$ is minuscule. 

To complete the proof in the minuscule case, we show that $wt(\xelt) = wt(\yelt)$ in $R$ means that $\xelt = \yelt$. 
We will apply \InitialSplittingTheorem\ to $R$ with $J = I$ and $\nu = 0$. 
Let $\mathcal{S}_{I,0} := \{\melt\}$. 
Take $\tau$ to be the identity function on $R \setminus \{\melt\}$. 
For any $\xelt \not= \melt$ in $R$, freely pick an edge $\xelt \myarrow{k} \yelt$, and then set $\kappa(\xelt) := k$. 
Since $m_{\kappa(\xelt)}(\xelt) = -1$ and $\tau(\xelt) = \xelt$, we get the identity $wt(\tau(\xelt)) = wt(\xelt) - (1 - 0 + m_{\kappa(\xelt)}(\xelt))\alpha_{\kappa(\xelt)}$. 
By \InitialSplittingTheorem, $\WGF(R) = \chi_{_{\lambda}}$. 
Since $\Pi(\lambda) = W\lambda$, then for each $\mu \in \Pi(\lambda)$, we have $\mu = \sigma(\lambda)$ for some $\sigma \in W$.  
By \FinitenessCorollary, $d_{\lambda,\mu} = d_{\lambda,\sigma(\lambda)} = d_{\lambda,\lambda} = 1$. 
So if $\xelt$ and $\yelt$ in $R$ have the same weight $\mu$, then $d_{\lambda,\mu} = 1$ means that $\xelt$ must equal $\yelt$. 
Therefore $R \cong \Pi(R) = \Pi(\lambda) = R(\lambda)$. 

Now for the quasi-minuscule case, take $R$ satisfying the stated combinatorial conditions.  
As a preliminary observation, we demonstrate that for distinct colors $i$ and $j$, no $i$-component of length two in $R$ can have elements in common with a $j$-component of length two. 
Indeed, let us suppose that we have and $i$-component $\mathcal{C}_{i}$ depicted as $\xelt_{0} \myarrow{i} \xelt_{1} \myarrow{i} \xelt_{2}$ and a $j$-component $\mathcal{C}_{j}$ depicted as $\yelt_{0} \myarrow{j} \yelt_{1} \myarrow{j} \yelt_{2}$ in $R$. 
Since $\xelt_{1}$ has no incident edges of color $j$, then $\xelt_{1} \not\in \mathcal{C}_{j}$. 
If $\xelt_{0}$ is in $\mathcal{C}_{j}$, then since $\delta_{j}(\xelt_{0}) = 0$, we must have $\xelt_{0} = \yelt_{2}$. 
In particular, $\yelt_{1} < \xelt_{1}$, a strict inequality in the partially ordered set $R$. 
But, $wt(\xelt_{1}) = 0 = wt(\yelt_{1})$, so $\xelt_{1}$ and $\yelt_{1}$ must have the same rank by \WeightLemma.2. 
This contradiction means that $\xelt_{0} \not\in \mathcal{C}_{j}$. 
Similarly one can see that $\xelt_{2} \not\in \mathcal{C}_{j}$. 

Let us assume for the moment that $R$ has a unique maximal element $\melt$. 
Then as before, $\lambda := wt(\melt)$ is necessarily dominant. 
We will show that $\lambda$ is quasi-minuscule. 
Pick a coroot $\alpha^{\vee} \in \Phi^{\vee}$, so $\alpha^{\vee} = \sigma(\alpha_{i}^{\vee})$ for some simple coroot $\alpha_{i}$ and some Weyl group element $\sigma$. 
Then $\langle \lambda,\alpha^{\vee} \rangle = \langle \sigma^{-1}(\lambda),\alpha_{i}^{\vee} \rangle$. 
Since $\Pi(R) = \Pi(\lambda)$ by \PiTheorem, then $\sigma^{-1}(\lambda) = wt(\xelt)$ for some $\xelt \in R$. 
Similar to the minuscule case, as long as $l_{i}(\xelt) \leq 1$, we will have $\langle wt(\xelt),\alpha_{i}^{\vee} \rangle \in \{0,\pm 1\}$, and therefore $\langle \lambda,\alpha^{\vee} \rangle \in \{0,\pm 1\}$. 
If $l_{i}(\xelt) \geq 2$, then combinatorial condition (1) in the quasi-minuscule case implies that $l_{i}(\xelt) = 2$ and that $wt(\xelt)$ is one of  $-\alpha_{i}$, $0$, or $\alpha_{i}$, so that $\langle wt(\xelt),\alpha_{i}^{\vee} \rangle$ is one of $-2$, $0$, or $2$ (respectively). 
Thus in all circumstances $\langle \lambda,\alpha^{\vee} \rangle = \langle wt(\xelt),\alpha_{i}^{\vee} \rangle \in \{0, \pm 1, \pm 2\}$.  
Suppose now that $\langle \lambda,\alpha^{\vee} \rangle = 2$. 
The prior analysis shows this only happens when $wt(\xelt) = \alpha_{i}$, that is, only when $\sigma^{-1}(\lambda) = \alpha_{i}$, i.e.\ when $\lambda = \sigma(\alpha_{i})$. 
Of course, since $\alpha^{\vee} = \sigma(\alpha_{i}^{\vee})$, then $\alpha$ also has the same length as $\alpha_{i}$, hence $\alpha = \sigma(\alpha_{i})$. 
Then $\alpha = \lambda$. 
So, there is a unique root $\beta \in \Phi$ such that $\langle \lambda,\beta^{\vee} \rangle = 2$, namely $\beta = \lambda$. 
So $\lambda$ meets the conditions of (QM1) of \MinQuasiMinProposition, therefore $\lambda$ is quasi-minuscule. 

Following the pattern of the proof in the minuscule case, we apply \InitialSplittingTheorem\ to $R$ with $J = I$ and $\nu = 0$. 
Let $\mathcal{S}_{I,0} := \{\melt\}$. 
Let $\xelt \not= \melt$ in $R$. 
If $l_{i}(\xelt) \leq 1$ for all $i \in I$, then set $\tau(\xelt) := \xelt$; freely pick an edge $\xelt \myarrow{k} \yelt$, and then set $\kappa(\xelt) := k$. 
On the other hand, if $l_{i}(\xelt) = 2$ for some $i \in I$, then $\comp_{i}(\xelt)$ can be depicted as $\xelt_{0} \myarrow{i} \xelt_{1} \myarrow{i} \xelt_{2}$. 
Now by our preliminary observation, none of these elements can be part of a $j$-component of length two if $j \not= i$. 
In this case, we take $\tau(\xelt_{0}) := \xelt_{1}$, $\tau(\xelt_{1}) := \xelt_{0}$, and (as long as $\xelt_{2} \not= \melt$) $\tau(\xelt_{2}) := \xelt_{2}$. 
Moreover, we let $\kappa(\xelt_{0}) := i$ and $\kappa(\xelt_{1}) := i$; as long as $\xelt_{2} \not= \melt$, freely pick an edge $\xelt_{2} \myarrow{k} \yelt$, and then set $\kappa(\xelt_{2}) := k$. 
Clearly $\tau$ is a bijection. 
To check the key identity in the hypothesis of \InitialSplittingTheorem, we consider cases. 
One can see that for any $\xelt \not= \melt$ such that $\tau(\xelt) = \xelt$ we will have  
$m_{\kappa(\xelt)}(\xelt) = -1$, and so we get the identity $wt(\tau(\xelt)) = wt(\xelt) - (1 - 0 + m_{\kappa(\xelt)}(\xelt))\alpha_{\kappa(\xelt)}$. 
The other case we must check is when $\tau(\xelt) = \yelt$ and $\tau(\yelt) = \xelt$ when $\xelt \myarrow{i} \yelt$ is an edge in an $i$-component of length two. 
Then $m_{i}(\xelt) = -2$ and $m_{i}(\yelt) = 0$. 
So we get $wt(\tau(\xelt)) = wt(\yelt) = wt(\xelt) + \alpha_{i} = wt(\xelt) - (1 + m_{i}(\xelt))\alpha_{i}$ and $wt(\tau(\yelt)) = wt(\xelt) = wt(\yelt) - \alpha_{i} = wt(\yelt) - (1 + m_{i}(\yelt))\alpha_{i}$, which is what we needed to verify. 
By \InitialSplittingTheorem, $\WGF(R) = \chi_{_{\lambda}}$. 

Now $\Pi(\lambda) = W\lambda \disjointunion \{0\}$ by \MinQuasiMinProposition.  
Let $\mu \in \Pi(\lambda)$, so $\mu = \sigma(\lambda)$ for some $\sigma \in W$.  
By \FinitenessCorollary, $d_{\lambda,\mu} = d_{\lambda,\sigma(\lambda)} = d_{\lambda,\lambda} = 1$. 
So if $\xelt$ and $\yelt$ in $R$ have the same weight $\mu$, then $\xelt$ must equal $\yelt$. 
Also, for each short simple root $\alpha_{i}$, there is exactly one $\xelt \in R$ such that $wt(\xelt) = \alpha_{i}$. 
Then $\langle wt(\xelt),\alpha_{i}^{\vee} \rangle = 2$, so $\rho_{i}(\xelt) = 2$ and $\delta_{i}(\xelt) = 0$. 
In particular, $\xelt = \xelt_{2}$ in a chain $\xelt_{0} \myarrow{i} \xelt_{1} \myarrow{i} \xelt_{2}$. 
By \QuasiMinLemma, $d_{\lambda,0}$ is the number of short simple roots, so the above correspondence accounts for all $i$-components of length two in $R$ ($i \in I$). 
One can see now that $R \cong R(\lambda)$. 

So, to complete the proof we must show that $R$ has a unique maximal element. 
Indeed suppose that we have $\selt \leftarrow \relt \rightarrow \telt$ for some $\relt$, $\selt$, and $\telt$ in $R$. 
Since $R$ is fibrous, we have $\selt \mybackarrow{i} \relt \myarrow{j} \telt$ for some colors $i \not = j$ in $I$. 
We have four cases to consider: (0) $\langle \alpha_{i},\alpha_{j}^{\vee} \rangle = 0 = \langle \alpha_{j},\alpha_{i}^{\vee} \rangle$, (1) $\langle \alpha_{i},\alpha_{j}^{\vee} \rangle = -1 = \langle \alpha_{j},\alpha_{i}^{\vee} \rangle$, (2) $\{\langle \alpha_{i},\alpha_{j}^{\vee} \rangle, \langle \alpha_{j},\alpha_{i}^{\vee} \rangle\} = \{-1,-2\}$, and (3) $\{\langle \alpha_{i},\alpha_{j}^{\vee} \rangle, \langle \alpha_{j},\alpha_{i}^{\vee} \rangle\} = \{-1,-3\}$. 
In cases (0) and (1), the hypotheses of the theorem statement guarantee that the $\{i,j\}$-component of $R$ containing $\relt$, $\selt$, and $\telt$ has a unique maximal element. 
In particular, $\selt$ and $\telt$ have a common upper bound in $R$. 
Now suppose we are in case (2), and without loss of generality assume that $\langle \alpha_{j},\alpha_{i}^{\vee} \rangle = -2$. 
Then $m_{i}(\telt) = \langle wt(\relt),\alpha_{i}^{\vee} \rangle + \langle \alpha_{j},\alpha_{i}^{\vee} \rangle = m_{i}(\relt) - 2 \leq -2$.  
Since $\comp_{i}(\telt)$ has length at most two, we get $m_{i}(\telt) \geq -2$. 
Therefore, $m_{i}(\telt) = -2$, from which it follows that $m_{i}(\relt) = 0$. 
This means that $\comp_{i}(\relt)$ must have length two. 
Since $\comp_{i}(\telt)$ intersects $\comp_{j}(\telt)$, then by our preliminary observation $\comp_{j}(\telt)$ cannot have length two. 
So, $m_{j}(\telt) = 1$ and $m_{j}(\relt) = -1$. 
Then $m_{j}(\selt) = \langle wt(\relt),\alpha_{j}^{\vee} \rangle + \langle \alpha_{i},\alpha_{j}^{\vee} \rangle = m_{j}(\relt) - 1 = -2$. 
So $\comp_{j}(\selt)$ has length two, hence (by our preliminary observation) $\comp_{i}(\selt)$ has length one. 
But we have already shown that $\comp_{i}(\relt) = \comp_{i}(\selt)$ must have length two. 
This contradiction means that when $\{\langle \alpha_{i},\alpha_{j}^{\vee} \rangle, \langle \alpha_{j},\alpha_{i}^{\vee} \rangle\} = \{-1,-2\}$, we cannot have two edges of colors $i$ and $j$ respectively meeting in a ``vee'' like $\selt \mybackarrow{i} \relt \myarrow{j} \telt$. Similar reasoning rules out such a structure in case (3) when $\{\langle \alpha_{i},\alpha_{j}^{\vee} \rangle, \langle \alpha_{j},\alpha_{i}^{\vee} \rangle\} = \{-1,-3\}$.
It follows from this reasoning that any two elements of $R$ share a common upper bound. 
Finiteness of $R$ means, then, that $R$ has a unique maximal element.\hfill\QED

\noindent 
{\bf \MinQuasiMinRemark}\ \ In fact, the sufficient conditions stated in \MinQuasiMinTheorem\ can be weakened as follows.  
In the minuscule case, {\sl (2)} can be replaced by: Whenever $\langle \alpha_{i},\alpha_{j}^{\vee} \rangle = 0 = \langle \alpha_{j},\alpha_{i}^{\vee} \rangle$, then no $\{i,j\}$-component of $R$ is isomorphic to any member of the family of edge-colored posets whose first two members are depicted in \TangledFigure.1. 
In the quasi-minuscule case, {\sl (2)} can be replaced by this same statement and {\sl (3)} can be replaced by:  whenever $\langle \alpha_{i},\alpha_{j}^{\vee} \rangle = -1 = \langle \alpha_{j},\alpha_{i}^{\vee} \rangle$, then no $\{i,j\}$-component of $R$ is isomorphic to any member of the family of edge-colored posets whose first two members are depicted in \TangledFigure.2. 
However, the proof of \MinQuasiMinTheorem\ given above uses the simpler and stronger hypotheses of the theorem statement.\hfill\QED

{\bf Supporting graphs as splitting posets.} 
{\em For the remainder of this section we assume Weyl's 
character formula 
(\WeylsTheorem).  No results from later sections depend on the results we develop in the remainder of this section.} 

This paragraph follows \cite{DonSupp} and \cite{Hum}. 
Associate to the root system $\Phi$ with its given choice of simple 
roots the rank $n$ complex semisimple 
Lie algebra $\mathfrak{g}$ 
with Chevalley generators $\{\myqx_{i},\myqy_{i},\myqh_{i}\}_{i \in 
I}$ satisfying the Serre relations, cf.\ \S 18 of \cite{Hum}.  
For the remainder of this discussion  of supporting 
graphs, $V$ denotes a finite-dimensional 
(f.d.\ for short) $\mathfrak{g}$-module.  
For any $\mu \in \Lambda$, $V_{\mu} = \{v \in V\, |\, \myqh_{i}.v = 
\langle \mu,\alpha_{i}^{\vee} \rangle v  
\mbox{ for all } i \in I\}$ is the $\mu$-weight space for 
$V$.  We have $V = \bigoplus_{\mu \in \Lambda} V_{\mu}$, and a weight 
basis is any basis for $V$ that respects this decomposition.  
Finite-dimensional $\mathfrak{g}$-modules are completely reducible, and 
the irreducible f.d.\  
modules are indexed by dominant weights.  
In particular, if 
$V(\lambda)$ is an irreducible f.d.\ 
$\mathfrak{g}$-module 
corresponding to dominant weight $\lambda$, then there is a ``highest'' 
weight vector $v_{\lambda}$ (unique up to scalar multiple) such that 
$\myqx_{i}.v_{\lambda} = 0$ for all $i \in I$.  For such $V(\lambda)$, 
we have $\mathrm{char}(V(\lambda)) := \sum_{\mu \in 
\Lambda}(\dim V(\lambda)_{\mu})e^{\mu} = \chi_{_{\lambda}}$, the 
latter equality by WCF. 
For the generic $\mathfrak{g}$-module $V$, 
$\mathrm{char}(V) := 
\sum_{\mu \in \Lambda}\, (\dim V_{\mu})e^{\mu} \in 
\mathbb{Z}[\Lambda]^{W}$ is a Weyl symmetric function, and 
$\mathrm{char}(V) = 
\sum_{i=1}^{r}\chi_{_{\lambda_{i}}}$ if and only if $V \cong 
V(\lambda_{1}) \oplus \cdots \oplus V(\lambda_{r})$ for irreducible 
f.d.\ 
$\mathfrak{g}$-modules $V(\lambda_{i})$.   
Given any weight basis $\{v_{\telt}\}_{\telt \in R}$ for $V$ (with basis 
vectors indexed by a set $R$), build an edge-colored directed graph 
with elements of $R$ as vertices and with 
edges $\selt \myarrow{i} \telt$ if $\myqX_{\telt,\selt} \not= 0$ 
or $\myqY_{\selt,\telt} \not= 0$ when we write $\myqx_{i}.v_{\selt} = 
\sum_{\uelt \in R}\myqX_{\uelt,\selt}v_{\uelt}$ and 
$\myqy_{i}.v_{\telt} = 
\sum_{\relt \in R}\myqY_{\relt,\telt}v_{\relt}$. Call $R$ a 
{\em supporting 
graph} for $V$. 

For the statement of the next proposition, we need the following notion. 
Fix a weight basis for the $\mathfrak{g}$-module $V$ and a numbering $\mu_{1},\ldots,\mu_{k}$ of the weights for which $\dim V_{\mu_i} > 0$. 
The set of all weight bases for $V$ can then be identified with $GL(\dim V_{\mu_1},\mathbb{C}) \times \cdots \times GL(\dim V_{\mu_k},\mathbb{C})$, where the latter is viewed only as a Cartesian product of sets: for any $k$-tuple $(P_1, \ldots , P_k) \in GL(\dim V_{\mu_1},\mathbb{C}) \times \cdots \times GL(\dim V_{\mu_k},\mathbb{C})$, view each $P_j$ as a change of basis matrix from the fixed basis for $V_{\mu_j}$ to some new basis. 
We say that a set of weight bases for $V$ comprises {\em almost all} weight bases if the corresponding subset of $GL(\dim V_{\mu_1},\mathbb{C}) \times \cdots \times GL(\dim V_{\mu_k},\mathbb{C})$ is Zariski open. 
One part of the next proposition is that almost all weight bases for an irreducible f.d.\ $\mathfrak{g}$-module $V(\lambda)$ have $U(\lambda)$ as their supporting graph. 
This observation is due to Proctor \cite{PrPersonal}. 

\noindent 
{\bf \SupportingGraphProp}\ \ {\sl Any supporting graph $R$ for an f.d.\ $\mathfrak{g}$-module $V$ is a splitting poset for the Weyl symmetric function $\mathrm{char}(V)$. 
If $V(\lambda)$ is an irreducible f.d.\ $\mathfrak{g}$-module with (dominant) highest weight $\lambda$, then any supporting graph for $V(\lambda)$ is a connected splitting poset for $\chi_{_{\lambda}}$, has a unique maximal element $\mathbf{m}$, and $wt(\mathbf{m}) = \lambda$. 
Moreover, the unique maximal splitting poset $U(\lambda)$ is the supporting graph for almost all weight bases for $V(\lambda)$.}   

{\em Proof.} Any supporting graph $R$ for $V$ is ${\Phi}$-structured by Lemmas 3.1.A and 3.2.A of \cite{DonSupp}.  
Now $\displaystyle \WGF(R) = \sum_{\telt \in 
R}e^{wt(\telt)} = \sum_{\mu \in \Lambda}\left(\sum_{\telt: wt(\telt) = 
\mu}e^{\mu}\right) = \sum_{\mu \in \Lambda}(\dim V_{\mu})e^{\mu} = 
\mathrm{char}(V)$, which is a Weyl symmetric function by the paragraph 
preceding the proposition statement.  
So $R$ is a splitting poset for $\mathrm{char}(V)$. Now say $V = 
V(\lambda)$.  Lemma 3.1.F of \cite{DonSupp} guarantees that the 
supporting graph $R$ is connected and has a unique maximal element 
$\mathbf{m}$ with $wt(\mathbf{m}) = \lambda$.  

Before we prove the assertion that almost all weight bases for $V(\lambda)$ have supporting graph $U(\lambda)$, we make the following simple observation. 
Let $V$ and $W$ be arbitrary complex vector spaces of dimensions $p$ and $q$ respectively, and suppose $T: V \longrightarrow W$ is linear and injective. 
Then one can fix bases for $V$ and $W$ with respect to which the representing matrix for $T$ is $\left[\begin{array}{c} \mathcal{I}\\ \hline \mathcal{O}\end{array}\right]$, where $\mathcal{I}$ is the $p \times p$ identity matrix and $\mathcal{O}$ is the $(q-p) \times p$ zero matrix. 
All pairs of bases for $V$ and $W$ can be identified with $GL(p,\mathbb{C}) \times GL(q,\mathbb{C})$ by associating to each pair of bases $(\mathcal{B},\mathcal{C})$ the pair $(P_{\mathcal{B}},Q_{\mathcal{C}})$, where $P_{\mathcal{B}} \in GL(p,\mathbb{C})$ and $Q_{\mathcal{C}} \in GL(q,\mathbb{C})$ are the obvious change of basis matrices. 
So the representing matrix for $T$ with respect to the pair of bases $(\mathcal{B},\mathcal{C})$ is $Q_{\mathcal{C}} \left[\begin{array}{c} \mathcal{I}\\ \hline \mathcal{O}\end{array}\right] P_{\mathcal{B}}$.  
Then those pairs of bases for which the representing matrix for $T$ has no zero entries correspond to a Zariski open subset of $GL(p,\mathbb{C}) \times GL(q,\mathbb{C})$. 
A similar statement holds when $T$ is surjective. 

Now fix a weight basis for $V := V(\lambda)$, and let $R$ denote its supporting graph. 
Then by our work so far we know that $\Pi(R) = \Pi(\lambda)$ is precisely the set of weights $\mu$ for which $\dim V_{\mu} > 0$. 
In particular, every such weight is part of an edge in $\Pi(\lambda)$. 
Say $\mu \myarrow{i} \nu$ is such an edge. 
Let $\mathfrak{g}_i$ be the Lie subalgebra of $\mathfrak{g}$ generated by $\{\myqx_i, \myqy_i, \myqh_i\}$, so $\mathfrak{g}_i$ is isomorphic to $\mathfrak{sl}(2,\mathbb{C})$. 
Let $\myqX_i$ and $\myqY_i$ be the transformations $V \longrightarrow V$ corresponding to the actions of $\myqx_i$ and $\myqy_i$ respectively. 
Let $\tilde{\myqX_i}: V_{\mu} \longrightarrow V_{\nu}$ be the restriction of $\myqX_i$ to $V_{\mu}$ and $\tilde{\myqY_i}: V_{\nu} \longrightarrow V_{\mu}$ be the restriction of $\myqY_i$ to $V_{\nu}$. 
A straightforward analysis of the subspace $\bigoplus_{j \in \mathbb{Z}}V(\lambda)_{\mu+j\alpha_i}$ as a $\mathfrak{g}_i$-submodule shows that one of $\tilde{\myqX_i}$ or $\tilde{\myqY_i}$ is injective and the other is surjective.  
(See for example the proof of Lemma 2.1 of \cite{DonSupp}.) 
By the previous paragraph, then, the set of pairs of bases for $V_{\mu}$ and $V_{\nu}$ for which the representing matrices for $\tilde{\myqX_i}$ and $\tilde{\myqY_i}$ have no zero entries comprises a Zariski open subset of $GL(\dim V_{\mu},\mathbb{C}) \times GL(\dim V_{\nu},\mathbb{C})$, which we can view as a Zariski open subset of $GL(\dim V_{\mu_1},\mathbb{C}) \times \cdots \times GL(\dim V_{\mu_k},\mathbb{C})$. 
Any such choice of bases for $V_{\mu}$ and $V_{\nu}$ yields the maximum number of edges between the sets $wt^{-1}(\mu)$ and $wt^{-1}(\nu)$ in the supporting graph. 
Repeat this argument for all edges in $\Pi(\lambda)$. 
The intersection of all of the resulting Zariski open subsets is Zariski open in $GL(\dim V_{\mu_1},\mathbb{C}) \times \cdots \times GL(\dim V_{\mu_k},\mathbb{C})$. 
Any corresponding weight basis has $U(\lambda)$ as its supporting graph.\hfill\QED 

\noindent 
{\bf \SupportingGraphRemarks}\ \ Some further comments on supporting graphs: 

(A) One view of \SupportingGraphProp\ is:  
A prerequisite for a poset $R$ to be the supporting 
graph for some weight basis for a 
$V$ is that it be a splitting poset for 
$\mathrm{char}(V)$.  

(B) In \cite{DonSupp} it is observed that any 
connected supporting graph is rank symmetric, rank unimodal, and also 
``strongly Sperner,'' i.e.\ for each positive integer $k$, 
the largest 
union of  $k$ antichains of $R$ is no larger than the largest union 
of  $k$ ranks. See Proposition 3.11 of \cite{DonSupp}, which 
applies Proctor's ``Peck poset theorem'' \cite{PrPeck}.\hfill\QED 

{\bf Weyl bialternants with unique splitting posets.} 
The next results   
classify those Weyl bialternants which have unique splitting 
posets, where ``unique'' is understood to mean 
up to isomorphism of edge-colored posets.  The first result 
concerns irreducible root systems, the second  
concerns reducible root systems.  
These results extend a classification result originally due to Howe 
\cite{Howe} (Theorem 4.6.3) and rederived in \cite{DonSupp} (Theorem 6.7). 
That classification result says that all weight spaces are one-dimensional 
for an irreducible finite-dimensional 
representation of a simple Lie algebra if and only if it is 
one of the irreducible representations (identified by type of the Lie 
algebra and highest weight) specified in \UniqueSplittingPropOne.  
For several other equivalent statements and a uniform construction of 
these representations, see \S 6 of \cite{DonSupp}.  
See \cite{Panyushev} for discussion of and further development of 
properties of algebras associated with these representations. 
For these representations, the 
unique splitting poset for the associated Weyl bialternant 
is a distributive lattice, cf.\ Corollary 6.8 
of \cite{DonSupp}. For explicit descriptions of these 
well-known distributive 
lattices, see for example \S 4 of \cite{DonSupp} for the $\myA_{n}$ 
cases;  
\S 5 of \cite{DonSupp} for the $\myC_{n}$ cases; \cite{DLP1} for the 
$(\myB_{n},\omega_1)$ 
and $(\myG_{2},\omega_1)$ cases; and \cite{PrEur}, 
\cite{WildbergerAdv}, and \cite{WildbergerEur} 
for the remaining (minuscule) cases. 

\noindent 
{\bf \UniqueSplittingPropOne}\ \ {\sl Let $\chi_{_{\lambda}}$ be 
a Weyl bialternant for an irreducible root system} 
$\Phi${\sl , with} $\lambda \not= 0${\sl . 
There is a unique splitting poset for 
$\chi_{_{\lambda}}$ if and only if the pair} $(\Phi,\lambda)$ {\sl is one 
of:}  
$(\myA_{n},m\omega_{1})$ ($m \geq 1$){\sl ;} 
$(\myA_{n},m\omega_{n})$ ($m \geq 1$){\sl ;} 
$(\myA_{n},\omega_{k})$ ($1 \leq k \leq n$){\sl ;} 
$(\myB_{n},\omega_{1})${\sl ;} $(\myB_{n},\omega_{n})${\sl ;} 
$(\myC_{n},\omega_{1})${\sl ;} $(\myC_{2},\omega_{2})${\sl ;} 
$(\myC_{3},\omega_{3})${\sl ;} 
$(\myD_{n},\omega_{1})${\sl ;} $(\myD_{n},\omega_{n-1})${\sl ;} 
$(\myD_{n},\omega_{n})${\sl ;} $(\myE_{6},\omega_{1})${\sl ;} 
$(\myE_{6},\omega_{6})${\sl ;} $(\myE_{7},\omega_{7})${\sl ; or} 
$(\myG_{2},\omega_{1})${\sl .} 

{\em Proof.}  We will consider the 
associated simple Lie algebra $\mathfrak{g}$ with irreducible 
representation $V$ having highest weight $\lambda$. 
First, suppose the Weyl bialternant 
$\chi_{_{\lambda}}$ has a unique splitting 
poset as in the proposition statement.  By \SupportingGraphProp, 
any supporting 
graph for $V$ is a splitting poset for $\chi_{_{\lambda}}$, and 
therefore $V$ has (up to isomorphism of edge-colored posets) only one 
supporting graph.  Apply Proposition 6.3/Theorem 6.7 of 
\cite{DonSupp} to see that 
$(\Phi,\lambda)$ must be one of the pairs listed in the proposition 
statement. 

Conversely, say $(\Phi,\lambda)$ is one of the indicated pairs.  It is 
well-known (see for example \S 6 of \cite{DonSupp}) that all weight 
spaces of $V$ are one-dimensional.  
Then for the unique maximal splitting poset, we have $U(\lambda) 
\cong \Pi(\lambda)$, an isomorphism of edge-colored directed graphs. 
But, since the $i$-components of $\Pi(\lambda)$ are chains, then by 
\MinimalProp, $U(\lambda)$ is an edge-minimal splitting poset for 
$\chi_{_{\lambda}}$.  
Since no proper subgraph of $U(\lambda)$ can be a splitting poset for $\chi_{_{\lambda}}$, it follows from \UniqueMaximalProp\ that $U(\lambda)$ is the only splitting poset for $\chi_{_{\lambda}}$.\hfill\QED

For $1 \leq i \leq k$, let $\Phi_{i}$ be an irreducible root 
system with corresponding lattice of weights $\Lambda_{i}$.  
As in \S \WeylSection, suppose for a root system $\Phi$ we have 
$\Phi = \Phi_{1} \disjointunion \cdots 
\disjointunion \Phi_{k}$, so the 
corresponding lattice of weights is $\Lambda = 
\Lambda_{1} \oplus \cdots \oplus \Lambda_{k}$.  For a 
dominant weight $\lambda$, write $\lambda = 
\lambda_{1}+\cdots+\lambda_{k}$ 
where $\lambda_{i}$ is dominant in $\Lambda_{i}$ for $1 \leq i \leq 
k$. Then 
the $\Phi$-Weyl bialternant $\chi^{\Phi}_{_{\lambda}}$ 
is $\chi^{\Phi_{1}}_{_{\lambda_{1}}}\cdots
\chi^{\Phi_{k}}_{_{\lambda_{k}}}$, 
where $\chi^{\Phi_{i}}_{_{\lambda_{i}}}$ is the $\Phi_{i}$-Weyl 
bialternant corresponding to the pair $(\Phi_{i},\lambda_{i})$. 

\noindent 
{\bf \UniqueSplittingPropTwo}\ \ {\sl Keep the above notation.  
(1) Suppose $\chi^{\Phi}_{_{\lambda}}$ 
has a unique splitting poset $R$.  Then for each $1 \leq i \leq k$, 
either 
$\lambda_{i}=0$ (in which case we let $R_{i}$ be the one-element 
poset) or $(\Phi_{i},\lambda_{i})$ is one of the pairs from 
\UniqueSplittingPropOne\ (in which case we let $R_{i}$ be the 
unique splitting poset 
associated to the pair $(\Phi_{i},\lambda_{i})$).   Moreover,  
$R \cong R_{1} \times \cdots \times R_{k}$, an isomorphism of 
edge-colored posets.  (2) Conversely, suppose that for each 
$1 \leq i \leq k$, either 
$\lambda_{i}=0$ (in which case we let $R_{i}$ be the one-element 
poset) or $(\Phi_{i},\lambda_{i})$ is one of the pairs from 
\UniqueSplittingPropOne\ (in which case we let $R_{i}$ be the 
unique splitting poset 
associated to the pair $(\Phi_{i},\lambda_{i})$). Then 
$R_{1} \times \cdots \times R_{k}$ is the unique splitting poset for 
$\chi^{\Phi}_{_{\lambda}}$.} 

{\em Proof.}  For {\sl (1)}, for each $1 \leq i \leq k$, let 
$Q_{i}$ be a splitting poset for the $\Phi_{i}$-Weyl bialternant 
$\chi^{\Phi_{i}}_{_{\lambda_{i}}}$.  Then by \OperationsLemma.2, 
$Q_{1} \times 
\cdots \times Q_{k}$ is a splitting poset for 
$\chi^{\Phi}_{_{\lambda}}$.  
Since $\chi^{\Phi}_{_{\lambda}}$ is assumed to have a unique 
splitting poset 
$R$, then we must have 
$R \cong Q_{1} \times 
\cdots \times Q_{k}$, an isomorphism of edge-colored posets as claimed. 
Moreover, each $\chi^{\Phi_{i}}_{_{\lambda_{i}}}$ 
must have a unique splitting poset, otherwise we can construct a 
splitting poset $Q_{1}' \times \cdots \times Q_{k}'$ for 
$\chi^{\Phi}_{_{\lambda}}$ distinct from $R \cong Q_{1} \times 
\cdots \times Q_{k}$.  Then \UniqueSplittingPropOne\ applies, so 
that either $\lambda_{i} = 0$ or $(\Phi_{i},\lambda_{i})$ is one of 
the pairs from \UniqueSplittingPropOne.  Note that $Q_{i}$ is a 
one-element poset if $\lambda_{i} = 0$. 
For {\sl (2)}, we use reasoning from the proof 
of \UniqueSplittingPropOne.  Observe that an irreducible 
representation with highest 
weight $\lambda$ for the associated semisimple Lie algebra will have 
one-dimensional weight spaces.  The argument used in the second 
paragraph of the proof of \UniqueSplittingPropOne\ applies in this 
case since it nowhere depends on $\Phi$ being irreducible.  We 
conclude that $\Pi(\lambda)$ is the unique splitting poset for 
$\chi^{\Phi}_{_{\lambda}}$.  By \OperationsLemma.2, $R_{1} 
\times \cdots \times 
R_{k}$ is a splitting poset for $\chi^{\Phi}_{_{\lambda}}$, 
and hence it is 
the unique splitting poset.\hfill\QED

\newpage
\noindent 
{\Large \bf \S \FibrousSection.\ Crystal products of fibrous posets and vertex coloring.} 

The Weyl character formula has long been an advantageous point of contact between combinatorics and Lie representation theory.  
The development of crystal graphs in the past several decades has provided another fruitful point of contact.  
In the next two sections, we develop a theory of crystal graphs more or less from scratch.  
Our approach is largely modelled on Stembridge's groundbreaking paper \cite{Stem}.  
However, the computational details are quite different in places and will allow us to generalize the notion of crystal graph in subsequent work, where (for example) the ``fibrous'' hypothesis used here will not be necessary.   

Our purpose in this section is to exhibit those aspects of our set-up that are independent of root system or Weyl symmetric function results, although these results are mainly meaningful as tools for our eventual use in constructing splitting posets.  
However, it is possible that such results could be helpful in studying characters of non-reductive Lie groups such as odd symplectic groups, and it is possible that extracting out some of the purely combinatorial aspects of the theory might help provide further combinatorial insight into the nature of the splitting problem.  
The most important idea of this section is vertex coloring, which will play a large role in this and other related papers. 

{\bf Crystal products of fibrous posets.}  
For the remainder of this section, $R_{1}$ and $R_{2}$ denote fibrous posets with edges colored by the same set $I$.  The product poset $R_{1} \times R_{2}$ need not be fibrous (consider a product of two length one chains of the same color, for example).  Our immediate aim is to define a different kind of product of fibrous posets whose resulting edge-colored poset is also fibrous.  For the remainder of this section, $R_{1} \times R_{2}$ is to be considered a set product only.  

We form a new edge-colored directed graph $\Rcrystal$ from the elements of $R_{1} \times R_{2}$ as follows.  If $\xxelt$ and $\yyelt$ are in $R_{1} \times R_{2}$, write $\xxelt \myarrow{i} \yyelt$ if and only if either 
\begin{enumerate}
\item[(1)] $\xelt_{1} \myarrow{i} \yelt_{1}$ in $R_{1}$ (so $\delta_{i}(\xelt_{1}) > 0$), $\xelt_{2} = \yelt_{2}$ in $R_{2}$, $\rho_{i}(\xelt_{1}) \geq \delta_{i}(\xelt_{2})$, or 
\item[(2)] $\xelt_{1} = \yelt_{1}$ in $R_{1}$, $\xelt_{2} \myarrow{i} \yelt_{2}$ in $R_{2}$, and $\rho_{i}(\xelt_{1}) < \delta_{i}(\xelt_{2})$.
\end{enumerate}
Equivalently we could say that $\xxelt \myarrow{i} \yyelt$ in $\Rcrystal$ if and only if either 
\begin{enumerate}
\item[(1$'$)] $\xelt_{1} \myarrow{i} \yelt_{1}$ in $R_{1}$, $\xelt_{2} = \yelt_{2}$ in $R_{2}$, and $\rho_{i}(\yelt_{1}) > \delta_{i}(\yelt_{2})$, or 
\item[(2$'$)] $\xelt_{1} = \yelt_{1}$ in $R_{1}$ (so $\rho_{i}(\yelt_{2}) > 0$), $\xelt_{2} \myarrow{i} \yelt_{2}$ in $R_{2}$, $\rho_{i}(\yelt_{1}) \leq \delta_{i}(\yelt_{2})$.
\end{enumerate}
We call $\Rcrystal$ the {\em crystal product} of fibrous posets $R_{1}$ and $R_{2}$. The next lemma is stated for the record, and its proof is routine. 

\noindent 
{\bf \FibrousProductLemma}\ \ {\sl The edge-colored directed graph $R_{1} \otimes R_{2}$ is the edge-colored Hasse diagram for a ranked poset, and it is fibrous.}\hfill\QED 

Movement ``up'' or ``down'' within any fibrous poset $R$ can be described by ``raising'' and ``lowering'' operators, as follows.  In our computations, the symbol $\theta$ will act as something of a trash bin.  For each $i \in I$, we define a {\em raising} operator $\widetilde{E}_{i}: R \longrightarrow R \disjointunion \{\theta\}$ and a {\em lowering} operator $\widetilde{F}_{i}: R \longrightarrow R \disjointunion \{\theta\}$ by the rules: For all $\selt \in R$, $\widetilde{E}_{i}(\selt) = \telt$ if $\selt \myarrow{i} \telt$ in $R$ and $\widetilde{E}_{i}(\selt) = \theta$ if $\selt$ is $i$-maximal.  Similarly, for all $\telt \in R$, $\widetilde{F}_{i}(\telt) = \selt$ if $\selt \myarrow{i} \telt$ in $R$ and $\widetilde{F}_{i}(\telt) = \theta$ if $\telt$ is $i$-minimal.  Observe that $\widetilde{E}_{i}$ restricts to a bijection $R \setminus \{\mbox{$i$-maximal elements}\} \longrightarrow R \setminus \{\mbox{$i$-minimal elements}\}$ and that $\widetilde{F}_{i}$ is its inverse. The effect of raising and lowering operators on crystal products follows immediately from our definition of edges in $\Rcrystal$, and we record this in the next lemma. 

\noindent 
{\bf \EFLemma}\ \ {\sl For each $i \in I$ and $\xxelt \in \Rcrystal$, we have}
\[\wE_{i}(\xelt_{1},\xelt_{2}) = \left\{\begin{array}{cl}
(\wE_{i}(\xelt_{1}),\xelt_{2}) & \hspace*{0.25in}\mbox{if } \rho_{i}(\xelt_{1}) \geq \delta_{i}(\xelt_{2}) \mbox{ and } \delta_{i}(\xelt_{1}) > 0\\
(\xelt_{1},\wE_{i}(\xelt_{2})) & \hspace*{0.25in}\mbox{if } \rho_{i}(\xelt_{1}) < \delta_{i}(\xelt_{2})\\
\theta & \hspace*{0.25in}\mbox{otherwise}
\end{array}\right.\]
{\sl and} 
\[\wF_{i}(\xelt_{1},\xelt_{2}) = \left\{\begin{array}{cl}
(\wF_{i}(\xelt_{1}),\xelt_{2}) & \hspace*{0.25in}\mbox{if } \rho_{i}(\xelt_{1}) > \delta_{i}(\xelt_{2})\\
(\xelt_{1},\wF_{i}(\xelt_{2})) & \hspace*{0.25in}\mbox{if } \rho_{i}(\xelt_{1}) \leq \delta_{i}(\xelt_{2}) \mbox{ and } \rho_{i}(\xelt_{2}) > 0\\
\theta & \hspace*{0.25in}\mbox{otherwise.}
\end{array}\right.\]

\vspace{-0.35in}\hfill\QED

The crystal product $R_{1} \otimes R_{2}$ need not be connected.  In fact, its potential to be {\em dis}connected is one of the crucial features of the theory. 

\noindent 
{\bf \CrystalExample}\ \ \CrystalProductFigures\ depict two crystal products.  
Each is the crystal product of a fibrous chain with itself.  
The four-element chain of \CrystalProductFigure.1 is ${\myC_{2}}$-structured, as are the connected components of the resulting crystal product. 
The seven-element chain of \CrystalProductFigure.2 is ${\myG_{2}}$-structured, as are the connected components of the crystal product of this chain with itself.\hfill\QED

\newcommand{\GtwoOneRowOneCol}[1]{\setlength{\unitlength}{0.7cm}
\begin{picture}(0,0)
\thinlines
\put(0,0){\line(1,0){0.4}}
\put(0,0.4){\line(1,0){0.4}}
\put(0,0){\line(0,1){0.4}}
\put(0.4,0){\line(0,1){0.4}}
\put(0.11,0.1){\tiny #1}
\end{picture}}
\newcommand{\GtwoOneRowTwoCol}[2]{\setlength{\unitlength}{0.7cm}
\begin{picture}(0,0)
\thinlines
\put(0,-0.05){\line(1,0){0.7}}
\put(0,0.3){\line(1,0){0.7}}
\put(0,-0.05){\line(0,1){0.35}}
\put(0.35,-0.05){\line(0,1){0.35}}
\put(0.7,-0.05){\line(0,1){0.35}}
\put(0.1,0.05){\tiny #1}
\put(0.45,0.05){\tiny #2}
\end{picture}}

\begin{figure}[t]
\begin{center}
{\CrystalProductFigure.1:} An example of a crystal product.  The factors of the crystal product, as well as the connected components of the resulting disjoint sum, are ${\myC_{2}}$-structured. 

\vspace*{0.1in}
\setlength{\unitlength}{0.7cm}
\begin{picture}(18,11)
\multiput(1,4)(0,1.5){4}{\circle*{0.22}}
\thicklines
\put(1,4.15){\color{Red}{\line(0,1){1.2}}}
\put(0.75,4.6){\color{Red}{\em 1}}
\put(0.1,3.75){\GtwoOneRowOneCol{4}}
\put(1,5.65){\color{Cyan}{\line(0,1){1.2}}}
\put(0.75,6.1){\color{Cyan}{\em 2}}
\put(0.1,5.25){\GtwoOneRowOneCol{3}}
\put(1,7.15){\color{Red}{\line(0,1){1.2}}}
\put(0.75,7.6){\color{Red}{\em 1}}
\put(0.1,6.75){\GtwoOneRowOneCol{2}}
\put(0.1,8.25){\GtwoOneRowOneCol{1}}
\multiput(3,4)(0,1.5){4}{\circle*{0.22}}
\thicklines
\put(3,4.15){\color{Red}{\line(0,1){1.2}}}
\put(2.75,4.6){\color{Red}{\em 1}}
\put(3.1,3.75){\GtwoOneRowOneCol{4}}
\put(3,5.65){\color{Cyan}{\line(0,1){1.2}}}
\put(2.75,6.1){\color{Cyan}{\em 2}}
\put(3.1,5.25){\GtwoOneRowOneCol{3}}
\put(3,7.15){\color{Red}{\line(0,1){1.2}}}
\put(2.75,7.6){\color{Red}{\em 1}}
\put(3.1,6.75){\GtwoOneRowOneCol{2}}
\put(3.1,8.25){\GtwoOneRowOneCol{1}}
\put(1.6,6.1){\Large $\bigotimes$}
\put(5.2,5.8){\Huge $\cong$}
\put(15.4,5.8){\Large $\bigoplus$}
\put(11.9,5.8){\Large $\bigoplus$}
%
\multiput(8,0)(1.5,1.5){4}{\multiput(0,6)(1.5,-1.5){4}{\circle*{0.22}}}
\put(11.35,1.35){\GtwoOneRowTwoCol{4}{4}}
\put(9.85,2.85){\GtwoOneRowTwoCol{4}{3}}
\put(8.35,4.35){\GtwoOneRowTwoCol{4}{2}}
\put(6.85,5.85){\GtwoOneRowTwoCol{4}{1}}
\put(14.15,2.85){\GtwoOneRowTwoCol{3}{4}}
\put(11.35,4.35){\GtwoOneRowTwoCol{3}{3}}
\put(9.85,5.85){\GtwoOneRowTwoCol{3}{2}}
\put(8.35,7.35){\GtwoOneRowTwoCol{3}{1}}
\put(15.65,4.35){\GtwoOneRowTwoCol{2}{4}}
\put(12.85,5.85){\GtwoOneRowTwoCol{2}{3}}
\put(11.35,7.35){\GtwoOneRowTwoCol{2}{2}}
\put(9.85,8.85){\GtwoOneRowTwoCol{2}{1}}
\put(17.15,5.85){\GtwoOneRowTwoCol{1}{4}}
\put(15.65,7.35){\GtwoOneRowTwoCol{1}{3}}
\put(14.15,8.85){\GtwoOneRowTwoCol{1}{2}}
\put(11.35,10.35){\GtwoOneRowTwoCol{1}{1}}
%
%
\put(12.4,1.6){\color{Red}{\line(-1,1){1.3}}}
\put(11.6,2.1){\color{Red}{\em 1}}
\put(10.9,3.1){\color{Cyan}{\line(-1,1){1.3}}}
\put(10.1,3.6){\color{Cyan}{\em 2}}
\put(9.4,4.6){\color{Red}{\line(-1,1){1.3}}}
\put(8.4,5.1){\color{Red}{\em 1}}
\put(12.4,4.6){\color{Cyan}{\line(-1,1){1.3}}}
\put(11.6,5.1){\color{Cyan}{\em 2}}
\put(15.4,4.6){\color{Red}{\line(-1,1){1.3}}}
\put(14.6,5.1){\color{Red}{\em 1}}
\put(12.4,7.6){\color{Red}{\line(-1,1){1.3}}}
\put(11.6,8.1){\color{Red}{\em 1}}
\put(15.4,7.6){\color{Cyan}{\line(-1,1){1.3}}}
\put(14.6,8.1){\color{Cyan}{\em 2}}
%
%
\put(14.1,3.1){\color{Cyan}{\line(1,1){1.3}}}
\put(14.5,3.5){\color{Cyan}{\em 2}}
\put(11.1,3.1){\color{Red}{\line(1,1){1.3}}}
\put(11.5,3.5){\color{Red}{\em 1}}
\put(14.1,6.1){\color{Red}{\line(1,1){1.3}}}
\put(14.5,6.5){\color{Red}{\em 1}}
\put(11.1,6.1){\color{Cyan}{\line(1,1){1.3}}}
\put(11.5,6.5){\color{Cyan}{\em 2}}
\put(8.1,6.1){\color{Red}{\line(1,1){1.3}}}
\put(8.5,6.5){\color{Red}{\em 1}}
\put(9.6,7.6){\color{Cyan}{\line(1,1){1.3}}}
\put(10,8){\color{Cyan}{\em 2}}
\put(11.1,9.1){\color{Red}{\line(1,1){1.3}}}
\put(11.5,9.5){\color{Red}{\em 1}}
\end{picture}

\vspace*{-0.5in}
\end{center}
\end{figure}

The crystal product interacts in a natural way with other operations on edge-colored posets. 
We record these results as follows. 
All of the indicated isomorphisms are natural. 

\begin{center}
{\CrystalProductFigure.2:} Another example of a crystal product.  The factors of the crystal product, as well as the connected components of the resulting disjoint sum, are ${\myG_{2}}$-structured. 

\vspace*{-0.2in}
\setlength{\unitlength}{0.7cm}
\begin{picture}(23,27)
\multiput(1,16)(0,1.5){7}{\circle*{0.22}}
\thicklines
\put(1,16.15){\color{Red}{\line(0,1){1.2}}}
\put(0.75,16.6){\color{Red}{\em 1}}
\put(0.1,15.75){\GtwoOneRowOneCol{7}}
\put(1,17.65){\color{Cyan}{\line(0,1){1.2}}}
\put(0.75,18.1){\color{Cyan}{\em 2}}
\put(0.1,17.25){\GtwoOneRowOneCol{6}}
\put(1,19.15){\color{Red}{\line(0,1){1.2}}}
\put(0.75,19.6){\color{Red}{\em 1}}
\put(0.1,18.75){\GtwoOneRowOneCol{5}}
\put(1,20.65){\color{Red}{\line(0,1){1.2}}}
\put(0.75,21.1){\color{Red}{\em 1}}
\put(0.1,20.25){\GtwoOneRowOneCol{4}}
\put(1,22.15){\color{Cyan}{\line(0,1){1.2}}}
\put(0.75,22.6){\color{Cyan}{\em 2}}
\put(0.1,21.75){\GtwoOneRowOneCol{3}}
\put(1,23.65){\color{Red}{\line(0,1){1.2}}}
\put(0.75,24.1){\color{Red}{\em 1}}
\put(0.1,23.25){\GtwoOneRowOneCol{2}}
\put(0.1,24.75){\GtwoOneRowOneCol{1}}
\put(1.7,20.3){\Large $\bigotimes$}
\multiput(3,16)(0,1.5){7}{\circle*{0.22}}
\thicklines
\put(3,16.15){\color{Red}{\line(0,1){1.2}}}
\put(2.75,16.6){\color{Red}{\em 1}}
\put(3.1,15.75){\GtwoOneRowOneCol{7}}
\put(3,17.65){\color{Cyan}{\line(0,1){1.2}}}
\put(2.75,18.1){\color{Cyan}{\em 2}}
\put(3.1,17.25){\GtwoOneRowOneCol{6}}
\put(3,19.15){\color{Red}{\line(0,1){1.2}}}
\put(2.75,19.6){\color{Red}{\em 1}}
\put(3.1,18.75){\GtwoOneRowOneCol{5}}
\put(3,20.65){\color{Red}{\line(0,1){1.2}}}
\put(2.75,21.1){\color{Red}{\em 1}}
\put(3.1,20.25){\GtwoOneRowOneCol{4}}
\put(3,22.15){\color{Cyan}{\line(0,1){1.2}}}
\put(2.75,22.6){\color{Cyan}{\em 2}}
\put(3.1,21.75){\GtwoOneRowOneCol{3}}
\put(3,23.65){\color{Red}{\line(0,1){1.2}}}
\put(2.75,24.1){\color{Red}{\em 1}}
\put(3.1,23.25){\GtwoOneRowOneCol{2}}
\put(3.1,24.75){\GtwoOneRowOneCol{1}}
\put(5.2,20.2){\Huge $\cong$}
\put(11.4,11.8){\Large $\bigoplus$}
\put(17.4,11.8){\Large $\bigoplus$}
\put(20.4,11.8){\Large $\bigoplus$}
%
\multiput(4,6)(1.5,1.5){7}{\multiput(0,6)(1.5,-1.5){7}{\circle*{0.22}}}
\put(11.85,2.85){\GtwoOneRowTwoCol{7}{7}}
\put(10.35,4.35){\GtwoOneRowTwoCol{7}{6}}
\put(8.85,5.85){\GtwoOneRowTwoCol{7}{5}}
\put(7.35,7.35){\GtwoOneRowTwoCol{7}{4}}
\put(5.85,8.85){\GtwoOneRowTwoCol{7}{3}}
\put(4.35,10.35){\GtwoOneRowTwoCol{7}{2}}
\put(2.85,11.85){\GtwoOneRowTwoCol{7}{1}}
\put(14.65,4.35){\GtwoOneRowTwoCol{6}{7}}
\put(11.85,5.85){\GtwoOneRowTwoCol{6}{6}}
\put(10.35,7.35){\GtwoOneRowTwoCol{6}{5}}
\put(8.85,8.85){\GtwoOneRowTwoCol{6}{4}}
\put(7.35,10.35){\GtwoOneRowTwoCol{6}{3}}
\put(5.85,11.85){\GtwoOneRowTwoCol{6}{2}}
\put(4.35,13.35){\GtwoOneRowTwoCol{6}{1}}
\put(16.15,5.85){\GtwoOneRowTwoCol{5}{7}}
\put(13.35,7.35){\GtwoOneRowTwoCol{5}{6}}
\put(11.85,8.85){\GtwoOneRowTwoCol{5}{5}}
\put(10.35,10.35){\GtwoOneRowTwoCol{5}{4}}
\put(8.85,11.85){\GtwoOneRowTwoCol{5}{3}}
\put(7.35,13.35){\GtwoOneRowTwoCol{5}{2}}
\put(5.85,14.85){\GtwoOneRowTwoCol{5}{1}}
\put(17.65,7.35){\GtwoOneRowTwoCol{4}{7}}
\put(14.85,8.85){\GtwoOneRowTwoCol{4}{6}}
\put(13.35,10.35){\GtwoOneRowTwoCol{4}{5}}
\put(13.15,11.85){\GtwoOneRowTwoCol{4}{4}}
\put(10.35,13.35){\GtwoOneRowTwoCol{4}{3}}
\put(8.85,14.85){\GtwoOneRowTwoCol{4}{2}}
\put(7.35,16.35){\GtwoOneRowTwoCol{4}{1}}
\put(19.15,8.85){\GtwoOneRowTwoCol{3}{7}}
\put(16.35,10.35){\GtwoOneRowTwoCol{3}{6}}
\put(14.85,11.85){\GtwoOneRowTwoCol{3}{5}}
\put(13.35,13.35){\GtwoOneRowTwoCol{3}{4}}
\put(11.85,14.85){\GtwoOneRowTwoCol{3}{3}}
\put(10.35,16.35){\GtwoOneRowTwoCol{3}{2}}
\put(8.85,17.85){\GtwoOneRowTwoCol{3}{1}}
\put(20.65,10.35){\GtwoOneRowTwoCol{2}{7}}
\put(19.15,11.85){\GtwoOneRowTwoCol{2}{6}}
\put(16.35,13.35){\GtwoOneRowTwoCol{2}{5}}
\put(14.85,14.85){\GtwoOneRowTwoCol{2}{4}}
\put(13.35,16.35){\GtwoOneRowTwoCol{2}{3}}
\put(11.85,17.85){\GtwoOneRowTwoCol{2}{2}}
\put(10.35,19.35){\GtwoOneRowTwoCol{2}{1}}
%
\put(21.45,12.25){\GtwoOneRowTwoCol{1}{7}}
\put(20.65,13.35){\GtwoOneRowTwoCol{1}{6}}
\put(19.15,14.85){\GtwoOneRowTwoCol{1}{5}}
\put(17.65,16.35){\GtwoOneRowTwoCol{1}{4}}
\put(16.15,17.85){\GtwoOneRowTwoCol{1}{3}}
\put(14.65,19.35){\GtwoOneRowTwoCol{1}{2}}
\put(11.85,20.85){\GtwoOneRowTwoCol{1}{1}}
%
%
%
\put(12.9,3.1){\color{Red}{\line(-1,1){1.3}}}
\put(12.1,3.6){\color{Red}{\em 1}}
\put(11.4,4.6){\color{Cyan}{\line(-1,1){1.3}}}
\put(10.6,5.1){\color{Cyan}{\em 2}}
\put(9.9,6.1){\color{Red}{\line(-1,1){1.3}}}
\put(9.1,6.6){\color{Red}{\em 1}}
\put(8.4,7.6){\color{Red}{\line(-1,1){1.3}}}
\put(7.6,8.1){\color{Red}{\em 1}}
\put(6.9,9.1){\color{Cyan}{\line(-1,1){1.3}}}
\put(6.1,9.6){\color{Cyan}{\em 2}}
\put(5.4,10.6){\color{Red}{\line(-1,1){1.3}}}
\put(4.4,11.1){\color{Red}{\em 1}}
%
\put(12.9,6.1){\color{Cyan}{\line(-1,1){1.3}}}
\put(12.1,6.6){\color{Cyan}{\em 2}}
\put(11.4,7.6){\color{Red}{\line(-1,1){1.3}}}
\put(10.6,8.1){\color{Red}{\em 1}}
\put(8.4,10.6){\color{Cyan}{\line(-1,1){1.3}}}
\put(7.6,11.1){\color{Cyan}{\em 2}}
%
\put(15.9,6.1){\color{Red}{\line(-1,1){1.3}}}
\put(15.1,6.6){\color{Red}{\em 1}}
\put(12.9,9.1){\color{Red}{\line(-1,1){1.3}}}
\put(12.1,9.6){\color{Red}{\em 1}}
\put(11.4,10.6){\color{Red}{\line(-1,1){1.3}}}
\put(10.6,11.1){\color{Red}{\em 1}}
\put(8.4,13.6){\color{Red}{\line(-1,1){1.3}}}
\put(7.6,14.1){\color{Red}{\em 1}}
%
\put(15.9,9.1){\color{Cyan}{\line(-1,1){1.3}}}
\put(15.1,9.6){\color{Cyan}{\em 2}}
\put(14.4,10.6){\color{Red}{\line(-1,1){1.3}}}
\put(13.6,11.1){\color{Red}{\em 1}}
\put(11.4,13.6){\color{Cyan}{\line(-1,1){1.3}}}
\put(10.6,14.1){\color{Cyan}{\em 2}}
%
\put(17.4,10.6){\color{Cyan}{\line(-1,1){1.3}}}
\put(16.6,11.1){\color{Cyan}{\em 2}}
\put(12.9,15.1){\color{Cyan}{\line(-1,1){1.3}}}
\put(12.1,15.6){\color{Cyan}{\em 2}}
%
\put(20.4,10.6){\color{Red}{\line(-1,1){1.3}}}
\put(19.6,11.1){\color{Red}{\em 1}}
\put(17.4,13.6){\color{Red}{\line(-1,1){1.3}}}
\put(16.6,14.1){\color{Red}{\em 1}}
\put(15.9,15.1){\color{Red}{\line(-1,1){1.3}}}
\put(15.1,15.6){\color{Red}{\em 1}}
\put(12.9,18.1){\color{Red}{\line(-1,1){1.3}}}
\put(12.1,18.6){\color{Red}{\em 1}}
%
\put(20.4,13.6){\color{Cyan}{\line(-1,1){1.3}}}
\put(19.6,14.1){\color{Cyan}{\em 2}}
\put(18.9,15.1){\color{Red}{\line(-1,1){1.3}}}
\put(18.1,15.6){\color{Red}{\em 1}}
\put(15.9,18.1){\color{Cyan}{\line(-1,1){1.3}}}
\put(15.1,18.6){\color{Cyan}{\em 2}}
%
%
\put(14.6,4.6){\color{Cyan}{\line(1,1){1.3}}}
\put(15,5){\color{Cyan}{\em 2}}
\put(17.6,7.6){\color{Red}{\line(1,1){1.3}}}
\put(18,8){\color{Red}{\em 1}}
\put(19.1,9.1){\color{Cyan}{\line(1,1){1.3}}}
\put(19.5,9.5){\color{Cyan}{\em 2}}
%
\put(11.6,4.6){\color{Red}{\line(1,1){1.3}}}
\put(12,5){\color{Red}{\em 1}}
\put(14.6,7.6){\color{Red}{\line(1,1){1.3}}}
\put(15,8){\color{Red}{\em 1}}
\put(16.1,9.1){\color{Red}{\line(1,1){1.3}}}
\put(16.5,9.5){\color{Red}{\em 1}}
\put(19.1,12.1){\color{Red}{\line(1,1){1.3}}}
\put(19.5,12.5){\color{Red}{\em 1}}
%
\put(11.6,7.6){\color{Cyan}{\line(1,1){1.3}}}
\put(12,8){\color{Cyan}{\em 2}}
\put(16.1,12.1){\color{Cyan}{\line(1,1){1.3}}}
\put(16.5,12.5){\color{Cyan}{\em 2}}
%
\put(10.1,9.1){\color{Cyan}{\line(1,1){1.3}}}
\put(10.5,9.5){\color{Cyan}{\em 2}}
\put(13.1,12.1){\color{Red}{\line(1,1){1.3}}}
\put(13.5,12.5){\color{Red}{\em 1}}
\put(14.6,13.6){\color{Cyan}{\line(1,1){1.3}}}
\put(15,14){\color{Cyan}{\em 2}}
%
\put(7.1,9.1){\color{Red}{\line(1,1){1.3}}}
\put(7.5,9.5){\color{Red}{\em 1}}
\put(10.1,12.1){\color{Red}{\line(1,1){1.3}}}
\put(10.5,12.5){\color{Red}{\em 1}}
\put(11.6,13.6){\color{Red}{\line(1,1){1.3}}}
\put(12,14){\color{Red}{\em 1}}
\put(14.6,16.6){\color{Red}{\line(1,1){1.3}}}
\put(15,17){\color{Red}{\em 1}}
%
\put(7.1,12.1){\color{Cyan}{\line(1,1){1.3}}}
\put(7.5,12.5){\color{Cyan}{\em 2}}
\put(10.1,15.1){\color{Red}{\line(1,1){1.3}}}
\put(10.5,15.5){\color{Red}{\em 1}}
\put(11.6,16.6){\color{Cyan}{\line(1,1){1.3}}}
\put(12,17){\color{Cyan}{\em 2}}
%
\put(4.1,12.1){\color{Red}{\line(1,1){1.3}}}
\put(4.5,12.5){\color{Red}{\em 1}}
\put(5.6,13.6){\color{Cyan}{\line(1,1){1.3}}}
\put(6,14){\color{Cyan}{\em 2}}
\put(7.1,15.1){\color{Red}{\line(1,1){1.3}}}
\put(7.5,15.5){\color{Red}{\em 1}}
\put(8.6,16.6){\color{Red}{\line(1,1){1.3}}}
\put(9,17){\color{Red}{\em 1}}
\put(10.1,18.1){\color{Cyan}{\line(1,1){1.3}}}
\put(10.5,18.5){\color{Cyan}{\em 2}}
\put(11.6,19.6){\color{Red}{\line(1,1){1.3}}}
\put(12,20){\color{Red}{\em 1}}
\end{picture}

\vspace*{-0.7in}
\end{center}

\noindent 
{\bf \AssociativeCrystalLemma}\ \ {\sl Let $R_{1}$, $R_{2}$, and $R_{3}$ be fibrous posets with edges colored by $I$, and let $\sigma: I \longrightarrow I'$ be some bijection of sets. 
Then} 
\begin{center}
$R_{1} \otimes (R_{2} \otimes R_{3}) \cong (\Rcrystal) \otimes R_{3}$,\\ 
$R_{1} \otimes (R_{2} \oplus R_{3}) \cong (R_{1} \otimes R_{2}) \oplus (R_{1} \otimes R_{3})$, $(R_{1} \oplus R_{2}) \otimes R_{3} \cong (R_{1} \otimes R_{3}) \oplus (R_{2} \otimes R_{3})$,\\ 
$(R_{1} \otimes R_{2})^{*} \cong R_{2}^{*} \otimes R_{1}^{*}$, and $(R_{1} \otimes R_{2})^{\sigma} \cong R_{1}^{\sigma} \otimes R_{2}^{\sigma}$. 
\end{center}
{\sl In particular, if $I$ indexes a set of simple roots for a root system $\Phi$ and if $\sigma = \sigma_{0}$, then $(R_{1} \otimes R_{2})^{\bowtie} \cong R_{2}^{\bowtie} \otimes R_{1}^{\bowtie}$.}\hfill\QED

\noindent 
{\bf \NoncommuteExample}\ \  Although the crystal product is well-behaved with respect to our usual poset operations, it is not, in general, commutative.  See \NoCommuteFigure\ for an example.\hfill\QED

\begin{figure}[h]
\begin{center}
{\NoCommuteFigure:} An example of fibrous posets that do not commute under the crystal product.\\
{\footnotesize (Specifically,  the product on the left results in exactly two maximal elements, namely $(\aelt,\xelt)$ and $(\aelt,\yelt)$,} 

\vspace*{-0.07in}
{\footnotesize while the product on the right results in three maximal elements, which are $(\xelt,\aelt)$, $(\yelt,\aelt)$, and $(\xelt,\celt)$.)}

\vspace*{0.1in}
\setlength{\unitlength}{0.7cm}
\begin{picture}(16,2)
\thicklines
\put(0,0){\circle*{0.22}} \put(1.5,1.5){\circle*{0.22}} \put(3,0){\circle*{0.22}}
\put(-0.45,0){\scriptsize $\belt$} \put(1.65,1.5){\scriptsize $\aelt$} \put(3.15,0){\scriptsize $\celt$}
\put(3.25,0.7){$\otimes$}
\put(4,1.5){\circle*{0.22}} \put(5.5,0){\circle*{0.22}} \put(7,1.5){\circle*{0.22}} 
\put(3.55,1.5){\scriptsize $\xelt$} \put(5.05,0){\scriptsize $\zelt$} \put(7.15,1.5){\scriptsize $\yelt$}
\put(7.5,0.5){\Large $\not\cong$}
\put(8.5,1.5){\circle*{0.22}} \put(10,0){\circle*{0.22}} \put(11.5,1.5){\circle*{0.22}} 
\put(8.05,1.5){\scriptsize $\xelt$} \put(9.55,0){\scriptsize $\zelt$} \put(11.65,1.5){\scriptsize $\yelt$}
\put(11.75,0.7){$\otimes$}
\put(12.5,0){\circle*{0.22}} \put(14,1.5){\circle*{0.22}} \put(15.5,0){\circle*{0.22}}
\put(12.05,0){\scriptsize $\belt$} \put(14.15,1.5){\scriptsize $\aelt$} \put(15.65,0){\scriptsize $\celt$}
%
%
\put(2.9,0.1){\color{Cyan}{\line(-1,1){1.3}}}
\put(2.1,0.6){\color{Cyan}{\em 2}}
\put(5.4,0.1){\color{Cyan}{\line(-1,1){1.3}}}
\put(4.6,0.6){\color{Cyan}{\em 2}}
%
%
\put(0.1,0.1){\color{Red}{\line(1,1){1.3}}}
\put(0.6,0.6){\color{Red}{\em 1}}
\put(5.6,0.1){\color{Green}{\line(1,1){1.3}}}
\put(6.1,0.6){\color{Green}{\em 3}}
%
%
\put(15.4,0.1){\color{Cyan}{\line(-1,1){1.3}}}
\put(14.6,0.6){\color{Cyan}{\em 2}}
\put(9.9,0.1){\color{Cyan}{\line(-1,1){1.3}}}
\put(9.1,0.6){\color{Cyan}{\em 2}}
%
%
\put(12.6,0.1){\color{Red}{\line(1,1){1.3}}}
\put(13.1,0.6){\color{Red}{\em 1}}
\put(10.1,0.1){\color{Green}{\line(1,1){1.3}}}
\put(10.6,0.6){\color{Green}{\em 3}}
\end{picture}

\vspace*{-0.2in}
\end{center}
\end{figure}

There are, however, some circumstances in which it is generally known that the crystal product is commutative.  
In particular, it is shown in \cite{HK2} that $A \otimes B \cong B \otimes A$ for crystal graphs $A$ and $B$.  (See also \CrystalOpsCorollary\ below.) 
So we ask: 

\noindent 
{\bf \CommuteProblem}\ \ Under what circumstances is the crystal product of fibrous posets commutative?\hfill\QED

\noindent 
{\bf \FibrousStructureLemma}\ \ {\sl Let $p \geq 2$, and say $R_1,\ldots,R_p$ are fibrous posets. 
Then for any $(\xelt_1,\ldots,\xelt_p) \in R_1 \otimes \cdots \otimes R_p$ and $i \in I$ we have:}
\begin{eqnarray*}
m_i(\xelt_1,\ldots,\xelt_p) & = & m_i(\xelt_1) + \cdots + m_i(\xelt_p)\\
\delta_i(\xelt_1,\ldots,\xelt_p) & = & \max_{1 \leq q \leq p}\bigg\{-(m_i(\xelt_1) + \cdots + m_i(\xelt_{q-1})) + \delta_i(\xelt_q)\bigg\}\\
\rho_i(\xelt_1,\ldots,\xelt_p) & = & \max_{1 \leq r \leq p}\bigg\{\rho_i(\xelt_r) + (m_i(\xelt_{r+1}) + \cdots + m_i(\xelt_{p}))\bigg\}
\end{eqnarray*}
{\sl Moreover, we have $\wE_{i}(\xelt_1,\ldots,\xelt_p) \not= \theta$ if and only if $\delta_{i}(\xelt_1,\ldots,\xelt_p) > 0$, in which case}
\[\wE_{i}(\xelt_1,\ldots,\xelt_p) = (\xelt_1,\ldots,\wE_{i}(\xelt_q),\ldots,\xelt_p),\]
{\sl where in the above formula for $\delta_{i}(\xelt_1,\ldots,\xelt_p)$, $q$ is the smallest index for which the max occurs.  Similarly, we have $\wF_{i}(\xelt_1,\ldots,\xelt_p) \not= \theta$ if and only if $\rho_{i}(\xelt_1,\ldots,\xelt_p) > 0$, in which case}
\[\wF_{i}(\xelt_1,\ldots,\xelt_p) = (\xelt_1,\ldots,\wF_{i}(\xelt_r),\ldots,\xelt_p),\]
{\sl where in the above formula for $\rho_{i}(\xelt_1,\ldots,\xelt_p)$, $r$ is the largest index for which the max occurs.} 

{\em Proof.} The proof is by induction on $p$. 
When $p=2$, the above formulas for $\wE_i\xxelt$ and $\wF_i\xxelt$ are consistent with \EFLemma, from which the above formulas for $\delta_i\xxelt$ and $\rho_i\xxelt$ are easily derived. 
We analyze $m_{i}\xxelt$ using cases. 
First, suppose that $\rho_i(\xelt_1) > \delta_i(\xelt_2)$. 
Then $\rho_i\xxelt = \rho_i(\xelt_1)-\delta_i(\xelt_2)+\rho_i(\xelt_2)$ and $\delta_i\xxelt = \delta_i(\xelt_1)$, and hence $m_i\xxelt = \rho_i\xxelt - \delta_i\xxelt = m_i(\xelt_1)+m_i(\xelt_2)$. 
Similar computations establish the latter equality in the cases that $\rho_i(\xelt_1) = \delta_i(\xelt_2)$ and $\rho_i(\xelt_1) < \delta_i(\xelt_2)$. 

For the induction hypothesis, suppose that for some $s \geq 2$, the lemma statement holds whenever $p \leq s$. 
Then take $p = s+1$. 
Let $R := R_1 \otimes \cdots \otimes R_{p-1}$, so that $R_1 \otimes \cdots \otimes R_{p-1} \otimes R_p \cong R \otimes R_p$. 
Apply the induction hypothesis to see that $m_i(\xelt_1,\ldots,\xelt_p) = m_i(\xelt_1,\ldots,\xelt_{p-1}) + m_i(\xelt_p) = m_1(\xelt_1) + \cdots + m_i(\xelt_{p-1}) + m_i(\xelt_p)$, as desired. 
Similarly,  
\begin{eqnarray*}
\delta_i(\xelt_1,\ldots,\xelt_p) & = & \max\{\delta_i(\xelt_1,\ldots,\xelt_{p-1}),-m_i(\xelt_1,\ldots,\xelt_{p-1})+\delta_i(\xelt_p)\}\\ 
 & = & \max\bigg\{\max_{1 \leq q \leq p-1}\{-(m_i(\xelt_1) + \cdots + m_i(\xelt_{q-1})) + \delta_i(\xelt_q)\}\, \mbox{\Large ,}\\
 & & \hspace*{2in}-m_i(\xelt_1,\ldots,\xelt_{p-1})+\delta_i(\xelt_p)\bigg\}\\
 & = & \max_{1 \leq q \leq p}\bigg\{-(m_i(\xelt_1) + \cdots + m_i(\xelt_{q-1})) + \delta_i(\xelt_q)\bigg\}
\end{eqnarray*}
Similarly see that $\displaystyle \rho_i(\xelt_1,\ldots,\xelt_p) = \max_{1 \leq r \leq p}\bigg\{\rho_i(\xelt_r) + (m_i(\xelt_{r+1}) + \cdots + m_i(\xelt_{p}))\bigg\}$. 
Assuming that $\delta_i(\xelt_1,\ldots,\xelt_p) > 0$, then by the induction hypothesis we have 
\[\wE_{i}(\xelt_1,\ldots,\xelt_p) = \left\{\begin{array}{cl}
(\wE_{i}(\xelt_1,\ldots,\xelt_{p-1}),\xelt_p) & \mbox{if } \delta_i(\xelt_1,\ldots,\xelt_{p-1}) \geq -m_i(\xelt_1,\ldots,\xelt_{p-1}) + \delta_i(\xelt_p)\\
(\xelt_1,\ldots,\xelt_{p-1},\wE_i(\xelt_p)) & \mbox{if } \delta_i(\xelt_1,\ldots,\xelt_{p-1}) < -m_i(\xelt_1,\ldots,\xelt_{p-1}) + \delta_i(\xelt_p).
\end{array}\right.\] 
In the circumstance that $\delta_i(\xelt_1,\ldots,\xelt_{p-1}) \geq -m_i(\xelt_1,\ldots,\xelt_{p-1}) + \delta_i(\xelt_p)$, then $\wE_{i}(\xelt_1,\ldots,\xelt_p) = (\xelt_1,\ldots,\wE_{i}(\xelt_q),\ldots,\xelt_{p-1},\xelt_p)$ where $q$ is smallest in $\{1,\ldots,p-1\}$ such that $-(m_i(\xelt_1)+\cdots+m_i(\xelt_{q-1}))+\delta_i(\xelt_q) = \delta_i(\xelt_1,\ldots,\xelt_{p-1})$. 
So the above formula for $\wE_{i}(\xelt_1,\ldots,\xelt_p)$ can be abbreviated as $(\xelt_1,\ldots,\wE_{i}(\xelt_q),\ldots,\xelt_{p})$ where $q$ is the smallest index in $\{1,\ldots,p\}$ such that $-(m_i(\xelt_1)+\cdots+m_i(\xelt_{q-1}))+\delta_i(\xelt_q) = \delta_i(\xelt_1,\ldots,\xelt_{p})$. 
One can similarly confirm the formula given in the lemma statement for $\wF_i(\xelt_1,\ldots,\xelt_p)$.\hfill\QED

\noindent 
{\bf \WGFCorollary}\ \ {\sl If $R_{1}$ and $R_{2}$ are fibrous posets with edges colored by a set $I$ that indexes a choice of simple roots for a root system $\Phi$, then} $\WGF(R_{1} \otimes R_{2}) = \WGF(R_{1})\, \WGF(R_{2})$.\hfill\QED 

The following lemma is patterned after the proof Lemma 3.1 of \cite{Stem} and will be helpful with some computations in the proof of \FibrousColoringLemma. 

\noindent 
{\bf \FibrousChainStructureLemma}\ \ {\sl In any $i$-component of $\Rcrystal$, there is a unique $\ttelt$ such that $\rho_{i}(\telt_{1}) = \delta_{i}(\telt_{2})$, in which case the entire $i$-component is}
\[(\telt_{1},\wF_{i}^{\rho_{i}(\telt_{2})}(\telt_{2})) \myarrow{i} \cdots \myarrow{i} \ttelt \myarrow{i} \cdots \myarrow{i} (\wE_{i}^{\delta_{i}(\telt_{1})}(\telt_{1}),\telt_{2}).\]
{\sl  For any $\xxelt$ in this $i$-component, we have $\delta_i\xxelt = \delta_i(\xelt_1) - \delta_i(\telt_2) + \delta_i(\xelt_2)$ and $\rho_i\xxelt = \rho_i(\xelt_1) - \rho_i(\telt_1) + \rho_i(\xelt_2)$. 
In particular, $\xxelt$ is $i$-maximal (respectively, $i$-minimal) in this $i$-component if and only if $\xelt_{1}$ is maximal in $\comp_{i}(\xelt_1)$ (respectively, $\xelt_{2}$ is minimal in $\comp_{i}(\xelt_2)$).} 

{\em Proof.} For $\xxelt \in \Rcrystal$, we analyze the quantity $\rho_i(\xelt_1) - \delta_i(\xelt_2)$ as follows. 
First, observe that if $\xxelt \myarrow{i} \yyelt$ in $\Rcrystal$, then $\rho_i(\xelt_1) - \delta_i(\xelt_2) + 1 = \rho_i(\yelt_1) - \delta_i(\yelt_2)$. 
Now $\rho_i\xxelt = \max\{\rho_i(\xelt_1)+\rho_i(\xelt_2)-\delta_i(\xelt_2),\rho_i(\xelt_2)\}$ by \FibrousStructureLemma. 
So when $\xxelt$ is $i$-minimal, then $\rho_i(\xelt_2) = 0$ and therefore the quantity $\rho_i(\xelt_1)-\delta_i(\xelt_2) \leq 0$. 
Also, $\delta_i\xxelt = \max\{\delta_i(\xelt_1),-\rho_i(\xelt_1)+\delta_i(\xelt_1)+\delta_i(\xelt_2)\}$, cf.\ \FibrousStructureLemma. 
So when $\xxelt$ is $i$-maximal, then $\delta_i(\xelt_1) = 0$ and therefore the quantity $-\rho_i(\xelt_1) + \delta_i(\xelt_2) \leq 0$, i.e.\ $\rho_i(\xelt_1) - \delta_i(\xelt_2) \geq 0$. 
This demonstrates the claimed existence of the special pair $\ttelt$. 
Since $\rho_i(\xelt_1) - \delta_i(\xelt_2) \geq 0$ for any $\xxelt \geq \ttelt$ in this $i$-component, then by \EFLemma, $\wE_{i}$ can be applied to $\ttelt$ precisely $\delta_i(\telt_1)$ times. 
Similarly from \EFLemma\ it follows that $\wF_{i}$ can be applied to $\ttelt$ precisely $\rho_i(\telt_2)$ times. 
So the entire $i$-component is $(\telt_{1},\wF_{i}^{\rho_{i}(\telt_{2})}(\telt_{2})) \myarrow{i} \cdots \myarrow{i} \ttelt \myarrow{i} \cdots \myarrow{i} (\wE_{i}^{\delta_{i}(\telt_{1})}(\telt_{1}),\telt_{2})$, as claimed. 
The last claim of the lemma statement follows directly from the preceding characterization of $\comp_{i}(\xelt_1,\xelt_2)$.

To see that $\delta_i\xxelt = \delta_i(\xelt_1) - \delta_i(\telt_2) + \delta_i(\xelt_2)$, we consider cases. 
If $\xelt_2 = \telt_2$, then the result is obvious. 
If $\xelt_2 \not= \telt_2$, then $\xelt_2 \leq \telt_2$ in $\comp_i(\xelt_2)$ and $\xelt_1 = \telt_1$.  
Now $\delta_i\xxelt$ is the number of steps from the maximal element of $\comp_{i}(\xelt_1)$ down to $\telt_1$ plus the number of steps from $\telt_2$ down to $\xelt_2$ in $\comp_{i}(\xelt_2)$. 
But this is precisely $\delta_i(\xelt_1) - \delta_i(\telt_2) + \delta_i(\xelt_2)$. 
Use similar case analysis to derive the given formula for $\rho_i\xxelt$. 
\hfill\QED

{\bf Toward a more refined view of fibrous posets.}  Now and for the remainder of this section fix $J \subseteq I$ and $\nu \in \Lambda^{+}_{\Phi_J}$.  We will use the pair $(J,\nu)$ to ``refine'' our view of a given fibrous poset $R$.  Note that in this section we will not make any use of $\nu$ as a weight; we only need the fact that it corresponds to a $J$-tuple $(\nu_{j})_{j \in J}$ of nonnegative integers, where $\nu = \sum_{j \in J}\nu_{j}\omega_{j}^{J}$.  The sets we define next along with the following lemma pertain to some arguments we will make about a certain kind of ``$(J,\nu)$ vertex coloring'' of $R$.  Define 
\[\mathcal{M}_{J,\nu}(R) := \{\melt \in R\, |\, \delta_{j}(\melt) \leq \nu_{j} \mbox{ for all } j \in J\}\] 
and for any $\xelt \in R$ set 
\[\mathcal{K}_{J,\nu}(\xelt) := \{j \in J\, |\, \delta_{j}(\xelt) > \nu_{j}\}.\]

{\bf Primary posets.} In the next section we will construct what we call ``crystalline'' splitting posets by taking crystal products of fibrous posets with the following special designation.  A {\em primary} poset is a fibrous poset $R$ with the property that if, for some $\xelt \in R$ and $i \in I$, the $i$-component $\comp_{i}(\xelt)$ is a chain $\xelt_0 \myarrow{i} \xelt_1 \myarrow{i} \ldots \myarrow{i} \xelt_{\mysmallerl}$ with $\myl \geq 2$,  
then for all $j \not= i$ and $0 \leq k \leq \myl-1$ we have $\delta_{j}(\xelt_k) = 0$ if $l_{j}(\xelt_k) \geq 2$. 
In other words, suppose we have an $i$- and a $j$-component of $R$ (with $i \not= j$), suppose both components have length at least two, and suppose they intersect. 
Then the intersection vertex is the maximal element of one of the two components. 

\noindent
{\bf \PrimaryExamples}\ \ 
Examples of primary posets include the minuscule and quasi-minuscule splitting posets of \S \SplittingSection\ (cf.\ \MinQuasiMinTheorem) and the unique connected edge-minimal splitting poset for $\chi_{_{\overline{\omega}}}$ when $\overline{\omega}$ is the highest root for an irreducible root system. 
The unique splitting poset for $\chi_{_{\omega_{3}}}^{\mytinyC_{3}}$ depicted in \PrimaryNotPrimaryFigure\ is a primary poset. 
However, the unique splitting poset for $\chi_{_{3\omega_{1}}}^{\mytinyA_{2}}$ depicted there is not a primary poset. 
Also, the $\myA_2$-structured poset of \NotSchurPositiveFigure\ is not primary. 
See \TangledFigure\ and \FundamentalFigure\ for more examples.\hfill\QED

\begin{figure}[t]
\begin{center}
{\PrimaryNotPrimaryFigure:} The $\myC_{3}$-structured fibrous poset on the left is primary,\\ but the $\myA_{2}$-structured fibrous poset on the right is not.

\vspace*{0.2in}
\setlength{\unitlength}{0.7cm}
\begin{picture}(5.5,13.5)
\thicklines
\multiput(3,3)(-1.5,1.5){3}{\circle*{0.22}}
\multiput(4.5,4.5)(-1.5,1.5){4}{\circle*{0.22}}
\multiput(4.5,7.5)(-1.5,1.5){3}{\circle*{0.22}}
\put(3,1.5){\circle*{0.22}}
\put(3,0){\circle*{0.22}}
\put(1.5,12){\circle*{0.22}}
\put(1.5,13.5){\circle*{0.22}}
%
%
\put(3,0.125){\color{Green}{\line(0,1){1.25}}}
\put(2.9,0.6){\color{Green}{\em 3}}
\put(3,1.625){\color{Cyan}{\line(0,1){1.25}}}
\put(2.9,2.1){\color{Cyan}{\em 2}}
\put(1.5,10.625){\color{Cyan}{\line(0,1){1.25}}}
\put(1.4,11.1){\color{Cyan}{\em 2}}
\put(1.5,12.125){\color{Green}{\line(0,1){1.25}}}
\put(1.4,12.6){\color{Green}{\em 3}}
%
%
\put(2.9,3.1){\color{Cyan}{\line(-1,1){1.3}}}
\put(2.1,3.6){\color{Cyan}{\em 2}}
\put(1.4,4.6){\color{Green}{\line(-1,1){1.3}}}
\put(0.6,5.1){\color{Green}{\em 3}}
\put(4.4,4.6){\color{Cyan}{\line(-1,1){1.3}}}
\put(3.6,5.1){\color{Cyan}{\em 2}}
\put(2.9,6.1){\color{Green}{\line(-1,1){1.3}}}
\put(2.1,6.6){\color{Green}{\em 3}}
\put(1.4,7.6){\color{Cyan}{\line(-1,1){1.3}}}
\put(0.6,8.1){\color{Cyan}{\em 2}}
\put(4.4,7.6){\color{Green}{\line(-1,1){1.3}}}
\put(3.6,8.1){\color{Green}{\em 3}}
\put(2.9,9.1){\color{Cyan}{\line(-1,1){1.3}}}
\put(2.1,9.6){\color{Cyan}{\em 2}}
%
%
\put(3.1,3.1){\color{Red}{\line(1,1){1.3}}}
\put(3.6,3.6){\color{Red}{\em 1}}
\put(1.6,4.6){\color{Red}{\line(1,1){1.3}}}
\put(2.1,5.1){\color{Red}{\em 1}}
\put(3.1,6.1){\color{Red}{\line(1,1){1.3}}}
\put(3.6,6.6){\color{Red}{\em 1}}
\put(0.1,6.1){\color{Red}{\line(1,1){1.3}}}
\put(0.6,6.6){\color{Red}{\em 1}}
\put(1.6,7.6){\color{Red}{\line(1,1){1.3}}}
\put(2.1,8.1){\color{Red}{\em 1}}
\put(0.1,9.1){\color{Red}{\line(1,1){1.3}}}
\put(0.6,9.6){\color{Red}{\em 1}}
\end{picture}
\hspace*{1in}
\begin{picture}(5.5,9.5)
\thicklines
\multiput(4.5,2)(-1.5,1.5){4}{\circle*{0.22}}
\multiput(4.5,5)(-1.5,1.5){3}{\circle*{0.22}}
\multiput(4.5,8)(-1.5,1.5){2}{\circle*{0.22}}
\multiput(4.5,11)(-1.5,1.5){1}{\circle*{0.22}}
%
%
\put(4.4,2,1){\color{Cyan}{\line(-1,1){1.3}}}
\put(3.6,2.6){\color{Cyan}{\em 2}}
\put(2.9,3.6){\color{Cyan}{\line(-1,1){1.3}}}
\put(2.1,4.1){\color{Cyan}{\em 2}}
\put(1.4,5.1){\color{Cyan}{\line(-1,1){1.3}}}
\put(0.6,5.6){\color{Cyan}{\em 2}}
\put(4.4,5.1){\color{Cyan}{\line(-1,1){1.3}}}
\put(3.6,5.6){\color{Cyan}{\em 2}}
\put(2.9,6.6){\color{Cyan}{\line(-1,1){1.3}}}
\put(2.1,7.1){\color{Cyan}{\em 2}}
\put(4.4,8.1){\color{Cyan}{\line(-1,1){1.3}}}
\put(3.6,8.6){\color{Cyan}{\em 2}}
%
%
\put(3.1,3.6){\color{Red}{\line(1,1){1.3}}}
\put(3.6,4.1){\color{Red}{\em 1}}
\put(1.6,5.1){\color{Red}{\line(1,1){1.3}}}
\put(2.1,5.6){\color{Red}{\em 1}}
\put(3.1,6.6){\color{Red}{\line(1,1){1.3}}}
\put(3.6,7.1){\color{Red}{\em 1}}
\put(0.1,6.6){\color{Red}{\line(1,1){1.3}}}
\put(0.6,7.1){\color{Red}{\em 1}}
\put(1.6,8.1){\color{Red}{\line(1,1){1.3}}}
\put(2.1,8.6){\color{Red}{\em 1}}
\put(3.1,9.6){\color{Red}{\line(1,1){1.3}}}
\put(3.6,10.1){\color{Red}{\em 1}}
\end{picture}
\end{center}
\end{figure}

The trichotomy observed in the next lemma is needed for the statement of \FibrousColoringLemma. 

\noindent 
{\bf \TidyProductLemma}\ \ {\sl Let $R$ be any disjoint sum of connected components of $\Rcrystal$, where $R_{2}$ is a primary poset.  
Say $\xxelt \in R \setminus \mathcal{M}_{J,\nu}(R)$, so $\mathcal{K}_{J,\nu}\xxelt \not= \emptyset$.  
Then exactly one of the following holds: 
(1) $\xelt_{1} \in R_1 \setminus \mathcal{M}_{J,\nu}(R_{1})$, 
(2) $\xelt_{1} \in \mathcal{M}_{J,\nu}(R_{1})$ and there is a unique color $j \in \mathcal{K}_{J,\nu}\xxelt$ such that 
$\max\{\rho_j(\xelt_2),l_j(\xelt_2)-\rho_j(\xelt_1)\} \geq 2$, 
or (3) $\xelt_{1} \in \mathcal{M}_{J,\nu}(R_{1})$ and for all $k \in \mathcal{K}_{J,\nu}\xxelt$ it is that case that 
$\max\{\rho_k(\xelt_2),l_k(\xelt_2)-\rho_k(\xelt_1)\} \leq 1$.} 

{\em Proof.}  Clearly the three conditions are mutually exclusive. 
They are exhaustive if we can show uniqueness of the color $j$ from condition (2). 
So suppose that $\xelt_{1}$ is in $\mathcal{M}_{J,\nu}(R_{1})$ and that for some $j \in \mathcal{K}_{J,\nu}\xxelt$ it is the case that 
$\max\{\rho_j(\xelt_2),l_j(\xelt_2)-\rho_j(\xelt_1)\} \geq 2$. 
Say $i \in \mathcal{K}_{J,\nu}\xxelt$ also has the property that $\max\{\rho_i(\xelt_2),l_i(\xelt_2)-\rho_i(\xelt_1)\} \geq 2$. 
Now $l_{j}(\xelt_2) = \rho_j(\xelt_2)+\delta_j(\xelt_2) \geq \max\{\rho_j(\xelt_2),l_j(\xelt_2)-\rho_j(\xelt_1)\} \geq 2$, and similarly $l_{i}(\xelt_2) \geq 2$. 
Since $i \in \mathcal{K}_{J,\nu}\xxelt$, then $\delta_{i}\xxelt = \max\{\delta_{i}(\xelt_1),-\rho_{i}(\xelt_1)+\delta_{i}(\xelt_1)+\delta_{i}(\xelt_2)\} > \nu_{i}$.  
Since $\xelt_{1}$ is in $\mathcal{M}_{J,\nu}(R_{1})$, then $\delta_{i}(\xelt_1) \leq \nu_{i}$. 
These latter two inequalities together require that $\delta_{i}(\xelt_2) > 0$, so $\xelt_2$ cannot be $i$-maximal. 
Similarly, $\xelt_2$ cannot be $j$-maximal. 
So $\xelt_2$ is neither $i$-maximal nor $j$-maximal within its $i$- and $j$-components, each of which has length at least two. 
Since $R_2$ is primary, we must therefore have $i=j$.\hfill\QED

{\bf Vertex coloring.} Much of our approach to splitting posets here and elsewhere will depend crucially on the problem of producing certain kinds of vertex-coloring functions.  Such functions allow for the cancelling of terms in computations related to $W$-symmetric functions as in \InitialSplittingTheorem.  For fibrous posets, we develop this idea as follows. 

If $R$ is fibrous, $\xelt \in R$, $i \in I$, and $p$ is a nonnegative integer, then let $U_{i}(\xelt,p)$ be the set of the topmost $1+p$ elements in the chain $\comp_{i}(\xelt)$, with the understanding that $U_{i}(\xelt,p) = \comp_{i}(\xelt)$ if $1+p > |\comp_{i}(\xelt)|$.  
A $(J,\nu)$-{\em coloring} of the vertices of $R$ is a function $\kappa: R \setminus \mathcal{M}_{J,\nu}(R) \longrightarrow J$ such that for all $\xelt \in  R \setminus \mathcal{M}_{J,\nu}(R)$, 
\[\{\yelt \in \comp_{j}(\xelt)\, |\, \yelt \in R \setminus \mathcal{M}_{J,\nu}(R) \mbox{ and } \kappa(\yelt) = j\} = \comp_{j}(\xelt) \setminus U_{j}(\xelt,\nu_j),\] where $j = \kappa(\xelt)$. 
In this case, note that $j \in \mathcal{K}_{J,\nu}(\xelt)$. 

When $R$ is primary, it is easy to identify a $(J,\nu)$-coloring of $R$. 
Note that the definition of primary poset requires that for any $\xelt \in R \setminus \mathcal{M}_{J,\nu}(R)$, the quantity $l_j(\xelt)$ can exceed one for at most one color $j \in \mathcal{K}_{J,\nu}(\xelt)$. 

\noindent 
{\bf \TidyColoringLemma}\ \ {\sl Let $R$ be primary.  
For any $\xelt \in R \setminus \mathcal{M}_{J,\nu}(R)$, let $\kappa(\xelt)$ be determined by}
\[\kappa(\xelt) = \left\{
\begin{array}{cl}
j & \rule[-7mm]{0mm}{12mm}\hspace*{0.25in}\parbox{4in}{if $j$ is the (necessarily unique) color in $\mathcal{K}_{J,\nu}(\xelt)$ such that 
$l_j(\xelt) \geq 2$}\\
k & \rule[-6mm]{0mm}{4mm}\hspace*{0.25in}\parbox{4in}{if for all $i \in \mathcal{K}_{J,\nu}(\xelt)$ it is that case that 
$l_i(\xelt) \leq 1$, in which case $k$ is a color freely chosen from $\mathcal{K}_{J,\nu}(\xelt)$} 
\end{array}
\right.\] 
{\sl Then $\kappa: R \setminus \mathcal{M}_{J,\nu}(R) \longrightarrow J$ is a $(J,\nu)$-coloring of $R$.} 

{\em Proof.} Let $\xelt \in R \setminus \mathcal{M}_{J,\nu}(R)$ and $j := \kappa(\xelt)$. 
Say $\yelt \in \comp_j(\xelt)$ with $\yelt \in R \setminus \mathcal{M}_{J,\nu}(R)$ and $\kappa(\yelt) = j$. 
Now $j \in \mathcal{K}_{J,\nu}(\yelt)$ means that $\delta_j(\yelt) > \nu_j$. 
Hence $\yelt \in \comp_j(\xelt) \setminus U_j(\xelt,\nu_j)$. 

Now suppose that $\yelt \in \comp_j(\xelt)$ with $\delta_j(\yelt) > \nu_j$. 
Then $\yelt \in R \setminus \mathcal{M}_{J,\nu}(R)$ and $j \in \mathcal{K}_{J,\nu}(\yelt)$. 
First, suppose $l_j(\xelt) \leq 1$. 
Since $j \in \mathcal{K}_{J,\nu}(\xelt)$, then $\delta_j(\xelt) > \nu_j \geq 0$. 
In particular, neither $\xelt$ nor $\yelt$ is the maximum element of the two-element chain $\comp_j(\xelt)$. 
Then we must have $\yelt = \xelt$ as the minimum element of this two-element chain. 
Therefore $\kappa(\yelt) = \kappa(\xelt) = j$. 
Second, suppose $l_j(\xelt) \geq 2$.  
Then $l_j(\yelt) \geq 2$ as well. 
Since $j \in \mathcal{K}_{J,\nu}(\yelt)$, then by the definition of $\kappa$ we must have $\kappa(\yelt) = j$.\hfill\QED

Our next result, which is the main result of this section, shows how we can obtain $(J,\nu)$-colorings by taking crystal products of primary posets. 

\noindent 
{\bf \FibrousColoringLemma}\ \ {\sl Let $R$ be any disjoint sum of connected components of $\Rcrystal$, where $R_{2}$ is a primary poset. Suppose $\kappa^{(1)}: R_{1} \setminus \mathcal{M}_{J,\nu}(R_{1}) \longrightarrow J$ is a $(J,\nu)$-coloring function on $R_{1}$.  Using the trichotomy of conditions from \TidyProductLemma, for any $\xxelt \in R \setminus \mathcal{M}_{J,\nu}(R)$ declare that} 
\[\kappa\xxelt = \left\{
\begin{array}{cl}
\kappa^{(1)}(\xelt_{1}) & \hspace*{0.25in}\mbox{if } \xelt_{1} \not\in \mathcal{M}_{J,\nu}(R_{1})\\
j & \rule[-3mm]{0mm}{12mm}\hspace*{0.25in}\parbox{4in}{if $\xelt_{1} \in \mathcal{M}_{J,\nu}(R_{1})$ and $j$ is the unique color in $\mathcal{K}_{J,\nu}\xxelt$ such that 
$\max\{\rho_j(\xelt_2),l_j(\xelt_2)-\rho_j(\xelt_1)\} \geq 2$}\\
k & \rule[-3mm]{0mm}{16mm}\hspace*{0.25in}\parbox{4in}{if $\xelt_{1} \in \mathcal{M}_{J,\nu}(R_{1})$ and for all $i \in \mathcal{K}_{J,\nu}\xxelt$ it is that case that 
$\max\{\rho_i(\xelt_2),l_i(\xelt_2)-\rho_i(\xelt_1)\} \leq 1$, in which case $k$ is a color freely chosen from $\mathcal{K}_{J,\nu}\xxelt$} 
\end{array}
\right.\]
{\sl Then  $\kappa: R \setminus \mathcal{M}_{J,\nu}(R) \longrightarrow J$ is a $(J,\nu)$-coloring of $R$.} 

{\em Proof.} Take $\xxelt \in R \setminus \mathcal{M}_{J,\nu}(R)$ with $j := \kappa\xxelt$,  and let $\ttelt$ be the special pair in $\comp_j\xxelt$ for which $\rho_j(\telt_1) = \delta_j(\telt_2)$ as identified in \FibrousChainStructureLemma. 
We must show that 
\[\{\yyelt \in \comp_{j}\xxelt\, |\, \yyelt \in R \setminus \mathcal{M}_{J,\nu}(R) \mbox{ and } \kappa\yyelt = j\}\hspace*{1in}\]

\vspace*{-0.35in}
\[\hspace*{2in} = \comp_{j}\xxelt \setminus U_{j}\left(\frac{}{}\xxelt,\nu_j\right).\]

For the first part of the proof we let $\yyelt \in \comp_{j}\xxelt$ with $\yyelt \in R \setminus \mathcal{M}_{J,\nu}(R)$ and $\kappa\yyelt = j$, and we show that $\delta_{j}\yyelt > \nu_j$. 
We consider two cases: (1) $\xelt_1 \in R_1 \setminus \mathcal{M}_{J,\nu}(R_1)$, and (2) $\xelt_1 \in \mathcal{M}_{J,\nu}(R_1)$. 
For case (1), $\kappa\xxelt = \kappa^{(1)}(\xelt_1) = j$, and since $\kappa^{(1)}$ is a $(J,\nu)$-coloring of $R_1$, it follows that $\nu_j < \delta_j(\xelt_1)$. 
Now $\delta_j(\xelt_1) \leq \delta_j(\telt_1)$. 
If $\yelt_1 \in \mathcal{M}_{J,\nu}(R_1)$, then $\delta_j(\yelt_1) \leq \nu_j$, so we must have $\yelt_1 \not= \telt_1$, and therefore $\yelt_2 = \telt_2$. 
Then $\delta_{j}\yyelt = \delta_{j}(\yelt_1) - \delta_{j}(\telt_2) + \delta_{j}(\yelt_2) = \delta_{j}(\yelt_1) \leq \nu_j$. 
But $\yelt_1 \in \mathcal{M}_{J,\nu}(R_1)$ also means $\kappa\yyelt = j$ for some $j \in \mathcal{K}_{J,\nu}\yyelt$, in which case $\delta_{j}\yyelt > \nu_j$. 
So we must have $\yelt_1 \in R_1 \setminus \mathcal{M}_{J,\nu}(R_1)$. 
Then $\kappa\yyelt = \kappa^{(1)}(\yelt_1) = j$, and since $\kappa^{(1)}$ is a $(J,\nu)$-coloring of $R_1$, it follows that $\nu_j < \delta_j(\yelt_1)$. 
So $\delta_{j}\yyelt = \delta_{j}(\yelt_1) - \delta_{j}(\telt_2) + \delta_{j}(\yelt_2) \geq \delta_{j}(\yelt_1) > \nu_j$. 
For case (2), it follows from the definition of $\kappa$ that $j \in \mathcal{K}_{J,\nu}\xxelt$, so that $\delta_{j}\xxelt = \delta_{j}(\xelt_1) - \delta_{j}(\telt_2) + \delta_{j}(\xelt_2) > \nu_j$. 
However, $\delta_{j}(\xelt_1) \leq \nu_j$, so it must be the case that $\xelt_2 \not= \telt_2$, and hence $\xelt_1 = \telt_1$. 
So, $\delta_{j}(\yelt_1) \leq \delta_{j}(\telt_1) = \delta_{j}(\xelt_1) \leq \nu_j$. 
Then it must be the case that $\yelt_1 \in \mathcal{M}_{J,\nu}(R_1)$. 
Then by the definition of $\kappa$, the fact that $\kappa\yyelt = j$ means that we must have $j \in \mathcal{K}_{J,\nu}\yyelt$, hence $\delta_{j}\yyelt > \nu_j$. 

In the other direction, we let $\yyelt \in \comp_{j}\xxelt$ with $\delta_{j}\yyelt > \nu_j$, and we show that $\kappa\yyelt = j$. 
We consider the same two cases as before: (1) $\xelt_1 \in R_1 \setminus \mathcal{M}_{J,\nu}(R_1)$, and (2) $\xelt_1 \in \mathcal{M}_{J,\nu}(R_1)$. 
As with case (1) from the previous paragraph, $\kappa\xxelt = \kappa^{(1)}(\xelt_1) = j$, and since $\kappa^{(1)}$ is a $(J,\nu)$-coloring of $R_1$, it follows that $\nu_j < \delta_j(\xelt_1) \leq \delta_j(\telt_1)$. 
If $\delta_j(\yelt_1) \leq \nu_j$, then we must have $\yelt_1 \not= \telt_1$, and therefore $\yelt_2 = \telt_2$. 
Then $\delta_{j}\yyelt = \delta_{j}(\yelt_1) - \delta_{j}(\telt_2) + \delta_{j}(\yelt_2) = \delta_{j}(\yelt_1) \leq \nu_j$, which contradicts the fact that $\delta_{j}\yyelt > \nu_j$. 
So in fact we must have $\delta_{j}(\yelt_1) > \nu_j$. 
Since $\yelt_1 \in \comp_{j}(\xelt_1) \setminus U_{j}(\xelt_1,\nu_j)$ and $\kappa^{(1)}$ is a $(J,\nu)$-coloring, then $\kappa^{(1)}(\yelt_1) = j$. 
And since $\yelt_1 \in R_1 \setminus \mathcal{M}_{J,\nu}(R_1)$, then $\kappa\yyelt = \kappa^{(1)}(\yelt_1) = j$. 
For case (2), we reason as before that $j \in \mathcal{K}_{J,\nu}\xxelt$, so that $\delta_{j}\xxelt = \delta_{j}(\xelt_1) - \delta_{j}(\telt_2) + \delta_{j}(\xelt_2) > \nu_j$. 
However, $\delta_{j}(\xelt_1) \leq \nu_j$, so it must be the case that $\xelt_2 \not= \telt_2$, and hence $\xelt_1 = \telt_1$. 
So, $\delta_{j}(\yelt_1) \leq \delta_{j}(\telt_1) = \delta_{j}(\xelt_1) \leq \nu_j$. 
But now $\delta_{j}\yyelt = \delta_{j}(\yelt_1) - \delta_{j}(\telt_2) + \delta_{j}(\yelt_2) > \nu_j$ means that we must have $\yelt_2 \not= \telt_2$, hence $\yelt_1 = \telt_1 = \xelt_1$. 
If $\yelt_2 = \xelt_2$ we are done, as we will have $\kappa\yyelt = \kappa\xxelt = j$. 
So suppose $\yelt_2 \not= \xelt_2$. 
Then $\xelt_2, \yelt_2, \telt_2$ are three distinct elements of $\comp_{j}(\xelt_2)$ in $R_2$ with $\xelt_2 < \telt_2$ and $\yelt_2 < \telt_2$. 
So $\rho_{j}(\telt_2) \geq 2$. 
Since $\yyelt$ is below $\ttelt$ in $\comp_{j}\xxelt$, then one can see from the proof of \FibrousChainStructureLemma\ that $\delta_j(\yelt_2)-\rho_j(\yelt_1)$ is the number of steps from $\yyelt$ up to $\ttelt$ in this chain. 
So, $\rho_j(\telt_2) = \rho_j(\yelt_2)+\delta_j(\yelt_2)-\rho_j(\yelt_1)$, and 
therefore $\max\{\rho_{j}(\yelt_2),l_j(\yelt_2)-\rho_{j}(\yelt_1)\} \geq 2$. 
It now follows from the definition of $\kappa$ that $\kappa\yyelt = j$.\hfill\QED

\newpage
\begin{center}
{\FundamentalFigure:} Some primary posets. 

{\footnotesize These posets are in fact crystal graphs in the sense of \S \CrystalSection.} 

\vspace*{-0.075in}
{\footnotesize In each picture, the pair of integers $(a,b)$ beside a vertex $\xelt$ indicates the weight $wt(\xelt) = a\omega_{1}+b\omega_{2}$.}

\setlength{\unitlength}{1.1cm}
\begin{picture}(14,18.0)
\put(0,0){\line(0,1){17.5}}
\put(8,0){\line(0,1){17.5}}
\put(14,0){\line(0,1){17.5}}
\put(0,0){\line(1,0){14}}
\put(0,13){\line(1,0){8}}
\put(8,12){\line(1,0){6}}
\put(0,17.5){\line(1,0){14}}
\put(1,16.5){\fbox{\Large $\myA_{1} \oplus \myA_{1}$}}
\put(0.5,13){\AoneAlphaWeights}
\put(3.5,13){\AoneBetaWeights}
\put(8.5,16.5){\fbox{\Large $\myA_{2}$}}
\put(8,12){\AtwoAlphaWeights}
\put(11,12){\AtwoBetaWeights}
\put(8.5,10.5){\fbox{\Large $\myC_{2}$}}
\put(8,4){\BtwoBetaWeights}
\put(11,3.5){\BtwoAlphaWeights}
\put(1,11.5){\fbox{\Large $\myG_{2}$}}
\put(0,2){\GtwoAlphaWeights}
\put(2,0){\GtwoBetaWeights}
\end{picture}
\end{center} 

\newpage
\noindent 
{\Large \bf \S \CrystallineSection.\ Crystalline splitting posets.}  

In this section, we use crystal products to produce a fibrous (and therefore edge-minimal) splitting poset $R(\lambda)$ for each $\chi_{_{\lambda}}$.   
These ``crystalline'' splitting posets will be shown to 
provide answers to the product decomposition and branching problems.  In the next section, we will see that these $R(\lambda)$'s are in fact the celebrated crystal graphs whose discovery was pioneered by Kashiwara and others.  
Our approach is inspired by Stembridge's paper \cite{Stem} in which he developed the notion of an ``admissible system'' as a kind of combinatorial and root system theoretic analog of the crystal graphs that are obtained from the representation theory of quantum groups.  

However, some differences with Stembridge's methodology afford some advantages to our approach.  
First, as we will see in \AdmissibleTheorem, data for a ``coherent timing pattern'' for an admissible system essentially prescribes, all at once, a $(J,\nu)$-coloring function for each pair $(J,\nu)$, whereas in some contexts it is possible within our framework to prescribe a $(J,\nu)$-coloring only for a particular pair or pairs $(J,\nu)$. 
Second, we are able to generalize our approach to non-fibrous posets using \InitialSplittingTheorem\, see e.g.\  \cite{DD} or Sections \CriteriaSection/\GaussianSection. 
And third, while the coloring functions found here are mostly obtained by repeated application of \FibrousColoringLemma\ to certain crystal products, in some circumstances we can explicitly identify vertex coloring functions and so obtain splitting results without any iterative machinery (again, \cite{DD} or Sections \CriteriaSection/\GaussianSection). 

The main theme of this section is to use crystal products of primary posets to construct refined splitting posets. 
This is stated as \MainColoringCorollary\ and is effected by \MainColoringTheorem\ together with \FibrousColoringLemma. 
In particular, \MainColoringTheorem\ applies \InitialSplittingTheorem\ to show how the $(J,\nu)$-splitting property is implied by a $(J,\nu)$-coloring of a fibrous poset. 
In \AdmissibleTheorem, we use \MainColoringTheorem\ to show how admissible systems are instances of crystalline splitting posets.

{\bf ${\Phi}$-structured fibrous posets.}  We begin by studying how some of the results of the previous section are influenced by the ${\Phi}$-structure requirement. 
The first result follows from  \FibrousStructureLemma. 

\noindent 
{\bf \CrystalProductStructureLemma}\ \ {\sl If $R_{1}$ and $R_{2}$ are ${\Phi}$-structured and fibrous posets, then $\Rcrystal$ is ${\Phi}$-structured as well.}\hfill\QED 

The next simple lemma is needed for the statement and proof of \MainColoringTheorem. 

\noindent
{\bf \FibrousInvariantLemma}\ \ {\sl Suppose $R$ is ${\Phi}$-structured and fibrous.  Let $J \subseteq I$ and $\nu \in \Lambda_{\Phi_J}^{+}$.  Then} $\WGF(R)|_{J}$ {\sl is $W_J$-invariant.  Also, if 
$\melt \in \mathcal{M}_{J,\nu}(R)$, then $\nu + wt^{J}(\melt) \in \Lambda_{\Phi_J}^{+}$, i.e.\ is $\Phi_{J}$-dominant.} 

{\em Proof.}  The first claim follows from \WInvariantLemma: since $R$ is fibrous, its $j$-components are necessarily rank symmetric for all $j \in J$, and since $R$ is ${\Phi}$-structured it is ${\Phi_J}$-structured as well. 
The last claim of the lemma is an easy calculation:   
For each $j \in J$, we have $\nu_{j}+m_{j}(\melt) = \nu_{j}+\rho_{j}(\melt)-\delta_{j}(\melt) \geq \nu_{j}-\delta_{j}(\melt) \geq 0$, the latter inequality by virtue of the fact that $\melt \in \mathcal{M}_{J,\nu}(R)$.\hfill\QED

\noindent 
{\bf \MainColoringTheorem}\ \ {\sl Suppose $R$ is ${\Phi}$-structured and fibrous.  Let $J \subseteq I$ and $\nu \in \Lambda_{\Phi_J}^{+}$, and suppose $R$ has a $(J,\nu)$-coloring function $\kappa: R \setminus \mathcal{M}_{J,\nu}(R) \longrightarrow J$.  Then $R$ is a $(J,\nu)$-splitting poset with}
\[\chi^{\Phi_J}_{_{\nu}} \cdot \WGF(R)|_{J} = \sum_{\melt \in \mathcal{M}_{J,\nu}(R)}\chi^{\Phi_J}_{_{\nu + wt^{J}(\melt)}}.\] 

{\em Proof.} 
To apply \InitialSplittingTheorem\ here, we identify a bijection $\tau: R \setminus \mathcal{M}_{J,\nu}(R) \longrightarrow R \setminus \mathcal{M}_{J,\nu}(R)$ that is compatible with the given vertex-coloring function $\kappa$. 
For any $\xelt \in R \setminus \mathcal{M}_{J,\nu}(R)$, define $\tau(\xelt)$ as follows: Let $j := \kappa(\xelt)$.  
Then let $\tau(\xelt)$ be the unique element of $\comp_{j}(\xelt)$ for which $\rho_j(\tau(\xelt)) = l_j(\xelt) - 1 - \nu_j - \rho_j(\xelt)$.  
Since $\kappa$ is a $(J,\nu)$-coloring function, then $\tau(\xelt) \in R \setminus \mathcal{M}_{J,\nu}(R)$. 
It is now easy to see that $\tau$ is an involution on $R \setminus \mathcal{M}_{J,\nu}(R)$. 
Since $R$ is ${\Phi}$-structured, then $wt^{J}(\tau(\xelt)) = wt^{J}(\xelt) - (1 + \nu_{\kappa(\xelt)} + m_{\kappa(\xelt)})\alpha_{\kappa(\xelt)}$. 
By \FibrousInvariantLemma, $\WGF(R)|_{J}$ is $W_J$-invariant, and if $\melt \in \mathcal{M}_{J,\nu}(R)$, then $\nu + wt^{J}(\melt) \in \Lambda_{\Phi_J}^{+}$. 
So we have met all the requirements of \InitialSplittingTheorem\ with $\mathcal{S}_{J,\nu}(R) = \mathcal{M}_{J,\nu}(R)$, so $R$ is a $(J,\nu)$-splitting poset with 
$\displaystyle \chi^{\Phi_J}_{_{\nu}} \cdot \WGF(R)|_{J} = \sum_{\melt \in \mathcal{M}_{J,\nu}(R)}\chi^{\Phi_J}_{_{\nu + wt^{J}(\melt)}}$.\hfill\QED

{\bf Crystalline splitting posets.} 
The preceding theorem is part of the inspiration for the following definition.  

\noindent 
{\bf \CrystallineDefinition}\ \ Say an ${\Phi}$-structured and fibrous poset $R$ is a {\em crystalline splitting poset} if for each $J \subseteq I$ and $\nu \in \Lambda_{\Phi_J}^{+}$, there is some subset $\mathcal{S}_{J,\nu}(R)$ of $R$ such that $\nu + wt^{J}(\selt) \in \Lambda_{\Phi_J}^{+}$ whenever $\selt \in \mathcal{S}_{J,\nu}(R)$ and 
\[\chi^{\Phi_J}_{_{\nu}} \cdot \WGF(R)|_{J} = \sum_{\selt \in \mathcal{S}_{J,\nu}(R)}\chi^{\Phi_J}_{_{\nu + wt^{J}(\selt)}},\] 
i.e.\ $R$ is a $(J,\nu)$-splitting poset.\hfill\QED 

Of course, a one-element poset $P := \{\melt\}$ trivially meets the conditions of the preceding definition, where we take each $\mathcal{S}_{J,\nu}(P)$ to be $\{\melt\}$. 
Since any crystalline splitting poset $R$ is fibrous, then $\chi := \WGF(R)$ is $W$-invariant. 
For such a Weyl symmetric function, we can regard the product decomposition and branching problems as solved, because we have $\displaystyle \chi^{\Phi_J}_{_{\nu}} \ \chi|_{_J} = \sum_{\selt \in \mathcal{S}_{J,\nu}(R)}\chi^{\Phi_J}_{_{\nu + wt^{J}(\selt)}}$ for each pair $(J, \nu)$. 

For each Weyl bialternant $\chi_{_\lambda}$, our goal is to identify a crystalline splitting poset $R(\lambda)$ such that $\WGF(R(\lambda)) = \chi_{_\lambda}$. 
This is accomplished in \CrystallineExistenceTheorem. 
The milieu for our approach is the following theorem, which we regard as the central result of \S \FibrousSection /\S \CrystallineSection. 
Notice that the hypotheses of \MainColoringCorollary\ do not invoke vertex coloring and that the end result is simply a statement about crystal products of certain kinds of fibrous posets.  

\noindent 
{\bf \MainColoringCorollary}\ \ {\sl Suppose $R$ is a disjoint sum of connected components of the crystal product $R_1 \otimes \cdots \otimes R_p$ of ${\Phi}$-structured and primary posets $R_{i}$ ($1 \leq i \leq p$).  Then $R$ is a crystalline splitting poset and for each $J \subseteq I$ and $\nu \in \Lambda_{\Phi_J}^{+}$ we have}
\[\chi^{\Phi_J}_{_{\nu}} \cdot \WGF(R)|_{J} = \sum_{\melt \in \mathcal{M}_{J,\nu}(R)}\chi^{\Phi_J}_{_{\nu + wt^{J}(\melt)}}\] 

{\em Proof.} Combine \TidyColoringLemma/\FibrousColoringLemma, \CrystalProductStructureLemmas, and  \MainColoringTheorem.\hfill\QED

Primary posets that are also ${\Phi}$-structured include, but certainly are not limited to, minuscule and quasi-minuscule splitting posets. 

\noindent 
{\bf \PrimaryProblem}\ \ Classify the connected, ${\Phi}$-structured, and primary posets. 
By the preceding theorem, of course, these are crystalline splitting posets.\hfill\QED 

\noindent 
{\bf \ProductProblem}\ \ Find explicit descriptions of $(J,\nu)$-coloring functions for various families of crystal products of $(J,\nu)$-colored and ${\Phi}$-structured posets. 
The motivation is that via \MainColoringTheorem, one would obtain more crystalline splitting posets, just as in \MainColoringCorollary. 
We do this in \S \GaussianSection\ for repeated crystal products of a minuscule splitting poset, and the result is another proof of the Proctor--Stanley theorem that minuscule posets are Gaussian (see \cite{PrEur}).\hfill\QED  

{\bf Constructing crystalline splitting posets.} 
One question that remains is how to use crystal products of ${\Phi}$-structured and primary posets to find, for each $\lambda \in \Lambda^{+}$, a crystalline splitting poset $R(\lambda)$ with $\WGF(R(\lambda)) = \chi_{_{\lambda}}$.  
The process we use here of taking crystal products of some especially nice graphs in order to obtain the desired crystal borrows heavily from both \cite{Stem} and \cite{Jan}. 

Initially, we will require $\Phi$ to be irreducible. 
Toward the end of this section, we account for the case that $\Phi$ is reducible. 
The main reason for the irreducible hypothesis at the outset is that, for reasons of convenience, our discussion of minuscule and quasi-minuscule splitting posets was limited to the circumstance that $\Phi$ is irreducible. 
In view of \MinQuasiMinTheorem\ and the definition of primary poset in \S \FibrousSection, we will refer to minuscule and quasi-minuscule splitting posets as {\em primary-plus} posets. 

The $\myA_2$- and $\myA_1 \oplus \myA_1$-splitting posets of \TangledFigure\ have the somewhat unfavorable property that they are fibrous and connected but with more than one maximal element.  In some sense, these splitting posets are ``tangled.''  Following \S 6 of \cite{Stem}, we say a fibrous poset $R$ is {\em strongly untangled} if for every $J \subseteq I$ and every $\xelt \in R$, then $\comp_{J}(\xelt)$ has precisely one maximal element.  Note that for fibrous posets, an element is maximal if and only if it is prominent. 

\noindent 
{\bf \StembridgeUntangledTheorem\ (Stembridge)}\ \ {\sl Let $\Phi$ be irreducible. 
Then any crystal product of primary-plus posets is strongly untangled.} 

{\em Proof.} This follows from Theorem 6.4 of \cite{Stem}.  
While the statement of that theorem regards the minuscule and quasi-minuscule splitting posets that are used in the crystal product to be ``admissible systems,'' the ``coherent timing pattern'' plays no role in the proof.  
The proof depends only on the ${\Phi}$-structure and fibrous properties of these posets as well as an understanding of the weights of the poset elements.\hfill\QED

The following corollary is automatic, cf.\ Remark 6.7.(a) of \cite{Stem}. 

\noindent 
{\bf \StembridgeUntangledCorollary}\ \ {\sl Let $\Phi$ be irreducible.  
For each $q \in \{1,\ldots,p\}$, suppose that $R_q$ is a disjoint sum of connected components from some crystal product of primary-plus posets.  
Then $R_1 \otimes \cdots \otimes R_p$ is strongly untangled.}\hfill\QED 

\noindent 
{\bf \ExistenceCorollary}\ \ {\sl Let $\Phi$ be irreducible. 
Suppose that the crystal product $R_1 \otimes \cdots \otimes R_p$ from \StembridgeUntangledCorollary\ has a maximal element whose (necessarily dominant) weight is $\lambda$. 
Let $R$ be the connected component containing this maximal element. 
Taking $\mathcal{S}_{J,\nu}(R) := \mathcal{M}_{J,\nu}(R)$ for each $J \subseteq I$ and $\nu \in \Lambda_{\Phi_J}^{+}$, then $R$ is a strongly untangled crystalline splitting poset with} $\WGF(R) = \chi_{_{\lambda}}$. 

{\em Proof.} By \MainColoringCorollary, $R$ is a crystalline splitting poset.  
Since $R$ is strongly untangled by \StembridgeUntangledCorollary, its maximal element is unique.  
Together with the fact that $R$ is $(I,0)$-splitting, it follows that $\WGF(R) = \chi_{_{\lambda}}$.\hfill\QED 

The preceding result resolves the problem of constructing crystalline splitting posets as long as there exists such a crystal product having a maximal element of weight $\lambda$. 
The proof of the following existence theorem can be found in \S 7 of \cite{Stem}, see specifically Theorem 7.6.  
The result is actually constructive.  

\noindent 
{\bf \StembridgeExistenceTheorem\ (Stembridge)}\ \ {\sl For any nonzero $\lambda \in \Lambda^{+}$, there is a positive integer $p$ and weights $\mu_1,\ldots,\mu_p$ such that (1) $\lambda = \mu_{1} + \cdots + \mu_{p}$, (2) $\mu_{1} + \cdots + \mu_{q}$ is dominant for each $1 \leq q \leq p$, and (3) either all the weights $\{\mu_{1},\ldots,\mu_{p}\}$ are minuscule or all are quasi-minuscule.}\hfill\QED 

Fix once and for all a total order on the set $\Omega$ of all minuscule and quasi-minuscule weights.  The {\em lexicographically minimal $\Omega$-expression} for $\lambda \in \Lambda^{+}$ is an expression $\mu_{1}+\cdots+\mu_{p}$ for $\lambda$ as in the preceding theorem so that the ``word'' $(\mu_{1},\ldots,\mu_{p})$ in the alphabet $\Omega$ is lexicographically smallest amongst all such expressions. 
The point of this is to identify one such expression for $\lambda$ for use in the following construction. 
We can now resolve as follows the problem of constructing crystalline splitting posets which are splitting posets for the Weyl bialternants. 

\noindent 
{\bf \CrystallineExistenceTheorem}\ \ {\sl Let $\Phi$ be irreducible. 
Let $\lambda \in \Lambda^{+}$ be nonzero with the lexicographically minimal $\Omega$-expression $\lambda = \mu_{1} + \cdots + \mu_{p}$. 
For each $1 \leq q \leq p$, let $\widehat{\mu_{q}}$ be the unique dominant weight in the $W$-orbit of the minuscule or quasi-minuscule weight $\mu_{q}$. 
Let $\xelt_{q}$ be an element of the minuscule or quasi-minuscule splitting poset $R(\widehat{\mu_{q}})$ such that $wt(\xelt_{q}) = \mu_{q}$. 
Then $(\xelt_{1},\ldots,\xelt_{p}) \in R(\widehat{\mu_{1}}) \otimes \cdots \otimes R(\widehat{\mu_{p}})$ is maximal with weight $\lambda$. 
Moreover, the connected component $R(\lambda)$ of $R(\widehat{\mu_{1}}) \otimes \cdots \otimes R(\widehat{\mu_{p}})$ that contains $(\xelt_{1},\ldots,\xelt_{p})$ is a strongly untangled crystalline splitting poset when we take $\mathcal{S}_{J,\nu}(R(\lambda)) = \mathcal{M}_{J,\nu}(R(\lambda))$ for each $J \subseteq I$ and each $\nu \in \Lambda_{\Phi_J}^{+}$, and} $\WGF(R(\lambda)) = \chi_{_{\lambda}}$.

{\em Proof.} First we argue that $(\xelt_{1},\ldots,\xelt_{p})$ is maximal in $R(\widehat{\mu_{1}}) \otimes \cdots \otimes R(\widehat{\mu_{p}})$. 
Well, by \FibrousStructureLemma, for any $i \in I$ we have $\delta_{i}(\xelt_{1},\ldots,\xelt_{p}) = \max\{-(m_i(\xelt_1) + \cdots + m_i(\xelt_{q-1}))+\delta_{i}(\xelt_q)\}$, where the max is taken over all $1 \leq q \leq p$. 
Let $q$ be the leftmost index where this max is attained. 
We aim to show that $-(m_i(\xelt_1) + \cdots + m_i(\xelt_{q-1}))+\delta_{i}(\xelt_q) = 0$. 
Of course, $\delta_{i}(\xelt_{1},\ldots,\xelt_{p}) \geq 0$, so it suffices to show that $-(m_i(\xelt_1) + \cdots + m_i(\xelt_{q-1}))+\delta_{i}(\xelt_q) \leq 0$. 
Suppose that $l_{i}(\xelt_q) \leq 1$ in $R(\widehat{\mu_{q}})$. 
If $\delta_{i}(\xelt_q) = 0$, we are done, since then $-(m_i(\xelt_1) + \cdots + m_i(\xelt_{q-1}))+\delta_{i}(\xelt_q) = -(m_i(\xelt_1) + \cdots + m_i(\xelt_{q-1})) \leq 0$. 
If $\delta_{i}(\xelt_q) = 1$, then $m_i(\xelt_q) = -1$, so $m_i(\xelt_1) + \cdots + m_i(\xelt_{q1}) + m_i(\xelt_q) \geq 0$ means than $-(m_i(\xelt_1) + \cdots + m_i(\xelt_{q-1})) \leq -1$. 
It follows that $-(m_i(\xelt_1) + \cdots + m_i(\xelt_{q-1}))+\delta_{i}(\xelt_q) \leq 0$. 
Now suppose that $l_{i}(\xelt_q) = 2$. 
We cannot have $\delta_{i}(\xelt_q) = 1$ in this case, otherwise $\mu_q = 0$, which is not minuscule or quasi-minuscule. 
So, we have $\delta_{i}(\xelt_q) = 0$ or $\delta_{i}(\xelt_q) = 2$. 
These cases are argued exactly as before, with the conclusion that $-(m_i(\xelt_1) + \cdots + m_i(\xelt_{q-1}))+\delta_{i}(\xelt_q) \leq 0$. 
Thus, $(\xelt_{1},\ldots,\xelt_{p})$ is maximal as claimed. 

Since $wt(\xelt_{1},\ldots,\xelt_{p}) = wt(\xelt_{1}) + \cdots wt(\xelt_{p})$ in the crystal product, then $wt(\xelt_{1},\ldots,\xelt_{p}) = \lambda$. 
The remaining claims follow from \ExistenceCorollary.\hfill\QED

In fact, the construction of the crystalline splitting poset $R(\lambda)$ from the previous theorem does {\em not} depend on the $\Omega$-expression used to represent $\lambda$, although this is not obvious. 
For more discussion on this point, see \CrystallineProblem.A below and \UniquenessCorollary\ in the next section. 

So, we now have one way to resolve the problem of constructing crystalline splitting posets when $\Phi$ is irreducible.  
Now suppose $\Phi$ is reducible.  
Write $I = I_{1} \disjointunion \cdots \disjointunion I_{k}$, where each $I_{j}$ is the index set for an irreducible root subsystem $\Phi_{I_j}$ of $\Phi$ so that $\Phi = \Phi_{I_1} \disjointunion \cdots \disjointunion \Phi_{I_k}$. 
Now take $\lambda \in \Lambda^{+}$. 
We can write $\lambda$ uniquely as  $\lambda_{1} + \cdots + \lambda_{k}$ with each $\lambda_{j}$ dominant in $\Lambda_{j} = \Lambda_{\Phi_{I_j}}$, since $\Lambda = \Lambda_{\Phi_{I_1}} \oplus \cdots \oplus \Lambda_{\Phi_{I_k}}$. 
It is a consequence of the definitions that the crystal product $R(\lambda_1) \otimes \cdots \otimes R(\lambda_k)$ coincides with the Cartesian product $R(\lambda_1) \times \cdots \times R(\lambda_k)$.  
Let $R(\lambda) := R(\lambda_1) \otimes \cdots \otimes R(\lambda_k) = R(\lambda_1) \times \cdots \times R(\lambda_k)$. 
The proof of the following result is routine and is left as an exercise. 

\noindent 
{\bf \CrystallineReducibleProposition}\ \ {\sl Take $\Phi$, $\lambda$, and $R(\lambda)$ as in the preceding paragraph. 
Let $J \subseteq I$ and $\nu \in \Lambda_{\Phi_J}^{+}$. 
Then $\mathcal{M}_{J,\nu}(R(\lambda)) = \mathcal{M}_{J,\nu}(R(\lambda_1)) \times \cdots \times \mathcal{M}_{J,\nu}(R(\lambda_k))$. 
Taking $\mathcal{S}_{J,\nu}(R(\lambda)) := \mathcal{M}_{J,\nu}(R(\lambda))$, then $R(\lambda)$ is a strongly untangled crystalline splitting poset with} $\WGF(R(\lambda)) = \chi_{_{\lambda}}$.\hfill\QED

{\bf When do product decomposition coefficients coincide with $\mathbf{\Phi}$-Kostka numbers?} 
The answer is: Almost always.  
We demonstrate this in the next theorem, whose proof uses reasoning about crystalline splitting posets. 
For each $i, k \in I$, let $u_i^{(k)}$ denote the length of the longest $i$-component in $\Pi(\omega_k)$. 

\noindent 
{\bf \EasyProductTheorem}\ \ {\sl Take $\lambda = \sum_{i \in I} a_{i}\omega_{i}$ in $\Lambda^{+}$. 
Let  $J \subseteq I$, and say $\nu = \sum_{j \in J} \nu_{j}\omega_{j}^{J}$ is in $\Lambda_{\Phi_J}^{+}$.  
Then $\mathcal{M}_{J,\nu}(R(\lambda)) = R(\lambda)$ if and only if $\sum_{k \in I}a_{k}u_{j}^{(k)} \leq \nu_j$ for all $j \in J$. 
If $J = I$, then $\displaystyle \chi_{_{\nu}} \chi_{_{\lambda}} = \sum_{\mu \in \Pi(\lambda)} c_{\nu+\mu} \chi_{_{\nu+\mu}}$ with each $\nu+\mu \in \Lambda^{+}$ and each $c_{\nu+\mu}$ being some nonzero integer if and only if $\displaystyle \chi_{_{\nu}} \chi_{_{\lambda}} = \sum_{\mu \in \Pi(\lambda)} d_{\lambda,\mu} \chi_{_{\nu+\mu}}$ with each $\nu+\mu \in \Lambda^{+}$ if and only if $\sum_{k \in I}a_{k}u_{i}^{(k)} \leq \nu_i$ for all $i \in I$.}

Before we prove the theorem, first a simple lemma. 

\noindent 
{\bf \CrystalProdWeightDiagramLemma}\ \ {\sl With $\lambda = \sum_{i \in I} a_{i}\omega_{i}$ in $\Lambda^{+}$, let $S(\lambda) := R(\omega_1)^{\times a_1} \times \cdots \times R(\omega_n)^{\times a_n}$, where $R(\omega_k)^{\times a_k}$ means obviously a Cartesian product of $R(\omega_k)$ with itself $a_k$ times.  Similarly let $T(\lambda) := R(\omega_1)^{\otimes a_1} \otimes \cdots \otimes R(\omega_n)^{\otimes a_n}$.  Now $S(\lambda)$ is ${\Phi}$-structured by \WeightLemma.3, and $T(\lambda)$ is ${\Phi}$-structured by \CrystalProductStructureLemma. Then we have the following identities of (generalized) weight diagrams:} $\Pi(R(\lambda)) = \Pi(S(\lambda)) = \Pi(T(\lambda)) = \Pi(\lambda).$ 

{\em Proof.} By \PiTheorem, we have $\Pi(S(\lambda)) = \Pi(\lambda) = \Pi(R(\lambda))$, with the latter equality because $\WGF(R(\lambda)) = \chi_{_{\lambda}}$.  
By \WGFCorollary, $\displaystyle \WGF(T(\lambda)) = \prod_{k=1}^{n} \chi_{_{\omega_{k}}}^{a_k} = \psi_{\lambda}$, while $\psi_{\lambda} = \WGF(S(\lambda))$ by \WGFLemma.  
The fact that $\WGF(T(\lambda)) = \WGF(S(\lambda))$ means that $\Pi(T(\lambda)) = \Pi(S(\lambda))$.\hfill\QED

{\em Proof of \EasyProductTheorem.} By \CrystalProdWeightDiagramLemma, $\Pi(\lambda) = \Pi(R(\lambda)) = \Pi(R(\omega_1)^{\otimes a_1} \otimes \cdots \otimes R(\omega_n)^{\otimes a_n}) = \Pi(R(\omega_1)^{\times a_1} \times \cdots \times R(\omega_n)^{\times a_n})$.  
Now the longest $i$-component of each $R(\omega_k)$ is $u_{i}^{(k)}$, since $R(\omega_k)$ is fibrous and $\Pi(R(\omega_k)) = \Pi(\omega_k)$.  
Since the longest $i$-component in $R(\omega_1)^{\times a_1} \times \cdots \times R(\omega_n)^{\times a_n}$ is a product of the longest $i$-components in each $R(\omega_k)$, this longest $i$-component has length $\sum_{k \in I} a_{k}u_{i}^{(k)}$.  
So the longest $i$-component in $\Pi(\lambda)$ has length $\sum_{k \in I} a_{k}u_{i}^{(k)}$.  Since $R(\lambda)$ is fibrous and $\Pi(R(\lambda)) = \Pi(\lambda)$, then the longest $i$-component of $R(\lambda)$ has length $\sum_{k \in I} a_{k}u_{i}^{(k)}$. 
This is enough to establish the first equivalence of the theorem statement.  

Now let $J=I$. 
We aim to prove the three assertions from the last sentence of the theorem statement are equivalent by showing that each is equivalent to the assertion that $\mathcal{M}_{I,\nu}(R(\lambda)) = R(\lambda)$. 
Given what was demonstrated in the previous paragraph, we only have to show this for the first two of the three assertions. 
To that end, assume that $\displaystyle \chi_{_{\nu}} \chi_{_{\lambda}} = \sum_{\mu \in \Pi(\lambda)} c_{\nu+\mu} \chi_{_{\nu+\mu}}$ with each $\nu+\mu \in \Lambda^{+}$ and $c_{\nu+\mu} \not= 0$. 
We will show that $\mathcal{M}_{I,\nu}(R(\lambda)) = R(\lambda)$. 
Otherwise, we can choose $\xelt \in R(\lambda) \setminus \mathcal{M}_{I,\nu}(R(\lambda))$ and $j \in I$ such that $\delta_{j}(\xelt)$ is as large as possible. 
Therefore $\delta_{j}(\xelt) > \nu_{j}$. 
But now $|\{\melt \in \mathcal{M}_{I,\nu}(R(\lambda))\, |\, wt(\melt) = wt(\xelt)\}|$ equals $c_{\nu+wt(\xelt)}$, which is nonzero. 
Therefore, there must be some $\yelt$ with $wt(\yelt) = wt(\xelt)$ and $\delta_{j}(\yelt) \leq \nu_{j}$. 
In particular, $m_{j}(\yelt) = m_{j}(\xelt)$, and hence $\rho_{j}(\yelt) - \delta_{j}(\yelt) = \rho_{j}(\xelt) - \delta_{j}(\xelt)$. 
Since $\delta_{j}(\xelt)$ is as large as possible, then $\rho_{i}(\xelt) = 0$.  
So $\rho_{j}(\yelt) - \delta_{j}(\yelt) = -\delta_{j}(\xelt)$, which means that $\rho_{j}(\yelt) = \delta_{j}(\yelt) - \delta_{j}(\xelt)$. 
That is, $\delta_{j}(\yelt) - \delta_{j}(\xelt) \geq 0$, i.e.\ $\delta_{j}(\yelt) \geq \delta_{j}(\xelt)$. 
This contradicts the fact that $\delta_{j}(\yelt) \leq \nu_{j} < \delta_{j}(\xelt)$. 
So, we must have $\mathcal{M}_{I,\nu}(R(\lambda)) = R(\lambda)$. 

Assume now that $\mathcal{M}_{I,\nu}(R(\lambda)) = R(\lambda)$.  
Then 
$\displaystyle \chi_{_{\nu}} \cdot \WGF(R(\lambda)) = \sum_{\melt \in \mathcal{M}_{I,\nu}(R(\lambda))}\chi_{_{\nu + wt(\melt)}}$ implies that each coefficient $c_{\nu+\mu}$ in the expansion $\displaystyle \chi_{_{\nu}} \chi_{_{\lambda}} = \sum_{\mu \in \Pi(\lambda)} c_{\nu+\mu} \chi_{_{\nu+\mu}}$ is equal to $|\{\melt \in R(\lambda)\, |\, wt(\melt) = \mu\}|$, which is obviously just $d_{\lambda,\mu}$. 
Therefore $\displaystyle \chi_{_{\nu}} \chi_{_{\lambda}} = \sum_{\mu \in \Pi(\lambda)} d_{\lambda,\mu} \chi_{_{\nu+\mu}}$ with each $\nu+\mu \in \Lambda^{+}$. 
Of course, since each of these $d_{\lambda,\mu}$'s is nonzero, then assuming $\displaystyle \chi_{_{\nu}} \chi_{_{\lambda}} = \sum_{\mu \in \Pi(\lambda)} d_{\lambda,\mu} \chi_{_{\nu+\mu}}$ with each $\nu+\mu \in \Lambda^{+}$ we can immediately conclude that $\displaystyle \chi_{_{\nu}} \chi_{_{\lambda}} = \sum_{\mu \in \Pi(\lambda)} c_{\nu+\mu} \chi_{_{\nu+\mu}}$ with each $\nu+\mu \in \Lambda^{+}$ and $c_{\nu+\mu} \not= 0$. 
\hfill\QED 

For irreducible root systems, we can easily determine the $u_{i}^{(k)}$'s.   

\noindent 
{\bf \DataProp}\ \ {\sl Let $\Phi$ be an irreducible root system. (1) Let $i, k \in I$.  Then $\alpha^{\vee} \in \Phi^{\vee}$ is the highest coroot in the $W$-orbit of $\alpha_{i}^{\vee}$ if and only if $\alpha \in \Phi$ is the highest root in the $W$-orbit of $\alpha_{i}$. For such an $\alpha$, we have}
\[u_{i}^{(k)} = \max\{\langle \nu,\alpha^{\vee} \rangle\, |\, \nu \in \Pi(\omega_{k}) \cap \Lambda^{+}\}.\] 
{\sl (2) We have the following data for the vector $u^{(k)} := (u_{1}^{(k)}, \ldots , u_{n}^{(k)})$ with $k \in I$}:\\ 
\hspace*{0.5in}{\sl When} $\Phi = \myA_n$, {\sl we have} $u^{(k)} = (1, \ldots , 1)$.\\ 
\hspace*{0.5in}{\sl When} $\Phi = \myB_n$, {\sl we have} $u^{(k)} = \rule[-8mm]{0mm}{22mm}\left\{
\begin{array}{cl}
(1, \ldots , 1, 2) & \hspace*{0.1in}\mbox{if } k=1\\
(2, \ldots , 2) & \hspace*{0.1in}\mbox{if } 2 \leq k \leq n-1\\
(1, \ldots , 1) & \hspace*{0.1in}\mbox{if } k=n
\end{array}\right.$\\
\hspace*{0.5in}{\sl When} $\Phi = \myC_n$, {\sl we have} $u^{(k)} = \rule[-5mm]{0mm}{16mm}\left\{
\begin{array}{cl}
(1, \ldots , 1) & \hspace*{0.1in}\mbox{if } k=1\\
(2, \ldots , 2, 1) & \hspace*{0.1in}\mbox{if } 2 \leq k \leq n
\end{array}\right.$\\
\hspace*{0.5in}{\sl When} $\Phi = \myD_n$, {\sl we have} $u^{(k)} = \rule[-5mm]{0mm}{16mm}\left\{
\begin{array}{cl}
(1, \ldots , 1) & \hspace*{0.1in}\mbox{if } k=1, n-1, n\\
(2, \ldots , 2) & \hspace*{0.1in}\mbox{if } 2 \leq k \leq n-2
\end{array}\right.$\\
\hspace*{0.5in}{\sl When} $\Phi = \myE_6$, {\sl we have} $u^{(k)} = \rule[-8mm]{0mm}{22mm}\left\{
\begin{array}{cl}
(1, 1, 1, 1, 1, 1) & \hspace*{0.1in}\mbox{if } k=1, 6\\
(2, 2, 2, 2, 2, 2) & \hspace*{0.1in}\mbox{if } k=2\\
(3, 3, 3, 3, 3, 3) & \hspace*{0.1in}\mbox{if } k=3, 4, 5
\end{array}\right.$\\
\hspace*{0.5in}{\sl When} $\Phi = \myE_7$, {\sl we have} $u^{(k)} = \rule[-10mm]{0mm}{26mm}\left\{
\begin{array}{cl}
(1, 1, 1, 1, 1, 1, 1) & \hspace*{0.1in}\mbox{if } k=7\\
(2, 2, 2, 2, 2, 2, 2) & \hspace*{0.1in}\mbox{if } k=1, 2, 6\\
(3, 3, 3, 3, 3, 3, 3) & \hspace*{0.1in}\mbox{if } k=3, 5\\
(4, 4, 4, 4, 4, 4, 4) & \hspace*{0.1in}\mbox{if } k=4
\end{array}\right.$\\
\hspace*{0.5in}{\sl When} $\Phi = \myE_8$, {\sl we have} $u^{(k)} = \rule[-14mm]{0mm}{34mm}\left\{
\begin{array}{cl}
(2, 2, 2, 2, 2, 2, 2, 2) & \hspace*{0.1in}\mbox{if } k=1, 8\\
(3, 3, 3, 3, 3, 3, 3, 3) & \hspace*{0.1in}\mbox{if } k=2, 7\\
(4, 4, 4, 4, 4, 4, 4, 4) & \hspace*{0.1in}\mbox{if } k=3, 6\\
(5, 5, 5, 5, 5, 5, 5, 5) & \hspace*{0.1in}\mbox{if } k=5\\
(6, 6, 6, 6, 6, 6, 6, 6) & \hspace*{0.1in}\mbox{if } k=4
\end{array}\right.$\\
\hspace*{0.5in}{\sl When} $\Phi = \myF_4$, {\sl we have} $u^{(k)} = \rule[-10mm]{0mm}{26mm}\left\{
\begin{array}{cl}
(2, 2, 2, 2) & \hspace*{0.1in}\mbox{if } k=1\\
(3, 3, 4, 4) & \hspace*{0.1in}\mbox{if } k=2\\
(2, 2, 3, 3) & \hspace*{0.1in}\mbox{if } k=3\\
(1, 1, 2, 2) & \hspace*{0.1in}\mbox{if } k=4
\end{array}\right.$\\
\hspace*{0.5in}{\sl When} $\Phi = \myG_2$, {\sl we have} $u^{(k)} = \rule[-8mm]{0mm}{18mm}\left\{
\begin{array}{cl}
(2, 1) & \hspace*{0.1in}\mbox{if } k=1\\
(3, 2) & \hspace*{0.1in}\mbox{if } k=2
\end{array}\right.$ 

{\em Proof.} 
The first statement of part {\sl (1)} of the proposition uses the facts that for irreducible root systems, (a) all roots of a given length are $W$-conjugate, (b) a root $\beta$ is short (respectively long) if and only if the coroot $\beta^{\vee}$ is long (resp.\ short), and (c) $\sum_{j \in I} k_{j}\alpha_{j}$ is a root in $\Phi$ if and only if $\sum_{j \in I} k_{j}\alpha_{j}^{\vee}$ is a root in $\Phi^{\vee}$. 
So suppose $\alpha$ is the highest root in the $W$-orbit of some simple root $\alpha_i$. 
If we write $\alpha = \sum_{j \in I} k_{j}\alpha_{j}$, it follows from fact (c) that $\alpha^{\vee} = \sum_{j \in I} k_{j}\alpha_{j}^{\vee} \in \Phi^{\vee}$ is highest in its $W$-orbit. 
Moreover, by fact (b), $\alpha^{\vee}$ and $\alpha_{i}^{\vee}$ have the same length, and so are $W$-conjugate by fact (a). 
That is, $\alpha^{\vee}$ is the highest coroot in the $W$-orbit of $\alpha_{i}^{\vee}$. 
The proof of the converse is entirely similar. 

Now for any $\mu \in \Pi(\omega_k)$, there is a $\sigma \in W$ such that $\nu := \sigma(\mu) \in \Pi(\omega_k) \cap \Lambda^{+}$. 
Then $\langle \mu,\alpha_{i}^{\vee} \rangle = \langle \nu,\sigma(\alpha_{i}^{\vee}) \rangle$. 
With $\alpha^{\vee}$ as the highest root in the $W$-orbit of $\alpha_{i}^{\vee}$, then dominance of $\nu$ implies that for any other $\beta^{\vee} = \tau(\alpha_{i}^{\vee})$ in this same $W$-orbit ($\tau \in W$), we have $\langle \nu,\beta^{\vee} \rangle \leq \langle \nu,\alpha^{\vee} \rangle$.  Therefore, 
\[u_{i}^{(k)} := \max\{\langle \mu,\alpha_{i}^{\vee} \rangle \, |\, \mu \in \Pi(\omega_{k})\} = \max\{\langle \nu,\alpha^{\vee} \rangle\, |\, \nu \in \Pi(\omega_{k}) \cap \Lambda^{+}\}.\] 

We will see in \S \CrystalSection\ that the crystal graph $\mathcal{G}_{\mbox{\tiny crystal}}(\omega_{k})$ associated with a fundamental weight $\omega_{k}$ is connected and ${\Phi}$-structured and that its unique maximal element has weight $\omega_{k}$. 
So by \PiTheorem, $\Pi(\mathcal{G}_{\mbox{\tiny crystal}}(\omega_{k})) = \Pi(\omega_{k})$. 
For the $\myA_{n}$---$\myD_{n}$ cases of part {\sl (2)} of the proposition statement, one can therefore obtain the data for $u^{(k)}$ using the explicit descriptions given in \cite{KN} of the crystal graph $\mathcal{G}_{\mbox{\tiny crystal}}(\omega_{k})$. 
For the $\myE_{6}$, $\myE_{7}$, $\myE_{8}$, $\myF_{4}$, and $\myG_{2}$ cases, one can obtain each $u^{(k)}$ by using the formula from part {\sl (1)} of the proposition. 
For explicit root and weight data for these computations, see \cite{BMP}.\hfill\QED

{\bf Further consequences, questions, and comments.}  
We close this section by enumerating some questions and comments about the foregoing construction of crystalline splitting posets. 

\noindent 
{\bf \CrystallineProblems}

(A) Suppose that in some crystal product of primary-plus posets, $R'(\lambda)$ is the connected component of some maximal element of weight $\lambda$. 
\ExistenceCorollary\ shows that $R'(\lambda)$ is a strongly untangled crystalline splitting poset with $\WGF(R'(\lambda)) = \chi_{_{\lambda}}$. 
It follows rather easily from crystal base theory that necessarily $R(\lambda) \cong R'(\lambda)$ as edge-colored directed graphs, see \UniquenessCorollary\ below. 
One could ask for a combinatorial proof of this uniqueness result.  

(B) The question in (A) is closely related to the question of providing a combinatorial algorithm showing that $R(\lambda) \otimes R(\nu) \cong R(\nu) \otimes R(\lambda)$ for any $\nu, \lambda \in \Lambda^{+}$.  See \cite{LenartTwo} for an answer to this question for the cases other than $\myE_{8}$, $\myF_4$, and $\myG_2$. 

(C) It is natural to ask under what circumstances the crystal product of fibrous posets is strongly untangled, cf. Question 6.3 of \cite{Stem}.  In particular, to what extent is \StembridgeUntangledTheorem\ a purely poset-theoretic result and to what extent does it require the ${\Phi}$-structure property or other root system related properties?\hfill\QED 

\noindent 
{\bf \CrystallineRemarks}

(A) It is natural to ask what relationship $R(\lambda)^{\bowtie}$ has to $R(\lambda)$.  In fact, it is another easy consequence of crystal base theory that $R(\lambda)^{\bowtie} \cong R(\lambda)$ as edge-colored directed graphs, see \CrystalOpsCorollary.  In \cite{Lenart}, Lenart provides a combinatorial algorithm for producing such an isomorphism.  His algorithm generalizes the ``evacuation'' involution on semistandard tableaux. 

(B) Positivity: From the fact that $R(\lambda)$ is a crystalline splitting poset with $\WGF(R(\lambda)) = \chi_{_{\lambda}}$, it follows that the product $\chi_{_{\nu}} \chi_{_{\lambda}}$ is a positive integer linear combination of some Weyl bialternants, i.e.\ the product decomposition coefficients $c_{\mu}$ in the expansion $\chi_{_{\nu}} \chi_{_{\lambda}} = \sum_{\mu \in \Lambda^{+}} c_{\mu}\chi_{_{\mu}}$ are nonnegative integers. 
Indeed, $c_{\mu} = |\{\melt \in \mathcal{M}_{I,\nu}(R(\lambda))\, |\, wt(\melt) = \mu\}|$. 
For $J \subseteq I$, it also follows that $\chi_{_{\lambda}}|_{_J}$ is a positive integer linear combination of some $\Phi_{J}$-Weyl bialternants, i.e.\ the restriction coefficients $b_{\mu}$ in the expansion $\chi_{_{\lambda}}|_{_J} = \sum_{\mu \in \Lambda_{\Phi_J}^{+}} b_{\mu}\chi^{\Phi_J}_{_{\mu}}$ are nonnegative integers. 
In fact, we have $b_{\mu} = |\{\melt \in \mathcal{M}_{J,0}(R(\lambda))\, |\, wt^{J}(\melt) = \mu\}|$. 

(C) We can make {\em any} fibrous splitting poset $R$ for $\chi_{_{\lambda}}$ into a crystalline splitting poset as follows.  
Take any weight-preserving bijection $\phi: R(\lambda) \longrightarrow R$. 
Now let $J \subseteq I$ and $\nu \in \Lambda_{\Phi_J}^{+}$. 
Then define $\mathcal{S}_{J,\nu}(R) := \phi(\mathcal{M}_{J,\nu}(R(\lambda))$. 
Since $\WGF(R)|_J = \WGF(R(\lambda))|_J$, it follows that $R$ is a crystalline splitting poset. 
More generally, even if $R$ is not fibrous, we can use this same procedure to identify sets $\mathcal{S}_{J,\nu}(R)$ so that $R$ is $(J,\nu)$-splitting for each pair $(J,\nu)$. 
How interesting this result is in any particular case will depend on how explicit and natural the bijection $\phi$ is, how nice the sets $\mathcal{S}_{J,\nu}(R)$ are, etc.  

(D) Another natural question is when the product $\chi_{_{\nu}} \chi_{_{\lambda}}$ will be ``multiplicity-free'' in the sense that when we write $\chi_{_{\nu}} \chi_{_{\lambda}} = \sum_{\mu \in \Lambda^{+}} c_{\mu}\chi_{_{\mu}}$, each $c_{\mu} \in \{0,1\}$. 
This problem is combinatorial in the sense that the given product of Weyl bialternants is multiplicity-free if and only if for all distinct $\melt_1, \melt_2 \in \mathcal{M}_{I,\nu}(R(\lambda))$ we have $wt(\melt_1) \not= wt(\melt_2)$. 
Similarly, one can ask, for a given $J \subset I$, when the branching $\chi_{_{\lambda}}|_{_J}$ will be multiplicity-free. 
This property is equivalent to having $wt^{J}(\melt_1) \not= wt^{J}(\melt_2)$ whenever $\melt_1$ and $\melt_2$  are distinct in $\mathcal{M}_{J,0}(R(\lambda))$. 
This ``mutliplicity-free calculus'' is studied thoroughly in \cite{StemMultFree}, with explicit answers.\hfill\QED

\newpage
\renewcommand{\thefootnote}{1}\vspace{1ex}
\noindent 
{\Large \bf \S \CrystalSection.\ Stembridge's Admissible systems and the crystal bases/graphs of Kashiwara {\em et al}.} 

For any dominant weight $\lambda$, Kashiwara's crystal base on the one hand and a particular admissible system of Stembridge on the other yield certain fibrous edge-colored directed graphs, notated here as  $\mathcal{G}_{\mbox{\tiny crystal}}(\lambda)$ and $\mathcal{G}_{\mbox{\tiny adm-syst}}(X(\lambda))$ respectively.\footnote{By uniqueness of the crystal base -- see \cite{Kash2} or Ch.\ 5 of \cite{HK} -- there is only one crystal graph $\mathcal{G}_{\mbox{\tiny crystal}}(\lambda)$ associated with the dominant weight $\lambda$. 
This graph is sometimes called a ``regular'' or ``normal'' crystal graph, to distinguish it from other edge-colored directed graphs which share some of the combinatorial properties of regular crystal graphs but do not arise from crystal bases.}   These graphs encode much of the information about the related irreducible representation of the associated semisimple Lie algebra or quantized enveloping algebra.  So, they are obvious candidates for (refined) splitting posets for $\chi_{_{\lambda}}$.  

Our main objective here is to demonstrate directly that these graphs are crystalline splitting posets in the sense of the previous section, see \MainCrystalTheorems.  
This will serve to provide an introduction to crystal graphs/admissible systems for the interested but perhaps uninitiated combinatorial reader and 
to enter these and related results into the record via the (mostly) self-contained vertex-coloring proofs we provide.  
In particular, we will use crystal graphs to see how certain combinatorial statements about our $R(\lambda)$ crystalline splitting posets have algebraic justifications, and we will pose some of these algebraic results as combinatorial problems. 
We will also see exactly how the Weyl character formula is a corollary of the proof given here that crystal graphs are splitting posets. 

One could add a third graph to this collection, namely the edge-colored directed graph associated with Littelmann's path model, which we denote $\mathcal{G}_{\mbox{\tiny path}}(\lambda)$.  The Isomorphism Theorem of \cite{LitPaths} guarantees that the path model is, in a precise sense, independent of the choice of an initial path, and therefore the path model also yields only one graph $\mathcal{G}_{\mbox{\tiny path}}(\lambda)$; see also Theorem 6.5 of \cite{Lit4}.  That these graphs are crystalline splitting posets follows from the fact, proved by Joseph \cite{Joseph} and Kashiwara \cite{Kash3} and stated as Theorem 11.1 in \cite{Lit4}, that $\mathcal{G}_{\mbox{\tiny path}}(\lambda)$ and $\mathcal{G}_{\mbox{\tiny crystal}}(\lambda)$ are isomorphic edge-colored directed graphs.  It also follows from the fact that $\mathcal{G}_{\mbox{\tiny path}}(\lambda)$ arises as $\mathcal{G}_{\mbox{\tiny adm-syst}}(X(\lambda))$ for some admissible system $X(\lambda)$, cf.\ \S 8 of \cite{Stem}. See also \LittelmannRemark.B below. 

{\bf Admissible systems as crystalline splitting posets.} 
Stembridge's admissible systems are characterized axiomatically in \cite{Stem}.  An admissible system is a 4-tuple $(X,\mu,\Stemdel,\{\widetilde{F}_{i}|i\in{I}\})$, where $X$ is a set of objects, $\mu$ and $\Stemdel$ are maps $X \rightarrow \Lambda$, and each $\widetilde{F}_{i}$ is a bijection between two subsets of $X$. 
We sometimes use $X$ by itself to denote the admissible system.  
We will require that $X$ be finite. 
Below we see that $\Stemdel$ of the admissible system is a nonpositive-valued version of our nonnegative-valued depth-measuring quantity $\delta$. 
Denote by $\widetilde{E}_{i}$ the inverse function $\widetilde{F}_{i}^{-1}$. For reasons that will be clear shortly, we think of $\widetilde{E}_{i}$ and $\widetilde{F}_{i}$ as ``raising'' and ``lowering'' operators respectively; we use the tilde in denoting these operators because of their similarities with Kashiwara's raising and lowering operators in the context of crystal graphs. 
We let $\varepsilon(\xelt) := \mu(\xelt)-\Stemdel(\xelt)$ for all $\xelt \in X$. Also, for each $i \in I$ set $\mu(\xelt,i) := \langle \mu(\xelt),\alpha_{i}^{\vee} \rangle$, $\Stemdel(\xelt,i) := \langle \Stemdel(\xelt),\alpha_{i}^{\vee} \rangle$, $\varepsilon(\xelt,i) := \langle \varepsilon(\xelt),\alpha_{i}^{\vee} \rangle$.  
Write $\xelt \preceq_{i} \yelt$ if $\xelt = \widetilde{F}_{i}^{k}(\yelt)$ for some $k \geq 0$.  

Assuming $X$ is finite, the system is required to satisfy four axioms: (A1) $\Stemdel(\xelt) \in -\Lambda^{+}$ and  $\varepsilon(\xelt) \in \Lambda^{+}$; (A2) $\widetilde{F}_{i}$ is a bijection $\{\xelt \in X\, |\, \varepsilon(\xelt,i) > 0\} \longrightarrow \{\xelt \in X\, |\, \Stemdel(\xelt,i) < 0\}$; (A3) $\mu(\widetilde{F}_{i}(\xelt)) = \mu(\xelt) - \alpha_{i}$, $\Stemdel(\widetilde{F}_{i}(\xelt),i) = \Stemdel(\xelt,i)-1$, and therefore $\varepsilon(\widetilde{F}_{i}(\xelt),i) = \varepsilon(\xelt,i)-1$; and (A4) the system has a {\em coherent timing pattern}, which is a real-valued function $t(\xelt,i)$ defined for all $\xelt \in X$ and $i \in I$ for which $\Stemdel(\xelt,i) < 0$ and satisfying the following properties: (i) If $t(\xelt,i)$ and $\widetilde{F}_{i}(\xelt)$ are defined then $t(\xelt,i) \leq t(\widetilde{F}_{i}(\xelt),i)$, and (ii) for all $j \not= i$, all integers $\delta < 0$, and all $t \leq t(\xelt,i)$, there is $\yelt \in X$ with $\yelt \succeq_{j} \xelt$ and $\Stemdel(\yelt,i) = \delta$ and $t(\yelt,i) = t$ if and only if there is $\yelt' \in X$ with $\yelt' \succeq_{j} \widetilde{F}_{i}(\xelt)$ and $\Stemdel(\yelt',i)= \Stemdel$ and $t(\yelt',i) = t$.  

It can be seen as follows that underlying the admissible system $(X,\mu,\Stemdel,\{\widetilde{F}_{i}|i\in{I}\})$ is a ranked poset $\ASG(X)$ with edges colored by the set $I$. Vertices of $\ASG(X)$ are the elements of $X$.  For $\selt, \telt \in X$ and $i \in I$, write $\selt \myarrow{i} \telt$ if and only if $\widetilde{E}_{i}(\selt) = \telt$ if and only if $\selt = \widetilde{F}_{i}(\telt)$.  The fact that $\mu(\selt) + \alpha_{i} = \mu(\telt)$ whenever $\selt \myarrow{i} \telt$ (axiom (A3)) shows that the resulting graph is acyclic (no directed cycles), that the transitive closure of these edge relations defines a partial order, that the given edges are covering relations with respect to this partial order, and that the resulting poset $\ASG(X)$ is ranked.  By (A2), we see that the $i$-components of $\ASG(X)$ are chains.  Now let $\xelt \in X$ and $i \in I$.  Then $\mathbf{comp}_{i}(\xelt)$ in $\ASG(X)$ is a chain $\xelt_{0} \myarrow{i} \xelt_{1} \myarrow{i} \cdots \myarrow{i} \xelt_{p}$, where $\xelt = \xelt_{r}$ for some $0 \leq r \leq p$.  Since $\xelt_{0}$ is minimal in the chain $\mathbf{comp}_{i}(\xelt)$, then by (A2) we have $\varepsilon(\xelt_{0},i) \leq 0$.  But by (A2) we also see that $\varepsilon(\xelt_{1},i) > 0$.  Since by (A1) we must have $\varepsilon(\xelt_{q},i) \in \mathbb{Z}$ for each $0 \leq q \leq p$, and since by (A3) we must have $\varepsilon(\xelt_{q-1},i) = \varepsilon(\xelt_{q},i)-1$ for each $1 \leq q \leq p$, then we conclude that $\varepsilon(\xelt_{0},i) = 0$ and $\varepsilon(\xelt_{1},i) = 1$. Further application of (A3) shows that for $0 \leq q \leq p$, $\varepsilon(\xelt_{q},i) = \rho_{i}(\xelt_{q}) = q$ and $\delta(\xelt_{q},i) = \rho_{i}(\xelt_{q})-l_{i}(\xelt_{q}) = q-p$, and hence that $m_{i}(\xelt) = \rho_{i}(\xelt) - (l_{i}(\xelt)-\rho_{i}(\xelt)) = \varepsilon(\xelt,i)+\Stemdel(\xelt,i) = \mu(\xelt,i) = \langle \mu(\xelt),\alpha_{i}^{\vee} \rangle$.  In view of \WeightLemma, we now have $\langle wt(\xelt),\alpha_{i}^{\vee} \rangle = \langle \mu(\xelt),\alpha_{i}^{\vee} \rangle$ for each $i \in I$, and hence $wt(\xelt) = \mu(\xelt)$.  So $\ASG(X)$ is the Hasse diagram for an ${\Phi}$-structured poset. It is now evident that in the edge-colored poset $\ASG(X)$, for all $i \in I$ and $\xelt \in X$ we have $\delta_{i}(\xelt) = - \Stemdel(\xelt,i)$. 

The following is Lemma 2.3 of \cite{Stem} and is needed for our next theorem. 

\noindent 
{\bf \StemLemma\ (Stembridge)}\ \ {\sl For $j \not= k$ in $I$, suppose $\xelt' \preceq_{j} \xelt$, $d$ is a positive integer, and $T \leq t(\xelt,j)$. Then there is a $\yelt \succeq_{k} \xelt$ such that $\delta_{k}(\yelt) = d$ and $t(\yelt,k) = T$ if and only if there is a $\yelt' \succeq_{k} \xelt'$ such that $\delta_{k}(\yelt') = d$ and $t(\yelt',k) = T$.}\hfill\QED 

The next result reworks Theorem 2.4 of \cite{Stem} by showing exactly how an admissible system $X$ can be realized as a crystalline splitting poset. 
For any $\xelt \in \ASG(X) \setminus \mathcal{M}_{J,\nu}(\ASG(X))$ and $i \in \mathcal{K}_{J,\nu}(\xelt)$, we let $\xelt^{(i)}$ denote the unique element of $\comp_{i}(\xelt)$ for which $\delta_{i}(\xelt^{(i)}) = \nu_{i}+1$, when $\nu_i \in \{0,1,\ldots,l_{i}(\xelt)\}$. 

\noindent 
{\bf \AdmissibleTheorem}\ \ {\sl For any (finite) admissible system $(X,\mu,\delta,\{\widetilde{F}_{i}|i\in{I}\})$, the edge-colored directed graph $\ASG(X)$ is the fibrous edge-colored Hasse diagram of an ${\Phi}$-structured poset.  Fix a total ordering of $I$.  Then for each $J \subseteq I$ and $\nu \in \Lambda_{\Phi_J}^{+}$, the function $\kappa: \ASG(X) \setminus \mathcal{M}_{J,\nu}(\ASG(X)) \longrightarrow J$ given by}
\[\kappa(\xelt) = \min\bigg\{i \in \mathcal{K}_{J,\nu}(\xelt)\, :\, t(\xelt^{(i)},i) = \min_{k \in \mathcal{K}_{J,\nu}(\xelt)}\{t(\xelt^{(k)},k)\}\bigg\}\]
{\sl is a $(J,\nu)$-coloring of $\ASG(X)$.  
In particular, it follows that $\ASG(X)$ is a crystalline splitting poset and that for each such $(J,\nu)$ pair we have}
\[\chi^{\Phi_J}_{_{\nu}} \cdot \WGF(\ASG(X))|_{J} = \sum_{\melt \in \mathcal{M}_{J,\nu}(\ASG(X))}\chi^{\Phi_J}_{_{\nu + wt^{J}(\melt)}}.\]

{\em Proof.} That $\ASG(X)$ is fibrous and ${\Phi}$-structured follows from the paragraphs preceding the theorem statement.  Once we show that the function $\kappa$ is a $(J,\nu)$-coloring of $\ASG(X)$, then we will be done, since the identity of $W_{J}$-symmetric functions as well as the crystalline claim of the theorem statement follow from \MainColoringTheorem.  So, if $\telt \in \ASG(X) \setminus \mathcal{M}_{J,\nu}(\ASG(X))$, then we need to show that 
\[\{\selt \in \comp_{j}(\telt)\, :\, \selt \not\in \mathcal{M}_{J,\nu}(\ASG(X)) \mbox{ and } \kappa(\selt) = j\} = \comp_{j}(\telt) \setminus U_{j}(\telt,\nu_j).\] 

To this end, let $\xelt \in \comp_{j}(\telt) \setminus U_{j}(\telt,\nu_j)$.  Then in this case we  automatically get $\delta_{j}(\xelt) > \nu_{j}$, so $\xelt \not\in \mathcal{M}_{J,\nu}(\ASG(X))$. Now suppose that $t(\xelt^{(k)},k) \leq t(\xelt^{(j)},j)$ for some $k \not= j$.  Now $\xelt^{(k)} \succeq_{k} \xelt$.  Let $d := \delta_{k}(\xelt^{(k)}) = \nu_{k}+1 > 0$ and $T := t(\xelt^{(k)},k)$.  Since $\xelt$ and $\telt$ are in the chain $\comp_{j}(\telt)$, then we must have one of $\xelt \preceq_{j} \telt$ or $\xelt \succeq_{j} \telt$.  If $\xelt \preceq_{j} \telt$, then observe that $T = t(\xelt^{(k)},k) \leq t(\xelt^{(j)},j) = t(\telt^{(j)},j) \leq t(\telt,j)$.  On the other hand, if $\xelt \succeq_{j} \telt$, then $T = t(\xelt^{(k)},k) \leq t(\xelt^{(j)},j) \leq t(\xelt,j)$.  Either way, by \StemLemma, we get a $\yelt \succeq_{k} \telt$ with $\delta_{k}(\yelt) = d = \nu_{k}+1$ and $t(\yelt,k) = T = t(\xelt^{(k)},k)$.  Then $\yelt = \telt^{(k)}$.  In particular, $t(\telt^{(k)},k) = t(\xelt^{(k)},k) \leq t(\xelt^{(j)},j) = t(\telt^{(j)},j)$. The fact that $j = \kappa(\telt)$ says that $t(\telt^{(j)},j) \leq t(\telt^{(k)},k)$, and combining this with the inequality of the preceding sentence gives us equality all the way through: $t(\telt^{(k)},k) = t(\xelt^{(k)},k) = t(\xelt^{(j)},j) = t(\telt^{(j)},j)$.  Also, since now $t(\telt^{(k)},k)  = t(\telt^{(j)},j)$, then by definition of $\kappa$, it must be the case that $j < k$ in the fixed ordering of $J$. Thus we have shown that if $t(\xelt^{(k)},k) \leq t(\xelt^{(j)},j)$ for some $k \not= j$, then $t(\xelt^{(k)},k) = t(\xelt^{(j)},j)$ and $j < k$.  Therefore, $\kappa(\xelt) = j$, so $\xelt \in \{\selt \in \comp_{j}(\telt)\, :\, \selt \not\in \mathcal{M}_{J,\nu}(\ASG(X)) \mbox{ and } \kappa(\selt) = j\}$. 

On the other hand, suppose now that $\xelt \in \{\selt \in \comp_{j}(\telt)\, :\, \selt \not\in \mathcal{M}_{J,\nu}(\ASG(X)) \mbox{ and } \kappa(\selt) = j\}$.  Then automatically $j \in \mathcal{K}_{J,\nu}(\xelt)$ (cf.\ \S \FibrousSection), hence $\delta_{j}(\xelt) > \nu_{j}$, and hence $\xelt \in \comp_{j}(\telt) \setminus U_{j}(\telt,\nu_j)$.\hfill\QED

\noindent 
{\bf \LittelmannRemarks} 

(A) Starting with the lexicographically minimal $\Omega$-expression for a given dominant weight $\lambda$, Stembridge can build an admissible system, here denoted $X(\lambda)$, using a crystal-type product of ``thin'' admissible systems whose underlying graphs are minuscule and quasi-minuscule splitting posets. 
The result is a strongly untangled system with a unique maximal element of dominant weight $\lambda$. 
It follows (by \AdmissibleTheorem, for example) that $\ASG(X(\lambda))$ is a crystalline splitting poset and that $\WGF(\ASG(X(\lambda))) = \chi_{_{\lambda}}$. 
A different $\Omega$-expression for $\lambda$ could conceivably result in a different graph $\ASG(X'(\lambda))$. 
As an algebraic consequence of crystal base theory, we argue in \UniquenessCorollary\ that in fact $\ASG(X'(\lambda)) \cong \ASG(X(\lambda)) \cong R(\lambda)$. 

(B) In \S 8 of \cite{Stem}, Stembridge shows how to realize a particular version of Littelmann's path model as an admissible system.  
In particular, for each element of $\mathcal{G}_{\mbox{\tiny path}}(\lambda)$, he provides in the proof of his Theorem 8.3 an explicit formula for the timing pattern $t(\xelt,i)$.   
By \AdmissibleTheorem, this gives an explicit $(J,\nu)$-coloring of $\mathcal{G}_{\mbox{\tiny path}}(\lambda)$ for each $J \subseteq I$ and $\nu \in \Lambda_{\Phi_J}^{+}$.  
That is, in this way we can view $\mathcal{G}_{\mbox{\tiny path}}(\lambda)$ as a crystalline splitting poset.\hfill\QED

{\bf Crystal graphs as crystalline splitting posets.}  
In \cite{Kash1} and \cite{Kash2}, Kashiwara introduced the notion of a crystal base and its associated crystal graph for a module of the quantized enveloping algebra $\mathcal{U}_{q}(\mathfrak{g})$ associated to a semisimple Lie algebra $\mathfrak{g}$. 
See \cite{Jan} or \cite{HK} for readable expositions of these and related ideas. 
Roughly, the crystal graph is an edge-colored directed graph that depicts a basis for the module as $q \rightarrow 0$.  
One representation-theoretic use for crystal graphs is in decomposing $\mathcal{U}_{q}(\mathfrak{g})$-modules and $\mathfrak{g}$-modules: The connected components of a crystal graph are in one-to-one correspondence with the irreducible modules in the decomposition of any given module.  
In \cite{KN}, crystal graphs are explicitly constructed for the $\myA_{n}$, $\myB_{n}$, $\myC_{n}$, and $\myD_{n}$ cases using tableaux. 

In \cite{Kash2}, Kashiwara establishes the existence and uniqueness of the crystal base -- and therefore of the crystal graph -- associated with an irreducible module (see also \cite{Jan} Ch.\ 9, \cite{HK} Ch.\ 5). 
Existence of the crystal base can be established by a product construction similar to the product construction of $R(\lambda)$ in \S \CrystallineSection. 
Uniqueness follows from a certain kind of complete reducibility for crystal bases and a Schur-type lemma for $\mathcal{U}_{q}(\mathfrak{g})$-modules. 
We will not carry out those arguments here. 
We refer the interested reader (that elusive creature) to the original sources, or \cite{Jan} or \cite{HK}. 

Our main purpose here is to record in \MainCrystalGraphTheorem\ that, for a given Weyl bialternant $\chi_{_{\lambda}}$, the associated crystal graph is the crystalline splitting poset $R(\lambda)$.  
This has several consequences. 
First, because irreducible $\mathfrak{g}$-modules are $q=1$ instances of irreducible $\mathcal{U}_{q}(\mathfrak{g})$-modules, then we get Weyl's character formula as an immediate corollary (\WeylCorollary). 
Second, we can say how certain operations on $\mathcal{U}_{q}(\mathfrak{g})$-modules are nicely behaved, as are the corresponding operations on crystal bases and crystal graphs (\CrystalOpsCorollary).  
Third, uniqueness of the crystal base allows us to conclude in \UniquenessCorollary\ that the construction of $R(\lambda)$ in \S \CrystallineSection\ is independent of the $\Omega$-expression used to represent $\lambda$, cf.\ \CrystallineProblem.A. 

The following synopsis of definitions and standard results from the theory of crystal bases largely follows Ch.\ 4, 5, and 9 of \cite{Jan} and is used to set up the sequence of results that lead to \MainCrystalGraphTheorem\ and its corollaries. 
We have $\mathfrak{g}$ as the aforementioned rank $n$ complex 
semisimple Lie algebra 
associated to the root system $\Phi$. Let $q$ be an indeterminate
over $\mathbb{Q}$. 
For $1 \leq i \leq n$, set $q_{i} :=
q^{\langle{\alpha_{i}},{\alpha_{i}}\rangle/2}$. For any 
nonnegative integer $m$, 
set $[m]_{i} := {\displaystyle
\frac{q_{i}^{m}-q_{i}^{-m}}{q_{i}-q_{i}^{-1}}}$, 
and set $[m]_{i}! := 
[m]_{i}[m-1]_{i}\cdots[1]_{i}$ (where an empty product is 1).  
For $0 \leq j \leq m$, the quantity $\qbin{m}{j}{i} := 
\frac{[m]_{i}!}{[j]_{i}![m-j]_{i}!}$ is a $q_{i}$-binomial coefficient. 
The {\em quantized enveloping algebra} 
$\mathcal{U} := \mathcal{U}_{q}(\mathfrak{g})$ is the
$\mathbb{Q}(q)$-algebra with generators $\{e_{i}, f_{i}, k_{i},
k^{-1}_{i}\}_{i=1}^{n}$ satisfying the quantum versions of the
Serre relations:   
\begin{center}
\parbox{4in}{
(QS1) $k_{i}k^{-1}_{i} = 1 = k^{-1}_{i}k_{i}$ and $k_{i}k_{j} =
k_{j}k_{i}$}\\ %
\vspace*{0.1in}
\parbox{4in}{
(QS2) $e_{i}f_{j} - f_{j}e_{i} = \delta_{i,j}{\displaystyle
\frac{k_{i}-k^{-1}_{i}}{q_{i}-q^{-1}_{i}}}$}\\ %
\vspace*{0.1in}
\parbox{4in}{
(QS3) \parbox[t]{3.5in}{$k_{i}e_{j}k^{-1}_{i} =
q^{\innprod{\alpha_{i}}{\alpha_{j}}}e_{j}$\\ %
$k_{i}f_{j}k^{-1}_{i} =
q^{-\innprod{\alpha_{i}}{\alpha_{j}}}f_{j}$}}\\ %
\vspace*{0.1in}
\parbox{4in}{ (QS$_{ij}^{+}$) $\displaystyle
\sum_{r=0}^{1-\innprod{\alpha_{j}}{\alpha_{i}^{\vee}}}
(-1)^{r}\qbin{1-\innprod{\alpha_{j}}{\alpha_{i}^{\vee}}}{r}{i}
e^{1-\innprod{\alpha_{j}}
{\alpha_{i}^{\vee}}-r}_{i}e_{j}e^{r}_{i} = 0$}\\ %
\vspace*{0.1in}
\parbox{4in}{
(QS$_{ij}^{-}$) $\displaystyle
\sum_{r=0}^{1-\innprod{\alpha_{j}}{\alpha_{i}^{\vee}}}
(-1)^{r}\qbin{1-\innprod{\alpha_{j}}{\alpha_{i}^{\vee}}}{r}{i}
f^{1-\innprod{\alpha_{j}}{\alpha_{i}^{\vee}}-r}_{i}f_{j}f^{r}_{i}
= 0$.}
\end{center}
These relations are $q$-deformations of the relations defining 
the universal enveloping algebra $\mathcal{U}(\mathfrak{g})$. 
For $\emptyset \not= J \subseteq I$, let $\mathcal{U}_{J}$ be the 
subalgebra of $\mathcal{U}$ generated by 
$\{e_j,f_j,k_j,k_{j}^{-1}\}_{j \in J}$.  

\renewcommand{\thefootnote}{2}Given a $\mathbb{Q}(q)$-vector space and 
$\mathcal{U}$-module 
$V$, we use $E_i$, $F_i$, and $K_{i}$ for
the images of $e_i$, $f_i$, and $k_{i}$ in $\mathfrak{gl}(V)$. 
The $\mathcal{U}$-module $V$ is {\em integrable} if $V =
\bigoplus_{\mu \in \Lambda} V_{\mu}$, where
$V_{\mu} = \{v \in V : K_{i}(v) = 
q_{i}^{\innprod{\mu}{\alpha_{i}^{\vee}}}v\
\mbox{for all}\ 1 \leq i \leq n\}$ is the $\mu$-weight
space of $V$.  For the remainder of this discussion  of crystal 
bases/crystal graphs, $V$ denotes an integrable finite-dimensional 
(i.f.d.\ for short) $\mathcal{U}$-module, and $\mbox{char}(V) := 
\sum_{\mu \in \Lambda}(\dim V_{\mu})e^{\mu}$.  

Certain operations on such $\mathcal{U}$-modules, particularly taking tensor products and duals, are effected by a Hopf algebra structure on $\mathcal{U}$, cf.\ \S 9.13 of \cite{Jan}. 
In particular, the category of i.f.d.\ $\mathcal{U}$-modules is closed under $\oplus$, $\otimes$, and $^{*}$, cf.\ \S 5.3 of \cite{Jan}. 
These operations combine in the expected ways: 
\begin{center}
$V_{1} \otimes (V_{2} \otimes V_{3}) \cong (V_{1} \otimes V_{2}) \otimes V_{3}$,\\ 
$V_{1} \otimes (V_{2} \oplus V_{3}) \cong (V_{1} \otimes V_{2}) \oplus (V_{1} \otimes V_{3})$, $(V_{1} \oplus V_{2}) \otimes V_{3} \cong (V_{1} \otimes V_{3}) \oplus (V_{2} \otimes V_{3})$,\\ 
$(V_{1} \oplus V_{2})^{*} \cong V_{1}^{*} \oplus V_{2}^{*}$, and $(V_{1} \otimes V_{2})^{*} \cong V_{2}^{*} \otimes V_{1}^{*}$, 
\end{center}
although it is not apparent yet that $\otimes$ is commutative. 

Moreover, i.f.d.\ $\mathcal{U}$-modules are completely reducible (\S 5.17 \cite{Jan})\footnote{Some of the theory of i.f.d.\ $\mathcal{U}$ -moldules is developed in Ch.\ 5 of \cite{Jan} using the usual $q=1$ theory of finite-dimensional $\mathfrak{g}$-modules over $\mathbb{C}$. 
In particular, the proof of Theorem 5.15 of \cite{Jan}, upon which this complete reducibility result depends, invokes Weyl's character formula, but not in a crucial way. 
When WCF is invoked, the issue in the argument is to deduce that two particular i.f.d.\ $\mathcal{U}$-modules $\widetilde{L}(\lambda)$ and $L(\lambda)$ have the same dimension by considering their $q=1$ counterpart, which in both cases is an irreducible finite-dimensional $\mathfrak{g}$-module $\overline{V}(\lambda)$ with highest weight $\lambda$, cf.\ \S 5.14 of \cite{Jan}. 
But in this case we can just use \S 5.13.3 of \cite{Jan} to conclude that $\displaystyle \dim \widetilde{L}(\lambda) = \sum_{\mu \in \Lambda}\dim \widetilde{L}(\lambda)_{\mu} = \sum_{\mu \in \Lambda}\dim \overline{V}(\lambda)_\mu = \sum_{\mu \in \Lambda}\dim L(\lambda)_{\mu} = \dim L(\lambda)$.}, and the irreducible i.f.d.\  $\mathcal{U}$-modules are indexed by dominant weights (\S 5.10 \cite{Jan}). 
The latter means in particular that if $V(\lambda)$ is an irreducible i.f.d.\ $\mathcal{U}$-module corresponding to dominant weight $\lambda$, then there is a ``highest'' weight vector $v_{\lambda}$ (unique up to scalar multiple) such that $E_{i}(v_{\lambda}) = 0$ for all $i \in I$.  
Throughout the remaining discussion, the notation $V(\lambda)$ refers to some generic irreducible i.f.d.\ $\mathcal{U}$-module with highest weight $\lambda$. 

Given an i.f.d.\ $\mathcal{U}$-module $V$ associated to a $\mathbb{Q}(q)$-algebra homomorphism $\phi: \mathcal{U} \longrightarrow \mathfrak{gl}(V)$, we can define a new $\mathcal{U}$-module structure on $V$ via a new $\mathbb{Q}(q)$-algebra homomorphism $\phi^{\bowtie}: \mathcal{U} \longrightarrow \mathfrak{gl}(V)$ which acts on generators for $\mathcal{U}$ as follows: 
For all $i \in I$ and $v \in V$, declare that $\phi^{\bowtie}(e_{i})(v) := \phi(f_{\sigma_{_{0}}(i)})(v)$, $\phi^{\bowtie}(f_{i})(v) := \phi(e_{\sigma_{_{0}}(i)})(v)$, $\phi^{\bowtie}(k_{i})(v) := \phi(k^{-1}_{\sigma_{_{0}}(i)})(v)$, and $\phi^{\bowtie}(k^{-1}_{i})(v) := \phi(k_{\sigma_{_{0}}(i)})(v)$, where $\sigma_{0}$ is the permutation of $I$ associated with the longest Weyl group element $w_0$, cf.\ \S \WeylSection. 
It is easy to see that the images $\{\phi^{\bowtie}(e_{i}), \phi^{\bowtie}(f_{i}), \phi^{\bowtie}(k_{i}), \phi^{\bowtie}(k^{-1}_{i})\}$ in $\mathfrak{gl}(V)$ satisfy the quantum Serre relations. 
We refer to this new $\mathcal{U}$-module by the notation $V^{\bowtie}$ and call this the {\em bow tie} of $V$. 
It is evident that the category of i.f.d.\ $\mathcal{U}$-modules is closed under $^{\bowtie}$. 

For any dominant weight $\lambda$, let $\overline{V}(\lambda)$ be an irreducible f.d.\ $\mathfrak{g}$-module with highest weight $\lambda$.  It follows from \S 5.13 of \cite{Jan} that $\mbox{char}(\overline{V}(\lambda)) = \mbox{char}(\overline{V}(\lambda))$. 
Let $\Pi(\overline{V}(\lambda))$ be the set of weights $\mu \in \Lambda$ with $\dim V(\lambda)_{\mu} > 0$, and similarly define $\Pi(V(\lambda))$. 
It follows from Proposition 21.3 of \cite{Hum} that $\Pi(\overline{V}(\lambda)) = \Pi(\lambda)$ (an equality of sets), so therefore $\Pi(V(\lambda)) = \Pi(\lambda)$ as well.

\renewcommand{\thefootnote}{3}Fix $i \in I$.  Kashiwara's ``raising'' and ``lowering'' operators $\widetilde{E}_{i}$ and $\widetilde{F}_{i}$ are defined on $V$ as follows. 
An analysis of the decomposition of $V$ as a $\mathcal{U}_{i}$-module shows that for each $\mu \in \Lambda$ and each $v \in V_{\mu}$, there exist unique vectors\footnote{Let $\{u_{1},\ldots,u_{s}\}$ be the collection of $\mathcal{U}_{i}$-highest vectors in the decomposition of $V$ as a $\mathcal{U}_{i}$-module. 
So for $1 \leq p \leq s$ there is a  nonnegative integer $\mysmalll_{p}$ for which  $K_{i}(u_{p}) = q_{i}^{\mysmallerl_{p}}u_{p}$, and $\{\frac{1}{[j]_{i}!}F_{i}^{j}(u_{p})\}_{1 \leq p \leq s, 0 \leq j \leq \mysmallerl_{p}}$ is a weight basis for $V$.  
Then there exist unique scalars $c_{p} \in \mathbb{Q}(q)$ for which $v = \sum c_{p}\frac{1}{[j_{p}]_{i}!}F_{i}^{j_{p}}(u_{p})$, where the sum is over all indices $p$ for which $1 \leq p \leq s$, $\mysmalll_{p}-2j_{p}=\langle \mu,\alpha_{i}^{\vee} \rangle$, and $0 \leq j_{p} \leq 
\mysmalll_{p}$.  
Then for any $j \geq \max\{0,-\langle \mu,\alpha_{i}^{\vee} \rangle\}$, we take $v_{j} = \sum c_{p}u_{p}$, where the sum is over all indices $p$ for which $\mysmalll_{p}-2j=\langle \mu,\alpha_{i}^{\vee} \rangle$.} $\{v_{j}\}_{j \geq \max\{0,-\langle \mu,\alpha_{i}^{\vee} \rangle\}}$ 
such that $v_{j} \in V_{\mu+j\alpha_{i}}\cap \ker{E_i}$ and 
\[v = \sum_{j \geq \max\{0,-\langle \mu,\alpha_{i}^{\vee} \rangle\}}\frac{1}{[j]_{i}!}F_{i}^{j}(v_{j}).\] 
Then define $\widetilde{E}_{i},\widetilde{F}_{i}: V_{\mu} \longrightarrow V$ by the rules 
\[\widetilde{E}_{i}(v) := \sum_{j \geq \max\{1,-\langle \mu,\alpha_{i}^{\vee} \rangle\}} \frac{1}{[j-1]_{i}!}F_{i}^{j-1}(v_{j}) \hspace*{0.2in}\mbox{and}\hspace*{0.2in} \widetilde{F}_{i}(v) := \sum_{j \geq \max\{0,-\langle \mu,\alpha_{i}^{\vee} \rangle\}} \frac{1}{[j+1]_{i}!}F_{i}^{j+1}(v_{j}).\] 
It is not hard to see that each of $\widetilde{E}_{i}$ and $\widetilde{F}_{i}$ is linear, and since $V$ is integrable we may extend these to linear transformations $V \longrightarrow V$. 

Let $A$ be the local ring of
rational functions in $\mathbb{Q}(q)$ well-defined at $q=0$: $A$ is a 
principal ideal domain with unique maximal ideal $qA$ and residue 
field $A/qA \cong \mathbb{Q}$.  A
{\em crystal base} for $V$ is a pair $\CB{V}{B}$, with $\mathcal{V}$ 
an $A$-submodule of $V$ and with the induced action of 
each $\widetilde{E}_{i}$ and $\widetilde{F}_{i}$ 
on $\mathcal{V}$ and $\mathcal{V}/q\mathcal{V}$, 
satisfying:
\begin{center}
\parbox{5.5in}{
(CB1) $\mathcal{V}$ is a finitely generated $A$-module   
and generates $V$ as a $\mathbb{Q}(q)$-vector space.}\\ %
\parbox{5.5in}{
(CB2) $\mathcal{B}$ is a basis for the $\mathbb{Q}$-vector space
$\mathcal{V}/q\mathcal{V}$.}\\ %
\parbox{5.5in}{
(CB3) $\mathcal{V} = \bigoplus_{\mu \in \Lambda}
\mathcal{V}_{\mu}$, where $\mathcal{V}_{\mu} := \mathcal{V} \cap
V_{\mu}$.}\\ %
\parbox{5.5in}{
(CB4) $\mathcal{B} = \disjointunion_{\mu \in \Lambda} 
\mathcal{B}_{\mu}$ (a disjoint union),
where $\mathcal{B}_{\mu} := \mathcal{B} \cap
(\mathcal{V}_{\mu}/q\mathcal{V}_{\mu})$.}\\ %
\parbox{5.5in}{
(CB5) $\widetilde{E}_{i}(\mathcal{V}) \subseteq \mathcal{V}$ and
$\widetilde{F}_{i}(\mathcal{V}) \subseteq \mathcal{V}$.}\\ %
\parbox{5.5in}{
(CB6) $\widetilde{E}_{i}(\mathcal{B}) \subseteq \mathcal{B} \cup \{0\}$
and $\widetilde{F}_{i}(\mathcal{B}) \subseteq \mathcal{B} 
\cup \{0\}$.}\\ %
\parbox{5.5in}{
(CB7) For $\selt$ and $\telt$ in $\mathcal{B}$ and $i \in I$, we have 
$\widetilde{E}_{i}(\selt) = \telt$ if and only if
$\widetilde{F}_{i}(\telt) = \selt$.}
\end{center}
We take the existence of crystal bases for granted, see for example Theorems 9.11/9.25 of \cite{Jan}. 

If $\CB{V}{B}$ is a crystal base for $V$, then its {\em crystal graph} $\mathcal{G}$ is the edge-colored directed graph whose vertices correspond to the elements of $\mathcal{B}$ and whose edges are defined by $\selt \myarrow{i} \telt$ if and only if $\widetilde{E}_{i}(\selt) = \telt$ if and only if $\selt = \widetilde{F}_{i}(\telt)$.  
For any $\xelt$ in $\mathcal{G}$, we define $wt_{\mbox{\tiny crystal}}(\xelt) := \mu$ (the {\em crystal weight} of $\xelt$) if and only if $\xelt \in \mathcal{B}_{\mu}$. 
It is not obvious at this point that any such $V$ has a unique (up to isomorphism) crystal graph $\mathcal{G}$. 
This is part of the content of the next result.

\noindent 
{\bf \ConnectedObservation}\ \ 
{\sl (1) If $V$ and $V'$ are irreducible i.f.d.\ $\mathcal{U}$-modules with highest weight $\lambda$, and if $\CB{V}{B}$ and $(\mathcal{V}',\mathcal{B}')$ are crystal bases with crystal graphs $\mathcal{G}$ and $\mathcal{G}'$ respectively, then $\mathcal{G} \cong \mathcal{G}'$.  Moreover, $\mathcal{G}$ has a unique maximal element $\melt$, and $wt_{\mbox{\tiny crystal}}(\melt) = \lambda$. 
(2) Suppose an i.f.d.\ $\mathcal{U}$-module $V$ decomposes as $V(\lambda_1) \oplus \cdots \oplus V(\lambda_{p})$, and suppose $\mathcal{G}$ is the crystal graph of some crystal base for $V$. 
Then $\mathcal{G} \cong \mathcal{G}_{1} \oplus \cdots \oplus \mathcal{G}_{p}$, where for each $1 \leq q \leq p$, $\mathcal{G}_{q}$ is the crystal graph for some crystal base for $V(\lambda_{q})$.} 

{\em Proof.} For (1), the fact that $\mathcal{G} \cong \mathcal{G}'$ follows from remarks at the end of \S 9.11 of \cite{Jan}. 
Uniqueness of the maximal element $\melt$ as well as the fact that $wt_{\mbox{\tiny crystal}}(\melt) = \lambda$ follow from Lemma 9.26 of \cite{Jan}. 
In view of (1), (2) then follows from Theorem 9.11 of \cite{Jan}.\hfill\QED  

From here on, then, we can speak of {\em the} crystal graph $\mathcal{G}(V)$ associated to any i.f.d.\ $\mathcal{U}$-module $V$. 
For an irreducible i.f.d.\ $\mathcal{U}$-module $V(\lambda)$, we denote by $\CG(\lambda)$ the associated crystal graph.  
Our main aim is to show in \MainCrystalGraphTheorem\ that $\CG(\lambda)$ 
coincides with the crystalline splitting poset $R(\lambda)$ of \S \CrystallineSection. 
To get there, we continue with some general observations about crystal graphs. 

\noindent 
{\bf \CharacterObservation}\ \ 
{\sl Let $(\mathcal{V},\mathcal{B})$ be a crystal base for an i.f.d.\ $\mathcal{U}$-module $V$, with $\mathcal{G}$ as the corresponding crystal graph. 
Then $|\mathcal{B}| = \dim_{\mathbb{Q}(q)}V$, and for each $\mu \in \Lambda$, we have we have $|\mathcal{B}_{\mu}| = \dim_{\mathbb{Q}(q)}V_{\mu}$. 
Moreover, $\mbox{char}(V) = \sum_{\xelt \in \mathcal{G}}e^{wt_{\mbox{\tiny crystal}}(\xelt)}$.}

{\em Proof.}  
As an $A$-submodule of $V$, $\mathcal{V}$ is necessarily torsion-free.  
So by the fundamental theorem of finitely-generated modules over a P.I.D., it follows that $\mathcal{V} \cong A^{d}$ for some positive integer $d$. 
Then, $\mathcal{V}/q\mathcal{V} \cong (A/qA)^{d} \cong \mathbb{Q}^{d}$.  
In extending scalars, one can use (CB1) to see that the linear transformation $\mathcal{V} \otimes_{A} \mathbb{Q}(q) \longrightarrow V$ induced by $v \otimes c \longmapsto cv$ is an isomorphism of $\mathbb{Q}(q)$-vector spaces. 
Then $\mathcal{V} \otimes_{A} \mathbb{Q}(q) \cong (\mathbb{Q}(q))^{d}$. 
Therefore, $d = |\mathcal{B}| = \dim_{A}\mathcal{V} = \dim_{\mathbb{Q}(q)}V$.  
Similarly, for each $\mu \in \Lambda$, we have $|\mathcal{B}_{\mu}| = \dim_{A}\mathcal{V}_{\mu} = \dim_{\mathbb{Q}(q)}V_{\mu}$.  
So $\sum_{\xelt \in \mathcal{G}}e^{wt_{\mbox{\tiny crystal}}(\xelt)} = \sum_{\mu \in \Lambda}(\dim V_{\mu})e^{\mu} = \mbox{char}(V)$.\hfill\QED  

\noindent
{\bf \CrystalOpsLemma}\ \ {\sl Given i.f.d.\ $\mathcal{U}$-modules $V$, $V_1$, and $V_2$. 
(1) Then $\mathcal{G}(V_{1} \oplus V_{2}) \cong \mathcal{G}(V_{1}) \oplus \mathcal{G}(V_{2})$, $\mathcal{G}(V^{*}) \cong \mathcal{G}(V)^{*}$, $\mathcal{G}(V^{\bowtie}) \cong \mathcal{G}(V)^{\bowtie}$, and $\mathcal{G}(V_{1} \otimes V_{2}) \cong \mathcal{G}(V_{1}) \otimes \mathcal{G}(V_{2})$. 
(2) Let $\emptyset \not= J \subseteq I$. 
Under the induced action, the $\mathcal{U}$-module $V$ is a $\mathcal{U}_{J}$-module, which we denote by $V^{J}$. 
A crystal base $\CB{V}{B}$ for $V$ is also a crystal base for $V^{J}$. 
Moreover, if $\mathcal{G}_{J}$ is the union of all $J$-components of $\mathcal{G}(V)$, then $\mathcal{G}_{J} \cong \mathcal{G}(V^{J})$ and the crystal weight $wt_{\mbox{\tiny crystal}}^{J}(\xelt)$ of any $\xelt \in \mathcal{G}_{J}$ is just $\sum_{j \in J}\langle wt_{\mbox{\tiny crystal}}(\xelt),\alpha_{j}^{\vee} \rangle\omega_{j}^{J}$.} 

{\em Proof.} 
The $\oplus$, $^{*}$, and $^{\bowtie}$ parts of (1) follow from the definitions.  
Below we work out the details for the $^{*}$ case.  
The $^{\bowtie}$ and $\oplus$ cases are similar, for the latter also see \S 9.4 of \cite{Jan}. 
For $\otimes$ we require a little more ammunition, in particular we appeal to Theorem 9.17 of \cite{Jan}. 

To show that $\mathcal{G}(V^{*}) \cong \mathcal{G}(V)^{*}$, we need to say how to produce a crystal base $(\mathcal{V},\mathcal{B})^{*} = (\mathcal{V}^{*},\mathcal{B}^{*})$ for the i.f.d.\ $\mathcal{U}$-module $V^{*}$. 
To that end, we let $\{w_1, \ldots , w_d\}$ be a weight basis for $V$ and $\{\gelt_1, \ldots , \gelt_d\}$ the corresponding basis for $V^{*}$, so $\gelt_j(w_k) = \delta_{jk}$. 
If $\varphi: V \longrightarrow V^{*}$ is the $\mathbb{Q}(q)$-vector space isomorphism induced by $v_j \mapsto \gelt_j$, then the restriction $\varphi|_{\mathcal{V}}$ is an $A$-module isomorphism $\mathcal{V} \longrightarrow \mathcal{V}^{*}$, where $\mathcal{V}^{*} := \varphi|_{\mathcal{V}}(\mathcal{V}) \subset V^{*}$. 
Since $\mathcal{V}$ is a finitely-generated $A$-submodule of $V$, then $\mathcal{V}^{*}$ is a finitely-generated $A$-submodule of $V^{*}$, cf.\ (CB1). 

Consider the following diagram: 
\[\begin{array}{ccc}
\mathcal{V} & \stackrel{\varphi|_{_{\mathcal{V}}}}{\longrightarrow} & \mathcal{V}^{*}\\
\setlength{\unitlength}{1cm}
\begin{picture}(0.2,0.8)
\put(0.1,.7){\vector(0,-1){0.7}}
\put(-0.3,.3){$\pi$}
\end{picture}
 & &\setlength{\unitlength}{1cm}
\begin{picture}(0.2,0.8)
\put(0.1,.7){\vector(0,-1){0.7}}
\put(0.3,.3){$\pi^{*}$}
\end{picture}
\\
\mathcal{V}/q\mathcal{V} & \stackrel{\tilde{\varphi}}{\longrightarrow} & \mathcal{V}^{*}/q\mathcal{V}^{*}  
\end{array}\]
where $\pi$ and $\pi^{*}$ are natural projection maps and $\tilde{\varphi}$ is induced. 
Since $\ker(\pi^{*} \circ \varphi|_{\mathcal{V}}) = \ker(\pi)$, then $\tilde{\varphi}$ is injective, and it is surjective since $\phi \circ \pi = \pi^{*} \circ \varphi|_{\mathcal{V}}$ is surjective. 
That is, $\tilde{\varphi}$ is a $\mathbb{Q}$-vector space isomorphism. 
We define $\mathcal{B}^{*} := \tilde{\varphi}(\mathcal{B})$. 
Then $\mathcal{B}^{*}$ is a $\mathbb{Q}$-basis for $\mathcal{V}^{*}/q\mathcal{V}^{*}$, cf.\ (CB2). 

The action of $\mathcal{U}$ on $V^{*}$ is prescribed as follows: for any $x \in \mathcal{U}$ and $\felt \in V^{*}$, the linear functional $x.\felt \in V^{*}$ is determined by $(x.\felt)(v) = \felt(S'(x).v)$, where $S'$ is the antipode of \S 9.13 of \cite{Jan}. 
Let $\mu_k$ be the weight of basis vector $w_k$, for any $1 \leq k \leq d$. 
We can see directly that $\{\gelt_1, \ldots , \gelt_d\}$ is a weight basis for $V^{*}$: $(k_i.\gelt_j)(w_k) = \gelt_j(S'(k_i).w_k) = \gelt_j(k_i^{-1}.w_k) = \gelt_j(q_i^{\langle -\mu_k,\alpha_i^{\vee}\rangle}w_k) = q_i^{\langle -\mu_j,\alpha_i^{\vee}\rangle}$. 
So, $k_i.\gelt_j = q_i^{\langle -\mu_j,\alpha_i^{\vee}\rangle}\gelt_j$. 
It follows that $V^{*} = \bigoplus_{\mu \in \Lambda} V^{*}_{\mu}$ with $V^{*}_{\mu} = \varphi(V_{-\mu})$. 
Using (CB3) and (CB4) for $V$ together with the maps $\varphi|_{\mathcal{V}}$ and $\tilde{\varphi}$, we immediately get (CB3) and (CB4) for $V^{*}$. 

For (2), consider $\emptyset \not= J \subseteq I$.  
Under the induced action, the $\mathcal{U}$-module $V$ is a $\mathcal{U}_{J}$-module, which we denote by $V^{J}$. 
One can see directly that $\CB{V}{B}$ is a crystal base for $V^{J}$: 
(CB1), (CB2), and (CB5-7) are automatic, and both (CB3) and (CB4) follow from the fact that for any $\mu \in \Lambda_{\Phi_J}$, the $\mu$-weight space $V^{J}_{\mu} = \bigoplus V_{\nu}$, where this sum is over all $\nu \in \Lambda$ for which  $\langle \nu,\alpha_{j}^{\vee} \rangle = \langle \mu,\alpha_{j}^{\vee} \rangle$ for all $j \in J$.  
So, the crystal graph $\mathcal{G}_{J}$ for $V^{J}$ is obtained from $\mathcal{G}$ by removing from $\mathcal{G}$ all edges of color $i$ with $i \not\in J$, but of course keeping all vertices of $\mathcal{G}$.  
Moreover, the crystal weight for any vertex $\xelt$ in $\mathcal{G}_{J}$ is therefore just $wt_{\mbox{\tiny crystal}}^{J}(\xelt) = \sum_{j \in J}\langle wt_{\mbox{\tiny crystal}}(\xelt),\alpha_{j}^{\vee} \rangle\omega_{j}^{J}$.\hfill\QED  

\noindent {\bf \CrystalProp}\ \ 
{\sl Let $V$ be an i.f.d.\ $\mathcal{U}$-module with associated crystal graph $\mathcal{G}$. 
For all $\xelt \in \mathcal{G}$, we have} $wt_{\mbox{\tiny crystal}}(\mathbf{\xelt}) = wt(\xelt)$. 
{\sl The crystal graph $\mathcal{G}$ is the Hasse diagram for a ranked poset, and as an edge-colored poset is fibrous, ${\Phi}$-structured, and strongly untangled. 
Moreover,} $\WGF(\mathcal{G}) = \mbox{char}(V)$. 

{\em Proof.}  Let $\CB{V}{B}$ be a crystal base for $V$. 
Now, from (QS3) it follows that for all $i \in I$ and $v \in V_{\mu}$, we have $E_{i}(v) \in V_{\mu+\alpha_{i}}$ and $F_{i}(v) \in V_{\mu-\alpha_{i}}$.  
From this and the definition of  $\widetilde{E}_{i}$, we see that $\selt \myarrow{i} \telt$ in $\mathcal{G}$ implies that 
\begin{equation}
wt_{\mbox{\tiny crystal}}(\selt) + \alpha_i =
wt_{\mbox{\tiny crystal}}(\telt).
\end{equation} 
From this it follows that $\mathcal{G}$ is acyclic (no directed cycles), that the transitive closure of the edge relations in $\mathcal{G}$ defines a partial order, that the given edge relations are precisely the covering relations relative to this partial order, and that as a poset $\mathcal{G}$ is ranked.  

Now take $\xelt \in \mathcal{G}$, and let $i \in I$. 
As in \CrystalOpsLemma.2, take $\{i\} \subseteq I$ and consider the subalgebra $\mathcal{U}_{i}$ acting on $V$. 
Then $\mathcal{G}_i$ is the crystal graph for the $\mathcal{U}_i$-module $V$.  
Then $\mathcal{C} := \comp_{i}(\xelt)$ is one of the connected components of $\mathcal{G}_{i}$, necessarily a chain by (CB6). 
That is, $\mathcal{G}$ is fibrous. 
The irreducible i.f.d.\ $\mathcal{U}_{i}$-module with crystal graph $\mathcal{C}$ must have dimension $l_{i}(\xelt)+1$, cf.\ \TwoObservations. 
Let $\melt_{i}$ be the maximal element of the chain $\mathcal{C}$, and let $\nu_i := wt_{\mbox{\tiny crystal}}(\melt_{i})$. 
From the explicit description of crystal bases for irreducible $\mathcal{U}_{i}$-modules given in Lemma 9.6 of \cite{Jan}, it follows that $wt_{\mbox{\tiny crystal}}^{\{i\}}(\melt_i) = l_{i}(\xelt)$.  
So by \CrystalOpsLemma.2, $l_{i}(\xelt) = \langle \nu_i,\alpha_{i}^{\vee} \rangle$. 
Since there are $l_{i}(\xelt)-\rho_{i}(\xelt)$ edges in the path $\xelt \myarrow{i} \cdots \myarrow{i} \melt_i$ from $\xelt$ to $\melt_i$ in the chain $\mathcal{C}$, then we have $wt _{\mbox{\tiny crystal}}(\xelt) = wt _{\mbox{\tiny crystal}}(\melt_i)-(l_{i}(\xelt)-\rho_{i}(\xelt))\alpha_{i} = \nu_i-(l_{i}(\xelt)-\rho_{i}(\xelt))\alpha_{i}$, by (7).  

So, $\displaystyle wt_{\mbox{\tiny crystal}}(\xelt) =\sum_{i=1}^{n}\innprod{wt_{\mbox{\tiny crystal}}(\xelt)}{\alpha_{i}^{\vee}}\omega_i = \sum_{i=1}^{n}\innprod{\nu_i-(l_{i}(\xelt)-\rho_{i}(\xelt))\alpha_{i}}{\alpha_{i}^{\vee}}\omega_i = \sum_{i=1}^{n}\Big(\innprod{\nu_i}{\alpha_{i}^{\vee}} - \innprod{(l_{i}(\xelt)-\rho_{i}(\xelt))\alpha_{i}}{\alpha_{i}^{\vee}}\Big)\omega_i = \sum_{i=1}^{n}\Big(l_{i}(\xelt) - 2(l_{i}(\xelt)-\rho_{i}(\xelt))\Big)\omega_i = \sum_{i=1}^{n}\Big(2\rho_{i}(\xelt) - l_{i}(\xelt)\Big)\omega_i = wt(\xelt)$. 
So $wt_{\mbox{\tiny crystal}}(\xelt) = wt(\xelt)$ for all $\xelt \in \mathcal{G}$.
It follows from \CharacterObservation\ that $\WGF(\mathcal{G}) = \mbox{char}(V)$.

Let $J \not= \emptyset$ be any subset of $I$, and consider the crystal graph $\mathcal{G}_{J}$ from \CrystalOpsLemma.2. 
It is a consequence of \ConnectedObservation\ that any connected component of $\mathcal{G}_{J}$ has a unique maximal element. 
So $\mathcal{G}$ is strongly untangled.\hfill\QED

\noindent 
{\bf \CharacterCorollary}\ \ {\sl Let $V$, $V_{1}$ and $V_{2}$ be i.f.d.\ $\mathcal{U}$-modules. 
Then $\mbox{char}(V_1 \oplus V_2) = \mbox{char}(V_1) + \mbox{char}(V_2)$, $\mbox{char}(V^{*}) = \mbox{char}(V)^{*}$, $\mbox{char}(V^{\bowtie}) = \mbox{char}(V)^{\bowtie}$, and $\mbox{char}(V_1 \otimes V_2) = \mbox{char}(V_1)\, \mbox{char}(V_2)$.} 

{\em Proof.} We have $\mbox{char}(V_1 \oplus V_2) = \WGF(\mathcal{G}(V_{1} \oplus V_{2})) \cong \WGF(\mathcal{G}(V_{1}) \oplus \mathcal{G}(V_{2})) = \WGF(\mathcal{G}(V_{1})) + \WGF(\mathcal{G}(V_{2})) = \mbox{char}(V_1) + \mbox{char}(V_2)$, applying \OpsProps\ and \WGFLemma.  
Argue similarly for the other operations, using \WGFCorollary\ for the $\otimes$ case.\hfill\QED

\noindent 
{\bf \CrystalMinQuasiMin}\ \ {\sl If $\lambda \in \Lambda^{+}$ is minuscule or quasi-minuscule, then $\CG(\lambda) \cong R(\lambda)$.}

{\em Proof.} Follows from the proof of Lemma 9.6.b of \cite{Jan}.\hfill\QED

\noindent 
{\bf \MainCrystalGraphTheorem}\ \ {\sl Let $\lambda$ be any dominant weight. 
Then $\CG(\lambda) \cong R(\lambda)$.} 

{\em Proof.} Take the lexicographically minimal $\Omega$-expression for $\lambda$ from \S \CrystallineSection, so $\lambda = \mu_1 + \cdots + \mu_p$ with each $\mu_q$ minuscule or quasi-minuscule and in the $W$-orbit of some dominant minuscule or quasi-minuscule $\widehat{\mu_q}$ ($1 \leq q \leq p$). 
Then consider the $\mathcal{U}$-module $V := V(\widehat{\mu_1}) \otimes \cdots \otimes V(\widehat{\mu_p})$. 
So $\mathcal{G}(V) \cong \CG(\widehat{\mu_1}) \otimes \cdots \otimes \CG(\widehat{\mu_p})$, by \CrystalOpsLemma. 
Now $\CG(\widehat{\mu_1}) \otimes \cdots \otimes \CG(\widehat{\mu_p}) \cong R(\widehat{\mu_1}) \otimes \cdots \otimes R(\widehat{\mu_q})$ by \CrystalMinQuasiMin. 
Also, $R(\lambda)$ is the connected component of $R(\widehat{\mu_1}) \otimes \cdots \otimes R(\widehat{\mu_q})$ identified in \CrystallineExistenceTheorem. 
Let $(\xelt_1,\ldots,\xelt_p)$ be the maximal element of $R(\lambda)$, with $wt(\xelt_1,\ldots,\xelt_p) = \lambda$. 
Then the corresponding element $\xelt$ of $\mathcal{G}(V)$ is therefore maximal, and its weight $wt(\xelt) = wt_{\mbox{\tiny crystal}}(\xelt)$ (cf.\ \CrystalProp) is therefore $\lambda$. 
By \ConnectedObservation, the connected component in $\mathcal{G}(V)$ of $\xelt$ is $\CG(\lambda)$. 
Then $\CG(\lambda) \cong R(\lambda)$. 
\hfill\QED

\noindent 
{\bf \WeylCorollary\ (Weyl's character formula)}\ \ {\sl For $\lambda \in \Lambda^{+}$, let $\overline{V}(\lambda)$ be an irreducible complex f.d.\ $\mathfrak{g}$-module with highest weight $\lambda$. 
Then} 
$\mbox{char}(\overline{V}(\lambda)) = \mbox{char}(V(\lambda)) = \chi_{_{\lambda}}$. 

{\em Proof.} By \S 5.13 of \cite{Jan}, $\mbox{char}(\overline{V}(\lambda)) = \mbox{char}(V(\lambda))$. 
It follows from \CrystalProp\ that $\mbox{char}(V(\lambda)) = \WGF(\CG(\lambda))$. 
By \MainCrystalGraphTheorem, $\WGF(\CG(\lambda)) = \WGF(R(\lambda))$. 
By \CrystallineExistenceTheorem, $\WGF(R(\lambda)) = \chi_{_{\lambda}}$.\hfill\QED 

\noindent 
{\bf \SymmetricCorollary}\ \ {\sl For any i.f.d.\ $\mathcal{U}$-module $V$, the character $\mbox{char}(V)$ is a Weyl symmetric function.} 

{\em Proof.} Write $V \cong V(\lambda_1) \oplus \cdots \oplus V(\lambda_p)$. 
Then by \CharacterAndWeylCorollaries, $\mbox{char}(V) = \mbox{char}(V(\lambda_1)) + \cdots + \mbox{char}(V(\lambda_p)) = \chi_{_{\lambda_1}} + \cdots + \chi_{_{\lambda_p}}$, which is a Weyl symmetric function.\hfill\QED

\noindent 
{\bf \CharacterInvariantCorollary}\ \ {\sl Let $V_1$ and $V_2$ be i.f.d.\ $\mathcal{U}$-modules. 
Then $V_1 \cong V_2$ if and only if $\mbox{char}(V_1) = \mbox{char}(V_2)$.}

{\em Proof.} We only need to show the converse. 
Suppose $\mbox{char}(V_1) = \mbox{char}(V_2)$. 
Write $V_1 \cong V(\lambda_1) \oplus \cdots \oplus V(\lambda_p)$ and $V_2 \cong V(\nu_1) \oplus \cdots \oplus V(\nu_q)$. 
Applying \CharacterCorollary\ and \WeylCorollary, we see that $\chi_{_{\lambda_1}} + \cdots + \chi_{_{\lambda_p}} = \mbox{char}(V_1) = \mbox{char}(V_2) = \chi_{_{\nu_1}} + \cdots + \chi_{_{\nu_r}}$. 
Since the Weyl bialternants are a basis for the ring of Weyl symmetric functions (cf.\ \BasisTheorem), then $p=r$ and with respect to a suitable re-ordering of the $\nu_q$'s, we have $\lambda_q = \nu_q$ for $1 \leq q \leq p$. 
This means that $V_1 \cong V_2$. 
\hfill\QED

\noindent 
{\bf \CrystalOpsCorollary}\ \ {\sl Let $V$, $V_{1}$ and $V_{2}$ be i.f.d.\ $\mathcal{U}$-modules. 
(1) Then $V^{\bowtie} \cong V$ and $V_1 \otimes V_2 \cong V_2 \otimes V_1$.  
(2) Moreover, $\mathcal{G}(V)^{\bowtie} \cong \mathcal{G}(V)$ and $\mathcal{G}(V_1) \otimes \mathcal{G}(V_2) \cong \mathcal{G}(V_2) \otimes \mathcal{G}(V_1)$.} 

{\em Proof.} For (1), note that $\mbox{char}(V^{\bowtie}) = \mbox{char}(V)^{\bowtie}$ by \CharacterCorollary. 
Since $\mbox{char}(V)$ is a Weyl symmetric function (\SymmetricCorollary), then by \StarProp\ we have $\mbox{char}(V)^{\bowtie} = \mbox{char}(V)$. 
Since we have $\mbox{char}(V^{\bowtie}) = \mbox{char}(V)$, then by \CharacterInvariantCorollary, $V^{\bowtie} \cong V$. 
A similar argument shows that $\mbox{char}(V_1 \otimes V_2) = \mbox{char}(V_2 \otimes V_1)$, from which it follows that $V_1 \otimes V_2 \cong V_2 \otimes V_1$. 
For part (2), note that $\mathcal{G}(V)^{\bowtie} \cong \mathcal{G}(V^{\bowtie})$ by \CrystalOpsLemma, and $\mathcal{G}(V^{\bowtie}) \cong \mathcal{G}(V)$ since $V^{\bowtie} \cong V$ by part (1). 
Therefore $\mathcal{G}(V)^{\bowtie} \cong \mathcal{G}(V)$. 
Similar reasoning shows that $\mathcal{G}(V_1) \otimes \mathcal{G}(V_2) \cong \mathcal{G}(V_2) \otimes \mathcal{G}(V_1)$. 
\hfill\QED

In \cite{Lenart}, Lenart uses the ``alcove path model'' to produce a combinatorially explicit isomorphism showing $\CG(\lambda)^{\bowtie} \cong \CG(\lambda)$.  
This isomorphism can be viewed as an involution of $\CG(\lambda)$.  
In the $\myA_{n}$ case, let $\lambda = \sum_{i \in I}a_{i}\omega_{i}$, and let $\mbox{shape}(\lambda)$ be the Ferrers board with $a_{n}$ columns of length $n$ followed by $a_{n-1}$ columns of length $n-1$ etc.  
Then $\CG(\lambda)$ is an ordering of the semistandard tableaux of $\mbox{shape}(\lambda)$ with entries from the set $\{1,2,\ldots,n+1\}$.  
In this case, Lenart's isomorphism $\CG(\lambda) \cong \CG(\lambda)^{\bowtie}$ is an involution of semistandard tableaux originally due to Sch\"{u}tzenberger, called evacuation. 

\noindent 
{\bf \UniquenessCorollary}\ \ {\sl (1) Suppose that $\mathcal{G}_{1},\ldots,\mathcal{G}_{p}$ are crystal graphs, or, alternatively, suppose that for each $q \in \{1,\ldots,p\}$, $\mathcal{G}_{q}$ is a disjoint sum of connected components from some crystal product of primary-plus posets. 
Suppose $R$ is a connected component of $\mathcal{G}_{1} \otimes \cdots \otimes \mathcal{G}_{p}$ with a maximal element $\melt$ of weight $\lambda := wt(\melt)$.  
Then $\lambda$ is dominant, and $R \cong R(\lambda)$. 
(2) Start with any $\Omega$-expression for a dominant $\lambda$ and build an admissible system $X'(\lambda)$ as in \LittelmannRemark.A; let $X(\lambda)$ be the admissible system resulting from the lexicographically minimal $\Omega$-expression for $\lambda$. 
Then $\ASG(X'(\lambda)) \cong \ASG(X(\lambda)) \cong R(\lambda)$.} 

{\em Proof.} 
For a fibrous and ${\Phi}$-structured poset, the weight of any maximal element is dominant. 
The proof of \MainCrystalGraphTheorem\ shows that every crystal graph $\CG(\lambda)$ can be realized as a connected component within some crystal product of minuscule or quasi-minuscule splitting posets. 
So the alternative hypothesis stated in part (1) of the corollary is equivalent to the initial hypothesis. 
\AssociativeCrystalLemma\ shows that for fibrous posets, $\otimes$ distributes over $\oplus$.  
So it suffices to assume that for each $1 \leq q \leq p$ we have $\mathcal{G}_{q} \cong \CG(\lambda_q)$ for some dominant $\lambda_q$. 
As a connected component of $\CG(\lambda_1) \otimes \cdots \otimes \CG(\lambda_q)$, 
$R$ is therefore isomorphic to $\CG(\lambda)$, cf.\ \ConnectedObservation. 
Both (1) and (2) now follow from \MainCrystalGraphTheorem. 
\hfill\QED

The main idea of the next result is to say how one can conclude that a poset is splitting if a subset of its edges forms a crystal graph. 
The various parts of the theorem consider various partial states of knowledge about certain crystal product subsets. 
If $S$ is any subset of a fibrous poset $\mathcal{G}$, say that $S$ is $\widetilde{E}${\em -preserving} if for all $\selt \in S$ and $i \in I$ we have $\widetilde{E}_{i}(\selt) \in S$ whenever $\widetilde{E}_{i}(\selt) \not= \theta$. 
Similarly define $\widetilde{F}${\em -preserving}. 

\noindent 
{\bf \ConfirmSplittingTheorem}\ \ {\sl Suppose that $\mathcal{G}_{1},\ldots,\mathcal{G}_{p}$ are crystal graphs, or, alternatively, suppose that for each $q \in \{1,\ldots,p\}$, $\mathcal{G}_{q}$ is a disjoint sum of connected components from some crystal product of primary-plus posets. 
Suppose also that $R$ is a subset of the elements of the crystal product $\mathcal{G}_{1} \otimes \cdots \otimes \mathcal{G}_{p}$.} 

{\sl (1E) Suppose that $R$ is $\widetilde{E}$-preserving and that $R$ contains exactly one minimal element from $\mathcal{G}_{1} \otimes \cdots \otimes \mathcal{G}_{p}$. 
Then $R$ is $\widetilde{F}$-preserving, contains a maximal element $\melt$ such that $\lambda := wt_{\mbox{\tiny crystal}}(\melt)$ is dominant, and the subgraph of $\mathcal{G}_{1} \otimes \cdots \otimes \mathcal{G}_{p}$ whose elements are from $R$ is isomorphic to $R(\lambda)$.} 

{\sl (1F) Suppose that $R$ is $\widetilde{F}$-preserving and that $R$ contains exactly one maximal element from $\mathcal{G}_{1} \otimes \cdots \otimes \mathcal{G}_{p}$. 
Then $R$ is $\widetilde{E}$-preserving, the weight $\lambda := wt_{\mbox{\tiny crystal}}(\melt)$ of the maximal element $\melt$ is dominant, and the subgraph of $\mathcal{G}_{1} \otimes \cdots \otimes \mathcal{G}_{p}$ whose elements are from $R$ is isomorphic to $R(\lambda)$.} 

{\sl (2) Suppose we have the hypotheses of either (1E) or (1F). 
In addition, suppose that $R$ is a ranked poset with edges colored by $I$ and that for all $\xelt \in R$, the weight $wt(\xelt)$ of $\xelt$ as an element of the edge-colored poset $R$ coincides with the crystal weight $wt_{\mbox{\tiny crystal}}(\mathbf{\xelt})$ of $\xelt$ as an element of the crystal product. 
Then} $\WGF(R) = \chi_{_{\lambda}}$. 

{\sl (3) Suppose in addition to (2) and (1E) or (1F) that $R$ is ${\Phi}$-structured. 
Then $R$ is a splitting poset for $\chi_{_{\lambda}}$.}

{\sl (4) Suppose that $R$ is a ranked poset with edges colored by $I$ and that in addition to (1E) (respectively, (1F)), we know that $\xelt \myarrow{i} \wE_{i}(\xelt)$ (respectively, $\wF_{i}(\xelt) \myarrow{i} \xelt$) is an edge in $R$ whenever $\wE_{i}(\xelt) \not= \theta$ (respectively, $\wF_{i}(\xelt) \not= \theta$) for some $\xelt \in R$ and $i \in I$. 
Then the subgraph of $R$ induced by taking all vertices and all edges of the form $\xelt \myarrow{i} \wE_{i}(\xelt)$ (respectively, $\wF_{i}(\xelt) \myarrow{i} \xelt$) is isomorphic to $R(\lambda)$, with $\lambda$ as in (1E) or (1F). 
Moreover, $R$ is connected and has $\melt$ from (1E) or (1F) as its unique maximal element.} 

{\sl (5) Suppose in addition to (4) that $R$ is ${\Phi}$-structured and that the weight $wt(\melt)$ of $\melt$ as an element of the edge-colored poset $R$ coincides with the crystal weight $wt_{\mbox{\tiny crystal}}(\mathbf{\melt})$ of $\melt$ as an element of the crystal product. 
Then for all $\xelt \in R$, $wt(\xelt) = wt_{\mbox{\tiny crystal}}(\mathbf{\xelt})$, and $R$ is a splitting poset for $\chi_{_{\lambda}}$. 
Moreover, each $i$-component of $R$ has an SCD as prescribed in \SCDLemma.} 

{\em Proof.} For (1E), the connected component in $\mathcal{G}_{1} \otimes \cdots \otimes \mathcal{G}_{p}$ of the minimal element $\nelt$ is isomorphic to some crystal graph $\CG(\lambda)$ by \ConnectedObservation, which is $R(\lambda)$ by \MainCrystalGraphTheorem. 
Let $S(\lambda)$ be the connected component of $\nelt$ in $\mathcal{G}_{1} \otimes \cdots \otimes \mathcal{G}_{p}$. 
Since $R(\lambda)$ is strongly untangled and $R(\lambda) \cong R(\lambda)^{\bowtie}$, then $S(\lambda)$ has a {\em unique} minimal element, namely $\nelt$. 
Therefore the elements of $\mathcal{G}_{1} \otimes \cdots \otimes \mathcal{G}_{p}$ that are reachable from $\nelt$ by sequences of $\wE_{i}$'s are precisely the elements of $S(\lambda)$. 
So $R = S(\lambda)$ set-wise. 
The remaining claims of (1E) now follow.  
The proof of (1F) is similar. 
For part (2), simply observe now that $\WGF(R) = \WGF(S(\lambda)) = \WGF(R(\lambda)) = \chi_{_{\lambda}}$. 
Part (3) follows from (2) and the definition of splitting. 
In view of (1E or F), the hypotheses of (4) tell us that $S(\lambda) \cong R(\lambda)$ is the subgraph of $R$ identified in the statement of part (4). 
In particular, $R$ is connected and has $\melt$ as its unique maximal element. 
For part (5), for any $\xelt \in R$, take a path $\xelt = \xelt_{0} \myarrow{i_1} \xelt_{1} \myarrow{i_2} \cdots \myarrow{i_r} \xelt_{r} = \melt$ in $S(\lambda)$ from $\xelt$ up to $\melt$. 
Then by \CrystalProp\ and the hypotheses of part (5), $wt_{\mbox{\tiny crystal}}(\xelt) = wt_{\mbox{\tiny crystal}}(\melt) - \alpha_{i_1} - \alpha_{i_2} - \cdots - \alpha_{i_r} = wt(\melt) - \alpha_{i_1} - \alpha_{i_2} - \cdots - \alpha_{i_r} = wt(\xelt)$. 
The splitting claim follows from part (3). 
The SCD claim follows from the cited lemma. 
\hfill\QED

We close this section with a daydream about symmetric chain decompositions, continuing a thought begun after \SCDLemma. 
The existence of an SCD is an open question for many very nicely described posets, for example the Gelfand-Tsetlin lattices of the form $L(\mathfrak{p},n+1)$, where the shape corresponding to $\mathfrak{p}$ is columnar. 
If one could realize $L(\mathfrak{p},n+1)$ as an $i$-component of a splitting poset $R$ meeting the criteria of \ConfirmSplittingTheorem.4/5, one would obtain an SCD. 
In fact, this can be done when the shape for $\mathfrak{p}$ is a column of length two, but of course an SCD in such a case is easy to see directly. 
This might not be feasible as an approach for other shapes $\mathfrak{p}$ or for other posets, but the idea does at least give a speculative algebraic and crystalline context for SCD problems.

\newpage
\noindent 
{\Large \bf \S \CriteriaSection.\ Splitting via vertex-coloring for non-fibrous posets.} 

Many of the candidate splitting posets to which we might apply \InitialSplittingTheorem\ (e.g.\ in \cite{DD} and \cite{DDgeneral}) have the property that all their $i$-components are products of chains.  
In these situations, some special versions of \InitialSplittingTheorem\ formulated below as \NewSplittingTheorem/\ChainProductVersions\ can be helpful.  
We apply these results in this section to re-obtain splitting distributive lattices for the $\myA_{n-1}$-Weyl bialternants; these are the famous Gelfand--Tsetlin lattices.  
We also use these results to obtain new splitting distributive lattices for the Weyl bialternants associated with multiples of certain ``end-node'' fundamental weights in the symplectic and orthogonal cases. 
In all these cases, the splitting distributive lattices are described explicitly as partial orders of Gelfand-type patterns. 

{\bf New splitting results for certain non-fibrous posets.}  
Before we state the new splitting results \NewSplittingTheorem/\ChainProductVersions, we briefly discuss some of their signature aspects.  
These results are methodological in that they render some of the technical criteria of \InitialSplittingTheorem\ in more purely combinatorial terms. 
Of course, our set $I$ indexes simple roots and fundamental weights associated with our fixed root system $\Phi$ and is to be viewed as a collection of colors; let $J$ be a subset of $I$.  
Let $R$ be a ranked poset with edges colored by $I$, and let $\mathcal{S}$ be a subset of the vertices of $R$.  
When a dominant weight $\nu \in \Lambda_{\Phi_J}^{+}$ is specified, we will most often take $\mathcal{S}$ to be a subset of $\{\xelt \in R\, |\, \delta_{j}(\xelt) \leq \nu_{j} \mbox{ for all } j \in J\}$, where the latter analogizes the set $\mathcal{M}_{J,\nu}$ for fibrous posets. 
Fix a set mapping (i.e.\ vertex-coloring function) $\kappa: R \setminus \mathcal{S} \longrightarrow J$. 
The {\em kindred set} of any $\xelt \in R \setminus \mathcal{S}$ is 
\[\myK(\xelt) := \{\yelt \in \comp_{\kappa(\xelt)}(\xelt)\, |\, \yelt \in R \setminus \mathcal{S} \mbox{ and } \kappa(\yelt) = \kappa(\xelt)\}.\]
In the subgraph of $R$ induced by the vertices of $\myK(\xelt)$, regard each connected component to be a subposet of $\comp_{\kappa(\xelt)}(\xelt)$ in the induced order. 
Then view $\myK(\xelt)$ as the disjoint sum of these subposets.

\noindent 
{\bf \NewSplittingTheorem}\ \ {\sl Let $J \subseteq I$, $\nu \in \Lambda_{\Phi_J}^{+}$, and write $\nu = \sum_{j \in J}\nu_{j}\omega_{j}^{J}$. 
Suppose $R$ is an ${\Phi_J}$-structured poset. 
Suppose also that} $\WGF(R)|_{J}$ {\sl is $W_J$-invariant. 
(In this setting, the latter is guaranteed if, for example, for all $j \in J$ the $j$-components of $R$ are rank symmetric, cf.\ \WInvariantLemma.) 
Let $\mathcal{S} = \mathcal{S}_{J,\nu}(R)$ be a subset of $\{\xelt \in R\, |\, \delta_{j}(\xelt) \leq \nu_{j} \mbox{ for all } j \in J\}$.
Suppose $\kappa: R \setminus \mathcal{S} \longrightarrow J$ is a vertex-coloring function such that for each $\xelt \in R \setminus \mathcal{S}$ with $k := \kappa(\xelt)$ we have 
(1) The restriction of the rank function $\rho_{k}: \comp_{k}(\xelt) \longrightarrow \{0,1,\cdots,l_{k}(\xelt)\}$ to} $\myK(\xelt)$ {\sl is a rank function with range $\{0,1,\cdots,l_{k}(\xelt)-(\nu_{k}+1)\}$,  
and (2)} $\myK(\xelt)$ {\sl has a symmetric chain decomposition with respect to the rank function} $\rho_{k}|_{\mytinyK(\xelt)}$. 
{\sl Then $\nu + wt^{J}(\selt) \in \Lambda_{\Phi_J}^{+}$ for all $\selt \in \mathcal{S}$ and $R$ is a $(J,\nu)$-splitting poset with}
\[\chi^{\Phi_J}_{_{\nu}} \cdot \WGF(R)|_{J} = \sum_{\selt \in \mathcal{S}_{J,\nu}(R)}\chi^{\Phi_J}_{_{\nu + wt^{J}(\selt)}}.\] 

{\em Proof.}  We demonstrate that the requirements of \InitialSplittingTheorem\ are met. 
First, note that if $\selt \in \mathcal{S}$ and $j \in J$, then $\langle \nu + wt^{J}(\selt), \alpha_{j}^{\vee} \rangle = \nu_{j} + \rho_{j}(\selt) - \delta_{j}(\selt) \geq \nu_{j} + 0 - \nu_{j} = 0$. 
So, $\nu + wt^{J}(\selt) \in \Lambda_{\Phi_J}^{+}$ for any $\selt \in \mathcal{S}$. 
Next, we construct a bijection $\tau: R \setminus \mathcal{S} \longrightarrow R \setminus \mathcal{S}$. 
For $\xelt \in R \setminus \mathcal{S}$, set $k := \kappa(\xelt)$. 
For a given SCD of $\myK(\xelt)$, consider the chain that contains $\xelt$. 
Within said chain, let $\xelt'$ be the unique element whose position is symmetric to $\xelt$ in the sense that $\rho_{k}(\xelt') = \delta_{k}(\xelt) - (1 + \nu_{k})$. 
Now, $wt^{J}(\xelt') = wt^{J}(\xelt) - (\rho_{k}(\xelt) - \rho_{k}(\xelt'))\alpha_{k} = wt^{J}(\xelt) - (\rho_{k}(\xelt) - \delta_{k}(\xelt) + 1 + \nu_{k})\alpha_{k} = wt^{J}(\xelt) - (1 + \nu_{k} + m_{k}(\xelt))\alpha_{k}$, where the first of the preceding equalities follows from the fact that $R$ is ${\Phi_J}$-structured. 
So, when we set $\tau(\xelt) := \xelt'$, we obtain a bijection that, together with our given vertex-coloring function $\kappa$, satisfies the requirements of \InitialSplittingTheorem.\hfill\QED

The following corollary is a consequence of the well-known fact that any chain product has an SCD.  
A classical and constructive version of this latter fact is Gansner's symmetric chain decomposition of chain products by parenthesization \cite{Gansner}\footnote{In v.\ 1 of this monograph from November 26, 2018, such an explicit construction of symmetric chains in a chain product was carried out in the proof of what was presented there as Theorem 8.1.}. 
Under the hypotheses of the following corollary, an explicit matching $\tau$ of vertices can be constructed using Gansner-prescribed symmetric chains of each $\myK(\xelt) \subseteq \comp_{\kappa(\xelt)}(\xelt)$ when $\comp_{\kappa(\xelt)}(\xelt)$ is a chain product.   
However, such explicitness is not required for the proof of the corollary.  

The following notation and nomenclature is needed for the statement and proof of \ChainProductVersion. 
Suppose $\mathcal{C}$ is a product $\mathcal{C}_{1} \times \cdots \times \mathcal{C}_{p}$ of chains $\mathcal{C}_{1}, \ldots ,\mathcal{C}_{p}$, and let $b$ be a positive integer. 
If for some $q \in \{1, 2, \ldots, p\}$ we have $\sum_{r = q+1}^{p}(|\mathcal{C}_{r}|-1) < b \leq \sum_{r=q}^{p}(|\mathcal{C}_{r}|-1)$, then the $b${\em -sub-block of} $\mathcal{C}$ is 
\[\left\{(\xelt_1,\ldots,\xelt_p) \in \mathcal{C}\, \left|\begin{array}{c} \\ \\ \end{array}\right.  \mbox{\parbox{4in}{$\xelt_{q+1},\ldots,\xelt_{p}$ are minimal in $\mathcal{C}_{q+1},\ldots,\mathcal{C}_{p}$ respectively and $\xelt_q$ has depth at least  $b - \sum_{r = q+1}^{p}(|\mathcal{C}_{r}|-1)$ in $\mathcal{C}_{q}$}}\right\}.\]
The $b$-sub-block of $\mathcal{C}$ is the empty set if $b > \sum_{r=1}^{p}(|\mathcal{C}_{r}|-1)$. 
If $\phi: \mathcal{C} \longrightarrow P$ is a poset isomorphism, then say $S \subseteq P$ is a $b${\em -sub-block of} $P$ if $\phi^{-1}(S)$ is the $b$-sub-block of $\mathcal{C}$. 
If the isomorphism $\phi$ is understood, we simply call $S$ a $b$-sub-block of $P$. 
When $b=0$, we call $S$ a {\em sub-face} of $P$. 

\noindent 
{\bf \ChainProductVersion.A}\ \ {\sl Let $J \subseteq I$, $\nu \in \Lambda_{\Phi_J}^{+}$, and write $\nu = \sum_{j \in J}\nu_{j}\omega_{j}^{J}$. 
Suppose $R$ is an ${\Phi_J}$-structured poset. 
Suppose also that} $\WGF(R)|_{J}$ {\sl is $W_J$-invariant. 
(In this setting, the latter is guaranteed if, for example, for all $j \in J$ the $j$-components of $R$ are rank symmetric, cf.\ \WInvariantLemma.) 
Let $\mathcal{S} = \mathcal{S}_{J,\nu}(R)$ be a subset of $\{\xelt \in R\, |\, \delta_{j}(\xelt) \leq \nu_{j} \mbox{ for all } j \in J\}$. 
Suppose $\kappa: R \setminus \mathcal{S} \longrightarrow J$ is a vertex-coloring function such that for each $\xelt \in R \setminus \mathcal{S}$ with $k := \kappa(\xelt)$ we have (1) $\comp_{k}(\xelt)$ is isomorphic to a product of chains and (2)} $\myK(\xelt)$ {\sl is a $(\nu_{k} + 1)$-sub-block of $\comp_{k}(\xelt)$. Then $\nu + wt^{J}(\selt) \in \Lambda_{\Phi_J}^{+}$ for all $\selt \in \mathcal{S}$ and $R$ is a $(J,\nu)$-splitting poset with}
\[\chi^{\Phi_J}_{_{\nu}} \cdot \WGF(R)|_{J} = \sum_{\selt \in \mathcal{S}_{J,\nu}(R)}\chi^{\Phi_J}_{_{\nu + wt^{J}(\selt)}}.\] 

\noindent 
{\bf \ChainProductVersion.B}\ \ {\sl Continuing with the hypotheses of} {\bf A}{\sl , suppose that $J = I$ and that $\nu = 0$, so that each} $\myK(\xelt)$ {\sl is a sub-face of $\comp_{k}(\xelt)$.  
Then $wt(\selt) = wt^{I}(\selt) \in \Lambda_{\Phi_{I}}^{+} = \Lambda_{\Phi}^{+}$ for all $\selt \in \mathcal{S}$ and $R$ is a $(I,0)$-splitting poset with}
\[\WGF(R)|_{I} = \WGF(R) = \sum_{\selt \in \mathcal{S}_{I,0}(R)}\chi^{\Phi}_{_{wt(\selt)}}.\]

{\em Proof.} Clearly {\bf B} is a special case of {\bf A}, so we only prove the latter. 
As in the proof of \NewSplittingTheorem, the fact that $\nu + wt^{J}(\selt) \in \Lambda_{\Phi_J}^{+}$ for all $\selt \in \mathcal{S}$ follows from the definitions. 
Now let $\xelt \in R \setminus \mathcal{S}$ with $k := \kappa(\xelt)$. 
Since $\myK(\xelt)$ is isomorphic to a sub-block of a product of chains, then $\myK(\xelt)$ itself is isomorphic to a product of chains. 
In fact, $\myK(\xelt)$ is easily seen to be an order ideal from $\comp_{k}(\xelt)$, so edges in $\comp_{k}(\xelt)$ are also edges in $\myK(\xelt)$. 
Now assume that $\comp_{k}(\xelt)$ is isomorphic to the chain product $\mathcal{C}$ in the paragraph preceding the corollary statement and that $\mathcal{K} \subseteq \mathcal{C}$ is isomorphic to $\myK(\xelt)$. 
If $\belt_{r}$ is the minimal element of the chain $\mathcal{C}_{r}$, then $(\belt_{1},\cdots,\belt_{p}) \in \mathcal{K}$. 
So, $\mathcal{K}$ contains the minimal element of $\mathcal{C}$. 
Set $b := \nu_{k}+1$. 
Say $\selt_{q}$ is an element of chain $\mathcal{C}_{q}$ with depth $b - \sum_{r = q+1}^{p}(|\mathcal{C}_{r}|-1)$, and, for $1 \leq i \leq q-1$, say $\telt_{i}$ is maximal in $\mathcal{C}_{i}$.  
Let $\selt \in \mathcal{C}$ be the $p$-tuple $(\telt_{1},\cdots,\telt_{q-1},\selt_{q},\belt_{q+1},\cdots,\belt_{p})$, which is clearly maximal in $\mathcal{K}$. 
Now, $\selt$ has rank $\sum_{r = 1}^{q-1}(|\mathcal{C}_{r}|-1) + \left[|\mathcal{C}_{q}|-1 - \left(b - \sum_{r = q+1}^{p}(|\mathcal{C}_{r}|-1)\right)\right] = \sum_{r = 1}^{p}(|\mathcal{C}_{r}|-1) - b = l_{k}(\xelt) - (\nu_{k}+1)$. 
Thus, the connected and ranked chain product $\myK(\xelt)$ contains the minimal element of $\comp_{k}(\xelt)$ and has length $l_{k}(\xelt) - (\nu_{k}+1)$. 
Then $\myK(\xelt)$ meets the criteria of \NewSplittingTheorem.\hfill\QED

{\bf Gelfand--Tsetlin lattices as splitting distributive lattices.} 
The well-known Gelfand--Tsetlin bases for the irreducible representations of the complex simple Lie algebra $\mathfrak{sl}(n,\mathbb{C})$ have diamond-colored  distributive lattice supporting graphs, see \cite{PrGZ}, \cite{DonSupp} and \cite{HL}. 
In this section we will produce these lattices from scratch and use \ChainProductVersion.B to conclude that they are splitting distributive lattices for the $\myA_{n-1}$-Weyl bialternants. 
This serves several purposes: First, it is an illustration of \ChainProductVersion.B.  
Second, we can often identify components of splitting posets as Gelfand--Tsetlin lattices 
(e.g.\ \cite{DD}, \cite{DDW}), 
and the description of Gelfand--Tsetlin lattices we offer here can be helpful in those circumstances. 
And third, we will use an almost identical approach in the next subsection to build the advertised symplectic and orthogonal splitting distributive lattices. 

Our root system in this subsection is of type $\myA_{n-1}$, and the index set $I = \{1, 2, \ldots , n-1\}$ gives the usual numbering of the simple roots. 
Fix a dominant weight $\lambda = \sum_{k \in I}a_{k}\omega_{k}$. 
We build Gelfand--Tseltin lattices from objects closely related to what are traditionally called ``Gelfand patterns'' (see for example \cite{HL}). 
A {\em special linear ideal pattern} of size $n-1$ and bounded by $\lambda$ is an array $(g_{i,j})_{1 \leq j \leq i \leq n}$ such that for $1 \leq j \leq n$ we have $g_{n,j} := \sum_{k=n+1-j}^{n-1}a_k$ (so $g_{n,1} = 0$) and satisfying the inequalities depicted in \SpecialPatternFigure. 
These are ``ideal patterns'' because if we set $b_{j} := \sum_{k=n+1-j}^{n-1}a_k$, then the pattern $(g_{i,j}-b_{j})_{1 \leq j \leq i \leq n-1}$ can be readily identified as data for an order ideal taken from the poset of join irreducibles for the associated Gelfand--Tsetlin lattice. 

We use these arrays as follows to build Gelfand--Tsetlin lattices.  
Let $L_{\mysmallA_{n-1}}(\lambda)$ be the collection of special linear ideal patterns of size $n-1$ and bounded by $\lambda$. 
For $\telt \in L_{\mysmallA_{n-1}}(\lambda)$, we use the notation $g_{i,j}(\telt)$ to refer to the $(i,j)$-entry of the corresponding array. 
Partially order the elements of $L_{\mysmallA_{n-1}}(\lambda)$ by componentwise comparison, so $\selt \leq \telt$ if and only if $g_{i,j}(\selt) \leq g_{i,j}(\telt)$ for all $1 \leq i \leq j \leq n$. 
It is routine to check that the partially ordered set $L_{\mysmallA_{n-1}}(\lambda)$ is a distributive lattice. 
Also, $\selt \rightarrow \telt$ in the Hasse diagram for $L_{\mysmallA_{n-1}}(\lambda)$ if and only if there is a pair $(i,j)$ such that $g_{p,q}(\selt) = g_{p,q}(\telt)$ when $(p,q) \not= (i,j)$ and $g_{i,j}(\selt) + 1 = g_{i,j}(\telt)$. 
In this case, we attach the color $i \in I$ to this edge of the Hasse diagram and write $\selt \myarrow{i} \telt$. 
One can see, then, that $L_{\mysmallA_{n-1}}(\lambda)$ is a diamond-colored distributive lattice called a {\em Gelfand--Tsetlin lattice}. 

\begin{figure}[ht] 
\begin{center}
\SpecialPatternFigure: Inequalities for a special linear ideal pattern.

\vspace*{0.25in}
\begin{tabular}{cccccccccc}
 & & & & & & & & & $g_{n,1}$ \\
 & & & & & & & & \begin{rotate}{45}$\geq$\end{rotate} & \\
 & & & & & & & $g_{n-1,1}$ & & \\
 & & & & & & \begin{rotate}{45}$\geq$\end{rotate} & & \begin{rotate}{315}$\!\!\!\!\leq$\end{rotate} & \\
 & & & & & $g_{n-2,1}$ & & & & $g_{n,2}$ \\
 & & & & \begin{rotate}{45}$\geq$\end{rotate} & & \begin{rotate}{315}$\!\!\!\!\leq$\end{rotate} & & \begin{rotate}{45}$\geq$\end{rotate} & \\
 & & & $\iddots$ & & $\cdot$ & $\cdot$ & $\cdot$ & & \\
 & \begin{rotate}{45}$\geq$\end{rotate} & & & & & & & & \\
$g_{1,1}$ & & & & & $\cdot$ & $\cdot$ & $\cdot$ & & \\
 & \begin{rotate}{315}$\!\!\!\!\leq$\end{rotate} & & & & & & & & \\
 & & & $\ddots$ & & $\cdot$ & $\cdot$ & $\cdot$ & & \\
 & & & & \begin{rotate}{315}$\!\!\!\!\leq$\end{rotate} & & \begin{rotate}{45}$\geq$\end{rotate} & & \begin{rotate}{315}$\!\!\!\!\leq$\end{rotate} & \\
 & & & & & $g_{n-2,n-2}$ & & & & $g_{n,n-1}$ \\
 & & & & & & \begin{rotate}{315}$\!\!\!\!\leq$\end{rotate} & & \begin{rotate}{45}$\geq$\end{rotate} & \\
 & & & & & & & $g_{n-1,n-1}$ & & \\
 & & & & & & & & \begin{rotate}{315}$\!\!\!\!\leq$\end{rotate} & \\
 & & & & & & & & & $g_{n,n}$ 
\end{tabular}
\end{center}
\end{figure} 

\noindent 
{\bf \GelfandProposition}\ \ {\sl Let $\lambda$ be as above.  For each $i \in I$, each $i$-component of} $L := L_{\mysmallA_{n-1}}(\lambda)$ {\sl is isomorphic to a product of chains. 
Moreover, $L$ is} ${\myA_{n-1}}${\sl -structured, and} $\WGF(L)$ {\sl is an} $\myA_{n-1}${\sl -Weyl symmetric function.} 

{\em Proof.}  
For any $\telt \in L$ and $i \in I$, it is clear that $\comp_{i}(\telt)$ is comprised precisely of those ideal patterns $\selt$ such that $g_{j,k}(\selt) = g_{j,k}(\telt)$ for all $j \not= i$. 
Moreover, the defining pattern of inequalities allows that changes can be made to positions $(i,k)$ and $(i,k')$ independently of one another. 
That is, $\comp_{i}(\telt) \cong \mathcal{C}_{1} \times \cdots \times \mathcal{C}_{i}$, where for each $1 \leq k \leq i$, $\mathcal{C}_{k}$ is a chain whose length is the difference of the maximum value that can replace $g_{i,k}(\telt)$ and the minimum such value. 
In fact, this quantity is $\min(g_{i-1,k}(\telt),g_{i+1,k+1}(\telt)) - \max(g_{i-1,k-1}(\telt),g_{i+1,k}(\telt))$, where $g_{p,q}(\telt)$ is omitted in the formula if the position $(p,q)$ is not part of the array. 

As in \cite{HL}, one can see that $L$ is isomorphic to the diamond-colored distributive lattice built from semistandard tableaux whose entries come from the set $\{1,2,\cdots,n\}$ and whose shape is a partition diagram (Ferrers diagram) with $a_{k}$ columns of length $k$. 
This isomorphism utilizes a natural one-to-one correspondence between special linear ideal patterns and semistandard tableaux. 
Therefore it follows from Theorem 4.1 of \cite{HL} that $L$ is ${\myA_{n-1}}$-structured. 
For later purposes, however, we present the proof details using the present notational conventions. 

A consequence of the first paragraph of the proof is that for each $\telt \in L$ and $i \in I$, we have 
\begin{eqnarray*}
m_{i}(\telt) & = & \sum_{k=1}^{i}\left[\rule[-3.5mm]{0mm}{8mm}\left(\rule[-2.5mm]{0mm}{6mm}g_{i,k}(\telt) - \min\left(\rule[-1.5mm]{0mm}{3.5mm}g_{i-1,k}(\telt),g_{i+1,k+1}(\telt)\right)\right) - \left(\rule[-2.5mm]{0mm}{6mm}\max\left(\rule[-1.5mm]{0mm}{3.5mm}g_{i-1,k-1}(\telt),g_{i+1,k}(\telt)\right)-g_{i,k}(\telt)\right)\right]\\
 & = & \sum_{k=1}^{i}\left[\rule[-2.5mm]{0mm}{6mm}2g_{i,k}(\telt) - \min\left(\rule[-1.5mm]{0mm}{3.5mm}g_{i-1,k}(\telt),g_{i+1,k+1}(\telt)\right) - \max\left(\rule[-1.5mm]{0mm}{3.5mm}g_{i-1,k-1}(\telt),g_{i+1,k}(\telt)\right)\right]. 
\end{eqnarray*}
To see that $L$ is $\myA_{n-1}$-structured, it suffices to check that when $\selt \myarrow{i} \telt$ in $L$, then $m_{j}(\selt) + \langle \alpha_{i},\alpha_{j}^{\vee} \rangle = m_{j}(\telt)$ for all $j \not= i$.  Let us say that $\selt \myarrow{i} \telt$ with $g_{i,q}(\telt) = g_{i,q}(\selt)+1$.  
First we analyze the case that $1 \leq i < n-1$ with $j=i+1$.  
In particular, $\langle \alpha_{i},\alpha_{j}^{\vee} \rangle = -1$.  
Then 
\begin{eqnarray*}
m_{j}(\telt) - m_{j}(\selt) & = & 2g_{j,q}(\telt) - \min\left(\rule[-1.5mm]{0mm}{3.5mm}g_{i,q}(\telt),g_{j+1,q+1}(\telt)\right) - \max\left(\rule[-1.5mm]{0mm}{3.5mm}g_{j,q-1}(\telt),g_{j+1,q}(\telt)\right)\\
& & -2g_{j,q}(\telt) + \min\left(\rule[-1.5mm]{0mm}{3.5mm}g_{i,q}(\telt)-1,g_{j+1,q+1}(\telt)\right) + \max\left(\rule[-1.5mm]{0mm}{3.5mm}g_{j,q-1}(\telt),g_{j+1,q}(\telt)\right)\\
& & +2g_{j,q+1}(\telt) - \min\left(\rule[-1.5mm]{0mm}{3.5mm}g_{i,q+1}(\telt),g_{j+1,q+2}(\telt)\right) - \max\left(\rule[-1.5mm]{0mm}{3.5mm}g_{i,q}(\telt),g_{j+1,q+1}(\telt)\right)\\
& & -2g_{j,q+1}(\telt) + \min\left(\rule[-1.5mm]{0mm}{3.5mm}g_{i,q+1}(\telt),g_{j+1,q+2}(\telt)\right) + \max\left(\rule[-1.5mm]{0mm}{3.5mm}g_{i,q}(\telt)-1,g_{j+1,q+1}(\telt)\right)\\
& = & - \min\left(\rule[-1.5mm]{0mm}{3.5mm}g_{i,q}(\telt),g_{j+1,q+1}(\telt)\right) + \min\left(\rule[-1.5mm]{0mm}{3.5mm}g_{i,q}(\telt)-1,g_{j+1,q+1}(\telt)\right)\\
& & - \max\left(\rule[-1.5mm]{0mm}{3.5mm}g_{i,q}(\telt),g_{j+1,q+1}(\telt)\right) + \max\left(\rule[-1.5mm]{0mm}{3.5mm}g_{i,q}(\telt)-1,g_{j+1,q+1}(\telt)\right)\\
& = & -1,
\end{eqnarray*}
where the concluding step follows by checking the cases (1) $g_{i,q}(\telt) > g_{j+1,q+1}(\telt)$ and (2) $g_{i,q}(\telt) \leq g_{j+1,q+1}(\telt)$. 
Analysis of the case $1 < i \leq n-1$ and $j=i-1$ is entirely similar. 
When $i$ and $j$ are ``distant,'' i.e.\ $|i-j| \geq 2$ so that $\langle \alpha_{i},\alpha_{j}^{\vee} \rangle = 0$, then the terms of the sum for $m_{j}(\telt)$ are exactly the same as the terms of the sum for $m_{j}(\selt)$, hence $m_{j}(\telt) = m_{j}(\selt)$. 

Thus, $L$ is $\myA_{n-1}$-structured.  It now follows from \WInvariantLemma\ that $\WGF(L)$ is an $\myA_{n-1}$-Weyl symmetric function.\hfill\QED

In order to develop a vertex-coloring scheme for the Gelfand--Tsetlin lattices, we require a total ordering of the positions $(i,j)$ for the special linear ideal patterns associated with a given $\myA_{n-1}$-dominant weight $\lambda$. 
To this end, fix such a special linear ideal pattern $\telt$, and say entry $g_{i,j}(\telt)$ is {\em slantwise-prior} to entry $g_{p,q}(\telt)$ if (i) $i-j < p-q$ or (ii) $i-j = p-q$ and $q < j$. 
This establishes a reading order of the array positions of $\telt$ such that the first position is $g_{n,n}(\telt)$, and from there we read along SE to NW diagonals (or ``slants'') to follow the total ordering of positions. 
We sometimes apply this slantwise nomenclature to the pairs $(i,j)$ that index the positions for special linear ideal patterns. 
Now let $\melt$ be the unique maximal element of $L_{\mysmallA_{n-1}}(\lambda)$. 
We say $\telt$ {\em can be maximized} in position $(i,j)$ if $g_{i,j}(\telt) < g_{i,j}(\melt)$ and when we replace entry $g_{i,j}(\telt)$ with $g_{i,j}(\melt)$ we get a valid special linear ideal pattern. 
Then say $(i,j)$ is the {\em slantwise-least maximizable position} for $\telt$ if $(i,j)$ is the smallest position in our slantwise reading such that $\telt$ can be maximized there. 

\noindent
{\bf \GelfandLemma}\ \ {\sl With $\lambda$ and} $L := L_{\mysmallA_{n-1}}(\lambda)$ {\sl as above, let $\melt$ denote the unique maximal element of $L$.  Let $\telt \not= \melt$ in $L$, and suppose $(i,j)$ is the slantwise-least maximizable position for $\telt$.  (1) Then $g_{p,q}(\telt) = g_{p,q}(\melt)$ for all positions $(p,q)$ that are slantwise-prior to $(i,j)$.  (2) If $\selt \in \comp_{i}(\telt)$ with $g_{i,j}(\selt) < g_{i,j}(\melt)$, then $(i,j)$ is the slantwise-least maximizable position for $\selt$. (3)  Suppose $\selt \in \comp_{i}(\telt)$ such that $\selt \not= \melt$ and the slantwise-least maximizable position for $\selt$ is some $(i,k)$. Then $k=j$.}

{\em Proof.}  For {\sl (1)}, let $(p,q)$ be the slantwise smallest position for which $g_{p,q}(\telt) < g_{p,q}(\melt)$.  
The inequalities defining special linear ideal patterns assert that $g_{p,q}(\telt) \leq g_{p-1,q}(\telt)$ and $g_{p,q}(\telt) \leq g_{p+1,q+1}(\telt)$, assuming $(p-1,q)$ and $(p+1,q+1)$ are valid positions in the array. 
Since both of the latter positions are slantwise-prior to $(p,q)$, it follows from our choice of $(p,q)$ that $g_{p-1,q}(\telt) = g_{p-1,q}(\melt)$ and $g_{p+1,q+1}(\telt) = g_{p+1,q+1}(\melt)$. 
Therefore $\telt$ can be maximized at position $(p,q)$, and it follows that $(p,q)$ is the slantwise-least maximizable position for $\telt$. 
Then $(i,j) = (p,q)$, from which the statement of {\sl (1)} follows. 

For {\sl (2)}, suppose the slantwise-least maximizable position for $\selt$ is $(p,q)$. 
Well, assuming positions $(i-1,j)$ and $(i+1,j+1)$ are valid, then $g_{i-1,j}(\selt) = g_{i-1,j}(\telt) = g_{i-1,j}(\melt)$ and $g_{i+1,j+1}(\selt) = g_{i+1,j+1}(\telt) = g_{i+1,j+1}(\melt)$. 
This means that $\selt$ can be maximized at position $(i,j)$. 
Thus, either $(i,j) = (p,q)$ or $(p,q)$ is slantwise-prior to $(i,j)$. 
Assume $(p,q)$ is slantwise-prior to $(i,j)$. 
Since $g_{r,s}(\selt) = g_{r,s}(\telt) = g_{r,s}(\melt)$ for all $(r,s)$ slantwise-prior to $(i,j)$ with $r \not= i$, then it follows that $p = i$ and $q > j$. 
Then in fact $(p-1,q-1)$ and $(p+1,q+1)$ are valid array positions, and we have the inequalities $g_{p-1,q-1}(\selt) \leq g_{p,q}(\selt) \leq g_{p+1,q+1}(\selt)$. 
But $g_{p-1,q-1}(\selt) = g_{p-1,q-1}(\telt) = g_{p-1,q-1}(\melt)$ and $g_{p+1,q+1}(\selt) = g_{p+1,q+1}(\telt) = g_{p+1,q+1}(\melt)$. 
And moreover, $g_{p-1,q-1}(\melt) = g_{p+1,q+1}(\melt) = g_{n,n-p+q}(\melt)$. 
This forces $g_{p,q}(\selt) = g_{p,q}(\melt)$, which contradicts the fact that $(p,q)$ is the slantwise-least maximizable position for $\selt$. 
Therefore $(p,q) = (i,j)$, as claimed. 

For {\sl (3)}, we show $(i,k) = (i,j)$ by ruling out the possibility that $(i,j)$ is slantwise-prior to $(i,k)$ and vice-versa.  
For our contradiction hypothesis, assume that $k \ne j$. 
Now, $k \ne j$ means that $i > 1$, and we already know that $i < n$. 
Let us suppose that $(i,k)$ is slantwise-prior to $(i,j)$.  
Then $(i-1,k-1)$ and $(i+1,k+1)$ are valid array positions, and we have $g_{i-1,k-1}(\telt) =  g_{i,k}(\telt) = g_{i+1,k+1}(\telt) = g_{n,n-i+k}(\melt)$, because, by part {\sl (1)}, $\telt$ attains the maximum value at each position slantwise-prior to $(i,j)$. 
Since $\selt \in \comp_{i}(\telt)$, then $g_{i-1,k-1}(\selt) = g_{i-1,k-1}(\telt)$ and $g_{i+1,k+1}(\selt) = g_{i+1,k+1}(\telt)$, which forces $g_{i-1,k-1}(\selt) =  g_{i,k}(\selt) = g_{i+1,k+1}(\selt) = g_{n,n-i+k}(\melt)$. 
This contradicts the fact that $(i,k)$ is the slantwise-least maximizable position for $\selt$. 
But now if $(i,j)$ is slantwise-prior to $(i,k)$, the same reasoning leads to a contradiction of the fact that $(i,j)$ is the slantwise-least maximizable position for $\telt$. 
Therefore, $k=j$, which completes the proof.\hfill\QED

\noindent
{\bf \GelfandColoring}\ \ {\sl With $\lambda$ and} $L := L_{\mysmallA_{n-1}}(\lambda)$ {\sl as above, let $\melt$ denote the unique maximal element of $L$. 
Let $\kappa: L \setminus \{\melt\} \longrightarrow I$ be defined as follows: For $\telt \not= \melt$ in $L$, say $\kappa(\telt) = i$ if  the slantwise-least maximizable position for $\telt$ is some $(i,j)$.  
Then $\{\selt \in \comp_{\kappa(\telt)}(\telt)\, |\, \selt \not= \melt \mbox{ and } \kappa(\selt) = \kappa(\telt)\}$ is a sub-face of $\comp_{\kappa(\telt)}(\telt)$.} 

{\em Proof.} From \GelfandLemma, particularly parts {\sl (2)} and {\sl (3)}, we get the following equality of sets: $\{\selt \in \comp_{i}(\telt)\, |\, \selt \not= \melt \mbox{ and } \kappa(\selt) = i\} = \{\selt \in \comp_{i}(\telt)\, |\, g_{i,j}(\selt) < g_{i,j}(\melt)\}$.
The latter set is clearly a sub-face of $\comp_{i}(\telt)$.\hfill\QED

From \GelfandPrelimProps, \ChainProductVersion.B, and \MainCorollary, it follows immediately that the Gelfand--Tsetlin lattices are splitting distributive lattices for the $\myA_{n-1}$-Weyl bialternants and that their rank generating functions have nice quotient-of-product expressions.  That is: 

\noindent 
{\bf \GelfandTheorem}\ \ {\sl With $\lambda$ as above, then} $L_{\mysmallA_{n-1}}(\lambda)$ {\sl is a rank symmetric and rank unimodal splitting distributive lattice for the} $\myA_{n-1}${\sl -Weyl bialternant} $\chi_{_{\lambda}}^{\mytinyA_{n-1}}$.  {\sl Moreover,} 
\[\RGF(L_{\mysmallA_{n-1}}(\lambda),q) = \frac{\mbox{$\displaystyle \prod_{\alpha \in \Phi^{+}}$}(1-q^{\langle \lambda+\varrho,\alpha^{\vee} \rangle})}{\mbox{$\displaystyle \prod_{\alpha \in \Phi^{+}}$}(1-q^{\langle \varrho,\alpha^{\vee} \rangle})}.\]

\vspace*{-0.4in}\hfill\QED

The $\myA_{n-1}$-Weyl bialternants are versions of Schur functions;  
this connection is developed in detail in \CaseAExample. 
For a more explicit version of the $\RGF$ polynomial above, see Theorem 10.6.3 of \cite{DonDistributive}.

{\bf Splitting distributive lattices for certain symplectic and orthogonal Weyl bialternants.} 
Analogizing the previous subsection, we will use some Gelfand-type patterns to build splitting distributive lattices for certain $\myB_n$-, $\myC_n$-, and $\myD_n$-Weyl bialternants. 
We index the nodes of the Dynkin diagram in each case using $I = \{1,2,\ldots,n\}$ as in \cite{Hum}. 
The dominant weights associated with these Weyl bilalternants are multiples of 
certain end-node fundamental weights. 
\ApplicationTheorem\ below is an application of \ChainProductVersion.B to these cases. 

Fix positive integers $m$ and $n$. 
An {\em odd orthogonal ideal pattern} of size $n$ and bounded by $m$ is an array $(c_{i,j})_{1 \leq j \leq i \leq n}$ of nonnegative integers with $c_{n,n} \leq m$ and satisfying the inequalities pictured in \OddOrthPatternFigure. 
These are ``ideal'' in the sense that the numbers of the array identify a particular order ideal taken from the product of a chain with $m$ elements and a staircase-shaped poset.  
This latter staircase-shaped poset is, in fact, the minuscule poset associated with the ``spin-node'' fundamental weight $\omega_n$ associated with the type $\myB_n$ root system, see \cite{PrEur}.  

\begin{figure}[ht] 
\begin{center}
\OddOrthPatternFigure: Inequalities for an odd orthogonal ideal pattern.

\vspace*{0.25in}
\begin{tabular}{cccccccccc}
 & & & & & & & & & $c_{n,1}$ \\
 & & & & & & & & \begin{rotate}{45}$\geq$\end{rotate} & \\
 & & & & & & & $c_{n-1,1}$ & & \\
 & & & & & & \begin{rotate}{45}$\geq$\end{rotate} & & \begin{rotate}{315}$\!\!\!\!\leq$\end{rotate} & \\
 & & & & & $c_{n-2,1}$ & & & & $c_{n,2}$ \\
 & & & & \begin{rotate}{45}$\geq$\end{rotate} & & \begin{rotate}{315}$\!\!\!\!\leq$\end{rotate} & & \begin{rotate}{45}$\geq$\end{rotate} & \\
 & & & $\iddots$ & & $\cdot$ & $\cdot$ & $\cdot$ & & \\
 & \begin{rotate}{45}$\geq$\end{rotate} & & & & & & & & \\
$c_{1,1}$ & & & & & $\cdot$ & $\cdot$ & $\cdot$ & & \\
 & \begin{rotate}{315}$\!\!\!\!\leq$\end{rotate} & & & & & & & & \\
 & & & $\ddots$ & & $\cdot$ & $\cdot$ & $\cdot$ & & \\
 & & & & \begin{rotate}{315}$\!\!\!\!\leq$\end{rotate} & & \begin{rotate}{45}$\geq$\end{rotate} & & \begin{rotate}{315}$\!\!\!\!\leq$\end{rotate} & \\
 & & & & & $c_{n-2,n-2}$ & & & & $c_{n,n-1}$ \\
 & & & & & & \begin{rotate}{315}$\!\!\!\!\leq$\end{rotate} & & \begin{rotate}{45}$\geq$\end{rotate} & \\
 & & & & & & & $c_{n-1,n-1}$ & & \\
 & & & & & & & & \begin{rotate}{315}$\!\!\!\!\leq$\end{rotate} & \\
 & & & & & & & & & $c_{n,n}$ 
\end{tabular}
\end{center}
\end{figure}

We use these arrays as follows to build splitting distributive lattices for certain odd orthogonal Weyl bialternants.  
Let $L_{\mysmallB_n}(m\omega_n)$ be the collection of odd orthogonal ideal patterns of size $n$ and bounded by $m$. 
For $\telt \in L_{\mysmallB_n}(m\omega_n)$, we use the notation $c_{i,j}(\telt)$ to refer to the $(i,j)$-entry of the corresponding array. 
Partially order the elements of $L_{\mysmallB_n}(m\omega_n)$ by componentwise comparison, so $\selt \leq \telt$ if and only if $c_{i,j}(\selt) \leq c_{i,j}(\telt)$ for all $1 \leq i \leq j \leq n$. 
It is routine to check that the partially ordered set $L_{\mysmallB_n}(m\omega_n)$ is a distributive lattice. 
Then $\selt \rightarrow \telt$ in the Hasse diagram for $L_{\mysmallB_n}(m\omega_n)$ if and only if there is a pair $(i,j)$ such that $c_{p,q}(\selt) = c_{p,q}(\telt)$ when $(p,q) \not= (i,j)$ and $c_{i,j}(\selt) + 1 = c_{i,j}(\telt)$. 
In this case, we attach the color $i$ to this edge of the Hasse diagram and write $\selt \myarrow{i} \telt$. 
In this way we realize $L_{\mysmallB_n}(m\omega_n)$ as a diamond-colored distributive lattice.

Similarly, we define a {\em symplectic ideal pattern} of size $n$ and bounded by $m$ to be an array $(c_{i,j})_{1 \leq j \leq i \leq n}$ of nonnegative integers with $c_{n,n} \leq m$ and satisfying the inequalities pictured in \SymplecticPatternFigure. 
If we partially order the collection $L_{\mysmallC_n}(m\omega_n)$ as in the odd orthogonal case, the result is a diamond-colored distributive lattice. 
In the symplectic case, we refer to the fundamental weight $\omega_n$ as the ``Catalan fundamental weight'' because the dimension of the corresponding fundamental representation of the symplectic Lie algebra is well-known to be a Catalan number. 

\begin{figure}[ht] 
\begin{center}
\SymplecticPatternFigure: Inequalities for a symplectic ideal pattern.

\vspace*{0.25in}
\begin{tabular}{cccccccccc}
 & & & & & & & & & $2c_{n,1}$ \\
 & & & & & & & & \begin{rotate}{45}$\geq$\end{rotate} & \\
 & & & & & & & $c_{n-1,1}$ & & \\
 & & & & & & \begin{rotate}{45}$\geq$\end{rotate} & & \begin{rotate}{315}$\!\!\!\!\leq$\end{rotate} & \\
 & & & & & $c_{n-2,1}$ & & & & $2c_{n,2}$ \\
 & & & & \begin{rotate}{45}$\geq$\end{rotate} & & \begin{rotate}{315}$\!\!\!\!\leq$\end{rotate} & & \begin{rotate}{45}$\geq$\end{rotate} & \\
 & & & $\iddots$ & & $\cdot$ & $\cdot$ & $\cdot$ & & \\
 & \begin{rotate}{45}$\geq$\end{rotate} & & & & & & & & \\
$c_{1,1}$ & & & & & $\cdot$ & $\cdot$ & $\cdot$ & & \\
 & \begin{rotate}{315}$\!\!\!\!\leq$\end{rotate} & & & & & & & & \\
 & & & $\ddots$ & & $\cdot$ & $\cdot$ & $\cdot$ & & \\
 & & & & \begin{rotate}{315}$\!\!\!\!\leq$\end{rotate} & & \begin{rotate}{45}$\geq$\end{rotate} & & \begin{rotate}{315}$\!\!\!\!\leq$\end{rotate} & \\
 & & & & & $c_{n-2,n-2}$ & & & & $2c_{n,n-1}$ \\
 & & & & & & \begin{rotate}{315}$\!\!\!\!\leq$\end{rotate} & & \begin{rotate}{45}$\geq$\end{rotate} & \\
 & & & & & & & $c_{n-1,n-1}$ & & \\
 & & & & & & & & \begin{rotate}{315}$\!\!\!\!\leq$\end{rotate} & \\
 & & & & & & & & & $2c_{n,n}$ 
\end{tabular}
\end{center}
\end{figure} 

When $n$ is even, we define an {\em even orthogonal ideal pattern} of size $n$ and bounded by $m$ to be an array $(c_{i,j})_{1 \leq j \leq i \leq n-2} \cup (c_{n-1,j})_{1 \leq j \leq \frac{n}{2}} \cup (c_{n,j})_{1 \leq j \leq \frac{n}{2}-1}$ of nonnegative integers with $c_{n-1,\frac{n}{2}} \leq m$ and satisfying the inequalities pictured in \EvenOrthPatternFigure.E. 
When $n$ is odd, we define an {\em even orthogonal ideal pattern} of size $n$ and bounded by $m$ to be an array $(c_{i,j})_{1 \leq j \leq i \leq n-2} \cup (c_{n-1,j})_{1 \leq j \leq \frac{n-1}{2}} \cup (c_{n,j})_{1 \leq j \leq \frac{n-1}{2}}$ of nonnegative integers with $c_{n,\frac{n-1}{2}} \leq m$ and satisfying the inequalities pictured in \EvenOrthPatternFigure.O. 
We use $L_{\mysmallD_n}(m\omega_{n-1})$ to denote the collection of even orthogonal ideal patterns of size $n$ and bounded by $m$ and partially ordered as in the odd orthogonal and symplectic cases. 
The result, as in the preceding cases, is a diamond-colored distributive lattice. 
An alternative array of inequalities for even orthogonal ideal patterns replaces each $c_{n-1,k}$ with $c_{n,k}$ and vice-versa. 
The resulting diamond-colored distributive lattice of such arrays is denoted $L_{\mysmallD_n}(m\omega_{n})$. 
Another way to realize $L_{\mysmallD_n}(m\omega_{n})$ is by recoloring the edges of $L_{\mysmallD_n}(m\omega_{n-1})$ by exchanging colors $n-1$ and $n$. 
However, whenever we refer to even orthogonal ideal patterns in $L_{\mysmallD_n}(m\omega_{n})$, we will mean those patterns resulting from replacing each $c_{n-1,k}$ with $c_{n,k}$ and vice-versa. 
In the $\myD_n$ case, we refer to $\omega_{n-1}$ and $\omega_{n}$ as the ``spin-node fundamental weights,'' for obvious reasons. 

We will show that these symplectic and orthogonal lattices are splitting distributive lattices for the corresponding Weyl bialternants.  
The next result establishes hypotheses that are needed in order to invoke \ChainProductVersion.B. 

\noindent 
{\bf \StructureProp}\ \ {\sl Let $\Phi$ be one of the root systems} $\{\myB_n, \myC_n, \myD_n\}$ {\sl and let $L$ be one of the edge-colored distributive lattices} $\{L_{\mysmallB_n}(m\omega_n),L_{\mysmallC_n}(m\omega_n),L_{\mysmallD_n}(m\omega_{n-1}),L_{\mysmallD_n}(m\omega_{n})\}$. 
{\sl For each $i \in I$, each $i$-component of $L$ is isomorphic to a product of chains. 
Moreover, $L$ is ${\Phi}$-structured, and} $\WGF(L)$ {\sl is a $\Phi$-Weyl symmetric function.}  

{\em Proof.} If we show that $L$ is ${\Phi}$-structured and that for each $i \in I$, each $i$-component of $L$ is isomorphic to a product of chains, then it follows from \WInvariantLemma\ that $\WGF(L)$ is a $\Phi$-Weyl symmetric function. 
For any $\telt \in L$ and $i \in I$, it is clear that $\comp_{i}(\telt)$ is comprised precisely of those ideal patterns $\selt$ such that $c_{j,k}(\selt) = c_{j,k}(\telt)$ for all $j \not= i$. 
Moreover, the defining pattern of inequalities allows that changes can be made to positions $(i,k)$ and $(i,k')$ independently of one another. 
That is, $\comp_{i}(\telt) \cong \mathcal{C}_{1} \times \cdots \times \mathcal{C}_{i}$, where for each $1 \leq k \leq i$, $\mathcal{C}_{k}$ is a chain whose length is the difference of the maximum value that can replace $c_{i,k}(\telt)$ and the minimum such value. 

To establish the ${\Phi}$-structure property, one can use reasoning similar to the proof of \GelfandProposition\ to work out formulas for $m_{i}(\telt)$. 
In all of the formulas that are stated next, any $c_{p,q}(\telt)$ is omitted in the formula if the position $(p,q)$ is not part of the array. 
For case $\Phi = \myB_n$ and $1 \leq i \leq n-1$ we have: 
\[m_{i}(\telt) = \sum_{k=1}^{i}[2c_{i,k}(\telt) - \min(c_{i-1,k}(\telt),c_{i+1,k+1}(\telt)) - \max(c_{i-1,k-1}(\telt),c_{i+1,k}(\telt))],\]
while for $i = n$ we have: 
\[m_{n}(\telt) = \sum_{k=1}^{n}[2c_{n,k}(\telt) - c_{n-1,k}(\telt) - c_{n-1,k-1}(\telt)].\]
For case $\Phi = \myC_n$ and $1 \leq i \leq n-2$ we have: 
\[m_{i}(\telt) = \sum_{k=1}^{i}[2c_{i,k}(\telt) - \min(c_{i-1,k}(\telt),c_{i+1,k+1}(\telt)) - \max(c_{i-1,k-1}(\telt),c_{i+1,k}(\telt))],\]
while for $i = n-1$ we have: 
\[m_{n-1}(\telt) = \sum_{k=1}^{n-1}[2c_{n-1,k}(\telt) - \min(c_{n-2,k}(\telt),2c_{n,k+1}(\telt)) - \max(c_{n-2,k-1}(\telt),2c_{n,k}(\telt))],\] 
and for $i = n$ we have: 
\[m_{n}(\telt) = \sum_{k=1}^{n}\left[2c_{n,k}(\telt) - \left\lfloor \frac{c_{n-1,k}(\telt)}{2} \right\rfloor - \left\lceil \frac{c_{n-1,k-1}(\telt)}{2} \right\rceil\right],\]
with $c_{n-1,n}(\telt) := 2k$ and $c_{n-1,0}(\telt) := 0$. 
For case $\Phi = \myD_n$ and $L = L_{\mysmallD_n}(m\omega_{n-1})$, then for $1 \leq i \leq n-3$ we have: 
\[m_{i}(\telt) = \sum_{k=1}^{i}[2c_{i,k}(\telt) - \min(c_{i-1,k}(\telt),c_{i+1,k+1}(\telt)) - \max(c_{i-1,k-1}(\telt),c_{i+1,k}(\telt))],\]
while for $i = n-2$ we have: 
\begin{eqnarray*}
m_{n-2}(\telt) & = & \sum_{k=1}^{n-2}\left[2c_{n-2,k}(\telt) - \min\left(c_{n-3,k}(\telt), \left\{\begin{array}{c}c_{n,\frac{k+1}{2}}\mbox{ if }k\mbox{ is odd} \\ c_{n-1,\frac{k+2}{2}}\mbox{ if }k\mbox{ is even}\end{array}\right\}\right)\right.\\ 
& &  \left. \hspace*{1.5in} - \max\left(c_{n-3,k-1}(\telt), \left\{\begin{array}{c}c_{n,\frac{k+1}{2}}\mbox{ if }k\mbox{ is odd} \\ c_{n-1,\frac{k}{2}}\mbox{ if }k\mbox{ is even}\end{array}\right\}\right)\right],\end{eqnarray*}
while for $i = n-1$ we have: 
\[m_{n-1}(\telt) = \sum_{k=1}^{\lfloor \frac{n}{2} \rfloor}[2c_{n-1,k}(\telt) - c_{n-2,2k-2}(\telt) - c_{n-2,2k-1}(\telt)],\]
and for $i = n$ we have: 
\[m_{n}(\telt) = \sum_{k=1}^{\lfloor \frac{n-1}{2} \rfloor}[2c_{n,k}(\telt) - c_{n-2,2k-1}(\telt) - c_{n-2,2k}(\telt)].\] 
The formulas for $L = L_{\mysmallD_n}(m\omega_{n})$ are the same except that $m_{n-1}(\telt)$ and $m_{n}(\telt)$ are switched. 

To see that $L$ is ${\Phi}$-structured, it suffices to check that when $\selt \myarrow{i} \telt$ in $L$, then $m_{j}(\selt) + \langle \alpha_{i},\alpha_{j}^{\vee} \rangle = m_{j}(\telt)$ for all $j \not= i$. 
To begin, let $\Phi \in \{\myB_{n},\myC_{n}\}$, and let $J := \{1,2,\cdots,n-1\} \subset \{1,2,\cdots,n-1,n\} = I$. 
Any given $J$-component of $L$ is formed by fixing all $c_{n,k}$'s and allowing all $c_{p,q}$'s to vary as long as $1 \leq p \leq n-1$. 
Therefore each such component is a Gelfand--Tsetlin lattice and, by \GelfandProposition, is $\myA_{n-1}$-structured. 
So, to confirm that $L$ is ${\Phi}$-structured, we must check that $m_{j}(\selt) + \langle \alpha_{i},\alpha_{j}^{\vee} \rangle = m_{j}(\telt)$ when $\selt \myarrow{i} \telt$ and when exactly one of $i$ or $j$ is equal to $n$. 
As in the proof of  \GelfandProposition, it is easy to see that $m_{j}(\selt) = m_{j}(\telt)$ when $i$ and $j$ are ``distant,'' i.e.\ when $|i-j| \geq 2$ and $\langle \alpha_{i},\alpha_{j}^{\vee} \rangle = 0$. 
So suppose $\Phi = \myB_{n}$, $i=n-1$, and $j=n$. 
Let's say $\selt \mylongarrow{n-1} \telt$ with $c_{n-1,q}(\telt) = c_{n-1,q}(\selt)+1$. 
Then 
{\small \begin{eqnarray*}
m_{n}(\telt) - m_{n}(\selt) & = & 2c_{n,q}(\telt) - c_{n-1,q}(\telt) - c_{n-1,q-1}(\telt) + 2c_{n,q+1}(\telt) - c_{n-1,q+1}(\telt) - c_{n-1,q}(\telt)\\ 
& & - 2c_{n,q}(\telt) + c_{n-1,q}(\telt) -1 + c_{n-1,q-1}(\telt) - 2c_{n,q+1}(\telt) + c_{n-1,q+1}(\telt) + c_{n-1,q}(\telt) - 1\\
& = & -2.
\end{eqnarray*}}Next take 
$\Phi = \myB_{n}$, $i=n$, $j=n-1$, and $\selt \myarrow{n} \telt$ with $c_{n,q}(\telt) = c_{n,q}(\selt)+1$. 
Then 
{\small \begin{eqnarray*}
m_{n-1}(\telt) - m_{n-1}(\selt) & = & 2c_{n-1,q}(\telt) - \min(c_{n-2,q}(\telt),c_{n,q+1}(\telt)) - \max(c_{n-2,q-1}(\telt),c_{n,q}(\telt))\\
& & -2c_{n-1,q}(\telt) + \min(c_{n-2,q}(\telt),c_{n,q+1}(\telt)) + \max(c_{n-2,q-1}(\telt),c_{n,q}(\telt)-1)\\
& & 2c_{n-1,q-1}(\telt) - \min(c_{n-2,q-1}(\telt),c_{n,q}(\telt)) - \max(c_{n-2,q-2}(\telt),c_{n,q-1}(\telt))\\
& & -2c_{n-1,q-1}(\telt) + \min(c_{n-2,q-1}(\telt),c_{n,q}(\telt)-1) + \max(c_{n-2,q-2}(\telt),c_{n,q-1}(\telt))\\
& = & -1.
\end{eqnarray*}}The latter equality follows by assessing the cases (1) $c_{n,q}(\telt) > c_{n-2,q-1}(\telt)$ and (2) $c_{n,q}(\telt) \leq c_{n-2,q-1}(\telt)$. 
Similarly 
see that when $\Phi = \myC_{n}$, then $m_{n-1}(\selt) - 2 = m_{n-1}(\telt)$ when $\selt \myarrow{n} \telt$ and $m_{n}(\selt) - 1 = m_{n}(\telt)$ when $\selt \mylongarrow{n-1} \telt$. 
This completes our demonstration that $L$ is ${\Phi}$-structured when $\Phi \in \{\myB_{n},\myC_{n}\}$.

Now take $\Phi = \myD_{n}$.  
Curiously, our demonstration in this case that $L$ is ${\Phi}$-structured is {\em easier} than when $\Phi \in \{\myB_{n},\myC_{n}\}$. 
Indeed, if we let $J := \{1,2,\cdots,n-2,n-1\}$ as in the previous paragraph, then one can see that each $J$-component of $L$ is a Gelfand--Tsetlin lattice, and therefore $L$ is $\myA_{n-1}$-structured, where $\myA_{n-1}$ is regarded as the root subsystem spanned by $\{\alpha_{j}\}_{j \in J}$.  
On the other hand, if we let $J' := \{1,2,\cdots,n-2,n\}$, then it is apparent that each $J'$-component of $L$ is also a Gelfand--Tsetlin lattice. 
So, $L$ is $\myA_{n-1}$-structured, where this time $\myA_{n-1}$ is regarded as the root subsystem spanned by $\{\alpha_{j}\}_{j \in J'}$. 
In particular, to show that $L$ is $\myD_{n}$-structured, it only remains to be shown that $m_{n-1}(\selt) = m_{n-1}(\telt)$ when $\selt \myarrow{n} \telt$ and that $m_{n}(\selt) = m_{n}(\telt)$ when $\selt \mylongarrow{n-1} \telt$. 
But these assertions follow readily from the formulas above for $m_{n-1}$ and $m_{n}$. 
\hfill\QED

In order to apply \ChainProductVersion.B, we need an appropriate vertex-coloring function for each of our symplectic and orthogonal lattices. 
Since the arrays defining the symplectic and orthogonal ideal patterns have the same general shape as the arrays for the special linear ideal patterns, then we will use the same slantwise reading order --- from SE to NW along successive diagonals, starting from the bottom of a given array --- in order to determine our vertex-coloring function, cf.\ \GelfandLemma\ and the paragraph preceding that lemma. 
Here is the analog of \GelfandLemma\ for our symplectic and orthogonal cases. 
The reasoning of the proof of \GelfandLemma\ carries over straightforwardly. 

\noindent
{\bf \ColorLemma}\ \ {\sl Keep the set-up of \StructureProp. 
Let $\melt$ denote the unique maximal element of $L$.  Let $\telt \not= \melt$ in $L$, and suppose $(i,j)$ is the slantwise-least maximizable position for $\telt$.  (1) Then $c_{p,q}(\telt) = c_{p,q}(\melt)$ for all positions $(p,q)$ that are slantwise-prior to $(i,j)$.  (2) If $\selt \in \comp_{i}(\telt)$ with $c_{i,j}(\selt) < c_{i,j}(\melt)$, then $(i,j)$ is the slantwise-least maximizable position for $\selt$. (3)  Suppose $\selt \in \comp_{i}(\telt)$ such that $\selt \not= \melt$ and the slantwise-least maximizable position for $\selt$ is some $(i,k)$. Then $k=j$.}

{\em Proof.}  Entirely similar to the proof of \GelfandLemma.\hfill\QED

Now we define the requisite vertex-coloring function. 
Take \[L \in \{L_{\mysmallB_n}(m\omega_n), L_{\mysmallC_n}(m\omega_n), L_{\mysmallD_n}(m\omega_{n-1}), L_{\mysmallD_n}(m\omega_{n})\}\] with $\melt$ maximal in $L$. 
Define $\kappa: L \setminus \{\melt\} \longrightarrow I$ by the following rule: For $\telt \ne \melt$ in $L$, such that $(i,j)$ is the slantwise-least maximizable position for $\telt$, declare that $\kappa(\telt) = i$. 

\noindent 
{\bf \ColorProp}\ \ {\sl Keep the notation of the preceding paragraph, with $\telt \in L \setminus \{\melt\}$. 
Then $\{\selt \in \comp_{\kappa(\telt)}(\telt)\, |\, \selt \not= \melt \mbox{ and } \kappa(\selt) = \kappa(\telt)\}$ is a sub-face of $\comp_{\kappa(\telt)}(\telt)$.} 

{\em Proof.} For $\telt \not= \melt$, take $i := \kappa(\telt)$, and let $(i,j)$ be the slantwise-least maximizable position for $\telt$. 
As in the proof of \GelfandColoring, we get the following equality of sets by applying \ColorLemma: $\{\selt \in \comp_{i}(\telt)\, |\, \selt \not= \melt \mbox{ and } \kappa(\selt) = i\} = \{\selt \in \comp_{i}(\telt)\, |\, c_{i,j}(\selt) < c_{i,j}(\melt)\}$.
The latter set is clearly a sub-face of $\comp_{i}(\telt)$.\hfill\QED

From \PrelimProps, \ChainProductVersion.B, and \MainCorollary, it follows immediately that the symplectic and orthogonal lattices constructed above are splitting distributive lattices for symplectic and orthogonal Weyl bialternants respectively associated with multiples of the Catalan and spin-node fundamental weights and that their rank generating functions have nice quotient-of-product expressions.  In particular, we have: 

\noindent 
{\bf \ApplicationTheorem}\ \ {\sl Let $\Phi$ be one of} $\myB_n$, $\myC_n$, {\sl or} $\myD_n$, {\sl and let $L_{\Phi}(\lambda)$ be one of} $L_{\mysmallB_n}(m\omega_n)$, $L_{\mysmallC_n}(m\omega_n)$, $L_{\mysmallD_n}(m\omega_{n-1})$, {\sl or} $L_{\mysmallD_n}(m\omega_{n})$ {\sl with $\lambda \in \{m\omega_{n-1}, m\omega_{n}\}$ as appropriate. 
Then $L_{\Phi}(\lambda)$ is a rank symmetric and rank unimodal splitting distributive lattice for the $\Phi$-Weyl bialternant $\chi_{_{\lambda}}^{\Phi}$.  Moreover,} 
\[\RGF(L_{\Phi}(\lambda),q) = \frac{\mbox{$\displaystyle \prod_{\alpha \in \Phi^{+}}$}(1-q^{\langle \lambda+\varrho,\alpha^{\vee} \rangle})}{\mbox{$\displaystyle \prod_{\alpha \in \Phi^{+}}$}(1-q^{\langle \varrho,\alpha^{\vee} \rangle})}.\]

\vspace*{-0.4in}\hfill\QED

For a more explicit version of the $\RGF$ polynomials above, see Theorem 10.6.3 of \cite{DonDistributive}. 
Only the symplectic results of this theorem are brand new. 
The preceding formulas for the rank generating functions of $L_{\mysmallB_n}(m\omega_n)$, $L_{\mysmallD_n}(m\omega_{n-1})$, and $L_{\mysmallD_n}(m\omega_{n})$ can be viewed as instances of the Bender--Knuth (ex-)conjecture. 
For some further discussion of this result, see \cite{PrEur} or \cite{Fischer}. 
So in the orthogonal cases we have, in effect, a new proof of this result by way of splitting posets, although our approach is very similar in philosophy to Proctor's in \cite{PrEur}. 
The above expression for the rank generating function in the symplectic case could be viewed as a new symplectic analog of Bender--Knuth. 

\begin{figure}[th] 
\begin{center}
\EvenOrthPatternFigure.E: Inequalities for an even orthogonal ideal pattern, with $n$ even.

\vspace*{0.25in}
\begin{tabular}{ccccccccccccc}
 & & & & & & & & & & & & $c_{n-1,1}$ \\
 & & & & & & & & & & & \begin{rotate}{45}$\geq$\end{rotate} & \\
 & & & & & & & & & & $c_{n-2,1}$ & & \\
 & & & & & & & & & \begin{rotate}{45}$\geq$\end{rotate} & & \begin{rotate}{315}$\!\!\!\!\leq$\end{rotate} & \\
 & & & & & & & & $c_{n-3,1}$ & & & & $c_{n,1}$ \\
 & & & & & & & \begin{rotate}{45}$\geq$\end{rotate} & & \begin{rotate}{315}$\!\!\!\!\leq$\end{rotate} & & \begin{rotate}{45}$\geq$\end{rotate} & \\
 & & & & & & $c_{n-4,1}$ & & & & $c_{n-2,2}$ & & \\
 & & & & & \begin{rotate}{45}$\geq$\end{rotate} & & \begin{rotate}{315}$\!\!\!\!\leq$\end{rotate} & & \begin{rotate}{45}$\geq$\end{rotate} & & \begin{rotate}{315}$\!\!\!\!\leq$\end{rotate} & \\
 & & & & $c_{n-5,1}$ & & & & $c_{n-3,2}$ & & & & $c_{n-1,2}$ \\
 & & & \begin{rotate}{45}$\geq$\end{rotate} & & \begin{rotate}{315}$\!\!\!\!\leq$\end{rotate} & & \begin{rotate}{45}$\geq$\end{rotate} & & \begin{rotate}{315}$\!\!\!\!\leq$\end{rotate} & & \begin{rotate}{45}$\geq$\end{rotate} & \\
 & & $\ \iddots\ \ $ & & & & & $\cdot$ & $\cdot$ & $\cdot$ & & &  \\
 & \begin{rotate}{45}$\geq$\end{rotate} & & & & & & & & & & & \\
$\ c_{1,1}\ \ $ & & & & & & & $\cdot$ & $\cdot$ & $\cdot$ & & & \\
 & \begin{rotate}{315}$\!\!\!\!\leq$\end{rotate} & & & & & & & & & & \\
 & & $\ \ddots\ \ $ & & & & &  $\cdot$ & $\cdot$ & $\cdot$ & & & \\
 & & & \begin{rotate}{315}$\!\!\!\!\leq$\end{rotate} & & \begin{rotate}{45}$\geq$\end{rotate} & & \begin{rotate}{315}$\!\!\!\!\leq$\end{rotate} & & \begin{rotate}{45}$\geq$\end{rotate} & & \begin{rotate}{315}$\!\!\!\!\leq$\end{rotate} & \\
 & & & & $c_{n-5,n-5}$ & & & & $c_{n-3,n-4}$ & & & & $c_{n-1,\frac{n}{2}-1}$ \\
 & & & & & \begin{rotate}{315}$\!\!\!\!\leq$\end{rotate} & & \begin{rotate}{45}$\geq$\end{rotate} & & \begin{rotate}{315}$\!\!\!\!\leq$\end{rotate} & & \begin{rotate}{45}$\geq$\end{rotate} & \\
 & & & & & & $c_{n-4,n-4}$ & & & & $c_{n-2,n-3}$ & & \\
 & & & & & & & \begin{rotate}{315}$\!\!\!\!\leq$\end{rotate} & & \begin{rotate}{45}$\geq$\end{rotate} & & \begin{rotate}{315}$\!\!\!\!\leq$\end{rotate} & \\
 & & & & & & & & $c_{n-3,n-3}$ & & & & $c_{n,\frac{n}{2}-1}$ \\
 & & & & & & & & & \begin{rotate}{315}$\!\!\!\!\leq$\end{rotate} & & \begin{rotate}{45}$\geq$\end{rotate} & \\
 & & & & & & & & & & $c_{n-2,n-2}$ & & \\
 & & & & & & & & & & & \begin{rotate}{315}$\!\!\!\!\leq$\end{rotate} & \\
 & & & & & & & & & & & & $c_{n-1,\frac{n}{2}}$ 
\end{tabular}
\end{center}
\end{figure} 

\begin{figure}[th] 
\begin{center}
\EvenOrthPatternFigure.O: Inequalities for an even orthogonal ideal pattern, with $n$ odd.

\vspace*{0.25in}
\begin{tabular}{ccccccccccccc}
 & & & & & & & & & & & & $c_{n-1,1}$ \\
 & & & & & & & & & & & \begin{rotate}{45}$\geq$\end{rotate} & \\
 & & & & & & & & & & $c_{n-2,1}$ & & \\
 & & & & & & & & & \begin{rotate}{45}$\geq$\end{rotate} & & \begin{rotate}{315}$\!\!\!\!\leq$\end{rotate} & \\
 & & & & & & & & $c_{n-3,1}$ & & & & $c_{n,1}$ \\
 & & & & & & & \begin{rotate}{45}$\geq$\end{rotate} & & \begin{rotate}{315}$\!\!\!\!\leq$\end{rotate} & & \begin{rotate}{45}$\geq$\end{rotate} & \\
 & & & & & & $c_{n-4,1}$ & & & & $c_{n-2,2}$ & & \\
 & & & & & \begin{rotate}{45}$\geq$\end{rotate} & & \begin{rotate}{315}$\!\!\!\!\leq$\end{rotate} & & \begin{rotate}{45}$\geq$\end{rotate} & & \begin{rotate}{315}$\!\!\!\!\leq$\end{rotate} & \\
 & & & & $c_{n-5,1}$ & & & & $c_{n-3,2}$ & & & & $c_{n-1,2}$ \\
 & & & \begin{rotate}{45}$\geq$\end{rotate} & & \begin{rotate}{315}$\!\!\!\!\leq$\end{rotate} & & \begin{rotate}{45}$\geq$\end{rotate} & & \begin{rotate}{315}$\!\!\!\!\leq$\end{rotate} & & \begin{rotate}{45}$\geq$\end{rotate} & \\
 & & $\ \iddots\ \ $ & & & & & $\cdot$ & $\cdot$ & $\cdot$ & & &  \\
 & \begin{rotate}{45}$\geq$\end{rotate} & & & & & & & & & & & \\
$\ c_{1,1}\ \ $ & & & & & & & $\cdot$ & $\cdot$ & $\cdot$ & & & \\
 & \begin{rotate}{315}$\!\!\!\!\leq$\end{rotate} & & & & & & & & & & \\
 & & $\ \ddots\ \ $ & & & & &  $\cdot$ & $\cdot$ & $\cdot$ & & & \\
 & & & \begin{rotate}{315}$\!\!\!\!\leq$\end{rotate} & & \begin{rotate}{45}$\geq$\end{rotate} & & \begin{rotate}{315}$\!\!\!\!\leq$\end{rotate} & & \begin{rotate}{45}$\geq$\end{rotate} & & \begin{rotate}{315}$\!\!\!\!\leq$\end{rotate} & \\
 & & & & $c_{n-5,n-5}$ & & & & $c_{n-3,n-4}$ & & & & $c_{n,\frac{n-1}{2}-1}$ \\
 & & & & & \begin{rotate}{315}$\!\!\!\!\leq$\end{rotate} & & \begin{rotate}{45}$\geq$\end{rotate} & & \begin{rotate}{315}$\!\!\!\!\leq$\end{rotate} & & \begin{rotate}{45}$\geq$\end{rotate} & \\
 & & & & & & $c_{n-4,n-4}$ & & & & $c_{n-2,n-3}$ & & \\
 & & & & & & & \begin{rotate}{315}$\!\!\!\!\leq$\end{rotate} & & \begin{rotate}{45}$\geq$\end{rotate} & & \begin{rotate}{315}$\!\!\!\!\leq$\end{rotate} & \\
 & & & & & & & & $c_{n-3,n-3}$ & & & & $c_{n-1,\frac{n-1}{2}}$ \\
 & & & & & & & & & \begin{rotate}{315}$\!\!\!\!\leq$\end{rotate} & & \begin{rotate}{45}$\geq$\end{rotate} & \\
 & & & & & & & & & & $c_{n-2,n-2}$ & & \\
 & & & & & & & & & & & \begin{rotate}{315}$\!\!\!\!\leq$\end{rotate} & \\
 & & & & & & & & & & & & $c_{n,\frac{n-1}{2}}$ 
\end{tabular}
\end{center}
\end{figure}

\clearpage
\noindent 
{\Large \bf \S \GaussianSection.\ Two new proofs that minuscule posets are Gaussian.} 

In \S \SplittingSection, and particularly \MinQuasiMinTheorem, we considered the minuscule splitting posets, which are the unique splitting posets affiliated with minuscule dominant weights. 
It is easy to see (for example, by \MinQuasiMinProposition) that a minuscule dominant weight must be a fundamental weight. 
In Theorem 12.6 of \cite{DonDistributive}, it is shown, using general principles and independent of the classification of minuscule dominant weights by root system type, that any minuscule splitting poset is a diamond-colored distributive lattice. 
Henceforth, we call a minuscule splitting poset a `minuscule splitting distributive lattice'. 
We may therefore consider its associated vertex-colored poset of join irreducibles, which in \cite{DonDistributive} is called a `minuscule compression poset'. 
Again, using general principles, it is shown in Theorem 12.12 of \cite{DonDistributive} that any minuscule compression poset is itself a distributive lattice and a full-length sublattice of some product of two (non-empty) chains, although this fact will not be needed here. 
Minuscule compression posets have many other remarkable properties, some of which we re-explore here.  

Fix for the remainder of this section an irreducible root system $\Phi$ with simple roots indexed by a set $I$ with $n$ elements and a minuscule dominant weight $\omega$. 
The affiliated minuscule splitting distributive lattice is $R_{\Phi}(\omega)$. 
Let $P_{\Phi}(\omega) := \mathbf{j}_{color}(R_{\Phi}(\omega))$ denote its associated minuscule compression poset.  
Throughout, $\mysmallindexM$ denotes a positive integer and $\mathcal{C}_{\mysmallerindexM}$ denotes a chain with $\mysmallindexM$ elements.  
The product poset $P_{\Phi}(\omega) \times C_{\mysmallerindexM}$ is naturally vertex-colored by declaring that the color of any element $(v,k) \in P_{\Phi}(\omega) \times C_{\mysmallerindexM}$ is the same as the color of $v$.  
In \cite{DonDistributive}, the vertex-colored poset $P_{\Phi}(\omega) \times C_{\mysmallerindexM}$ is called a {\em minuscule prism}. 
We let $L_{\Phi}(\mysmallindexM\omega) := \mathbf{J}_{color}(P_{\Phi}(\omega) \times C_{\mysmallerindexM})$, the lattice of order ideals from the minuscule prism $P_{\Phi}(\omega) \times C_{\mysmallerindexM}$. 

{\bf Gaussian posets.} 
In \cite{StanThesis}, Stanley defined a poset $P$ to be {\em Gaussian} if there exists a positive integer $r$ and positive integers $\hbar_{1}$, \ldots, $\hbar_{r}$ such that for every positive integer $\mysmallindexM$, 
\[\RGF\left(\rule[-1.5mm]{0mm}{5mm}\mathbf{J}(P \times C_{\mysmallerindexM}),q\right) = \mbox{\large $\displaystyle \prod_{i=1}^{r}$}\left(\frac{1-q^{\mysmallerindexM+\hbar_{i}}}{1-q^{\hbar_{i}}}\right),\] 
where $C_{\mysmallerindexM}$ denotes a chain with $\mysmallindexM$ elements. 
In \cite{PrEur}, Proctor, in collaboration with Stanley, demonstrated the following: 

\noindent
{\small {\bf The Proctor--Stanley Gaussian Poset Theorem}\ \  {\sl All minuscule compression posets are Gaussian.}}

\noindent Proctor's and Stanley's proof of this result crucially utilized Seshadri's standard monomial theory \cite{Seshadri}. 
In the nomenclature of this manuscript, their key idea was to recognize that the rank generating function $\RGF\left(\rule[-1.5mm]{0mm}{5mm}\mathbf{J}(P \times C_{\mysmallerindexM}),q\right)$ can be viewed as a specialization of a particular Weyl bialternant (see \GaussianFromWGF, below). 
Whether minuscule compression posets are the only connected Gaussian posets is a longstanding open question. 
For some related discussion, see \cite{Green}.

{\bf Two new approaches.} 
In this section, we use methodologies developed in this monograph to provide what we believe are two new proofs of the Proctor--Stanley Gaussian Poset Theorem. 
First, we will take a crystalline splitting poset perspective that uses ideas from \S \FibrousSection\ and \S \CrystallineSection.  
And second, we will utilize the vertex-coloring method of \S \CriteriaSection.
With both approaches, we proceed from the assumption that minuscule splitting distributive lattices have already been obtained as in \S \SplittingSection, so properties of these lattices and their minuscule compression posets will be developed within this context. 
A different conceptual starting point might begin by defining minuscule compression posets as `dominant minuscule heaps,' 
as in \cite{StemMinuscule}, or by characterizing them combinatorially, as in \cite{StrayerThesis} and \cite{StrayerArxiv}; one would then derive properties of these posets and their minuscule splitting distributive lattices from these respective starting points. 

There are at least the following commonalities to both of our new proofs of the Proctor--Stanley Gaussian Poset Theorem. 
Both arguments are general and do not depend on the classification of irreducible root systems. 
In both cases, Gaussian-ness of $P_{\Phi}(\omega)$ is concluded by demonstrating that each $L_{\Phi}(\mysmallindexM\omega)$ is a splitting distributive lattice for $\chi_{_{\mysmallestindexM\omega}}^{\Phi}$ and then applying the following lemma.  

\noindent 
{\bf \GaussianFromWGF}\ \ {\sl Let $R$ be any splitting poset for $\chi_{_{\mysmallestindexM\omega}}^{\Phi}$. 
Then} $\displaystyle \RGF(R,q) = \mbox{\large $\displaystyle \prod_{i=1}^{r}$}\left(\frac{1-q^{\mysmallerindexM+\hbar_{i}}}{1-q^{\hbar_{i}}}\right)$ {\sl where $\hbar_{i} := \langle \varrho,\beta_{i}^{\vee} \rangle$ and $\{\beta_{i}\}_{i=1}^{r}$ is a listing of exactly those positive roots $\beta \in \Phi^{+}$ for which $\langle \omega,\beta^{\vee} \rangle \not= 0$.}

{\em Proof.} Let $\beta$ be a generic positive root in $\Phi^{+}$. 
It follows from definitions that the quantity $\langle \lambda,\beta^{\vee} \rangle$ is nonnegative whenever $\lambda$ is dominant. 
By \MinQuasiMinProposition, we have $\langle \omega,\beta^{\vee} \rangle \in \{-1,0,1\}$. 
The two preceding facts combine to give us $\langle \omega,\beta^{\vee} \rangle \in \{0,1\}$. 
So, when $\beta=\beta_{i}$, then $\langle \mysmallindexM\omega+\varrho,\beta_{i}^{\vee} \rangle = \mysmallindexM+\hbar_{i}$.  
The result now follows from \MainCorollary.\hfill\QED

Moreover, both of our approaches to proving The Proctor--Stanley Gaussian Poset Theorem share the following result. 

\noindent
{\bf \GaussianStructureProp}\ \ {\sl The diamond-colored distributive lattice $L_{\Phi}(\mysmallindexM\omega)$ is $\Phi$-structured and its monochromatic components are chain products.}

{\em Proof.} This is a special case of Theorem 13.2 of \cite{DonDistributive}.\hfill\QED

Our first approach to proving that $L_{\Phi}(\mysmallindexM\omega)$ is a splitting poset for $\chi_{_{\mysmallestindexM\omega}}^{\Phi}$ will produce a weight-preserving bijection $R_{\Phi}(\mysmallindexM\omega)\ \widetilde{\mbox{\Large $\longrightarrow$}}\ L_{\Phi}(m\omega)$, where the crystalline splitting poset $R_{\Phi}(\mysmallindexM\omega)$ is understood to be the connected component of the maximal element of $R_{\Phi}(\omega)^{\otimes \mysmallerindexM}$. 
Our second approach will demonstrate that $L_{\Phi}(\mysmallindexM\omega)$ satisfies the main hypotheses of \ChainProductVersion, i.e.\ that $L_{\Phi}(\mysmallindexM\omega)$ is $\Phi$-structured, has a unique maximal element with weight $\mysmallindexM\omega$, and has kindred sets that are sub-faces of their respective monochromatic components. 

{\bf Some further aspects of our set-up.} 
From here on, we mostly suppress the super/subscript `$\Phi$.' 
Now, an element $\xelt$ in $L(\mysmallindexM\omega)$ is an order ideal from $P(\omega) \times C_{\mysmallerindexM}$ but can be identified as a nonnegative integer array $(g_{v}(\xelt))_{v \in P(\omega)}$ in the following way. 
The positions of the array are the elements of $P(\omega)$, and the array entry at position $v \in P(\omega)$ is $g_{v}(\xelt) = \left|\left\{k \in \{0,1,\ldots,\mysmallindexM\}\, \rule[-2.25mm]{0.1mm}{6mm}\, (v,k) \in \xelt\right\}\right|$. 
Observe that an array $(h_{v})_{v \in P(\omega)}$ issues from an order ideal from $P(\omega) \times C_{\mysmallerindexM}$ if and only if $h_{v} \in \{0,1,\ldots,\mysmallindexM\}$ for each $v \in P(\omega)$ and $h_{u} \geq h_{v}$ whenever $u \leq v$ in $P(\omega)$, in which case we call $(h_{v})_{v \in P(\omega)}$ an {\em ideal array}. 
The unique maximal element of $L(\mysmallindexM\omega)$ is $\melt := (\mysmallindexM)_{v \in P(\omega)}$.  

Many elementary combinatorial features of minuscule compression posets are developed in \S 12 of \cite{DonDistributive}. 
Amongst these key properties are the following. 
The minuscule compression poset $P(\omega)$ has a unique maximal element whose color is the same as the color of the fundamental weight $\omega$, cf.\ Lemma 12.7 of \cite{DonDistributive}. 
Now suppose $u$ and $v$ are vertices of $P(\omega)$ with vertex colors $i$ and $j$ respectively. 
See Lemma 12.8 of \cite{DonDistributive} for a proof of the following two conditional statements.  
If $u$ and $v$ are incomparable elements of $P(\omega)$, then $i \not= j$ and $M_{ij} = 0 = M_{ji}$; and 
If $u \rightarrow v$ is a covering relation in $P(\omega)$, then $i \not= j$ and the Cartan matrix entries $M_{ij}$ and $M_{ji}$ are both nonzero. 
From the latter conditional, together with the uniqueness of the maximal element of $P(\omega)$, it follows that the maximal element of $L(\mysmallindexM\omega)$ has weight $\mysmallindexM\omega$. 

{\bf First new proof.} 
For our first new proof of the Proctor--Stanley Gaussian Poset Theorem, we produce a weight-preserving bijection $R(\mysmallindexM\omega)\ \widetilde{\mbox{\Large $\longrightarrow$}}\ L(\mysmallindexM\omega)$. 
Since $R(\omega)$ is naturally isomorphic to $L(\omega)$, we regard their elements as the same. 
Let $\mathbf{max}$ denote the unique maximal element of $L(\omega) = R(\omega)$, so $g_{v}(\mathbf{max}) = 1$ for each $v \in P(\omega)$. 
A maximal element in $R(\omega)^{\otimes \mysmallerindexM}$ is the $\mysmallindexM$-tuple $\melt' = (\mathbf{max},\ldots,\mathbf{max})$. 
By \FibrousStructureLemma, $wt(\melt') = \mysmallindexM\omega$. 
Note that each $(\xelt_{1},\xelt_{2},\ldots,\xelt_{\mysmallerindexM}) \in R(\omega)^{\otimes \mysmallerindexM}$ is an $\mysmallindexM$-tuple of $\xelt_{i}$'s where $\xelt_{i} \in R(\omega)$ and that $R(\mysmallindexM\omega) \hookrightarrow R(\omega)^{\otimes \mysmallerindexM}$ is the connected component of $\melt'$. 
Let us define a subset $\mathcal{I}(\mysmallindexM\omega)$ of $R(\omega)^{\otimes \mysmallerindexM}$ by the rule that $(\xelt_{1},\ldots,\xelt_{\mysmallerindexM}) \in R(\omega)^{\otimes \mysmallerindexM}$ is in $\mathcal{I}(\mysmallindexM\omega)$ if and only if for each $j \in \{1,\ldots,\mysmallindexM-1\}$ and $v \in P(\omega)$ we have $g_{v}(\xelt_{j}) \leq g_{v}(\xelt_{j+1})$, i.e.\ $\xelt_{j} \leq \xelt_{j+1}$ in $R(\omega)$.  
Observe that $\melt' \in \mathcal{I}(\mysmallindexM\omega)$.

\noindent
{\bf \CrystalProductIdeal}\ \ {\sl Keeping the above notation, $\mathcal{I}(\mysmallindexM\omega) = R(\mysmallindexM\omega)$, an equality of sets.}

{\em Proof.} Suppose $(\yelt_{1},\ldots,\yelt_{\mysmallerindexM}) \in \mathcal{I}(\mysmallindexM\omega)$. 
Now suppose $(\xelt_{1},\ldots,\xelt_{\mysmallerindexM}) \myarrow{i} (\yelt_{1},\ldots,\yelt_{\mysmallerindexM})$ in $R(\omega)^{\otimes \mysmallerindexM}$, in which case we have $\rho_{i}(\yelt_{1},\ldots,\yelt_{\mysmallerindexM}) > 0$. 
We wish to show that $(\xelt_{1},\ldots,\xelt_{\mysmallerindexM}) \in \mathcal{I}(\mysmallindexM\omega)$. 
By \FibrousStructureLemma, we have  
\[(\xelt_{1},\ldots,\xelt_{\mysmallerindexM}) = \wF_{i}(\yelt_{1},\ldots,\yelt_{\mysmallerindexM}) = (\yelt_{1},\ldots,\wF_{i}(\yelt_r),\ldots,\yelt_{\mysmallerindexM}),\] where $r$ is the largest index for which the max occurs in the formula 
\[\rho_i(\yelt_{1},\ldots,\yelt_{\mysmallerindexM}) = \max_{1 \leq p \leq \mysmallerindexM}\bigg\{\rho_{i}(\yelt_p) + \left(\rule[-1.5mm]{0mm}{5mm}m_{i}(\yelt_{p+1}) + \cdots + m_{i}(\yelt_{\mysmallerindexM})\right)\bigg\}.\]  
Since $\xelt_{r} \myarrow{i} \yelt_{r}$ in $L(\omega)$, then for some $v \in P(\omega)$ of color $i$ we must have $g_{v}(\xelt_{r}) = 0 = g_{v}(\yelt_{r})-1$ while $g_{u}(\xelt_{r}) = g_{u}(\yelt_{r})$ for all $u \ne v$. 
So, to see that $(\xelt_{1},\ldots,\xelt_{\mysmallerindexM})$ meets the requirement for membership in $\mathcal{I}(\mysmallindexM\omega)$, it suffices to check that $g_{v}(\yelt_{r-1}) \leq g_{v}(\yelt_{r}) - 1$, i.e.\ that $g_{v}(\yelt_{r-1}) = 0$. 

Since $L(\omega)$ is a minuscule splitting distributive lattice and $\xelt_{r} \myarrow{i} \yelt_{r}$, then $\delta_{i}(\yelt_{r}) = 0$ and $m_{i}(\yelt_{r}) = \rho_{i}(\yelt_{r})$. 
Now suppose that $g_{v}(\yelt_{r-1}) = 1$. 
This hypothesis, together with the fact that $\yelt_{r-1} \leq \yelt_{r}$ in $L(\omega)$, affords that $\rho_{i}(\yelt_{r-1}) = 1$. 
Therefore, $\rho_{i}(\yelt_{r-1})+m_{i}(\yelt_{r}) + m_{i}(\yelt_{r+1}) + \cdots + m_{i}(\yelt_{\mysmallerindexM}) = 1 + \rho_{i}(\yelt_{r}) +  m_{i}(\yelt_{r+1}) + \cdots + m_{i}(\yelt_{\mysmallerindexM}) > \rho_{i}(\yelt_{r}) + m_{i}(\yelt_{r+1}) + \cdots + m_{i}(\yelt_{\mysmallerindexM})$, which contradicts that maximality of $\rho_{i}(\yelt_{r}) + m_{i}(\yelt_{r+1}) + \cdots + m_{i}(\yelt_{\mysmallerindexM})$. 
So, we must have $g_{v}(\yelt_{r-1}) = 0$. 
This is what we needed to show in order to conclude that $(\xelt_{1},\ldots,\xelt_{\mysmallerindexM}) \in \mathcal{I}(\mysmallindexM\omega)$. 
(Note that an entirely similar argument can be used to see that if $(\xelt_{1},\ldots,\xelt_{\mysmallerindexM}) \in \mathcal{I}(\mysmallindexM\omega)$ and if $(\xelt_{1},\ldots,\xelt_{\mysmallerindexM}) \myarrow{i} (\yelt_{1},\ldots,\yelt_{\mysmallerindexM})$ for some $(\yelt_{1},\ldots,\yelt_{\mysmallerindexM}) \in R(\omega)^{\otimes \mysmallerindexM}$, then $(\yelt_{1},\ldots,\yelt_{\mysmallerindexM}) \in \mathcal{I}(\mysmallindexM\omega)$.) 
Since $\melt' \in \mathcal{I}(\mysmallindexM\omega)$, and since each element of $R(\mysmallindexM\omega)$ can be obtained by applying some sequence of lowering operators to $\melt'$, we conclude that $R(\mysmallindexM\omega) \subseteq \mathcal{I}(\mysmallindexM\omega)$, an inclusion of sets. 

To show that $\mathcal{I}(\mysmallindexM\omega) \subseteq R(\mysmallindexM\omega)$, it suffices to show that for any $\xelt = (\xelt_{1},\ldots,\xelt_{\mysmallerindexM}) \in \mathcal{I}(\mysmallindexM\omega)$, there is a sequence (possibly empty) of raising operators that can be applied to $\xelt$ in order to obtain $\melt'$. 
We induct on the depth of elements $\xelt = (\xelt_{1},\ldots,\xelt_{\mysmallerindexM}) \in \mathcal{I}(\mysmallindexM\omega)$ as measured by $\delta(\xelt) := \sum_{j=1}^{\mysmallerindexM}\delta^{(\omega)}(\xelt_{j})$, where `$\delta^{(\omega)}$' is the depth function for $L(\omega)$. 
Of course, $\mathbf{max}$ is the unique element of $L(\omega)$ with zero depth. 
Suppose there is some $\xelt = (\xelt_{1},\ldots,\xelt_{\mysmallerindexM}) \in \mathcal{I}(\mysmallindexM\omega)$ such that $\delta(\xelt) = 0$. 
Then $\delta^{(\omega)}(\xelt_{j}) = 0$ for each $j \in \{1,2,\ldots,\mysmallindexM\}$, so  $\xelt = (\xelt_{1},\ldots,\xelt_{\mysmallerindexM}) = (\mathbf{max},\ldots,\mathbf{max}) = \melt'$. 
So $\melt'$ can be obtained from $\xelt$ by applying the empty sequence of raising operators to $\xelt$. 

Now suppose that for some nonnegative integer $k$, it is the case that whenever $\yelt = (\yelt_{1},\ldots,\yelt_{\mysmallerindexM}) \in \mathcal{I}(\mysmallindexM\omega)$ has depth $\delta(\yelt) \leq k$, then there is a sequence of raising operators that can be applied to $\yelt$ in order to obtain $\melt'$. 
Next, suppose that $\xelt = (\xelt_{1},\ldots,\xelt_{\mysmallerindexM}) \in \mathcal{I}(\mysmallindexM\omega)$ has depth $\delta(\xelt) = k+1$. 
Since $\delta(\xelt) > 0$, then there exists some smallest index $q \in \{1,2,\ldots,\mysmallindexM\}$ such that $\xelt_{1} = \cdots = \xelt_{q-1} = \xelt_{q}$ and either $\xelt_{q} < \xelt_{q+1}$ with $q < \mysmallindexM$ or $\xelt_{q} < \mathbf{max}$ with $q = \mysmallindexM$. 
Let $\xelt'_{q}$ be an element of $L(\omega)$ such that $\xelt_{q} \myarrow{i} \xelt'_{q}$ and such that $\xelt'_{q} \leq \xelt_{q+1}$ (if $q < \mysmallindexM$) or $\xelt'_{q} \leq \mathbf{max}$ (if $q = \mysmallindexM$). 
Now, $\xelt_{1} = \cdots = \xelt_{q}$ implies that $\delta_{i}(\xelt_{j}) = 1$ and $\rho_{i}(\xelt_{j}) = 0$ for each $j \in \{1,\ldots,q\}$, since $L(\omega)$ is a minuscule splitting distributive lattice. 
That is, $m_{i}(\xelt_{1}) = \cdots = m_{i}(\xelt_{q-1}) = -1$. 
Then, $-(m_{i}(\xelt_{1}) + \cdots + m_{i}(\xelt_{q-1}))+\delta_{i}(\xelt_{q}) = q > 0$, so by \FibrousStructureLemma, $\delta_{i}(\xelt) > 0$ when we regard $\xelt$ as an element of $R(\mysmallindexM\omega)$ and $\widetilde{E}_{i}(\xelt) = (\xelt_{1},\ldots,\xelt_{q-1},\widetilde{E}_{i}(\xelt_{q}),\xelt_{q+1},\ldots,\xelt_{\mysmallerindexM}) = (\xelt_{1},\ldots,\xelt_{q-1},\xelt'_{q},\xelt_{q+1},\ldots,\xelt_{\mysmallerindexM}) =: \xelt'$. 
But, since $\xelt_{q} < \xelt'_{q} \leq \xelt_{q+1}$ (if $q < \mysmallindexM$) or $\xelt_{q} < \xelt'_{q} \leq \mathbf{max}$ (if $q = \mysmallindexM$), then $\xelt'$ meets the membership requirement for $\mathcal{I}(\mysmallindexM\omega)$ and depth $\delta(\xelt') = \delta(\xelt)-1 = k$. 
So the inductive hypothesis applies to $\xelt'$. 
That is, there is some sequence of raising operators that, when applied to $\xelt'$, yields $\melt'$, so if we extend said sequence by $\widetilde{E}_{i}$ we have a sequence of raising operators that, when applied to $\xelt$, yields $\melt'$.  
This completes the induction argument, and the proof.\hfill\QED

We utilize the foregoing identification of the sets $R(\mysmallindexM\omega)$ and $\mathcal{I}(\mysmallindexM\omega)$ in the statement and proof of the following result. 

\noindent
{\bf \WeightPreservingBijection}\ \ {\sl Let $\phi: R(\mysmallindexM\omega) \longrightarrow L(\mysmallindexM\omega)$ be the set mapping given by $\phi(\xelt_{1},\ldots,\xelt_{\mysmallerindexM}) := \left(\sum_{j=1}^{\mysmallerindexM}g_{v}(\xelt_{j})\right)_{v \in P(\omega)}$.  Then $\phi$ is a well-defined weight-preserving bijection of $I$-colored posets.}

{\em Proof.} By \CrystalProductIdeal, we may characterize $R(\mysmallindexM\omega)$ as consisting precisely of those $\mysmallindexM$-tuples $(\xelt_{1},\ldots,\xelt_{\mysmallerindexM})$ with each $\xelt_{j} \in L(\omega)$ and with $\xelt_{1} \leq \cdots \leq \xelt_{\mysmallerindexM}$. 
Given such an $\mysmallindexM$-tuple $\xelt = (\xelt_{1},\ldots,\xelt_{\mysmallerindexM})$, let $h_{v}(\xelt) := \sum_{j=1}^{\mysmallerindexM}g_{v}(\xelt_{j})$ for each $v \in P(\omega)$, so $(h_{v}(\xelt))_{v \in P(\omega)} = \phi(\xelt)$.  
The fact that $\xelt_{1} \leq \cdots \leq \xelt_{\mysmallerindexM}$ implies that $h_{u}(\xelt) \geq h_{v}(\xelt)$ whenever $u \leq v$ in $P(\omega)$. 
In particular, $\phi(\xelt)$ is in $L(\mysmallindexM\omega)$. 
Since we can invert $\phi$, then $\phi: R(\mysmallindexM\omega) \longrightarrow L(\mysmallindexM\omega)$ is a well-defined bijection. 

Now, there is a sequence of raising operators that can be applied to our $\xelt \in R(\mysmallindexM\omega)$ in order to obtain $\melt' = \phi(\melt)$. 
That is, $\xelt =: \xelt^{(0)} \myarrow{i_1} \xelt^{(1)} \myarrow{i_2} \cdots \myarrow{i_k} \xelt^{(k)} := \melt'$ for some sequence of edges $\xelt^{(j-1)} \myarrow{i_j} \xelt^{(j)}$ in $R(\mysmallindexM\omega)$. 
Clearly, $\xelt^{(j-1)} \myarrow{i_j} \xelt^{(j)}$ in $R(\mysmallindexM\omega)$ implies that $\phi(\xelt^{(j-1)}) \myarrow{i_j} \phi(\xelt^{(j)})$ in $L(\mysmallindexM\omega)$. 
We already know that $wt(\melt) = wt(\phi(\melt))$. 
By \GaussianStructureProp, $L(\mysmallindexM\omega)$ is $\Phi$-structured. 
Putting the preceding four sentences together, we see that $wt(\xelt) = wt(\melt) -\sum_{j=0}^{k}\alpha_{i_j} = wt(\phi(\melt)) -\sum_{j=0}^{k}\alpha_{i_j} = wt(\phi(\xelt))$.\hfill\QED

\noindent 
{\bf \WeightPreservingBijectionCorollary}\ \ {\sl The diamond-colored distributive lattice} $L(\mysmallindexM\omega)$ {\sl is a splitting distributive lattice for} $\chi_{_{\mysmallestindexM\omega}}${\sl , so by \GaussianFromWGF,} $P(\mysmallindexM\omega)${\sl is Gaussian.} 

{\em Proof.} By \GaussianStructureProp, $L(\mysmallindexM\omega)$ is $\Phi$-structured, and by \CrystallineExistenceTheorem, we have $\WGF(R(\mysmallindexM\omega)) = \chi_{_{\mysmallestindexM\omega}}$. 
So, $L(\mysmallindexM\omega)$ is a splitting distributive lattice for $\chi_{_{\mysmallestindexM\omega}}$. 
By \WeightPreservingBijection, we get $\WGF(L(\mysmallindexM\omega)) = \WGF(R(\mysmallindexM\omega))$.\hfill\QED

{\bf Second new proof.} 
For our second new proof of the Proctor--Stanley Gaussian Poset Theorem, we design a vertex-coloring function $\kappa: L(\mysmallindexM\omega) \setminus \{\melt\} \longrightarrow I$ that meets the criteria of \ChainProductVersion.B. 
The key will be to demonstrate that, for each $\xelt \in L(\mysmallindexM\omega) \setminus \{\melt\}$, the kindred set $\myK(\xelt)$ is a sub-face of $\comp_{\kappa(\xelt)}(\xelt)$. 
To this end, fix a linear extension $\mathcal{T}: P(\omega) \longrightarrow \{1,2,\ldots,|P(\omega)|\}$, so $\mathcal{T}(u) \leq \mathcal{T}(v)$ whenever $u \leq v$ in $P(\omega)$. 
We regard $\mathcal{T}$ to be a total ordering of the elements of $P(\omega)$ and write $u \leq_{\mathcal{T}} v$ if and only if $\mathcal{T}(u) \leq \mathcal{T}(v)$. 
For any $\xelt \in L(\mysmallindexM\omega) \setminus \{\melt\}$, there exists a unique element $v_{\xelt} \in P(\omega)$ that is smallest in the total order $\mathcal{T}$ and has the property that $g_{v_{\xelt}}(\xelt) < \mysmallindexM$.
Declare $\kappa(\xelt)$ to be the color of the vertex $v_{\xelt}$. 
Call $v_{\xelt}$ the $P(\omega)$-{\em coloring vertex of} $\xelt$. 
Observe that if $v < v_{\xelt}$ in $P(\omega)$, then $v$ strictly precedes $v_{\xelt}$ in the total ordering $\mathcal{T}$ and therefore $g_{v}(\xelt) = \mysmallindexM$. 

\noindent 
{\bf \GaussianColoringTheorem}\ \ {\sl With the vertex-coloring function} $\kappa: L(\mysmallindexM\omega) \setminus \{\melt\} \longrightarrow I$ {\sl as described above, the kindred set} $\myK(\xelt)$ {\sl for any $\xelt \in L(\mysmallindexM\omega) \setminus \{\melt\}$ is a sub-face of $\comp_{\kappa(\xelt)}(\xelt)$.}

{\em Proof.} Let $k := \kappa(\xelt)$, and let $u_{1},\ldots,u_{r}$ be a numbering of the color $k$ vertices of $P(\omega)$ that respects the total order $\mathcal{T}$, i.e.\ $u_{p} \leq_{\mathcal{T}} u_{q}$ when $1 \leq p \leq q \leq r$. 
Since $v_{\xelt}$ must be one of these vertices, set $u_{s} := v_{\xelt}$. 
Let $\melt'$ and $\nelt'$ be, respectively, the unique maximal and unique minimal elements of $\comp_{k}(\xelt)$.
So, $\comp_{k}(\xelt)$ is isomorphic to the chain product $\prod_{q=1}^{r}\{g_{u_{q}}(\nelt'),g_{u_{q}}(\nelt')+1,\ldots,g_{u_{q}}(\melt')-1,g_{u_{q}}(\melt')\}$. 
Observe that $g_{u_{s}}(\nelt') < \mysmallindexM = g_{u_{s}}(\melt')$. 
We claim that $\myK(\xelt) = \{\yelt \in \comp_{k}(\xelt)\, |\, g_{u_{s}}(\yelt) < \mysmallindexM\}$ is a set equality, in which case it will follow immediately that $\myK(\xelt)$ is a sub-face of $\comp_{k}(\xelt)$. 

Suppose $\yelt \in \comp_{k}(\xelt)$ and $g_{u_{s}}(\yelt) < \mysmallindexM$. 
Since $\yelt \not= \melt$, then $\yelt$ has a $P(\omega)$-coloring vertex $v_{\yelt}$. 
The hypothesis `$g_{u_{s}}(\yelt) < \mysmallindexM$' means that $v_{\yelt} \leq_{\mathcal{T}} u_{s}$. 
Let $k' := \kappa(\yelt)$, which is the color of $v_{\yelt}$ in $P(\omega)$, and suppose $k' \not= k$. 
The latter hypothesis means that $v_{\yelt}$ strictly precedes $u_{s} = v_{\xelt}$ in the total order. 
But, since $\yelt \in \comp_{k}(\xelt)$, then $g_{v}(\xelt) = g_{v}(\yelt)$ for all $v \in P(\omega)$ having a color other than $k$. 
Therefore, $g_{v_{\yelt}}(\xelt) = g_{v_{\yelt}}(\yelt) < \mysmallindexM$, which implies that $v_{\xelt} \leq_{\mathcal{T}} v_{\yelt}$. 
We cannot simultaneously have $v_{\xelt} \leq_{\mathcal{T}} v_{\yelt}$ and $v_{\yelt} <_{\mathcal{T}} v_{\xelt}$, so our assumption that $k' \not= k$ does not hold. 
Then, $\kappa(\yelt) = k$, which suffices for us to conclude that $\yelt \in \myK(\xelt)$.

Suppose $\yelt \in \myK(\xelt)$, so $\yelt \in \comp_{k}(\xelt) \setminus \{\melt\}$ and $\kappa(\yelt) = k$. 
This means that the $P(\omega)$-coloring vertex $v_{\yelt}$ must be one of $\{u_{1},\ldots,u_{r}\}$, say $v_{\yelt} = u_{p}$. 
We wish to show that $g_{u_{s}}(\yelt) < \mysmallindexM$. 
Now, $g_{v}(\yelt) = \mysmallindexM$ for all $v < v_{\yelt}$ in $P(\omega)$. 
For any $q \in \{1,2,\ldots,p-1\}$, there is some $v \in P(\omega)$ with color other $k$ and such that $u_{q} < v < u_{p}$; it follows that $g_{u_{q}}(\xelt) = \mysmallindexM$ for each such $q$. 
Therefore $p \leq s$. 
So $g_{u_{s}}(\yelt) \leq g_{u_{p}}(\yelt) = g_{v_{\yelt}}(\yelt) < \mysmallindexM$.\hfill\QED

\noindent 
{\bf \GaussianColoringCorollary}\ \ {\sl The diamond-colored distributive lattice} $L(\mysmallindexM\omega)$ {\sl is a splitting distributive lattice for} $\chi_{_{\mysmallestindexM\omega}}${\sl , so by \GaussianFromWGF,} $P(\mysmallindexM\omega)${\sl is Gaussian.}

{\em Proof.} By \GaussianStructureProp, $L(\mysmallindexM\omega)$ is $\Phi$-structured and has monochromatic components that are isomorphic to chain products. 
By \GaussianColoringTheorem, we can apply \ChainProductVersion.B to conclude that $L(\mysmallindexM\omega)$ is a splitting distributive lattice for $\chi_{_{\mysmallestindexM\omega}}$.\hfill\QED


\newpage
%
\renewcommand{\refname}{\Large \bf References}
\renewcommand{\baselinestretch}{1.1}

\end{document}